\def\date{July 11, 2006}

\input amssym.def
\input amssym.tex

\def\item#1{\vskip1.3pt\hang\textindent {\rm #1}}


\newskip\litemindent
\litemindent=0.7cm  
\def\Litem#1#2{\par\noindent\hangindent#1\litemindent
\hbox to #1\litemindent{\hfill\hbox to \litemindent
{\ninerm #2 \hfill}}\ignorespaces}
\def\litem{\Litem1}
\def\litemitem{\Litem2}

\tolerance=300
\pretolerance=200
\hfuzz=1pt
\vfuzz=1pt

\hoffset=0in
\voffset=0.5in

\hsize=5.8 true in 
\vsize=9.2 true in
\parindent=25pt
\mathsurround=1pt
\parskip=1pt plus .25pt minus .25pt
\normallineskiplimit=.99pt

\countdef\revised=100
\mathchardef\emptyset="001F 
\chardef\ss="19
\def\3{\ss}
\def\anf{$\lower1.2ex\hbox{"}$}
\def\frac#1#2{{#1 \over #2}}
\def\>{>\!\!>}
\def\<{<\!\!<}

\def\into{\hookrightarrow}
\def\onto{\to\mskip-14mu\to} 
\def\ssssarr{\hbox to 15pt{\rightarrowfill}}
\def\sssarr{\hbox to 20pt{\rightarrowfill}}
\def\ssarr{\hbox to 30pt{\rightarrowfill}}
\def\sarr{\hbox to 40pt{\rightarrowfill}}
\def\arr{\hbox to 60pt{\rightarrowfill}}
\def\larr{\hbox to 60pt{\leftarrowfill}}
\def\Arr{\hbox to 80pt{\rightarrowfill}}
\def\mapdown#1{\Big\downarrow\rlap{$\vcenter{\hbox{$\scriptstyle#1$}}$}}

\def\sssmapright#1{\smash{\mathop{\sssarr}\limits^{#1}}}
\def\ssmapright#1{\smash{\mathop{\ssarr}\limits^{#1}}}
\def\smapright#1{\smash{\mathop{\sarr}\limits^{#1}}}

\def\mapup#1{\Big\uparrow\rlap{$\vcenter{\hbox{$\scriptstyle#1$}}$}}
\def\lmapup#1{\llap{$\vcenter{\hbox{$\scriptstyle#1$}}$}\Big\uparrow}

\def\indlim{{\displaystyle \lim_{\longrightarrow}}\ }
\def\prolim{{\displaystyle \lim_{\longleftarrow}}\ }

\def\ad{\mathop{\rm ad}\nolimits}

\def\Ad{\mathop{\rm Ad}\nolimits}

\def\Aut{\mathop{\rm Aut}\nolimits}

\def\Idem{\mathop{\rm Idem}\nolimits}

\def\der{\mathop{\rm der}\nolimits}
\def\det{\mathop{\rm det}\nolimits}

\def\Diff{\mathop{\rm Diff}\nolimits}

\def\Exp{\mathop{\rm Exp}\nolimits}

\def\evol{\mathop{\rm evol}\nolimits}
\def\ev{\mathop{\rm ev}\nolimits}

\def\End{\mathop{\rm End}\nolimits}
\def\Ext{\mathop{\rm Ext}\nolimits}

\def\G{\mathop{\rm G{}}\nolimits}

\def\Gau{\mathop{\rm Gau}\nolimits}
\def\GL{\mathop{\rm GL}\nolimits}

\def\hats #1{\hat{\hat{\hbox{$#1$}}}}

\def\Hom{\mathop{\rm Hom}\nolimits}%
\def\id{\mathop{\rm id}\nolimits} 
\def\im{\mathop{\rm im}\nolimits}

\def\Im{\mathop{\rm Im}\nolimits}
\def\inf{\mathop{\rm inf}\nolimits}

\def\Inn{\mathop{\rm Inn}\nolimits}
\def\Int{\mathop{\rm int}\nolimits}

\def\OO{\mathop{\rm O{}}\nolimits}

\def\PU{\mathop{\rm PU}\nolimits}
\def\per{\mathop{\rm per}\nolimits}

\def\Rot{\mathop{\rm Rot}\nolimits}

\def\sgn{\mathop{\rm sgn}\nolimits}
\def\SL{\mathop{\rm SL}\nolimits}

\def\span{\mathop{\rm span}\nolimits}
\def\Skew{\mathop{\rm Skew}\nolimits}

\def\Spec{\mathop{\rm Spec}\nolimits}

\def\SU{\mathop{\rm SU}\nolimits}
\def\sup{\mathop{\rm sup}\nolimits}
\def\supp{\mathop{\rm supp}\nolimits}

\def\tr{\mathop{\rm tr}\nolimits}
\def\Tr{\mathop{\rm Tr}\nolimits}
\def\trile{\trianglelefteq}

\def\UU{\mathop{\rm U{}}\nolimits}

\def\0{{\bf 0}}
\def\1{{\bf 1}}

\def\a{{\frak a}}
\def\aff{{\frak {aff}}}

\def\gau{{\frak {gau}}}

\def\e{{\frak e}}

\def\g{{\frak g}}
\def\gl{{\frak {gl}}}
\def\h{{\frak h}}

\def\k{{\frak k}}

\def\n{{\frak n}}
\def\oo{{\frak o}}

\def\q{{\frak q}}
\def\r{{\frak r}}
\def\s{{\frak s}}

\def\ham{{\frak {ham}}}

\def\su{{\frak {su}}}

\def\sL{{\frak {sl}}}

\def\uu{{\frak u}}

\def\z{{\frak z}}

\def\L{\mathop{\bf L{}}\nolimits}

\def\C{{{\Bbb C}{\mskip+1mu}}} 
\def\K{{{\Bbb K}{\mskip+2mu}}} 

\def\R{{\Bbb R}} 
\def\Z{{\Bbb Z}} 
\def\N{{\Bbb N}} 
 
\def\E{{\Bbb E}} 
\def\F{{\Bbb F}} 
\def\K{{\Bbb K}}

\def\P{{\Bbb P}} 
\def\Q{{\Bbb Q}} 
\def\V{{\Bbb V}} 
\def\SS{{\Bbb S}} 
\def\T{{\Bbb T}} 

\def\:{\colon}  
\def\.{{\cdot}}
\def\|{\Vert}
\def\bsk{\bigskip}

\def\giantskip{\vskip2\bigskipamount}
\def\gsk{\giantskip}
\def \la {\langle}

\def\msk{\medskip}
\def \ra {\rangle}
\def \res {\!\mid\!\!}

\def\bbr{\bigbreak}
\def\giantbreak{\par \ifdim\lastskip<2\bigskipamount \removelastskip
         \penalty-400 \giantskip\fi}

\def\nin{\noindent}
\def\cen{\centerline}
\def\pagebreak{\vskip 0pt plus 0.0001fil\break}
\def\linebreak{\break}

\def\hat{\widehat}

\def\derat#1{{d \over dt} \hbox{\vrule width0.5pt 
                height 5mm depth 3mm${{}\atop{{}\atop{\scriptstyle t=#1}}}$}}

\def\eps{\varepsilon}
\def\epsilon{\varepsilon}
\def\eset{\emptyset}

\def\nin{\noindent}
\def\oline{\overline}

\def\pder#1,#2,#3 { {\partial #1 \over \partial #2}(#3)}
\def\pde#1,#2 { {\partial #1 \over \partial #2}}
\def\phi{\varphi}


\def\subeq{\subseteq}

\def\Rarrow{\Rightarrow}

\def\tilde{\widetilde}

\font\ninerm=cmr9
\font\eightrm=cmr8

\font\eightbf=cmbx8


\font\smc=cmcsc10
\font\bfone=cmbx10 scaled\magstep1 
\font\bftwo=cmbx10 scaled\magstep2 

\def\qed{{\unskip\nobreak\hfil\penalty50\hskip .001pt \hbox{}\nobreak\hfil
          \vrule height 1.2ex width 1.1ex depth -.1ex
           \parfillskip=0pt\finalhyphendemerits=0\medbreak}\rm}

\def\qeddis{\eqno{\vrule height 1.2ex width 1.1ex depth -.1ex} $$
                   \medbreak\rm}

\def\Lemma #1. {\bigbreak\vskip-\parskip\noindent{\bf Lemma #1.}\quad\it}

\def\Sublemma #1. {\bigbreak\vskip-\parskip\noindent{\bf Sublemma #1.}\quad\it}

\def\Proposition #1. {\bigbreak\vskip-\parskip\noindent{\bf Proposition #1.}
\quad\it}

\def\Corollary #1. {\bigbreak\vskip-\parskip\nin{\bf Corollary #1.}
\quad\it}

\def\Theorem #1. {\bigbreak\vskip-\parskip\noindent{\bf Theorem #1.}
\quad\it}

\def\Definition #1. {\rm\bigbreak\vskip-\parskip\noindent
{\bf Definition #1.}
\quad}

\def\Remark #1. {\rm\bigbreak\vskip-\parskip\noindent{\bf Remark #1.}\quad}

\def\Example #1. {\rm\bigbreak\vskip-\parskip\noindent{\bf Example #1.}\quad}
\def\Examples #1. {\rm\bigbreak\vskip-\parskip\noindent{\bf Examples #1.}\quad}

\def\Problems #1. {\bigbreak\vskip-\parskip\noindent{\bf Problems #1.}\quad}
\def\Problem #1. {\bigbreak\vskip-\parskip\noindent{\bf Problem #1.}\quad}
\def\Exercise #1. {\bigbreak\vskip-\parskip\noindent{\bf Exercise #1.}\quad}

\def\Conjecture #1. {\bigbreak\vskip-\parskip\noindent{\bf Conjecture #1.}\quad}

\def\Proof#1.{\rm\par\ifdim\lastskip<\bigskipamount\removelastskip\fi\smallskip
            \noindent {\bf Proof.}\quad}

\def\Axiom #1. {\bigbreak\vskip-\parskip\noindent{\bf Axiom #1.}\quad\it}

\def\Satz #1. {\bigbreak\vskip-\parskip\noindent{\bf Satz #1.}\quad\it}

\def\Korollar #1. {\bbr\vskip-\parskip\nin{\bf Korollar #1.} \quad\it}

\def\Folgerung #1. {\bbr\vskip-\parskip\nin{\bf Folgerung #1.} \quad\it}

\def\Folgerungen #1. {\bbr\vskip-\parskip\nin{\bf Folgerungen #1.} \quad\it}

\def\Bemerkung #1. {\rm\bigbreak\vskip-\parskip\noindent{\bf Bemerkung #1.}
\quad}

\def\Beispiel #1. {\rm\bigbreak\vskip-\parskip\noindent{\bf Beispiel #1.}\quad}
\def\Beispiele #1. {\rm\bigbreak\vskip-\parskip\noindent{\bf Beispiele #1.}\quad}
\def\Aufgabe #1. {\rm\bigbreak\vskip-\parskip\noindent{\bf Aufgabe #1.}\quad}
\def\Aufgaben #1. {\rm\bigbreak\vskip-\parskip\noindent{\bf Aufgabe #1.}\quad}

\def\Beweis#1. {\rm\par\ifdim\lastskip<\bigskipamount\removelastskip\fi
           \smallskip\noindent {\bf Beweis.}\quad}

\nopagenumbers

\def\date{\ifcase\month\or January\or February \or March\or April\or May
\or June\or July\or August\or September\or October\or November
\or December\fi\space\number\day, \number\year}

\def\title{Title ??}
\def\author{Author ??}

\def\thanks#1{\footnote*{\eightrm#1}}

\def\rightheadline{\hfil{\eightrm\title}\hfil\tenbf\folio}
\def\leftheadline{\tenbf\folio\hfil{\eightrm\author}\hfil}
\headline={\vbox{\line{\ifodd\pageno\rightheadline\else\leftheadline\fi}}}

\def\firstheadline{}
\def\firstfootline{\cen{\rm\folio}}

\def\seite #1 {\pageno #1
               \headline={\ifnum\pageno=#1 \firstheadline
               \else\ifodd\pageno\rightheadline\else\leftheadline\fi\fi}
               \footline={\ifnum\pageno=#1 \firstfootline\else{}\fi}}

\newdimen\dimenone
 \def\checkleftspace#1#2#3#4{
 \dimenone=\pagetotal
 \advance\dimenone by -\pageshrink   
 \ifdim\dimenone>\pagegoal          
   \else\dimenone=\pagetotal
        \advance\dimenone by \pagestretch
        \ifdim\dimenone<\pagegoal
          \dimenone=\pagetotal
          \advance\dimenone by#1         
          \setbox0=\vbox{#2\parskip=0pt                
                     \hyphenpenalty=10000
                     \rightskip=0pt plus 5em
                     \noindent#3 \vskip#4}    
        \advance\dimenone by\ht0
        \advance\dimenone by 3\baselineskip   
        \ifdim\dimenone>\pagegoal\vfill\eject\fi
          \else\eject\fi\fi}


\def\subheadline #1{\nin\bigbreak\vskip-\lastskip
      \checkleftspace{0.9cm}{\bf}{#1}{\medskipamount}
          \indent\vskip0.7cm\centerline{\bf #1}\medskip}
\def\subsection{\subheadline} 

\def\lsubheadline #1 #2{\bigbreak\vskip-\lastskip
      \checkleftspace{0.9cm}{\bf}{#1}{\bigskipamount}
         \vbox{\vskip0.7cm}\cen{\bf #1}\msk \cen{\bf #2}\bsk}

\def\sectionheadline #1{\bigbreak\vskip-\lastskip
      \checkleftspace{1.1cm}{\bf}{#1}{\bigskipamount}
         \vbox{\vskip1.1cm}\cen{\bfone #1}\bsk}
\def\section{\sectionheadline} 

\def\lsectionheadline #1 #2{\bigbreak\vskip-\lastskip
      \checkleftspace{1.1cm}{\bf}{#1}{\bigskipamount}
         \vbox{\vskip1.1cm}\cen{\bfone #1}\msk \cen{\bfone #2}\bsk}

\def\lchapterheadline #1 #2{\bigbreak\vskip-\lastskip\indent\vskip3cm
                       \cen{\bftwo #1} \msk \cen{\bftwo #2} \gsk}
\def\llsectionheadline #1 #2 #3{\bigbreak\vskip-\lastskip\indent\vskip1.8cm
\cen{\bfone #1} \msk \cen{\bfone #2} \msk \cen{\bfone #3} \nobreak\bsk\nobreak}


\newtoks\literat
\def\[#1 #2\par{\literat={#2\unskip.}%
\hbox{\vtop{\hsize=.15\hsize\nin [#1]\hfill}
\vtop{\hsize=.82\hsize\nin\the\literat}}\par
\vskip.3\baselineskip}

\def\references{
\sectionheadline{\bf References}
\frenchspacing

\entries\par}

\mathchardef\emptyset="001F 
\def\address{Author: \tt$\backslash$def$\backslash$address$\{$??$\}$}

\def\abstract #1{{\narrower\baselineskip=10pt{\noindent
\eightbf Abstract.\quad \eightrm #1 }
\bigskip}}

\def\firstpage{\nin
{\obeylines \parindent 0pt }
\vskip2cm
\centerline{\bfone\title}
\gsk
\centerline{\bf\author}
\vskip1.5cm \rm}

\def\lastpage{\par\vbox{\vskip1cm\nin
\line{
\vtop{\hsize=.5\hsize{\parindent=0pt\baselineskip=10pt\nin\address}}
\hfill} }}

\def\Box #1 { \msk\par\nin 
\centerline{
\vbox{\offinterlineskip
\hrule
\hbox{\vrule\strut\hskip1ex\hfil{\smc#1}\hfill\hskip1ex}
\hrule}\vrule}\msk }

\def\adots{\mathinner{\mkern1mu\raise1pt\vbox{\kern7pt\hbox{.}}
                        \mkern2mu\raise4pt\hbox{.}
                        \mkern2mu\raise7pt\hbox{.}\mkern1mu}}


\pageno=1

\def\title{Towards a Lie theory of locally convex groups} 
\def\author{Karl-Hermann Neeb}
\def\leftheadline{\tenbf\folio\hfil{\tt japsurv.tex}\hfil\eightrm\date}

\def\Ham{\mathop{\rm Ham}\nolimits}

\def\Log{\mathop{\rm Log}\nolimits}
\def\fL{{\frak L}} 
\def\Gh{\mathop{\rm Gh}\nolimits}
\def\GaL{\mathop{\rm \Gamma L}\nolimits}
\def\Gf{\mathop{\rm Gf}\nolimits}
\def\Gs{\mathop{\rm Gs}\nolimits}
\def\Idem{\mathop{\rm Idem}\nolimits}
\def\Vir{\mathop{\rm Vir}\nolimits}
\def\Vol{\mathop{\rm Vol}\nolimits}
\def\evol{\mathop{\rm evol}\nolimits}
\def\Fl{\mathop{\rm Fl}\nolimits}

\def\gh{{\frak {gh}}}
\def\gs{{\frak {gs}}}
\def\gf{{\frak {gf}}}

\def\V{{{\Bbb V}}} 
\def\G{{{\Bbb G}}}

\firstpage 

\abstract{In this survey, we report on the state of the art of some of 
the fundamental problems in the Lie theory of Lie groups modeled 
on locally convex spaces, such as integrability of Lie algebras, 
integrability of Lie subalgebras to Lie subgroups, and integrability of 
Lie algebra extensions to Lie group extensions. We further 
describe how regularity or local exponentiality of a Lie group 
can be used to obtain quite satisfying answers to some of the 
fundamental problems. These results are illustrated by specialization 
to some specific classes of Lie groups, such as direct limit groups, 
linear Lie groups, groups of smooth maps and groups of diffeomorphisms. \hfill\linebreak
AMS Classification: 22E65, 22E15 \hfill\linebreak 
}

\sectionheadline{Introduction} 

\nin Symmetries play a decisive role in the natural sciences and 
throughout mathematics. Infinite-dimensional Lie theory 
deals with symmetries depending on infinitely many
parameters. Such symmetries may be studied on an
infinitesimal, local or global level, which amounts 
to studying Lie algebras, local Lie groups
and global Lie groups, respectively. Here the passage from the 
infinitesimal to the local level requires a smooth structure 
on the symmetry group (such as a Lie group structure as defined below), whereas the 
passage from the local to the global level is purely 
topological. 

Finite-dimensional Lie theory was created 
about 130 years ago by {\smc Marius Sophus Lie} and {\smc Friedrich Engel}, 
who showed that in
finite dimensions the local and the infinitesimal 
theory are essentially equivalent ([Lie80/95]). 
The differential geometric approach to 
finite-dimensional global Lie groups (as smooth or
analytic manifolds) is naturally complemented by the
theory of algebraic groups with which it interacts most 
fruitfully. A crucial point of the finite-dimensional theory
is that finiteness conditions permit to develop a powerful 
structure theory of finite-dimensional Lie groups in terms
of the Levi splitting and the fine structure of semisimple 
Lie groups ([Ho65], [Wa72]). 

A substantial part of the literature
on infinite-dimensional Lie
theory exclusively deals with the level
of Lie algebras, their structure, and their
representations (cf.\ [Ka90], [Neh96], 
[Su97], [AABGP97], [DiPe99], [ABG00]). 
However, 
only special classes of groups, such as Kac--Moody groups, can
be approached with success by purely algebraic methods ([KP83], [Ka85]); see also 
[Rod89] for an analytic approach to Kac--Moody groups. 
In mathematical physics,
the infinitesimal approach dealing mainly with Lie algebras 
and their representations is convenient for calculations, 
but a global analytic perspective is required 
to understand global phenomena (cf.\ [AI95], [Ot95], 
[CVLL98], [EMi99], [Go04], [Sch04]).
A similar statement applies to 
non-commutative geometry, throughout of which derivations 
and covariant derivatives are used. It 
would be interesting to understand how global symmetry groups 
and the associated geometry fit into the picture ([Co94], [GVF01]). 

In infinite dimensions, the passage from the infinitesimal to the local
 and from there to the global level is not possible
in general, whence Lie theory splits into three properly distinct levels.
It is a central point of this survey to explain some of 
the concepts and the results that can be used to translate 
between these three levels. 

We shall use a Lie group concept which is both simple and very general: 
A Lie group is a manifold, endowed with 
a group structure
such that multiplication and inversion are smooth maps. 
The main difference, compared
to the finite-dimensional theory, 
concerns the notion of a manifold:
The manifolds we consider 
are modeled on arbitrary locally convex spaces.
As we shall see later, it is natural to approach Lie groups
from such a general perspective,
because it leads to a unified treatment of
the basic aspects of the theory without unnatural restrictions 
on model spaces or the notion of a Lie group. Although we shall simply call them 
{\it Lie groups}, a more specific terminology is {\it locally convex Lie groups}. Depending on the type 
of the model spaces, we obtain in particular the classes of 
finite-dimensional, Banach--, Fr\'echet--, LF-- and Silva--Lie groups. 

\subheadline{The fundamental problems of Lie theory} 

As in finite dimensions, the Lie algebra $\L(G)$ of a Lie group 
$G$ is identified with the tangent space $T_\1(G)$, where the Lie bracket 
it obtained by identification with the space of left invariant vector fields. 
This turns $\L(G)$ into a locally convex (topological) Lie algebra. 
Associating, furthermore, to a morphism $\phi$ of Lie groups its tangent map 
$\L(\phi) := T_\1(\phi)$, 
we obtain the {\it Lie functor} from the category of (locally convex) Lie 
groups to the category of locally convex topological Lie algebras. 
The core of Lie theory now consists in determining how much information the Lie 
functor forgets and how much can be reconstructed. This leads to several 
integration problems such as: 
\litemindent=1.1cm 
\litem{(FP1)} When does a continuous homomorphism $\L(G) \to \L(H)$ between Lie algebras of 
connected Lie groups integrate to a (local) group homomorphism $G \to H$? 
\litem{(FP2)} {\bf Integrability Problem for subalgebras:} Which 
Lie subalgebras $\h$ of the Lie algebra $\L(G)$ of a Lie group $G$ 
correspond to Lie group morphisms $H \to G$ 
with $\L(H) = \h$? 
\litem{(FP3)} {\bf Integrability Problem for Lie algebras (LIE III):} 
For which locally convex Lie algebras $\g$ does a local/global 
Lie group $G$ with $\L(G) = \g$ exist? 
\litem{(FP4)} {\bf Integrability Problem for extensions:} 
When does an extension of the Lie algebra 
$\L(G)$ of a Lie group $G$ by the Lie algebra $\L(N)$ of a Lie group $N$ 
integrate to a Lie group extension of $G$ by $N$? 
\litem{(FP5)} {\bf Subgroup Problem:} Which subgroups of a Lie group $G$ carry natural 
Lie group structures? 
\litem{(FP6)} When does a Lie group have an exponential map $\exp_G \: \L(G) \to G$? 
\litem{(FP7)} {\bf Integrability Problem for smooth actions:} When 
does a homomorphism $\g \to {\cal V}(M)$ into the Lie algebra ${\cal V}(M)$ of 
smooth vector fields on a manifold $M$ integrate to a smooth action of a 
corresponding Lie group? 
\litem{(FP8)} {\bf Small Subgroup Problem:} Which Lie groups have identity 
neighborhoods containing no non-trivial subgroup? 
\litem{(FP9)} {\bf Locally Compact Subgroup Problem:} For which Lie groups 
are locally compact subgroups (finite-dimensional) Lie groups? 
\litem{(FP10)} {\bf Automatic Smoothness Problem:} When are 
continuous homomorphisms between Lie groups smooth? 
\litemindent=0.7cm 

An important tool in the finite-dimensional and Banach context is the exponential map, 
but as vector fields on locally convex manifolds need not possess integral 
curves, there is no general theorem that guarantees the existence of a (smooth) 
exponential map, i.e., a smooth function 
$$\exp_G \: \L(G) \to G,$$ 
for which the curves $\gamma_x(t) := \exp_G(tx)$ are homomorphisms 
$(\R,+) \to G$ with $\gamma_x'(0)= x$. Therefore the existence of an 
exponential function has to be treated as an additional requirement (cf.\ (FP6)). 
Even stronger 
is the requirement of {\it regularity}, meaning that, for each 
smooth curve $\xi \: [0,1] \to \L(G)$, the initial value problem 
$$ \gamma'(t) =\gamma(t).\xi(t) := T_\1(\lambda_{\gamma(t)})\xi(t), \quad \gamma(0) = \1 $$
has a solution $\gamma_\xi \:[0,1] \to G$ and that 
$\gamma_\xi(1)$ depends smoothly on $\xi$. 
Regularity is a natural assumption that provides a good deal of methods to 
pass from the infinitesimal to the global level. This regularity 
concept is due to {\smc Milnor} ([Mil84]). 
It weakens the $\mu$-regularity (in our terminology) introduced 
by {\smc Omori} et al. in [OMYK82/83a] (see [KYMO85] for a survey), 
but it is still strong enough for the essential 
Lie theoretic applications.  
Presently, we do not know of any Lie group modeled on a complete space 
which is not regular. 
For all major concrete classes discussed below,  
one can prove regularity, but there is 
no general theorem saying that each locally convex Lie group with a complete 
model space is regular or even that it has an exponential function. 
To prove or disprove such a theorem is another fundamental open problem of the theory. 
An assumption of a different nature than regularity, 
and which can be used to develop a profound 
Lie theory, is that $G$ is {\it locally exponential} in the sense that 
it has an exponential function which is a local diffeomorphism in $0$. 
Even stronger is the assumption that, in addition, $G$ is analytic and that 
the exponential function is an analytic local diffeomorphism in $0$. 
Groups with this property are called {\it BCH--Lie groups}, because the local 
multiplication in canonical local coordinates is given by the 
Baker--Campbell--Hausdorff series. This class contains in particular all Banach--Lie groups. 

\subheadline{Important classes of infinite-dimensional Lie groups} 

Each general theory lives from the concrete classes of objects it can be applied to. 
Therefore it is 
good to have certain major classes of Lie groups 
in mind to which the general theory should apply. Here we briefly describe four 
such classes:  

{\bf Linear Lie groups:} Loosely speaking, linear Lie groups are Lie groups 
of operators on locally convex spaces, but this has to be made more precise. 

Let $E$ be a locally convex space and ${\cal L}(E)$ the unital associative 
algebra of all continuous linear endomorphisms of $E$. Its unit group 
is the general linear group $\GL(E)$ of $E$. If  $E$ is not normable, 
there is no vector topology on ${\cal L}(E),$ 
for which the composition map is continuous ([Mais63, Satz 2]). 
In general, the group $\GL(E)$ is open for no vector topology 
on ${\cal L}(E)$, as follows from the observation that 
if the spectrum of the operator $D$ is unbounded, then $\1 + t D$ is not 
invertible for arbitrarily small values of $t$. 
Therefore we need a class of associative algebras which are better behaved 
than ${\cal L}(E)$ to define a natural class of linear Lie groups. 

The most natural class of associative algebras for 
infinite-dimensional Lie theory are {\it continuous inverse algebras} (CIAs for short), 
introduced in [Wae54a/b/c] by {\smc Waelbroeck} in the context of commutative spectral 
theory. A CIA is a unital associative locally convex algebra $A$ 
with continuous multiplication,  
for which the unit group $A^\times$ is open and the inversion 
$A^\times \to A, a \mapsto a^{-1}$ is a continuous map. As this implies the smoothness of 
the inversion map, $A^\times$ carries a natural Lie group structure. 
It is not hard to see that the CIA property is inherited by matrix algebras 
$M_n(A)$ over $A$ ([Swa62]), so that $\GL_n(A) = M_n(A)^\times$ also is a Lie group, 
and under mild completeness assumptions (sequential completeness) on $A$, the 
Lie group $A^\times$ is regular and locally exponential ([Gl02b], [GN06]). 
This leads to natural classes of Lie subgroups of CIAs and hence to a natural concept of a {\it linear Lie 
group}. 

{\bf Mapping groups and gauge groups:} There are many natural classes of groups 
of maps with values in Lie groups which can be endowed with Lie group 
structures. The most important cases are the following: 
If $M$ is a compact manifold and $K$ a Lie group (possibly infinite-dimensional) 
with Lie algebra $\L(K) = \k$,
then the group $C^\infty(M,K)$ always carries a natural Lie group structure 
such that $C^\infty(M,\k)$, endowed with the pointwise bracket, 
is its Lie algebra ([GG61], [Mil82/84]; see also [Mi80]). A prominent class of such groups 
are the smooth loop groups $C^\infty(\SS^1,K)$, which, for finite-dimensional simple 
groups $K$, are closely related to Kac--Moody groups ([PS86], [Mick87/89]). 

If, more generally, $q \: P \to M$ is a smooth principal bundle over a compact 
manifold $M$ with locally exponential structure group $K$, then its 
{\it gauge group} 
$$ \Gau(P) := \{ \phi \in \Diff(P) \: q\circ \phi = q, 
(\forall p \in P, k \in K)\ \phi(p.k) = \phi(p).k\} $$
also carries a natural Lie group structure. For trivial bundles, we 
obtain the mapping groups $C^\infty(M,K)$ as special cases. 
Natural generalizations are the groups $C^\infty_c(M,K)$ of compactly supported 
smooth maps on a $\sigma$-compact finite-dimensional manifold ([Mi80], [Gl02c]) 
and also Sobolev completions of the groups $C^\infty(M,K)$ ([Sch04]). 

{\bf Direct limit groups:} A quite natural method to obtain infinite-dimensional 
groups from finite-dimensional ones is to consider 
a sequence $(G_n)_{n \in\N}$ of 
finite-dimensional Lie groups and morphisms 
$\phi_n \: G_n \to G_{n+1}$, so that we can define a direct limit group 
$G := \indlim G_n$ whose representations correspond to compatible 
sequences of 
representations of the groups $G_n$. 
According to a recent theorem of {\smc Gl\"ockner} ([Gl05]), 
generalizing previous work of {\smc J.~Wolf} and his coauthors ([NRW91/93]),  
the direct limit group $G$ can always be endowed with a natural Lie group 
structure. Its Lie algebra $\L(G)$ is the countably dimensional 
direct limit space $\indlim \L(G_n)$, endowed with the finest locally convex 
topology. This provides an interesting class of infinite-dimensional Lie groups
which is still quite close to finite-dimensional groups and has a very rich 
representation theory ([DiPe99], [NRW99], [Wol05]). 

{\bf Groups of diffeomorphisms:} In a similar fashion as linear Lie groups 
arise as symmetry groups of linear structures, such as bilinear forms on 
modules of CIAs, Lie groups of diffeomorphisms arise as symmetry groups of 
geometric structures on manifolds, such as symplectic structures, contact 
structures or volume forms. 

A basic result is that, for any compact manifold $M$,  
the group $\Diff(M)$ can be turned 
into a Lie group modeled on the Fr\'echet space ${\cal V}(M)$ of smooth vector fields 
on $M$ (cf.\ [Les67], [Omo70], [EM69/70], [Gu77], [Mi80], [Ham82]). 

If $M$ is non-compact and finite-dimensional, but $\sigma$-compact, 
then there is no natural Lie group structure on $\Diff(M)$ such that smooth 
actions of Lie groups $G$ on $M$ correspond to Lie group homomorphisms 
$G \to \Diff(M)$. Nevertheless, it is possible to turn $\Diff(M)$ 
into a Lie group with Lie algebra 
${\cal V}_c(M)$, the Lie algebra of all smooth vector fields with compact support, 
endowed with the natural test function topology, turning it into an LF space 
(cf.\ [Mi80], [Mil82], [Gl02d]). 
If $M$ is compact, this yields the aforementioned 
Lie group structure on $\Diff(M)$, but if 
$M$ is not compact, then the corresponding topology on $\Diff(M)$ is so fine that 
the global flow generated by a vector field whose support is not compact does not 
lead to a continuous homomorphism $\R \to \Diff(M)$. 
For this Lie group structure, the normal subgroup 
$\Diff_c(M)$ of all diffeomorphisms  which coincide with $\id_M$ outside 
a compact set is an open subgroup. 

By {\it groups of diffeomorphisms} we mean groups of the type 
$\Diff_c(M)$, as well as natural subgroups defined as symmetry groups 
of geometric structures, such as groups of symplectomorphisms, groups of contact 
transformations and groups of volume preserving diffeomorphisms. 
Of a different nature, but also locally convex Lie groups, are 
groups of formal diffeomorphism as studied by 
{\smc Lewis} ([Lew39]), {\smc Sternberg} ([St61]) and {\smc Kuranishi} ([Kur59]), 
groups of germs of smooth and analytic diffeomorphisms of $\R^n$ fixing $0$ ([RK97], [Rob02]), 
and also germs of biholomorphic maps of $\C^n$ fixing $0$ ([Pis76/77/79]). 

As the discussion of these classes of examples shows, the 
concept of a locally convex Lie group subsumes quite different classes 
of Lie groups: Banach--Lie groups, 
groups of diffeomorphisms (modeled on Fr\'echet and LF spaces), 
groups of germs (modeled on Silva spaces) and formal groups (modeled 
on Fr\'echet spaces such as $\R^\N$). 

\msk 

In this survey article, we present our 
personal view of the current state of several aspects of the Lie theory of locally 
convex Lie groups. We shall focus on the general structures and concepts 
related to the fundamental problems (FP1)-(FP10) 
and on what can be said for the classes of Lie groups mentioned above. 

Due to limited space and time, we had to make choices, and as a result, 
we could not take up many interesting 
directions such as the modern theory of symmetries of differential equations as exposed 
in {\smc Olver}'s beautiful book [Olv93] and  
the fine structure and the geometry of specific groups of 
diffeomorphisms, such as the group $\Diff(M,\omega)$ of symplectomorphisms 
of a symplectic manifold $(M,\omega)$ ([Ban97], [MDS98], [Pol01] are recent  
textbooks on this topic). We do not go into 
(unitary) representation theory (cf.\ [AHMTT93], [Is96], [DP03], [Pic00a/b], [Ki05])  
and connections to physics, which are nicely described in 
the recent surveys of {\smc Goldin} [Go04] and {\smc Schmid} [Sch04]. 
Other topics are only mentioned very briefly, such as 
the ILB and ILH-theory of Lie groups of diffeomorphisms 
which plays an important role in geometric analysis (cf.\ [AK98], [EMi99]) 
and the group of invertible Fourier integral operators of order zero, 
whose Lie group property 
was the main goal of the series of papers [OMY80/81], [OMYK81/82/83a/b], completed in 
[MOKY85]. An alternative approach to these groups is described in 
[ARS84,86a/b]. More recently, very interesting results concerning diffeomorphism 
groups and Fourier integral operators on non-compact manifolds (with bounded geometry) 
have been obtained by {\smc Eichhorn} and {\smc Schmid} ([ES96/01]). 

\subheadline{Some history} 

To put the Lie theory of locally convex groups into proper perspective, we take a 
brief look at the historical development of infinite-dimensional Lie theory. 
Infinite-dimensional Lie algebras, such as 
Lie algebras of vector fields, where present in Lie theory right 
from the beginning, when {\smc Sophus Lie} started to study (local) Lie groups 
as groups ``generated'' by finite-dimensional Lie algebras 
of vector fields ([Lie80]). 
The general global theory 
of finite-dimensional Lie groups started to develop in the 
late 19th century, driven substantially by {\smc \'E.~Cartan}'s work 
on symmetric spaces ([CaE98]). 
The first exposition of a global theory, 
including the description of all connected groups with a given Lie algebra 
and analytic subgroups, was given by {\smc Mayer} and {\smc Thomas} 
([MaTh35]). After the combination 
with the structure theory of Lie algebras, the theory reached 
its mature form in the middle of the 20th century, which is exposed in 
the fundamental books of {\smc Chevalley} ([Ch46]) and {\smc Hochschild} ([Ho65]) 
(see also [Po39] for an early textbook situated on the borderline between topological groups 
and Lie groups). 

Already in the work of {\smc S.~Lie} infinite-dimensional groups show up as groups of 
(local) diffeomorphisms of open domains in $\R^n$ (cf.\ [Lie95]) 
and later {\smc \'E.~Cartan} undertook a more systematic study of certain types of infinite-dimensional 
Lie algebras, resp., groups of diffeomorphisms preserving geometric structures 
on a manifold, such a symplectic, contact or volume forms (cf.\ [CaE04]). 
On the other hand, the advent of Quantum Mechanics in the 1920s created 
a need to understand the structure of 
groups of operators on Hilbert spaces, which is a quite different class of 
infinite-dimensional groups (cf.\ [De32]). 

The first attempt to deal with infinite-dimensional groups as Lie groups, i.e., 
as a smooth manifolds,   
was undertaken by {\smc Birkhoff} in [Bir36/38], where he developed the local 
Lie theory of Banach--Lie groups, resp., Banach--Lie algebras (see also [MiE37] 
for first steps in extending Lie's theory of local transformation groups to the 
Banach setting). 
In particular, he proved that (locally) $C^1$-Banach--Lie groups 
admit exponential coordinates which leads to analytic Lie group structures, 
that continuous homomorphisms are analytic and  
that, for every Banach--Lie 
algebra, the BCH series defines an analytic local group structure. 
He also defines the Lie algebra of such a local group, 
derives its functoriality properties and establishes the correspondence 
between closed subalgebras/ideals and the corresponding local subgroups. 
Even product integrals, which play a central role in the modern theory, 
appear in his work as solutions of left invariant differential equations. 
The local theory of Banach--Lie groups was continued 
by {\smc Dynkin} ([Dy47/53]) who developed the algebraic theory of the 
BCH series further and by {\smc Laugwitz} ([Lau55/56]) who developed 
a differential geometric perspective, which is quite close in spirit to the 
theory of locally exponential Lie groups described in Section IV below. 
Put in modern terms, he uses the Maurer--Cartan form and integrability conditions 
on (partial) differential equations on Banach space, developed by {\smc Michal}
and {\smc Elconin} 
([MiE37], [MicA48]), to derive the existence of the local group structure 
from the Maurer--Cartan form, which in turn is obtained from the Lie bracket. 
In the finite-dimensional case, this strategy is due to {\smc F.~Schur} ([SchF90a]) 
and quite close to {\smc Lie}'s original approach. 
In [Lau55], {\smc Laugwitz} shows in particular that 
the center and any locally compact subgroup of a Banach--Lie group 
is a Banach--Lie subgroup. 
Formal Lie groups of infinitely many parameters were introduced by {\smc Ritt} a few years 
earlier ([Ri50]). 

The global theory of Banach--Lie groups started in the early 1960s with 
{\smc Maissen}'s paper 
[Mais62] which contains the first basic results on the Lie functor 
on the global level, such as the existence of integral subgroups for closed Lie 
subalgebras and the integrability of Lie algebra homomorphisms for simply connected 
groups. Later {\smc van Est} and {\smc Korthagen} 
studied the integrability problem for Banach--Lie algebras 
and found the first example 
of a non-integrable Banach--Lie algebra ([EK64]). Based on {\smc Kuiper}'s Theorem 
that the unitary group of an infinite-dimensional Hilbert space is contractible 
([Ku65]), simpler examples were given later by {\smc Douady} and {\smc Lazard} ([DL66]). 
Chapters 2 and 3 in {\smc Bourbaki}'s ``Lie groups and Lie algebras'' 
contain in particular the basic local theory of Banach--Lie groups and Lie algebras 
and also some global aspects ([Bou89]). 
Although most of the material in {\smc Hofmann}'s Tulane Lecture Notes 
([Hof68]), approaching the subject from a topological group perspective, 
was never published until recently ([HoMo98]), it was an important 
source of information for many people working on Banach--Lie theory 
(see also [Hof72/75]).

In the early 1970s, {\smc de la Harpe} extended {\smc \'E.~Cartan}'s classification of 
Riemannian symmetric spaces to Hilbert manifolds associated to a certain class of 
Hilbert--Lie algebras, called $L^*$-algebras, 
and studied different classes of operator groups related to Schatten ideals. 
Another context, where a structure-theoretic approach leads quite far is the
theory of bounded symmetric domains in Banach spaces and the related theory 
of (normed) symmetric spaces, developed by {\smc Kaup} and {\smc Upmeier} 
(cf.\ [Ka81/83a/b], [Up85]). 
For a more general approach to Banach symmetric spaces in the sense of {\smc Loos}
([Lo69]), extending the class of all finite-dimensional symmetric spaces, not 
only Riemannian ones,  we refer to [Ne02c] (cf.\ also [La99] for the corresponding 
basic differential geometry). In the context of symplectic geometry, resp., 
Hamiltonian flows, Banach manifolds were introduced by {\smc Marsden} ([Mar67]), and 
{\smc Weinstein} obtained a Darboux Theorem in this context ([Wei69]).
{\smc Schmid}'s monograph [Sch87] provides a nice introduction to 
infinite-dimensional Hamiltonian systems. 
For more recent results on Banach--K\"ahler manifolds and their connections 
to representation theory, we refer to 
([Ne04b], [Bel06]) and for Banach--Poisson manifolds to the recent work of 
{\smc Ratiu}, {\smc Odzijewicz} and {\smc Beltita} ([RO03/04], [BR05a/b]). 

Although {\smc Birkhoff} was already aware of the fact that his theory 
covered groups of operators on Banach spaces, but not groups of diffeomorphisms, 
it took 30 years until infinite-dimensional Lie groups modeled on 
(complete) locally convex spaces occurred for the first time,  
as an attempt to understand the Lie structure 
of the group $\Diff(M)$ of diffeomorphisms 
of a compact manifold $M$, in the work of {\smc Leslie} ([Les67]) and 
{\smc Omori} ([Omo70]). 
This theory was developed further by {\smc Omori} 
in the context of strong ILB--Lie groups (cf.\ [Omo74]). 
A large part of [Omo74] is devoted to the construction of a strong ILB--Lie group 
structure on various types of groups of diffeomorphisms. In the 1980s, this theory 
was refined substantially by imposing and proving additional regularity conditions 
for such groups ([OMYK82/83a], [KYMO85]). 
A different type of Lie group was studied by {\smc Pisanelli} in [Pis76/77/79], 
namely the group $\Gh_n(\C)$ of germs of biholomorphic maps 
of $\C^n$ fixing $0$. This group carries the structure of a Silva--Lie 
group and is one of the first non-Fr\'echet--Lie groups studied 
systematically in a Lie theoretic context. 
In [BCR81], {\smc Boseck}, {\smc Czichowski} and {\smc Rudolph} approach infinite-dimensional 
Lie groups from a topological group perspective. They use the same concept of an 
infinite-dimensional manifold as we do here, but a stronger Lie group concept. This 
leads them to a natural setting for mapping groups of non-compact manifolds 
modeled on spaces of rapidly decreasing functions. 

In his lecture notes [Mil84], {\smc Milnor} undertook the first attempt 
to develop a general theory of Lie groups modeled on complete locally convex spaces, 
which already contained important cornerstones of the theory. 
This paper and the earlier preprint [Mil82]  
had a strong influence on the development of the theory. 
Both contain precise formulations of several problems, 
some of which have been solved in the meantime and some of which are 
still open, as we shall see in more detail below (see also [Gl06b] for the state 
of the art on some of these problems). 

In the middle of the 1980s, groups of smooth maps, and in particular 
groups of smooth loops became popular because of their intimate connection 
with Kac--Moody theory and topology (cf.\ [PS86], [Mick87/89]). 
The interest in direct limits of finite-dimensional Lie groups grew in 
the 1990s (cf.\ [NRW91/93/94/99]). They show 
up naturally in the structure and representation theory of Lie algebras 
(cf.\ [Ne98/01b], [DiPe99], [NRW99], [NS01], [Wol05]). The general Lie theory of 
these groups was put into its definitive form by {\smc Gl\"ockner} in [Gl05].  

There are other, weaker, concepts of Lie groups, 
resp., infinite-dimensional manifolds. 
One is based on the ``convenient setting'' for global analysis developed 
by {\smc Fr\"olicher}, {\smc Kriegl} and {\smc Michor} ([FB66], [Mi84], [FK88] and [KM97]). 
In the context of Fr\'echet 
manifolds, this setting leads to Milnor's concept of a regular Lie group, but 
for more general 
model spaces, it provides a concept of a smooth map which does not imply 
continuity, hence leads to Lie groups which are not topological groups. 
Another approach, due to {\smc Souriau}, 
 is based on the concept of a diffeological space 
([So84/85], [DoIg85], [Los92]; see [HeMa02] for applications to diffeomorphism groups) 
which can be used to study spaces like quotients of $\R$ by 
non-discrete subgroups in a differential geometric context. On the one hand, it has the 
advantage that the category of diffeological spaces is cartesian closed 
and that any quotient of a diffeological space carries a natural diffeology. 
But on the other hand, this incredible freedom makes it harder to distinguish 
``regular'' objects from ``non-regular'' ones. 
Our discussion of smoothness of maps with values in diffeomorphism groups 
of (possibly infinite-dimensional) manifolds 
is inspired by the diffeological approach. We shall see in particular, 
that, to some extent, one can use differential methods to deal with 
groups with no Lie group structure, such as groups of diffeomorphisms 
of non-compact manifolds or groups of linear automorphisms of 
locally convex spaces, and that this provides a natural framework 
for a Lie theory of smooth actions on manifolds and smooth linear 
representations. 

There are other purposes, for which a Lie group structure on an infinite-dimensional 
group $G$ is indispensable. The most crucial one is that without the manifold 
structure, there is simply not enough structure available to pass from the 
infinitesimal level to the global level. For instance, to integrate 
abelian or central extensions of Lie algebras to corresponding group extensions, 
the manifold structure on the group is of crucial importance (cf.\ (FP4)). 
To deal with these extension problems, one is naturally lead to certain classes 
of closed differential $2$-forms on Lie groups, which in turn leads to 
infinite-dimensional symplectic geometry and Hamiltonian group actions. 
Although we do not know which coadjoint orbits of an infinite-dimensional Lie group 
carry manifold structures, for any such orbit, we have a 
natural Hamiltonian action of the group $G$ on itself with respect to a closed 
invariant $2$-form which in general is degenerate; so the reduction of 
free actions of infinite-dimensional Lie groups causes similar difficulties as 
singular reduction does in finite-dimensions; but still all the geometry is visible 
in the non-reduced system. It is our hope that this kind of 
symplectic geometry will ultimately lead to a more systematic 
``orbit method'' for infinite-dimensional Lie groups, in the sense that it 
paves the way to a better understanding of the unitary representations of 
infinite-dimensional Lie groups, based on symplectic geometry and Hamiltonian 
group actions (cf.\ [Ki05] for a recent survey on various aspects of the 
orbit method).

\msk 

{\bf Acknowledgement:} We are most grateful to {\smc H.\ Gl\"ockner} for 
many useful comments on earlier versions of the manuscript 
and for supplying immediate answers to some questions we thought still open. 
We also thank {\smc K.~H.\ Hofmann} for explaining many features of projective limits 
of Lie groups. 

\subheadline{Contents} 

The structure of this paper is as follows. Sections I-VI deal with the general theory 
and Sections VII-X explain how it applies to several classes of Lie groups. 
After explaining the basic concepts and some of the pitfalls of infinite-dimensional 
calculus in Section I, we discuss the basic concepts of Lie group theory in Section II. 
Section III is devoted to the concept of regularity, its relatives, and its applications 
to the fundamental problems. Section IV on locally exponential Lie groups is 
the longest one. This class of groups still displays many features 
of finite-dimensional, resp., Banach--Lie groups. It is also general 
enough to give a good impression of the difficulties arising beyond Banach spaces 
and for smooth Lie groups without any analytic structure. In this section one 
also finds several results concerning Lie group structures on closed subgroups 
and quotients by closed normal subgroups. Beyond the class of locally exponential 
Lie groups, hardly anything is known in this direction. 

In Sections V and VI, which are closely 
related, we discuss the integrability problem for abelian Lie algebra extensions (FP4) 
and the integrability problem for Lie algebras (FP3). In particular, we give a complete 
characterization of the integrable locally exponential Lie algebras and discuss 
integrability results for various other Lie algebras. The remaining sections are relatively 
short: Section VII describes some specifics of direct limits of Lie groups, 
Section VIII explains how the theory applies to linear Lie groups and Section~IX 
presents some results concerning infinite-dimensional Lie groups acting on finite-dimensional 
manifolds. We conclude this paper with a discussion of projective limits, a construction 
that leads beyond the class of Lie groups. It is amazing that projective limits of finite-dimensional 
Lie groups still permit a powerful structure theory ([HoMo06]) and it would be 
of some interest to develop a theory of projective limits of infinite-dimensional Lie groups. 
\msk

\litemindent0.9cm
\litem{I.} Locally convex manifolds \dotfill 8 
\litem{II.} Locally convex Lie groups \dotfill 19 
\litem{III.} Regularity  \dotfill 39
\litem{IV.} Locally exponential Lie groups \dotfill 46
\litem{V.} Extensions of Lie groups \dotfill 66
\litem{VI.} Integrability of locally convex Lie algebras \dotfill 81
\litem{VII.} Direct limits of Lie groups \dotfill 95 
\litem{VIII.} Linear Lie groups \dotfill 98
\litem{IX.} Lie transformation groups \dotfill 101 
\litem{X.} Projective limits of Lie groups  \dotfill 106 
\litemindent0.7cm

\subheadline{Notation} 

We write $\N := \{ 1,2,\ldots \}$ for the natural numbers and $\N_0 := \N \cup \{0\}$. 
Throughout this paper all vector spaces, algebras and Lie algebras are defined over the field 
$\K$, which is $\R$ or~$\C$. If $E$ is a real vector space, we write 
$E_\C := \C \otimes_\R E$ for its complexification, considered as a complex vector space. 
For two topological vector spaces $E,F$, we write 
${\cal L}(E,F)$ for the space of continuous linear operators 
$E \to F$ and put ${\cal L}(E) := {\cal L}(E,E)$. 
For $F = \K$, we write $E' := {\cal L}(E,\K)$ for the {\it topological dual space of $E$}. 

If $G$ is a group, we denote the identity element by  $\1$, and for $g \in G$, we write 
\par\nin $\lambda_g \: G \to G, x \mapsto gx$ for the
{\it left multiplication} by $g$, 
\par\nin $\rho_g \: G \to G, x \mapsto xg$ for the {\it right multiplication} by 
$g$, 
\par\nin $m_G \: G \times G \to G, (x,y) \mapsto xy$ for 
the {\it multiplication map}, and 
\par\nin $\eta_G \: G \to G, x \mapsto x^{-1}$ for the {\it inversion}. 

We always write $G_0$ for the connected component of the identity and, 
if $G$ is  connected, we write $q_G \: \tilde G \to G$ for the universal covering 
group. 

We call a manifold $M$  {\it $1$-connected} if it is connected and simply connected.

\subheadline{Index of notation and concepts} 

\nin LF-space, Silva space \dotfill Definition I.1.2

\nin (Split) submanifold \dotfill Definition I.3.5 

\nin Smoothly paracompact \dotfill Remark I.4.5 

\nin CIA (continuous inverse algebra) \dotfill Definition II.1.3

\nin $x_l(g) = g.x$ (left invariant vector fields) \dotfill Definition II.1.5 

\nin $\L(\phi) = T_\1(\phi)$ (Lie functor) \dotfill Definition II.1.7 

\nin $C^r_X(M,K)$ ($C^r$-maps supported in $X$, $0 \leq r \leq \infty$) 
\dotfill Definition II.2.7 

\nin $C^r_c(M,K)$ (compactly supported $C^r$-maps, $0 \leq r \leq \infty$) 
\dotfill Definition II.2.7 

\nin $\Fl^X_t$ (time $t$ flow of vector field $X$) \dotfill Example II.3.14 

\nin $\kappa_G$ (left Maurer--Cartan form of $G$) \dotfill Section II.4 

\nin $Z(G)$ (center of the group $G$) \dotfill Corollary~II.4.2 

\nin $\evol_G \: C^\infty([0,1],\g) \to G$ (evolution map) \dotfill Definition II.5.2 

\nin $\kappa_\g(x) = \int_0^1 e^{-t\ad x}\, dt$ (Maurer--Cartan form of $\g$) 
\dotfill Remark II.5.8 

\nin $\z(\g)$ (center of the Lie algebra $\g$) \dotfill Proposition II.5.11 

\nin $C^r_*(M,K)$ (base point preserving $C^r$-maps) \dotfill Proposition~II.6.3 

\nin $\L^d(H)$ (differential Lie algebra of subgroup $H$) \dotfill Proposition~II.6.3 

\nin $\g^{\rm op}$ and $G^{\rm op}$ (opposite Lie algebra/group) \dotfill Example II.3.14 

\nin $H^p_{\rm sing}(M,A)$ ($A$-valued singular cohomology) \dotfill Theorem III.1.9 

\nin $\Diff(M,\omega)$, ${\cal V}(M,\omega)$ \dotfill Theorem~III.3.1 

\nin BCH (Baker--Campbell--Hausdorff) Lie algebra/group \dotfill Definitions IV.1.5, 9 

\nin $\gau(P)$ (gauge Lie algebra of principal bundle $P$) \dotfill Theorem~IV.1.12 

\nin $\gau_c(P)$ (compactly supported gauge Lie algebra) \dotfill Theorem~IV.1.12 

\nin Pro-nilpotent Lie algebra/group \dotfill Example~IV.1.13 

\nin $\Gf_n(\K)$ (formal diffeomorphisms in dim.\ $n$) \dotfill Example IV.1.14 

\nin Locally exponential topological group  \dotfill Remark~IV.1.22 

\nin ${\frak L}(G) := \Hom_c(\R,G)$ (for a top. group $G$) \dotfill Definition IV.1.23 

\nin $\L^e(H)$ (exponential Lie algebra of subgroup $H$) \dotfill Lemma IV.3.1 

\nin Locally exponential Lie subgroup  \dotfill Definition~IV.3.2 

\nin Stable Lie subalgebra, resp., ideal ($e^{\ad x}\h = \h$ for $x \in\h$, resp., $\g$)
\dotfill Definition IV.4.1 

\nin Integral subgroup \dotfill Definition IV.4.7 

\nin $H^p_c(\g,\a)$ (continuous Lie algebra cohomology) \dotfill 
Definition~V.2.2 

\nin $\omega^{\rm eq}$ (equivariant $p$-form on $G$ with $\omega^{\rm eq}_\1 = \omega$) 
\dotfill Definition~V.2.3 

\nin $H^p_s(G,A)$ (locally smooth cohomology group) \dotfill 
Definition~V.2.5 

\nin $\per_\Omega \: \pi_k(G) \to E$ (period homo.\ of closed $E$-val.\ $k$-form) 
\dotfill Definition~V.2.12 

\nin $F_\omega \: \pi_1(G) \to H^1_c(\g,\a)$ (flux homo., $\omega \in Z^2_c(\g,\a)$) 
\dotfill Definition~V.2.12 

\nin Integrable/enlargible Lie algebra \dotfill Definition~VI.1.1 

\nin Generalized central extension \dotfill Definition~VI.1.3 

\nin $\gf_n(\R)$ (formal vector fields in dim.\ $n$) \dotfill Example VI.2.8 

\nin $\Gs_n(\R)$, $\gs_n(\R)$ (germs of smooth diffeomorphisms/vector fields) \dotfill 
Theorem~VI.2.9 

\nin $\Gh_n(\C)$, $\gh_n(\C)$ (germs of holomorphic diffeomorphisms/vector fields) \dotfill 
Example~VI.2.10 

\nin Linear Lie group \dotfill Definition~VIII.1 

\nin pro-Lie group/algebra \dotfill Section X.1

\sectionheadline{I. Locally convex manifolds} 

\nin In this section, we briefly explain the natural setup 
for manifolds modeled on locally convex spaces, 
vector fields and differential forms on these manifolds. 
An essential difference to the finite-dimensional, resp., the Banach setting is that 
we use a $C^k$-concept which on Banach spaces is slightly weaker than 
Fr\'echet differentiability, but implies $C^{k-1}$ in the Fr\'echet sense, so that we 
obtain the same class of smooth functions. 
The main point is that, for a non-normable locally convex space 
$E$, the space ${\cal L}(E,F)$ of continuous linear maps to some locally convex space $F$ 
does not carry any vector topology for which the evaluation map is continuous 
([Mais63]). Therefore it is more natural to develop calculus 
independently of any topology on spaces of linear maps and thus to deal instead  
with the differential of a function as a function of two arguments, not as an 
operator-valued function of one variable. 
One readily observes that once the 
Fundamental Theorem of Calculus is available, which is not in general the case beyond locally 
convex spaces, most basic calculus results can simply be reduced to the familiar 
finite-dimensional situation. This is done by restricting to finite-dimensional 
subspaces and composing with linear functionals, which separate the points 
due to the Hahn--Banach Theorems.

The first steps towards a calculus on locally convex spaces have been taken 
by {\smc Michal} (cf.\ [MicA38/40]), whose work was developed further by 
{\smc Bastiani} in [Bas64], so that the calculus we present below is named after 
Michal--Bastiani, and the $C^k$-concept is denoted 
$C^k_{\rm MB}$ accordingly (if there is any need to distinguish it from other $C^k$-concepts).  
{\smc Keller}'s comparative discussion of various notions  
of differentiability on topological vector spaces ([Ke74]) shows that the 
Michal--Bastiani calculus is the most natural one since it does not rely on 
convergence structures or topologies on spaces of linear maps. 
Streamlined discussions of the basic results of 
calculus, as we use it, can be found in [Mi80] and [Ham82]. 
In  [Gl02a], {\smc Gl\"ockner} treats real and complex analytic functions over not necessarily 
complete spaces, which presents some subtle difficulties. Beyond 
Fr\'echet spaces, it is more convenient to work 
with locally convex spaces which are not necessarily complete because quotients of complete 
non-metrizable locally convex spaces need not be complete (cf.\ [K\"o69], \S 31.6). 
Finally we mention that the MB-calculus can even be developed 
for topological vector spaces over general non-discrete 
topological fields (see [BGN04] for details). 

One of the earliest references for smooth manifolds modeled on (complete) 
locally convex spaces is {\smc Eell}'s paper [Ee58], but he uses a different smoothness concept, 
based on the topology of bounded convergence on the space of linear maps 
(cf.\ also [Bas64] and [FB66]). 
Lie groups in the context of MB-calculus 
show up for the first time in {\smc Leslie}'s paper on diffeomorphism groups of compact 
manifolds ([Les67]). 

\subheadline{I.1. Locally convex spaces} 

\Definition I.1.1. A topological vector space $E$ is said to 
be {\it locally convex} if each $0$-neighborhood in $E$ contains a convex one. 
Throughout, topological vector spaces $E$ are assumed to be Hausdorff. 
\qed

It is a standard result in functional analysis that local convexity 
is equivalent to the embeddability of $E$ into a product of normed spaces. 
This holds if and only if the topology can be defined by a family 
$(p_i)_{i \in I}$ of seminorms in the sense that a subset $U$ of $E$ is a 
$0$-neighborhood if and only if it contains a finite intersection of sets of the form 
$$ V(p_i, \eps_i) := \{x \in E \: p_i(x) < \eps_i\}, \quad i \in I, \eps_i > 0. $$

\Definition I.1.2. (a) A locally convex space $E$ is called a {\it Fr\'echet space} 
if there exists a sequence $\{p_n\: n \in \N\}$ of seminorms on $E$, 
such that the topology on $E$ is induced by the metric 
$$ d(x,y) := \sum_{n \in \N} 2^{-n} {p_n(x - y) \over 1 + p_n(x-y)}, $$
and the metric space $(E, d)$ is complete. 
Important examples of Fr\'echet spaces are Banach spaces, which are the ones 
where the topology is defined by a single (semi-)norm. 

(b)Let $E$ be a vector space which can be written as 
$E = \bigcup_{n = 1}^\infty E_n$, where $E_n \subeq E_{n+1}$ are subspaces of 
$E$, endowed with structures of locally convex 
spaces in such a way 
that the inclusion mappings $E_n \to E_{n+1}$ are continuous. 

Then we obtain a 
locally convex vector topology on $E$ by defining a seminorm $p$ on
$E$ to be continuous if and only if its restriction to all the
subspaces $E_n$ is continuous. 
We call $E$ the {\it inductive limit} (or {\it direct limit}) of the spaces $(E_n)_{n
\in \N}$. 

If all maps $E_n \into E_{n+1}$ are embeddings, 
we speak of a {\it strict inductive limit}. 
If, in addition, the spaces $E_n$ are Fr\'echet spaces, then 
each $E_n$ is closed in $E_{n+1}$ and $E$ is called an {\it LF space}. 
If the spaces $E_n$ are Banach spaces and 
the inclusion maps $E_n \to E_{n+1}$ are compact, then $E$ is called a {\it Silva space}. 
\qed 

\Examples I.1.3. To given an impression of the different types of locally convex spaces occurring 
below, we take a brief look at function spaces on the real line. 

(a) For $r \in \N_0$ and $a < b$, the spaces 
$C^r([a,b],\R)$ of $r$-times continuously differentiable functions on $[a,b]$ form a Banach space with 
respect to the norms 
$$ \|f\|_r := \sum_{k = 0}^r \|f^{(k)}\|_\infty $$ 
and $C^\infty([a,b],\R)$ is a Fr\'echet space with respect to the topology defined by 
the sequence $(\|\cdot\|_r)_{r \in \N_0}$ of norms. 

(b) For each fixed $r \in \N_0$, the space $C^r(\R,\R)$ is a Fr\'echet space with respect to the 
sequence of seminorms $p_n(f) := \|f\res_{[-n,n]}\|_r$. It is the projective limit of the 
Banach spaces $C^r([-n,n],\R)$. On $C^\infty(\R,\R)$ we also obtain a Fr\'echet space structure 
defined by the countable family of seminorms $p_{n,r}(f) := \|f\res_{[-n,n]}\|_r$, $n,r \in \N$. 

(c) For each $r \in \N_0 \cup \{\infty\}$, the space $C^r_c(\R,\R)$ of compactly supported 
$C^r$-functions on $\R$ is the union of the subspaces $C^r_{[-n,n]}(\R,\R)$ of all those 
functions supported by the interval $[-n,n]$. 
As a closed subspace of the Fr\'echet space $C^r([-n,n],\R)$, 
$C^r_{[-n,n]}(\R,\R)$ inherits a Fr\'echet space structure, so that we obtain on 
$C^r_c(\R,\R)$ the structure of an LF space. 

(d) For each $r > 0$, the space 
of all sequences $(a_n)_{n \in \N_0}$ for which 
$\sum_{n = 0}^\infty |a_n| r^n$ converges can be identified with the space 
$E_r$ of all functions $f \: [-r,r] \to \R$ which can be represented by a power series, uniformly convergent on $[-r,r]$. This is a Banach space with respect to the norm 
$\|f\|_r := \sum_{n = 0}^\infty {|f^{(n)}(0)| \over n!} r^n$. 
The direct limit space $E := \bigcup_{r > 0}^\infty E_r$ is the space of germs of analytic function 
in $0$. Since the inclusion maps $E_{1\over n} \to E_{1\over n+1}$ are compact operators, 
$E$ carries a natural Silva space structure. Note that the subspaces $E_r$, $r > 0$, 
are dense and not closed, so that $E$ is {\sl not} an LF space. 
\qed

The natural completeness requirement for calculus on locally convex spaces 
is the following: 

\Definition I.1.4. A locally convex space~$E$ is said to be 
{\it Mackey complete} if for each smooth curve $\xi \: [0,1] \to E$ there
exists a smooth curve $\eta \: [0,1] \to E$ with $\eta' = \xi$. 
\qed

For each continuous linear functional $\lambda \: E \to \R$ on a locally convex space $E$ 
and each continuous curve $\xi \: [0,1] \to E$, we have a continuous real-valued 
function $\lambda \circ \xi \: [0,1] \to \R$ which we may integrate 
to obtain a linear functional 
$$ I_\xi \: E' \to \R, \quad \lambda \mapsto \int_0^1 \lambda(\xi(t))\, dt, $$
called the {\it weak integral of $\xi$}. On the other hand, 
we have a natural embedding 
$$\eta_E \: E \to (E')^*, \quad \eta_E(x)(\lambda) := \lambda(x) $$
which is injective, because $E'$ separates the points of $E$ by the Hahn--Banach Theorem. 
Therefore Mackey completeness means that 
for each smooth curve $\xi$ the weak integral $I_\xi$ is represented by an 
element of $E$, i.e., contained in $\eta_E(E)$. If this is the case, 
 we simply write $\int_0^1 \xi(t)\, dt$ for the representing element of $E$. 
The curve $\eta(s) := \int_0^s \xi(t)\, dt$ then is differentiable and 
satisfies $\eta' = \xi$. 

For a more detailed discussion of Mackey completeness and equivalent conditions, we 
refer to [KM97, Th.~2.14], where it is shown in particular that 
integrals exist for Lipschitz curves and in particular for each $\eta \in C^1([0,1],E)$. 

\subheadline{I.2. Calculus on locally convex spaces} 

The following notion of $C^k$-maps is also known as $C^k_{MB}$ ($C^k$ in the 
Michal--Bastiani sense) ([MicA38/40], [Bas64]) or Keller's $C_c^k$-maps ([Ke74]). 
Its main advantage is that it does not refer to any topology on spaces 
of linear maps or any quasi-topology (cf.\ [Bas64]). 

\Definition I.2.1.  (a) Let $E$ and $F$ be locally convex spaces, $U
\subeq E$ open and $f \: U \to F$ a map. Then the {\it derivative
  of $f$ at $x$ in the direction $h$} is defined as 
$$ df(x)(h) := (D_h f)(x) := \derat0 f(x + t h) 
= \lim_{t \to 0} {1\over t}(f(x+th) -f(x)) $$
whenever it exists. The function $f$ is called {\it differentiable at
  $x$} if $df(x)(h)$ exists for all $h \in E$. It is called {\it
  continuously differentiable}, if it is differentiable at all
points of $U$ and 
$$ df \: U \times E \to F, \quad (x,h) \mapsto df(x)(h) $$
is a continuous map. It is called a {\it $C^k$-map}, $k \in \N \cup \{\infty\}$, 
if it is continuous, the iterated directional derivatives 
$$ d^{j}f(x)(h_1,\ldots, h_j)
:= (D_{h_j} \cdots D_{h_1}f)(x) $$
exist for all integers $j \leq k$, $x \in U$ and $h_1,\ldots, h_j \in E$, 
and all maps $d^j f \: U \times E^j \to F$ are continuous. 
As usual, $C^\infty$-maps are called {\it smooth}. 

  (b) If $E$ and $F$ are complex vector spaces, then a map $f$ is 
called {\it complex analytic} if it is continuous and for each 
$x \in U$ there exists a $0$-neighborhood $V$ with $x + V \subeq U$ and 
continuous homogeneous polynomials $\beta_k \: E \to F$ of degree $k$ 
such that for each $h \in V$ we have 
$$ f(x+h) = \sum_{k = 0}^\infty \beta_k(h), $$
as a pointwise limit ([BoSi71]). 

If $E$ and $F$ are real locally convex spaces, 
then we call $f$ {\it real analytic}, resp., $C^\omega$, 
if for each point $x \in U$ there exists an open neighborhood 
$V \subeq E_\C$ and a holomorphic map $f_\C \: V \to F_\C$ with 
$f_\C\res_{U \cap V} = f\res_{U \cap V}$ (cf.\ [Mil84]). 
The advantage of this definition, which differs from the one in [BoSi71], 
is that it works nicely for non-complete spaces, any analytic map is smooth, 
and the corresponding chain rule holds without any condition 
on the underlying spaces, which is the key to the definition of 
analytic manifolds (see [Gl02a] for details). 

The map $f$ is called {\it holomorphic} if it is $C^1$ and for each $x \in U$ the 
map $df(x) \: E \to F$ is complex linear (cf.\ [Mil84, p.\ 1027]). 
If $F$ is sequentially complete, then $f$ is holomorphic if and only if 
it is complex analytic (cf.\ [Gl02a], [BoSi71, Ths.~3.1, 6.4], 
[Mil82, Lemma~2.11]). 
\qed

\Remark I.2.2. If $E$ and $F$ are Banach spaces, then the 
Michal--Bastiani $C^1_{MB}$-concept from above is weaker than 
continuous Fr\'echet differentiability, which requires
that the map $x \mapsto df(x)$ is continuous with respect to the
operator norm (cf.\ [Mil82, Ex.~6.8]). 
Nevertheless, one can show that $C^{k+1}_{MB}$ implies $C^k$ in the sense of Fr\'echet differentiability, 
which in turn implies $C^k_{MB}$.  
Therefore the different $C^k$-concepts lead to the same class of smooth functions  
(cf.\ [Mil82, Lemma~2.10], [Ne01a, I.6 and I.7]). 
\qed

After clarifying the $C^k$-concept, we recall the precise statements of the most fundamental facts 
from calculus on locally convex spaces. 

\Proposition I.2.3. Let $E$ and $F$ be locally convex spaces, $U \subeq E$ an open 
subset, and $f \: U \to F$ a continuously differentiable function. 
\litem{(i)} For any $x \in U$, the map $df(x) \: E \to F$ is real linear 
and continuous. 
\litem{(ii)} {\rm(Fundamental Theorem of Calculus)} If $x + [0,1]h \subeq U$, then 
$$f(x + h) = f(x) + \int_0^1 df(x + t h)(h)\, dt.$$
In particular, $f$ is locally constant if and only if $df = 0$. 
\litem{(iii)} $f$ is continuous. 
\litem{(iv)} If $f$ is $C^n$, $n \geq 2$, then the functions 
 $d^{n}f(x)$, $x \in U$, are 
 symmetric $n$-linear maps. 
\litem{(v)} If $x + [0,1]h \subeq U$ and $f$ is $C^n$, then we have the Taylor Formula 
$$ \eqalign{ 
f(x + h) 
= f(x) + df(x)(h) + \ldots + &{1\over (n-1)!} d^{n-1}f(x)(h,\ldots,h) \cr
&\ \ \ \ + {1\over (n-1)!} \int_0^1 (1-t)^{n-1} d^{n}f(x+th)(h,\ldots,h)\, 
dt. \cr} $$
\litem{(vi)} {\rm(Chain Rule)} If, in addition, $Z$ is a locally convex space, 
$V \subeq F$ is open, and $f_1 \: U \to V$, $f_2 \: V \to Z$ are
$C^1$,  then $f_2 \circ f_1 \: U \to Z$ is $C^1$ with 
$$ d(f_2 \circ f_1)(x) = df_2\big(f_1(x)\big) \circ df_1(x)\quad \hbox{ for } \quad 
x \in U. $$
If $f_1$ and $f_2$ are $C^k$, $k \in \N \cup \{\infty\}$, the Chain Rule implies that 
$f_2 \circ f_1$ is also $C^k$. 
\qed

\Remark I.2.4. A continuous 
$k$-linear map $m \: E_1 \times \ldots \times E_k \to F$
is continuously differentiable with 
$$ dm(x)(h_1, \ldots, h_k) 
= m(h_1, x_2, \ldots, x_k) + 
\cdots + m(x_1, \ldots, x_{k-1}, h_k). $$
Inductively, one obtains that $m$ is smooth with 
$d^{k+1}m = 0$. 
\qed

\Example I.2.5. The following example shows that local convexity is crucial 
for the validity of the Fundamental Theorem of Calculus. 

Let $E$ denote the space of measurable functions $f \: [0,1]\to \R$ for which 
$$ |f| := \int_0^1 |f(x)|^{1\over 2}\, dx $$
is finite and identify functions that coincide on a 
set whose complement has measure zero. 
Then $d(f,g) := |f-g|$ defines a metric on $E$. 
We thus obtain a metric topological vector space $(E,d)$. 

For a subset $S \subeq [0,1]$, let $\chi_S$ denote its characteristic function.  
Consider the curve 
$$ \gamma \: [0,1] \to E, \quad \gamma(t) := \chi_{[0,t]}. $$
Then 
$|h^{-1}\big(\gamma(t+h) - \gamma(t)\big)| =   |h|^{-{1\over 2}} |h| \to 0$
for each $t \in [0,1]$ as $h \to 0$. Hence $\gamma$ is continuously differentiable 
with $\gamma' = 0$. 
Since $\gamma$ is not constant, the Fundamental Theorem of Calculus does not hold in 
$E$. 

The defect in this example 
is caused by the non-local convexity of $E$. In fact, one can even show that 
all continuous linear functionals on $E$ vanish. 
\qed

The preceding phenomenon could also be excluded by requiring that the topological 
vector spaces under consideration have the property that the continuous linear 
functionals separate the points, which is automatic for locally convex spaces. 
Another reason for working with locally convex spaces is that local convexity 
is also crucial for approximation arguments, more specifically to approximate 
continuous maps by smooth ones in the same homotopy class (cf.\ [Ne04c], 
[Wo05a]). Local convexity it also crucial for the continuous parameter-dependence of 
integrals which in turn goes into the proof of the Chain Rule. 

One frequently encounters situations where it is convenient to describe 
multilinear maps $m \: E_1 \times \cdots \times E_k \to F$ 
as continuous linear maps on the tensor product space 
$E_1 \otimes \cdots \otimes E_k$, endowed with a suitable topology. 
For locally convex spaces, there is a natural such topology, 
the {\it projective tensor topology}, and it has 
the nice property that projective tensor products are associative. 
That this is no longer true for more general topological vector spaces 
is one more reason to work in the locally convex setting (cf.\ [Gl04a]).

\Remark I.2.6. (Inverse Function Theorems) In the context 
of Banach spaces, one has an Inverse Function Theorem 
and also an Implicit Function Theorem (cf.\ [La99]). Such 
results cannot be expected in general 
for Fr\'echet spaces. One of the simplest examples 
demonstrating this fact arises from the algebra 
$A := C(\R,\R)$ of all continuous functions on $\R$, endowed with 
the topology of uniform convergence on compact subsets, turning 
$A$ into a Fr\'echet space on which the algebra multiplication is continuous. 
We have a smooth exponential map 
$$ \exp_A \: A \to A, \quad f \mapsto e^f $$
with $T_0(\exp_A) = \id_A$. Since the range of $\exp_A$ 
lies in the unit group $A^\times = C(\R, \R^\times)$, which 
apparently is not a neighborhood of the constant function $\1$, 
the Inverse Function Theorem fails in this case (cf.\ [Ee66, p.761]). 

In Example~II.5.9 below, we shall even encounter examples of exponential functions of 
Lie groups which, 
in spite of $T_\1(\exp_G) = \id_{\L(G)}$, are not a local diffeomorphism in $0$. 
In view of these examples, the usual Inverse Function Theorem cannot be generalized 
in any straightforward manner to arbitrary Fr\'echet spaces. 

Nevertheless, {\smc Gl\"ockner} ([Gl03a]) 
obtained a quite useful Implicit Function Theorem for maps of the type 
$f \: E \times F \to F$, where $F$ is a Banach space and $E$ is locally convex. 
These results have many interesting applications, even in the case where 
$F$ is finite-dimensional. Similar results have been achieved by {\smc Hiltunen} in [Hi99], 
but he uses a different notion of smoothness. 

A complementary Inverse Function Theorem is due to {\smc Nash} and {\smc Moser} (cf.~[Mo61] 
and [Ham82] for a nice exposition). 
This is a variant that can be applied to Fr\'echet spaces 
endowed with an additional structure, 
called a grading, and to smooth maps which are ``tame'' in the sense that they are 
compatible with the grading. 

Another variant based on compatibility with a projective limit of Banach spaces is 
the ILB-Implicit Function Theorem to be found in {\smc Omori}'s book ([Omo97]). 
\qed

\Remark I.2.7. (Non-complemented subspaces) Another 
serious pathology occurring already for Banach spaces is that 
a closed subspace $F$ of a locally convex space $E$ need not have a closed complement. 
A simple example is the subspace $F := c_0(\N,\R)$ of the Banach space 
$E := \ell^\infty(\N,\R)$ (cf.\ [Wer95, Satz IV.6.5] for an elementary proof).

This implies that if $q \: E \to E/F$ is a quotient map of locally convex spaces, 
 there need not be any continuous linear map 
$\sigma \: E/F \to E$ with $q \circ \sigma = \id_{E/F}$. If such a map $\sigma$ 
exists, then 
$$ F \times E/F \to E, \quad (x,y) \mapsto x + \sigma(y) $$
is a linear isomorphism of topological vector spaces, which implies that 
$\sigma(E/F)$ is a closed complement of $F$ in $E$. We then call the quotient map 
$q$, resp., the subspace $F$,  {\it topologically split}. If $E$ is a Fr\'echet space, 
then the Open Mapping Theorem implies that the 
existence of a closed complement for $F$ is equivalent to the existence of a splitting map 
$\sigma$. 

For Fr\'echet spaces, it is quite easy to find natural examples of non-splitting 
quotient maps: Let 
$E := C^\infty([0,1],\R)$ be the Fr\'echet space of smooth functions on 
the unit interval and 
$$ q \: E \to \R^\N, \quad q(f) = (f^{(n)}(0))_{n \in \N}. $$
In view of E.~Borel's Theorem, this map is surjective, hence a quotient map by 
the Open Mapping Theorem. Since every $0$-neighborhood in $\R^\N$ contains a 
non-trivial subspace, there is no continuous norm on $\R^\N$, hence there 
is no continuous linear cross section $\sigma \: \R^\N\to E$ for $q$ 
because the topology on $E$ is defined by a sequence of norms. 

If a continuous linear cross section $\sigma$ 
does not exist, then $q$ has no smooth local sections  
either, because for any such section $\sigma \: U \to E$, $U$ open in $E/F$, 
the differential of $\sigma$ in any point would be a continuous linear section of $q$. 
If $E$ is Fr\'echet, then $q$ has a continuous global 
section by Michael's Selection Theorem ([MicE59], [Bou87]), and the preceding argument shows 
that no such section is continuously differentiable.

For more detailed information on splitting conditions for extensions of 
Fr\'echet spaces, we refer to [Pala71] and [Vo87]. 
\qed

\subheadline{I.3. Smooth manifolds}  
  
Since the Chain Rule is valid for smooth maps between open subsets of 
locally convex spaces, we can define smooth manifolds as in the finite-dimensional case 
(see [Ee58] for one of the first occurrences of manifolds modeled on complete locally convex spaces). 

\Definition I.3.1. Let 
$M$ be a Hausdorff space and $E$ a locally convex space. 
An {\it $E$-chart} of $M$ is a pair $(\phi,U)$ of an open subset $U \subeq M$ and a homeomorphism 
$\phi \: U \to \phi(U) \subeq E$ onto an open subset $\phi(U)$ of $E$. 
For $k \in \N_0 \cup \{ \infty,\omega\}$, two $E$-charts 
$(\phi, U)$ and $(\psi, V)$ are said to be {\it $C^k$-compatible} if the maps 
$$ \psi \circ \phi^{-1} \res_{\phi(U \cap V)} \: \phi(U \cap V) \to \psi(U \cap V) $$
and $\phi \circ \psi^{-1}$ are $C^k$, where $k = \omega$ stands for analyticity. 
Since compositions of $C^k$-maps are $C^k$-maps, 
$C^k$-compatibility of $E$-charts is an equivalence relation. 
An {\it $E$-atlas of class $C^k$ of $M$} is a set ${\cal A} := 
\{ (\phi_i, U_i) \: i \in I\}$ of pairwise 
compatible $E$-charts of $M$ with $\bigcup_i U_i = M$. 
A {\it smooth/analytic $E$-structure on $M$} is a maximal 
$E$-atlas of class $C^\infty$/$C^\omega$, and a 
{\it smooth/analytic $E$-manifold} is a pair $(M,{\cal A})$, where ${\cal A}$ is a 
maximal smooth/analytic $E$-atlas 
on $M$. 

We call a manifold modeled on a locally convex, resp., Fr\'echet, resp., 
Banach space a {\it locally convex}, resp., {\it Fr\'echet}, resp., {\it Banach 
manifold}.
\qed

We do not make any further assumptions on the topology of smooth 
locally convex manifolds, such as regularity (as in [Mil84]) or paracompactness. 
But we impose the Hausdorff condition, an assumption not made in some textbooks 
(cf.\ [La99], [Pa57]). We refer to Example 6.9 in [Mil82] for a non-regular 
manifold. 

\Remark I.3.2. If $M_1, \ldots, M_n$ are smooth manifolds modeled on the spaces $E_i$, 
$i =1,\ldots, n$, then the product set $M := M_1 \times \ldots \times M_n$ 
carries a natural manifold structure with model space 
$E = \prod_{i =1}^n E_i$. 
\qed

Smooth maps between smooth manifolds are defined as usual. 

\Definition I.3.3. For $p \in M$, tangent vectors 
$v \in T_p(M)$ are defined as equivalence classes of 
smooth curves $\gamma \: ]{-\eps},\eps[ \to M$ with $\eps > 0$ and 
$\gamma(0) = p$, where the equivalence relation 
is given by $\gamma_1 \sim \gamma_2$ if 
$(\phi \circ \gamma_1)'(0) =  (\phi \circ \gamma_2)'(0)$ holds for 
a chart $(\phi,U)$ with $p \in U$. 
Then $T_p(M)$ carries a natural vector space structure such that for any 
$E$-chart $(\phi,U)$, the map 
$T_p(M) \to E, [\gamma] \mapsto (\phi \circ \gamma)'(0)$ is a linear isomorphism. 
We write 
$T(M) := \bigcup_{p \in M} T_p(M)$
for the {\it tangent bundle} of $M$. The map $\pi_{TM} \: T(M) \to M$ mapping 
elements of $T_p(M)$ to $p$ is called the {\it bundle projection}.

If $f \: M \to N$ is a smooth map between smooth manifolds, 
we obtain for each $p \in M$ a linear tangent map 
$$ T_p(f) \: T_p(M) \to T_{f(p)}(N), \quad [\gamma] \mapsto [f\circ \gamma], $$
and these maps combine to the tangent map $T(f) \: T(M) \to T(N)$. 
On the tangent bundle $T(M)$ 
we obtain for each $E$-chart $(\phi, U)$ of $M$ 
an $E \times E$-chart by 
$$ T(\phi) \: T(U) := \bigcup_{p \in U} T_p(M) \to T(\phi(U)) \cong \phi(U) \times E. $$
Endowing $T(M)$ with the topology for which $O \subeq T(M)$ is open if and only 
if for each $E$-chart $(\phi, U)$ of $M$ 
the set $T(\phi)(O \cap T(U))$ is open in 
$\phi(U) \times E$, we obtain on $T(M)$ the structure of an $E \times E$-manifold 
defined by the charts $(T(\phi), T(U))$, obtained from $E$-charts $(\phi,U)$ of $M$. 
This leads to an endofunctor $T$ on the 
category of smooth manifolds, preserving finite products (cf.\ Remark~I.3.2). 

If $f \: M \to V$ is a smooth map into a locally convex space, then 
$T(f) \: T(M) \to T(V) \cong V \times V$ is smooth, and can be written as 
$T(f) = (f,df)$, where $df \: T(M) \to V$ is called the {\it differential of $f$}. 
\qed

As a consequence of Proposition~I.2.3(ii), we have: 

\Proposition I.3.4. A smooth map $f \: M \to N$ is locally constant if and only if $T(f) =0$. 
\qed

\Definition I.3.5. Let $M$ be a smooth manifold modeled on the space $E$,  
and $N \subeq M$ a subset. 

(a) $N$ is called a {\it submanifold} of $M$ if there exists a closed subspace 
$F \subeq E$ and for each $n \in N$ there exists an $E$-chart $(\phi,U)$ of 
$M$ with $n \in U$ and $\phi(U \cap N) = \phi(U) \cap F$. 

(b) $N$ is called a {\it split submanifold} of $M$ if, in addition, 
there exists a subspace $G \subeq E$ for which the addition map 
$F \times G \to E, (f,g) \mapsto f + g$ is a topological isomorphism. 
\qed

\Definition I.3.6. A {\it (smooth) vector field} $X$ on $M$ is a smooth section 
of the tangent bundle $\pi_{TM} \: TM \to M$, i.e. a smooth map $X \: M \to TM$ 
with $\pi_{TM} \circ X = \id_M$. 
We write ${\cal V}(M)$ for the space of all vector fields
on $M$. If $f \in C^\infty(M,V)$ is a smooth function on $M$ with values 
in some locally convex space $V$ and 
$X \in {\cal V}(M)$, then we obtain a smooth function on $M$ via 
$$ X.f := df \circ X \: M  \to V. $$

For $X, Y \in {\cal V}(M)$, there exists a unique vector
field $ [X,Y] \in {\cal V}(M)$ determined by the
property that on each open subset $U \subeq M$ we have 
$$ [X,Y].f = X.(Y.f) - Y.(X.f) \leqno(1.3.1) $$
for all $f \in C^\infty(U,\R)$. We thus obtain on 
${\cal V}(M)$ the structure of a Lie algebra. 
\qed

\Remark I.3.7. If $M = U$ is an open subset of the locally convex space 
$E$, then $TU = U \times E$ with the bundle projection 
$\pi_{TU} \: U \times E \to U, (x,v) \mapsto x$. Each smooth 
vector field is of the form 
$X(x) = (x, \tilde X(x))$ for some smooth function $\tilde X \: U \to E$, 
and we may thus identify ${\cal V}(U)$ with the space $C^\infty(U,E)$. 
Then the Lie bracket satisfies 
$$ [X,Y]\,\tilde{}(p) = d\tilde Y(p)\tilde X(p) - d\tilde X(p)\tilde Y(p) \quad 
\hbox{ for each } \quad p \in U. 
\qeddis 

\Definition I.3.8. Let $M$ be a smooth $E$-manifold and $F$ a locally convex space. 
A {\it smooth vector bundle of type $F$} over $M$ is a pair $(\pi,\F,F)$, consisting of a 
smooth manifold $\F$, a smooth map 
$\pi \: \F \to M$ and a locally convex space $F$, with the following properties: 
\litem{(a)} For each $m \in M$, the fiber $\F_m := \pi^{-1}(m)$ carries a locally convex 
vector space structure isomorphic to $F$.
\litem{(b)} Each point $m \in M$ has an open neighborhood $U$ for which there exists a diffeomorphism 
$$ \phi_U \: \pi^{-1}(U) \to U \times F $$
with $\phi_U = (\pi\res_U, g_U)$, where $g_U \: \pi^{-1}(U) \to F$ is linear on each 
fiber $\F_m$, $m \in U$. 

We then call $U$ a {\it trivializing subset} of $M$ and $\phi_U$ a {\it bundle chart}. 
If $\phi_U$ and $\phi_V$ are two bundle charts  and $U \cap V \not=\eset$, then we obtain a 
diffeomorphism 
$$ \phi_U \circ \phi_V^{-1} \: (U \cap V) \times F \to (U \cap V) \times F $$
of the form $(x,v) \mapsto (x, g_{VU}(x)v)$. This leads to a map 
$$ g_{UV} \: U \cap V \to \GL(F) $$
for which it does not make sense to speak about smoothness because 
$\GL(F)$ is not a  Lie group if $F$ is not a Banach space. 
This is a major difference between the Banach 
and the locally convex context. 
Nevertheless, $g_{UV}$ is smooth in the sense that the map 
$$ \hat g_{UV} \: (U \cap V) \times F \to F \times F, \quad 
(x,v) \mapsto 
(g_{UV}(x)v, g_{UV}(x)^{-1}v) = (g_{UV}(x)v, g_{VU}(x)v) $$
is smooth (cf.\ Definition~II.3.1 below). 
\qed

Obviously, the tangent bundle $T(M)$ of a smooth (locally convex) manifold is 
an example of a vector bundle, but the cotangent bundle is more problematic: 

\Remark I.3.9. We define for each $E$-manifold $M$ 
the cotangent bundle by $T^*(M) := \bigcup_{m \in M} T_m(M)'$ and observe that, 
as a set, it carries a natural structure of a vector bundle over $M$, but to endow it with a smooth 
manifold structure we need a locally convex topology on the dual space $E'$ such that 
for each local diffeomorphism $f \: U \to E$, $U$ open in $E$, the map 
$U \times E' \to E', (x,\lambda) \mapsto \lambda \circ df(x)$ is smooth. 
If $E$ is a Banach space, then the norm topology on $E'$ has this property, 
but in general this property fails for non-Banach manifolds. 

Indeed, let $E$ be a locally convex space which is not normable 
and pick a non-zero $\alpha_0 \in E'$. We consider the smooth map 
$$ f \: E \to E, \quad x \mapsto x + \alpha_0(x)x = (1+ \alpha_0(x))x. $$
Then $df(x)v = (1 + \alpha_0(x))v + \alpha_0(v)x$ 
implies that $df(x) = (1+\alpha_0(x))\1 + \alpha_0 \otimes x$, which 
is invertible for $\alpha_0(x)\not\in \{-1, -{1\over 2}\}$. 
If $\phi \: ]{-{1\over 4}}, \infty[ \to \R$ is the inverse function 
of $\psi(x) = x + x^2$ on $]{-{1\over 2}}, \infty[$, then an easy calculation gives  
on $\{ y \in E \: \alpha_0(y) > -{1\over 4}\}$ the inverse function 
$f^{-1}(y) = (1+\phi(\alpha_0(y)))^{-1} \cdot y$. We conclude 
that $f$ is a local diffeomorphism on some $0$-neighborhood of~$E$. 

On the other hand, the map $U \times E' \to E', (x,\lambda) \mapsto \lambda \circ df(x)$ 
satisfies 
$$ \lambda \circ df(x) = (1 + \alpha_0(x)) \lambda + \lambda(x)\alpha_0. $$
Since the evaluation map $E' \times E \to \R$ is discontinuous in $0$ 
for any vector topology on $E'$ ([Mais63]), 
$f$ does not induce a continuous map 
on $T^*(E) \cong E \times E'$ for any such topology. Hence there is no natural 
smooth vector bundle structure on $T^*(M)$ if $E$ is not normable. 
\qed

In view of the difficulties caused by the cotangent bundle, 
we shall introduce differential forms directly, not as sections of a vector bundle.

\subheadline{I.4. Differential forms} 

Differential forms play a significant role throughout infinite-dimensional 
Lie theory. 
In the present section, we describe a natural approach to differential 
forms on manifolds modeled on locally convex spaces. A major difference 
to the finite-dimensional case is that in local charts there is no 
natural coordinate description of differential forms in terms of basic forms, 
that differential forms cannot be defined as the smooth sections of a natural 
vector bundle (Remark~I.3.9), and that, even for Banach manifolds,  
smooth partitions of unity need not exist, 
so that one has to be careful with localization arguments. 

In [KM97], one finds a discussion of various types of differential forms,
containing in particular those introduced below, which are also used by 
{\smc Beggs} in [Beg87].

\Definition I.4.1. (a) If $M$ 
is a differentiable manifold and $E$ a locally convex space, then an
{\it $E$-valued $p$-form} $\omega$ on $M$ is a function 
$\omega$ which associates to each $x \in M$ a $p$-linear 
alternating map $\omega_x \: T_x(M)^p \to E$ 
such that in local coordinates the
map 
$(x,v_1, \ldots, v_p) \mapsto \omega_x(v_1, \ldots, v_p)$
is smooth. We write $\Omega^p(M,E)$ for the space of $E$-valued $p$-forms
on $M$ and identify $\Omega^0(M,E)$ with the space $C^\infty(M,E)$ 
of smooth $E$-valued functions on $M$. 

(b) Let $E_1, E_2, E_3$ be locally convex spaces and 
$\beta \: E_1 \times E_2 \to E_3$ be a continuous bilinear map. 
Then the wedge product
$$ \Omega^p(M,E_1) \times \Omega^q(M,E_2) 
\to \Omega^{p+q}(M,E_3), \quad (\omega, \eta) \mapsto \omega \wedge \eta $$
is defined by $(\omega \wedge \eta)_x := \omega_x \wedge \eta_x$, where  
$$ \eqalign{ 
&(\omega_x \wedge \eta_x)(v_1,\ldots, v_{p+q}) 
:= \frac{1}{p!q!} \sum_{\sigma \in S_{p+q}} 
\sgn(\sigma) 
\beta\big(\omega_x(v_{\sigma(1)}, \ldots, v_{\sigma(p)}), 
\eta_x(v_{\sigma(p+1)}, \ldots, v_{\sigma(p+q)})\big). \cr} $$

Important special cases, where such wedge products are used, are: 
\litem{(1)} $\beta \: \R \times E \to E$ is the scalar multiplication of $E$. 
\litem{(2)} $\beta \: A \times A \to A$ is the multiplication of an associative 
algebra. 
\litem{(3)} $\beta \: \g \times \g \to \g$ is the Lie bracket of a Lie algebra. 
In this case, we also write $[\omega,\eta] := \omega \wedge \eta$. 
\qed

The definition of the exterior differential 
$d \: \Omega^p(M,E) \to \Omega^{p+1}(M,E)$
is a bit more subtle than in finite dimensions, where one usually uses 
local coordinates to define it in charts. 
Here the exterior differential is determined uniquely by the property that 
for each open subset $U \subeq M$ we have for $X_0, \ldots, X_p \in {\cal V}(U)$ 
in the space $C^\infty(U,E)$ the identity 
$$ \leqalignno{ 
(d\omega)(X_0, \ldots, X_p) 
&:= \sum_{i = 0}^p (-1)^i X_i.\omega(X_0, \ldots, \hat X_i, \ldots,X_p) \cr
&\ \ + \sum_{i < j} (-1)^{i+j} \omega([X_i, X_j], X_0, \ldots, \hat X_i,
\ldots, \hat X_j, \ldots, X_p).  \cr} $$
The main point is to show that in a point $x \in U$ 
the right hand side only depends 
on the values of the vector fields $X_i$ in~$x$. 
The exterior differential has the usual properties, such as 
$d^2 = 0$ and the compatibility with pullbacks:  $\phi^*(d\omega) = d(\phi^*\omega)$. 

Extending $d$ to a linear map on the space $\Omega(M,E) := \bigoplus_{p \in
\N_0} \Omega^p(M,E)$ of all $E$-valued differential forms on $M$, 
the relation $d^2 = 0$ implies that the space 
$$ Z^p_{\rm dR}(M,E) := \ker(d\res_{\Omega^p(M,E)}) $$
of {\it closed $p$-forms} contains the space 
$B^p_{\rm dR}(M,E) := d(\Omega^{p-1}(M,E))$
of {\it exact $p$-forms}, so that we may define the {\it $E$-valued de Rham cohomology space} 
by 
$$ H^p_{\rm dR}(M,E) := Z^p_{\rm dR}(M,E) / B^p_{\rm dR}(M,E). $$

For finite-dimensional manifolds, one usually defines the Lie derivative 
of a differential form in the direction of a vector field $X$ by using its local 
flow $t \mapsto \Fl^X_t$: 
$$ {\cal L}_X \omega := \derat0 (\Fl^X_{-t})^*\omega. $$
Since vector fields on infinite-dimensional manifold need not have a local flow 
(cf.\ Example~II.3.11 below), we introduce the Lie derivative more directly.

\Definition I.4.2. (a) For any smooth manifold $M$ and each locally convex space, 
we have a natural representation of 
the Lie algebra ${\cal V}(M)$ on the space $\Omega^p(M,E)$ of 
$E$-valued $p$-forms on $M$, given by the 
{\it Lie derivative}, which for $Y \in {\cal V}(M)$ is uniquely determined by 
$$({\cal L}_Y\omega)(X_1,\ldots, X_p) 
= Y.\omega(X_1, \ldots, X_p) - \sum_{j = 1}^p \omega(X_1, \ldots,
[Y, X_j], \ldots, X_p) $$
for $X_i \in {\cal V}(U)$, $U \subeq M$ open. 
Again one has to verify that the value of the right hand side in some $x \in M$ 
only depends on the values of the vector fields $X_i$ in $x$. 

(b) We further obtain for each $X \in {\cal V}(M)$ and $p \geq 1$ a unique linear map 
$$ i_X  \: \Omega^p(M,E) \to \Omega^{p-1}(M,E)  \quad \hbox{ with } \quad 
(i_X\omega)_x = i_{X(x)}\omega_x, $$
where 
$(i_v \omega_x)(v_1,\ldots, v_{p-1}) := \omega_x(v, v_1,\ldots, v_{p-1}).$
For $\omega \in \Omega^0(M,E) = C^\infty(M,E)$, we put $i_X \omega := 0$. 
\qed

\Proposition I.4.3. For $X, Y \in {\cal V}(M)$, we have on $\Omega(M,E)$ the 
Cartan formulas: 
$$ [{\cal L}_X, i_Y] = i_{[X,Y]}, \quad 
{\cal L}_X = d \circ i_X + i_X \circ d \quad \hbox{ and } \quad  
{\cal L}_X \circ d = d \circ {\cal L}_X. $$

\Remark I.4.4. Clearly integration of differential forms $\omega \in \Omega^p(M,E)$ 
only makes sense if $M$ is a $p$-dimensional compact oriented
 manifold (possibly with boundary) and $E$ is Mackey complete (Definition~I.1.4). 
We need the Mackey completeness to ensure  
that each smooth function $f \: Q \to E$ on a cube $Q := \prod_{i=1}^p 
[a_i, b_i] \subeq \R^p$ has an iterated integral 
$$ \int_Q f dx := \int_{a_1}^{b_1} \cdots 
\int_{a_p}^{b_p} f(x_1, \ldots, x_p)\, dx_p\cdots dx_1. 
\qeddis 

\Remark I.4.5.  (a) We call a smooth manifold $M$ {\it smoothly paracompact} if 
every open cover has a subordinated smooth partition of
unity. De Rham's Theorem holds for every 
smoothly paracompact manifold (cf.\ [Ee58], [KM97, Thm.\ 34.7], [Beg87]). 
Smoothly Hausdorff second countable
manifolds modeled on a smoothly regular space are smoothly paracompact 
([KM97, Cor.~27.4]). 
Typical examples of smoothly regular spaces are nuclear 
Fr\'echet spaces ([KM97,  Th.~16.10]). 

  (b) Examples of Banach spaces which are not
smoothly paracompact are $C([0,1],\R)$ and $\ell^1(\N,\R)$. 
On these spaces, there exists no non-zero smooth function supported in the 
unit ball ([KM97, 14.11]). 
\qed

\subheadline{I.5. The topology on spaces of smooth functions} 

In this subsection, we describe a natural topology on spaces of smooth maps which is 
derived from the compact open topology, the {\it compact open $C^r$-topology} 
(cf.\ [Mil82, Ex.~6.10] for a comparison of different topologies on spaces of 
smooth maps). Unfortunately, this topology has certain defects for functions 
on infinite-dimensional manifolds. 

\Definition I.5.1. (a) If $X$ and $Y$ are topological Hausdorff 
spaces, then the {\it compact open topology} 
on the space $C(X,Y)$ is defined as  the topology generated by the sets of the form 
$$ W(K,U) := \{ f \in C(X,Y) \: f(K) \subeq U \}, $$
where $K$ is a compact subset of $X$ and $U$ an open subset of $Y$. 
We write $C(X,Y)_c$ for the topological space obtained by endowing 
$C(X,Y)$ with the compact open topology. 

(b) If $G$ is a topological group and $X$ is Hausdorff, then $C(X,G)$ is a group with respect to the 
pointwise product. Then the compact open topology on $C(X,G)$ coincides 
with the topology of uniform convergence on compact subsets of $X$, for which 
the sets $W(K,U)$, $K \subeq X$ compact and $U \subeq G$ a $\1$-neighborhood, 
form a basis of $\1$-neighborhoods. In particular, 
$C(X,G)_c$ is a topological group. 

(c) If $Y$ is a locally convex space, then $C(X,Y)$ is a vector space with respect to 
the pointwise operations. In view of the preceding remark, the topology on $C(X,Y)_c$ 
is defined by the seminorms 
$$ p_K(f) := \sup \{ p(f(x)) \: x \in K \}, $$
where $K \subeq X$ is compact and $p$ is a continuous seminorm on $Y$. 
It follows in particular that $C(X,Y)_c$ is a locally convex space. 

(d) If $M$ and $N$ are smooth (possibly infinite-dimensional) manifolds, then every 
smooth map $f \: M \to N$ defines a sequence of smooth maps 
$T^k f \: T^k M \to T^k N$
on the iterated tangent bundles. We thus obtain for 
$r \in \N_0 \cup \{ \infty \}$ an embedding 
$$ C^r(M,N) \into \prod_{k=0}^r C(T^kM, T^kN)_c, $$
into a topological product space, that we use to define a topology on $C^r(M,N)$, 
called the {\it compact open $C^r$-topology}. For $r < \infty$, it suffices to consider 
the embedding $C^r(M,N) \into C(T^r(M), T^r(N))_c$. 
On the set $C^\infty(M,N)$, the compact open $C^\infty$-topology is the common 
refinement of all $C^r$-topologies for $r < \infty$. 
Since every compact subset of $M$ is contained in a finite union of 
chart domains, the topology on $C^r(M,N)$ is generated by 
sets of the form $W(K,U)$ in  $C(T^k(M), T^k(N))$, where 
$K$ lies in $T^k(U)$ for a chart $(\phi,U)$ of $M$. 

If $E$ is a locally convex space, then all spaces $C(T^k M, T^k E)$ are 
locally convex, by (c) above. Therefore the corresponding product topology is 
locally convex, and hence $C^\infty(M,E)$ is a locally convex space. 
If $M$ is finite-dimensional, for each chart $(\phi,U)$ of $M$, 
the topology on $C^\infty(U,E)$ 
coincides with the topology of uniform convergence 
of all partial derivatives on each compact subset of $U$. 
\qed

\Definition I.5.2. Since 
smooth vector fields are smooth functions $X \: M \to TM$, we have a natural embedding 
${\cal V}(M) \into C^\infty(M,TM)$, defining a topology 
on ${\cal V}(M)$. If $(\phi, U)$ is an $E$-chart of $M$, 
then $TU \cong U \times E,$
and smooth vector fields on $U$ correspond to smooth functions $U \to E$. 
This shows that, endowed with its natural topology, 
${\cal V}(M)$ is a locally convex space. 
\qed

\Remark I.5.3. As a consequence of Remark~I.3.7, 
the bracket on ${\cal V}(M)$ is continuous 
if $M$ is finite-dimensional. 

It is interesting to observe that, in general, the bracket on ${\cal V}(M)$ 
is not continuous if $M$ is infinite-dimensional. 
To see this, we assume that 
$M= U$ is an open subset of a locally convex space $E$ 
and consider the subalgebra 
$\aff(E) \cong E \rtimes \gl(E)$ 
of affine vector fields 
$X_{A,b}$ with $X_{A,b}(v) = Av + b$. It is easy to see that 
the natural topology on ${\cal V}(U)$ induces on $\aff(E)$ 
the product topology of the original topology on $E$ 
and the compact open topology on $\gl(E) \cong {\cal L}(E)_c$. 
In view of 
$$ [X_{A,b}, X_{A',b'}] = X_{[A',A], A'b - Ab'}, $$
it therefore suffices to show that the bilinear evaluation map 
${\cal L}(E)_c \times E \to E$ 
is not continuous if $\dim E = \infty$. 
Pick $0 \not = v \in E$ and embed $E'_c \into {\cal L}(E)_c$ by assigning to 
$\alpha \in E'$ the operator $v \otimes \alpha \: x \mapsto \alpha(x) v$. 
Hence it suffices to see that the evaluation map 
$$ \ev \: E'_c \times E \to \R, \quad (\alpha,v) \mapsto \alpha(v) $$ 
is not continuous. Basic neighborhoods of $(0,0)$ in $E'_c \times E$ are 
of the form $\hat K\times U_E$, where 
$U_E \subeq U$ is a $0$-neighborhood, $K \subeq E$ is compact, 
and $\hat K := \{ f \in E' \: (\forall k \in K)\, |f(k)| \leq 1\}$ is the polar set of~$K$. 
On $\hat K \times U_E$ the evaluation map is bounded if and only if 
$U_E$ is contained in some multiple of the bipolar $\hats K$, 
which, according to the Bipolar Theorem, coincides with the balanced convex closure 
of $K$, which is pre-compact ([Tr67, Prop.~7.11]). 
Then $\hats K$ is a pre-compact $0$-neighborhood in $E$, so that 
$E$ is finite-dimensional (cf.\ [Ru73, proof of Th.~1.22]). 
A similar argument shows that, if we endow $E'$ with the finer 
topology of uniform convergence on bounded subsets of $E$, 
then the evaluation map is continuous if and only if $E$ is normable, which is 
equivalent to the existence of a (weakly) bounded $0$-neighborhood ([Ru73]). 
\qed

\Remark I.5.4. The fact that for an infinite-dimensional locally convex space 
$E$ the evaluation map $\ev \: E'_c \times E \to \R$ is not continuous 
also causes trouble if one wants to associate to transformation groups 
corresponding continuous, resp.,  smooth representations on function spaces. 

A very simple example of a smooth group action is the translation 
action of $E$ on itself. The corresponding representation of $(E,+)$ 
on the space of smooth functions on $E$ is given by $(x.f)(y) := f(x+y)$. 
Clearly, the subspace of affine functions in $C^\infty(E,\R)$ is isomorphic to 
$\R \times E'_c$ as a locally convex space, and on this subspace the 
representation of $E$ is given by $x.(t,\alpha) = (t+\alpha(x), \alpha)$, 
which is discontinuous because $\ev(\alpha,x) = \alpha(x)$ is not 
continuous (Remark~I.5.3). In view of [Mais63], the 
same pathology occurs for any locally convex topology on $C^\infty(E,\R)$ if $E$ is not normable. 
\qed

\sectionheadline{II. Locally convex Lie groups}

\nin In this section, we give the definition of a locally convex Lie group.  
We explain how its Lie algebra and the corresponding Lie functor are defined 
and describe some basic properties. 
In our discussion of Lie groups, we essentially follow [Mil82/84], 
but, as for manifolds, we do not assume that the model space of a Lie group is complete ([Gl02a]). 
The natural strategy to endow groups with (infinite-dimensional) Lie group 
structures is to construct a chart around the identity in which the group 
operations are smooth. As we shall see in Subsection II.2, this suffices in many situations 
to specify a global Lie group structure. 

In Subsection~II.3, we discuss a smoothness concept for maps with values in diffeomorphism groups 
of locally convex manifolds. This specializes in particular to maps into  
general linear groups of locally convex spaces. 
The main point of this subsection is to obtain uniqueness results 
for solutions of certain ordinary differential equations on locally convex manifolds. 
In Subsection~II.4, we apply all this to smooth maps with values in Lie groups, 
where it shows in particular that morphisms of connected Lie groups are determined 
by their differential in $\1$. 

We conclude this section with some basic results on 
the behavior of the exponential function (Subsection~II.5), 
and a discussion of the concept of an 
initial Lie subgroup in Subsection~II.6.

\subheadline{II.1. Infinite-dimensional Lie groups and their Lie algebras} 

\Definition II.1.1. A {\it locally convex Lie group} $G$ is a 
locally convex manifold endowed with a group structure such that the multiplication 
map $m_G \: G \times G \to G$ and the inversion map $\eta_G \: G \to G$ 
are smooth. 

A morphism of Lie groups is a smooth group homomorphism. In the following,  
we call locally convex Lie groups simply {\it Lie groups}. 
\qed

\Example II.1.2. (Vector groups) Each locally convex space 
$E$ is an abelian Lie group with respect to addition and 
the obvious manifold structure. 
\qed

Vector groups $(E,+)$ form the most elementary Lie groups. 
The next natural class are unit groups of algebras. This leads us to the 
concept of a continuous inverse algebra, which came up in the 1950s 
(cf.\ [Wae54a/b] and [Wae71]): 

\Definition II.1.3. (a) A {\it locally convex algebra} is a locally convex space 
$A$, endowed with an associative continuous bilinear multiplication 
$A \times A \to A, (a,b) \mapsto ab$. 
A unital locally convex algebra $A$ is called a 
{\it continuous inverse algebra} (CIA for short) 
if its unit group $A^\times$ 
is open and the inversion is a continuous map $A^\times \to A, a \mapsto a^{-1}$. 

(b) If $A$ is a locally convex algebra which is 
not unital, then we obtain a monoid structure on $A$ 
by $x * y := x + y + xy$ for which $0$ is the identity element. 
In this case, we write $A^\times$ for the unit group of $(A,*)$ and say that 
$A$ is a {\it non-unital CIA} if $A^\times$ is open and the (quasi-)inversion map 
$\eta_A \: A^\times \to A$ is continuous. 

If $A_+ := A \times \K$ is the unital locally convex algebra with the 
multiplication $(x,t)(x',t') := (xx'+tx'+t'x,tt')$, then the map 
$(A,*) \to A \times \{1\}, a \mapsto (a,1)$ is an isomorphism of monoids, 
and it is easy to see that $A_+$ is a CIA if and only if $A$ is a (not necessarily unital) CIA. 
\qed

\Example II.1.4. Let $A$ be a continuous inverse algebra over $\K$ and $A^\times$ its unit group. 
As an open subset of $A$, the group $A^\times$ carries a natural manifold 
structure. The multiplication on $A$ is bilinear and continuous, 
hence a smooth map (Remark~I.2.4). 
Therefore the multiplication of $A^\times$ is smooth. 
One can further show quite directly that the continuity of the 
inversion $\eta_A \: A^\times \to A^\times$ implies that 
$d\eta_A(x)(y) = - x^{-1}yx^{-1}$ exists for each pair $(x,y)$, 
and this formula implies inductively that $\eta_A$ is smooth and 
hence that $A^\times$ is a Lie group. 

In some cases, it is also possible to obtain a Lie group structure on the unit group 
$A^\times$ of a unital locally convex algebra if $A^\times$ is not open 
(cf.\ Remark II.2.10 below). 
\qed

\Definition II.1.5. A vector field $X$ on the Lie group $G$ is called 
 {\it left invariant} if 
$$X \circ \lambda_g = T(\lambda_g) \circ X \: G \to T(G) $$ 
holds for each $g \in G$, i.e., $X$ is $\lambda_g$-related to itself for each $g \in G$. 
We write ${\cal V}(G)^l$ for the set of left invariant vector fields in 
${\cal V}(G)$. 
The left invariance of a vector field $X$ implies in particular that for each 
$g \in G$, we have $X(g) = g.X(\1)$, where $G \times T(G) \to T(G), (g,v) \mapsto g.v$ denotes 
the smooth action of $G$ on $T(G)$, induced by the left multiplication action of $G$ 
on itself. 
For each $x \in \g$, we have a unique left invariant vector field 
$x_l \in {\cal V}(G)^l$ defined by $x_l(g) := g.x$, and the map 
$$ \ev_\1 \: {\cal V}(G)^l \to T_\1(G), \quad X \mapsto X(\1) $$ 
is a linear bijection. If $X,Y$ are left invariant, then they are 
$\lambda_g$-related to themselves, and 
their Lie bracket $[X,Y]$ inherits this property. 
We thus obtain a unique Lie bracket $[\cdot, \cdot]$ on $T_\1(G)$ satisfying 
$$ [x,y]_l = [x_l, y_l] \quad \hbox{ for all} \quad x,y \in T_\1(G),  \leqno(2.1.1) $$
and from the formula for the bracket in local coordinates, 
it follows that it is continuous (cf.\ Remark~II.1.8 below). 
\qed

\Remark II.1.6. The tangent map $T(m_G) \: T(G \times G) \cong T(G) \times T(G) \to T(G)$ defines 
on the tangent bundle $T(G)$ of $G$ the structure of a Lie group with inversion map $T(\eta_G)$. 

In fact, let $\eps_G \: G  \to G, g \mapsto \1$, be the constant homomorphism. 
Then the group axioms for $G$ are encoded in the relations 
\litem{(1)} $m_G \circ (m_G \times \id_G) = m_G \circ (\id_G \times m_G)$ (associativity), 
\litem{(2)} $m_G \circ (\eta_G,\id_G) = m_G \circ (\id_G,\eta_G) = \eps_G$ 
(inversion), and 
\litem{(3)} $m_G \circ (\eps_G, \id_G) = m_G \circ (\id_G, \eps_G) = \id_G$ (unit element).
\par\nin Applying the functor $T$ to these relations, 
it follows that $T(m_G)$ defines a Lie group structure on 
$T(G)$ for which $T(\eta_G)$ is the inversion and $0_\1 \in T_\1(G)$ is the identity. 
\qed

\Definition II.1.7. {\rm(The Lie functor)} For a Lie group $G$, the locally convex Lie algebra 
$\L(G) := (T_\1(G), [\cdot, \cdot])$ is called {\it the Lie algebra of $G$}. 

To each morphism $\phi \: G \to H$ of Lie groups we further associate 
its tangent map 
$$\L(\phi) := T_\1(\phi) \: \L(G) \to \L(H),$$ 
and the usual argument with related vector fields implies that 
$\L(\phi)$ is a homomorphism of Lie algebras. 
\qed

This means that the assignments 
$G \mapsto \L(G)$ and $\phi \mapsto \L(\phi)$ define a functor 
$\L$ from the category of (locally convex) Lie groups to the category 
of locally convex Lie algebras. 
Since each functor maps isomorphisms to isomorphisms, 
for each isomorphism $\phi \: G \to H$ of Lie groups, 
the map $\L(\phi)$ is an isomorphism of locally convex Lie algebras. 

The following remark describes a convenient way to calculate the Lie algebra of a given group. 

\Remark II.1.8. For each chart $(\phi,U)$ of $G$ with $\1 \in U$ and $\phi(\1) = 0$, 
we identify $\g := T_\1(G)$ via the topological isomorphism $T_\1(\phi)$ 
with the corresponding model space. Then 
the second order Taylor expansion in $(0,0)$ of the multiplication 
$x * y := \phi(\phi^{-1}(x)\phi^{-1}(y))$ (cf.~Proposition~I.2.3) is of the form 
$$ x * y = x + y + b(x,y) + \hbox{higher order terms}, $$
where $b \: \g \times \g \to \g$ is a continuous bilinear map satisfying 
$$ [x,y] = b(x,y) - b(y,x). \leqno(2.1.2) $$
Using the chain rule for Taylor polynomials, it is easy to show that the 
second order Taylor polynomial of the commutator map 
$x*y*x^{-1}*y^{-1}$ is given by the Lie bracket: 
$$ x*y*x^{-1}*y^{-1} = [x,y] + \hbox{higher order terms}
\qeddis

We now take a look at the Lie algebras of the Lie groups from Examples~II.1.2/4. 
\Examples II.1.9. (a) If $G$ is an abelian Lie group, then the map $b \: \g \times \g \to \g$ 
in Remark~II.1.8 is symmetric, which implies that $\L(G)$ is abelian. 
This applies in particular to the additive Lie group $(E,+)$ of a locally 
convex space~$E$. 

(b) Let $A$ be a CIA. Then the map 
$\phi \: A^\times \to A, x \mapsto x - \1$
is a global chart of $A^\times$, satisfying $\phi(\1) = 0$. 
In this chart, the group multiplication is given by 
$$ x * y := \phi(\phi^{-1}(x)\phi^{-1}(y)) = (x +\1)(y+\1) - \1 
= x + y + xy. $$
In the terminology of Remark~II.1.8, we then have 
$b(x,y) = xy$ and therefore 
$\L(A^\times) = (A,[\cdot,\cdot])$, where 
$[x,y] =xy - yx$ is the commutator bracket on the associative algebra~$A$.
\qed

We conclude this subsection with the observation that the passage from groups 
to Lie algebras can also be established on the local level.  

\Definition II.1.10. (The Lie algebra of a local Lie group) 
There is a natural notion of a local Lie group. The corresponding algebraic 
concept is that of a local group: 
Let $G$ be a set and $D \subeq G \times G$ a subset 
on which we are given a map 
$$m_G  \: D \to G, \quad (x,y) \mapsto xy. $$ 
We say that the product $xy$ of two elements $x,y \in G$ 
is {\it defined} if $(x,y) \in D$. 
The quadruple $(G,D,m_G,\1)$, where $\1$ is a distinguished element of $G$, 
is called a {\it local group} if the following conditions are
satisfied: 
\litem{(1)} Suppose that $xy$ and $yz$ are defined. If $(xy)z$ or
$x(yz)$ is defined, then the other product is also defined and both
are equal. 
\litem{(2)} For each $x \in G$,  
the products $x\1$ and $\1x$ are defined and equal to $x$. 
\litem{(3)} For each $x \in G$, there exists a unique element $x^{-1}
\in G$ such that $xx^{-1}$ and $x^{-1} x$ are defined and $xx^{-1} = x^{-1}x = \1$.
\litem{(4)} If $xy$ is defined, then $y^{-1}x^{-1}$ is defined.

 
If $(G,D,m_G,\1)$ is a local group and, in addition, 
$G$ has a smooth manifold structure, 
$D$ is open, and the maps 
$$ m_G \: D \to G\quad \hbox{ and }  \quad \eta_G \: G \to G, x \mapsto x^{-1} $$
are smooth, then $G$, resp., $(G,D,m_G,\1)$ is called a {\it local Lie group}. 

Let $G$ be a local Lie group and $T_\1(G)$ its tangent space in $\1$. For 
each $x \in T_\1(G)$, we then obtain a {\it left invariant vector field} 
$x_l(g) := g.x := 0_g \cdot x$. 
The Lie bracket of two left invariant vector fields 
is left invariant and we thus obtain on $T_\1(G)$ a locally convex Lie algebra 
structure. We call $\L(G) := \L(G,D,m_G,\1) := (T_\1(G),[\cdot,\cdot])$ the 
Lie algebra of the local group $G$. For more details on local Lie groups, we refer to [GN06]. 
\qed

\Remark II.1.11. If $G$ is a Lie group and $U = U^{-1} \subeq G$ an open 
identity neighborhood, then $U$ carries a natural local Lie group 
structure with $D := \{ (x,y) \in U \times U \: xy \in U\}$ and
$m_U := m_G\res_D$. Clearly $U$ and $G$ have the same Lie algebras. 

Local groups of this type are called {\it enlargeable}. As we shall 
see in Example~VI.1.7 below, not all local Lie groups are enlargeable 
because not all Banach--Lie algebras are integrable (Example VI.1.16).
\qed 

\subheadline{II.2. From local data to global Lie groups} 

We now give the precise formulation of an elementary 
but extremely useful tool which helps to construct Lie group structures 
on groups containing a local Lie group. 
This theorem directly carries over from the finite-dimensional 
case, which can be found in [Ch46, \S 14, Prop.~2] or [Ti83, p.14]. 
In [GN06], it is our main method to endow groups with Lie group structures. 

\Theorem II.2.1. Let $G$ be a group and 
$U = U^{-1}$ a symmetric subset. 
We further assume that $U$ is a smooth manifold such that 
\litem{(L1)} there exists an open symmetric $\1$-neighborhood $V \subeq
U$ with $V \cdot V \subeq U$ such that the group multiplication 
$m_V \: V \times V \to U$ is smooth, 
\litem{(L2)} the inversion map $\eta_U \: U \to U, u \mapsto u^{-1}$ is
smooth, and 
\litem{(L3)} for each $g \in G$ there exists an open symmetric 
$\1$-neighborhood $U_g \subeq
U$ with $c_g(U_g) \subeq U$, and such that the conjugation map 
$c_g \: U_g \to U, x \mapsto gxg^{-1}$
is smooth. 

Then there exists a unique Lie group structure 
on $G$ for which there exists an open $\1$-neighborhood 
$U_0 \subeq U$ such that the inclusion map $U_0 \to G$ induces a 
diffeomorphism onto an open subset of $G$. 

If the group $G$ is generated by $V$, then condition {\rm(L3)} can be
omitted. 
\qed 

If $V$ is as above, then $D := \{ (x,y) \in V \times V \: xy \in V\}$ 
defines on $V$ the structure of a local Lie group, and the preceding 
theorem implies that the smooth structure of this local Lie group, together with 
the group structure of $G$, determines the global Lie group structure of $G$. 
The subtlety of condition (L3) is that it mixes local and global objects because 
it requires that each element of $G$ induces an isomorphism of the
 corresponding germ of local groups. The following corollary is a converse of 
Remark~II.1.11 (cf.\ [Swi65]). It is the central tool to pass from local to global 
subgroups of Lie groups. 

\Corollary II.2.2. Let $(U, D, m_U, \1)$ be a local Lie group, $G$ a group,  and 
$\eta \: U \to G$ an injective morphism of local groups. 
Then the subgroup $H := \la \eta(U)\ra \subeq G$ generated by $\eta(U)$ carries a unique 
Lie group structure for which $\eta$ is a diffeomorphism 
onto an open subset of~$H$. 
\qed

The preceding corollary shows in particular that if, in addition to the assumptions of 
Theorem~II.2.1, the group multiplication of $G$ restricts to a smooth map 
on the domain $D_U := \{ (x,y) \in U \times U \: xy \in U\}$, then the inclusion $U \into G$ is a diffeomorphism 
onto an open subset of $G$, endowed with the Lie group structure determined by~$U$. 

\Corollary II.2.3. Let $G$ be a group and $N \trile G$ a normal subgroup that 
carries a Lie group structure. Then there exists a Lie group structure on $G$ 
for which $N$ is an open subgroup if and only if for each $g \in G$,  
the restriction $c_g\res_N$ is a smooth automorphism of $N$. 
\qed

The preceding corollary is of particular interest for abelian groups. In this case, it leads 
for each Lie group structure on any subgroup $N \subeq G$ to a Lie group structure on $G$ 
for which $N$ is an open subgroup. 

The following corollary implies in particular 
that quotients of Lie groups by discrete normal subgroups are Lie groups. 

\Corollary II.2.4. Let $\phi \: G \to H$ be a covering of topological groups. 
If $G$ or $H$ is a Lie group, then the other group has a unique Lie group 
structure for which $\phi$ is a morphism of Lie groups which is a local diffeomorphism. 
\qed

\Remark II.2.5. (a) {\rm(Lie subgroups)} If $G$ is a Lie group with Lie algebra $\g$ 
and $H \subeq G$ is a submanifold which 
is a group, then $H$ inherits a Lie group structure from $G$. 
Moreover, there exists a closed subspace $\h \subeq \g \cong T_\1(G)$ 
and a chart $(\phi,U)$ of $G$ with $\1 \in U= U^{-1}$, $\phi(\1) = 0$ and 
$$ \phi(U \cap H) = \phi(U) \cap \h. $$
The local multiplication $x*y := \phi(\phi^{-1}(x)\phi^{-1}(y))$ on 
$$D := \{(x,y) \in \phi(U) \times \phi(U) \: \phi^{-1}(x)\phi^{-1}(y) \in U\}$$ 
then satisfies 
$$ x*y \in \h \quad \hbox{ for } \quad (x,y) \in D \cap (\h \times \h) \leqno(2.2.1) $$
and 
$$ x^{-1} \in \h \quad \hbox{ for } \quad x \in \h \cap \phi(U). \leqno(2.2.2) $$
In view of Remark~II.1.8, this implies that $\h$ is a closed Lie subalgebra of 
$\g$. 

If, conversely, $\h \subeq \g$ is a closed Lie subalgebra for which there 
is a chart $(\phi,U)$ as above, satisfying (2.2.1/2), then 
$\phi(U) \cap \h$ carries a local Lie group structure and we can apply Corollary~II.2.2 to 
the embedding $\phi^{-1} \: \phi(U) \cap \h \to G$, which leads to a Lie group 
structure on the subgroup $H := \la \phi^{-1}(\phi(U) \cap \h)\ra$ of $G$. 
We know already from the finite dimensional theory that, in general, this does 
not lead to a submanifold of $G$. 

(b) A weaker concept of a ``Lie subgroup'' is obtained by requiring only that 
$H \subeq G$ a subgroup, for  which there exists an identity neighborhood $U^H$ 
whose smooth arc-component $U^H_0$ of $\1$ is a submanifold of $G$ 
(cf.\ [KYMO85, p.45]). Then we can use Theorem~II.2.1 to  obtain a Lie group 
stucture on $H$ for which some identity neighborhood is diffeomorphic to 
an identity neighborhood in $U^H_0$. 
\qed

\Remark II.2.6. Since it also makes sense to consider manifolds without 
assuming that they are Hausdorff (cf.\ [Pa57], 
[La99]), it is worthwhile to observe that 
this does not lead to a larger class of Lie groups. 

In fact, let $G$ be a Lie group which is not necessarily Hausdorff. Then $G$  
is in particular a topological group which possesses an identity 
neighborhood $U$ homeomorphic to an open subset of a locally convex space. 
As $U$ is Hausdorff, and since the subgroup $\oline{\{\1\}}$ of $G$ coincides 
with the intersection of all $\1$-neighborhoods, the closedness of 
$\{\1\}$ in $U$ implies that $\{\1\}$ is a closed subgroup of $G$ and hence 
that $G$ is a Hausdorff topological group. 
\qed

To see how Theorem II.2.1 can be applied, we now take a closer look 
at groups of differentiable maps. First we introduce a natural topology on these groups. 

\Definition II.2.7. (Groups of differentiable maps as topological groups)  
Let $M$ be a smooth manifold (possibly in\-fi\-nite-di\-men\-sio\-nal), 
$K$ a Lie group with Lie algebra $\k$ and $r \in \N_0 \cup \{\infty\}$. 
We endow the group $G := C^r(M,K)$ with the compact open $C^r$-topology (Definition~I.5.1). 

We know already that the tangent bundle $TK$ of $K$ is a Lie group (Remark~II.1.6). 
Iterating this procedure, we obtain a Lie group structure on all higher tangent bundles $T^n K$. 
For each $n \in \N_0$, we thus obtain topological groups $C(T^n M, T^n K)_c$ 
by using the topology of uniform convergence on compact subsets of $T^n M$, 
which coincides with the compact open topology (Definition~I.5.1). 
We also observe that for two smooth maps $f_1, f_2 \: M \to K$,  
the functoriality of $T$ yields 
$$ T(f_1\cdot f_2) = T(m_G \circ (f_1 \times f_2)) 
=T(m_G) \circ (Tf_1 \times Tf_2) = Tf_1 \cdot Tf_2. $$
Therefore the inclusion map 
$$ C^r(M,K) \into \prod_{n = 0}^r  C(T^n M, T^n K)_c, \quad 
f \mapsto (T^n f)_{0 \leq n \leq r} $$
is a group homomorphism, so that the inverse image of the product 
topology from the right hand side is a 
group topology on $C^r(M,K)$. Hence 
the compact open $C^r$-topology turns $C^r(M,K)$ into a topological group, 
even if $M$ and $K$ are infinite-dimensional. 
\qed

We define the {\it support} of a Lie group-valued map 
$f \: M \to G$ by 
$$\supp(f) := \oline{\{x \in M \: f(x) \not=\1\}}, $$
for a closed subset $X \subeq M$ we put 
$$ C^r_X(M,K) := \{ f \in C^r(M,K) \: \supp(f) \subeq X \}, $$
and write $C^r_c(M,K)$ for the subgroup of compactly supported $C^r$-maps. 

\Theorem II.2.8. Let $K$ be a Lie group with Lie  algebra $\k$,  
$M$ a finite-dimensional manifold (possibly with boundary),  
and $r \in \N_0\cup \{\infty\}$. 

{\rm(a)} If $M$ is compact, then 
$C^r(M,K)$ carries a Lie group structure compatible with the compact open 
$C^r$-topology, 
and its Lie algebra is $C^r(M,\k)$, endowed with the pointwise bracket. 

{\rm(b)} If $M$ is $\sigma$-compact, then 
$C^r_c(M,\k)$, endowed with the locally convex direct limit topology of the spaces 
$C^r_X(M,\k)$, $X \subeq M$ compact, is a topological Lie algebra 
and $C^r_c(M,K)$ carries a natural Lie group structure with Lie algebra 
$C^r_c(M,\k)$.
 
\Proof. (Sketch) (a) Let $G := C^r(M,K)$ and $\g:= C^r(M,\k)$. 
The Lie group structure on $G$ can be constructed with Theorem~II.2.1 as follows. 
Let $\phi_K \: U_K \to\k$ be a chart of $K$. Then the set 
$U_G := \{ f \in G \: f(M) \subeq U_K\}$ is an open subset of 
$G$. Assume, in addition, that $\1 \in U_K$ and $\phi_K(\1)=0$. 
Then the map 
$$ \phi_G \: U_G \to \g, \quad f \mapsto \phi_K \circ f $$
defines a chart $(\phi_G, U_G)$ of $G$. To apply Theorem II.2.1, one 
has to verify that in this chart 
the inversion is a smooth map, that the multiplication map 
$$ D_G := \{ (f,g) \in U_G \times U_G \: fg \in U_G \} \to U_G $$
is smooth and that conjugation maps are smooth in some $\1$-neighborhood of $U_G$. 
For details, we refer to [Gl02c], resp., [GN06]. 

To calculate the Lie algebra of this group, we observe that for 
$m \in M$, we have for the multiplication in local coordinates 
$$ \eqalign{ (f *_G g)(m) 
&:= \phi_G\Big(\phi_G^{-1}(f)\phi_G^{-1}(g)\Big)(m) 
= \phi_K\big(\phi_K^{-1}(f(m))\phi_K^{-1}(g(m))\big) \cr
&= f(m) *_K g(m)= f(m) + g(m) + b_\k(f(m), g(m)) + \cdots. \cr} $$
In view of Remark~II.1.8, this implies that 
$b_\g(f,g)(m) = b_\k(f(m),g(m)),$
and hence that 
$$ [f,g](m) = b_\g(f,g)(m) - b_\g(g,f)(m) = b_\k(f(m),g(m))-b_\k(g(m),f(m)) = [f(m),g(m)].$$
Therefore $\L(C^r(M,K)) = C^r(M,\k)$, endowed 
with the pointwise defined Lie bracket. 

(b) is proved along the same lines. Note that it is not obvious that the Lie bracket 
on $C^r_c(M,\k)$ is continuous because it is a bilinear map. 
\qed

If $K$ is finite-dimensional, then the preceding Lie group construction 
can be found in {\smc Michor}'s book [Mi80], 
and also in [AHMTT93] (which basically deals with the topological level). 
In [BCR81], one finds interesting variants of groups of smooth maps on 
open subsets $U \subeq \R^n$ which are rapidly decreasing at the boundary with 
respect to certain weight functions. In particular, there is a Lie group 
${\cal S}(\R^n,K)$ whose Lie algebra is the space 
${\cal S}(\R^n,\k)$ of $\k$-valued Schwartz functions on $\R^n$.

\Remark II.2.9. (a) If $M$ is a non-compact finite-dimensional manifold, 
then one cannot expect the topological groups $C^r(M,K)$ to be 
Lie groups. A typical example arises for $M = \N$ (a $0$-dimensional manifold) 
and $K = \T := \R/\Z$. Then $C^r(M,K) \cong \T^\N$ is a compact topological 
group for which no $\1$-neighborhood is contractible, so that it carries 
no smooth manifold structure. 

(b) Non-linear maps on spaces 
of compactly supported functions such as 
$E := C^\infty_c(\R,\R)$ (Examples~I.1.3) require extreme caution. E.g., the map 
$$ f \: C^\infty_c(\R,\R) \to C^\infty_c(\R,\R), \quad \gamma \mapsto 
\gamma \circ \gamma - \gamma(0) $$
is smooth on each closed Fr\'echet subspace $E_n := C^\infty_{[-n,n]}(\R,\R)$, 
but it is discontinuous in $0$ ([Gl06a]). Therefore the LF space 
$E = \indlim E_n$ is a direct limit in the category of locally convex spaces, 
but not in the category of topological spaces. 
\qed

\Remark II.2.10. (a) Let $A$ be a commutative unital locally convex algebra 
with a smooth exponential function 
$$ \exp_A \: A \to A^\times, $$
i.e., $\exp_A \: (A,+) \to (A^\times, \cdot)$ is a group homomorphism with 
$T_0(\exp_A)= \id_A$. 

Then $\Gamma_A := \ker(\exp_A)$ is a closed subgroup of $A$ not containing any 
line. Suppose that $\Gamma_A$ is discrete. Then 
$N := A/\Gamma_A$ carries a natural Lie group structure (Corollary~II.2.4) 
and the exponential function factors through an injection $N \into A^\times$.
We may therefore use Corollary~II.2.3 to define a Lie group structure on the group 
$A^\times$ for which the identity component is $\exp_A(A) \cong N$.

(b) If $M$ is a $\sigma$-compact finite-dimensional manifold, then 
$A := C^\infty(M,\C)$ is a complex locally convex algebra with respect to the 
compact open $C^\infty$ topology, and 
$$ \exp_A \: A \to A^\times = C^\infty(M,\C^\times), \quad f \mapsto e^f $$
is a smooth exponential function. 

If $M$ is non-compact, then $A^\times$ is not open because 
for each unbounded function $f \: X \to \C$ the element 
$\1 + \lambda f$ is not invertible for $\lambda \in \C$ arbitrarily close to $0$. 
It follows that $A$ is a CIA if and only if $M$ is compact. 

The closed subgroup $\Gamma_A = \ker(\exp_A) = C^\infty(M,2\pi i \Z) \cong C^\infty(M,\Z)$ 
is discrete if and only if $M$ has only finitely many connected components. 
In this case, (a) implies that $A^\times$ carries a Lie group structure for which 
$\exp_A$ is a local diffeomorphism. 

A typical example is $M = \R$ and $A = C^\infty(\R,\C)$ with 
$$ A^\times = C^\infty(\R,\C^\times) 
\cong C^\infty_*(\R,\C^\times) \times \C^\times 
\cong C^\infty_*(\R,\C) \times \C^\times$$ as topological groups, where 
$C^\infty_*$ denotes functions mapping $0$ to $1$, resp., to $0$. 
For $M = \N$, we have $A \cong \C^\N$, and $\Gamma_A \cong \Z^\N$ is not discrete. 

If $M$ is connected, it is not hard to see that the map 
$$ \delta \: A^\times \cong C^\infty(M,\C^\times) \to Z^1_{\rm dR}(M,\C), \quad 
f \mapsto {df \over f} $$
induces a topological isomorphism of $A^\times/\C^\times$ onto the 
group $Z^1_{\rm dR}(M,\Z)$ of closed $1$-forms whose periods are contained in 
$2\pi i \Z$, and the arc-component $A^\times_a$ of the identity is mapped onto 
the set of exact $1$-forms. We conclude that, as topological groups, 
$$ \pi_0(A^\times) \cong A^\times/A^\times_a \cong H^1_{\rm dR}(M,\Z) 
\cong \Hom(\pi_1(M),\Z) \cong \Hom(H_1(M),\Z), $$
and this group is discrete if and only if 
$H^1_{\rm sing}(M,\Z) \cong \Hom(H_1(M),\Z)$ is finitely generated 
(cf.\ [NeWa06b]). This shows that the arc-component of the identity in 
$A^\times$ is open if and only if 
$H^1_{\rm sing}(M,\Z)$ is finitely generated. 

For $M := \C\setminus \N$, the group  
$H_1(M) \cong \Z^{(\N)}$ is of infinite rank, 
$H^1_{\rm sing}(M,\Z) \cong \Z^\N$ is not discrete, 
but $M$ is connected, so that $A^\times$ carries a Lie group structure whose underlying 
topology is finer than the original topology of $A^\times$ induced from $A$. 
\qed

\subheadline{II.3. Smoothness of maps into diffeomorphism groups} 

Although the notion of a smooth manifold provides us with a natural 
notion of a smooth map between such manifolds, it turns out to be 
convenient to have a notion of a smooth map of a manifold into spaces of smooth maps  
which do not carry a natural manifold structure. In this subsection,  
we discuss this notion of smoothness with an emphasis on  
maps with values in groups of diffeomorphisms of locally convex manifolds. 

\Definition II.3.1. Let $M$ be a smooth locally convex manifold and 
$\Diff(M)$ the group of diffeomorphisms of $M$. Further let $N$ be  a
smooth manifold. Although, in general, $\Diff(M)$ has no
natural Lie group structure, we call a map 
$\phi \: N \to \Diff(M)$ {\it smooth} if the map 
$$ \hat \phi \: N \times M \to M \times M, \quad 
(n,x) \mapsto (\phi(n)(x), \phi(n)^{-1}(x)) $$
is smooth. If $N$ is an interval in $\R$, we obtain in particular a notion of a 
smooth curve. 
\qed

To discuss derivatives of such smooth map, we take a closer look 
at the ``tangent bundle'' of $\Diff(M)$, which can be done without 
a Lie group structure on $\Diff(M)$ (which does not exist in a satisfactory fashion 
for non-compact $M$; cf.\ Theorem VI.2.6). 
We think of the set 
$$ T(\Diff(M)) := \{ X \in C^\infty(M,TM) \: \pi_{TM} \circ X \in
\Diff(M)\} $$
as the {\it tangent bundle of $\Diff(M)$}, with the map 
$$ \pi \: T(\Diff(M)) \to \Diff(M), \quad X \mapsto \pi_{TM}
\circ X $$
as the bundle projection, and 
$T_\phi(\Diff(M)) := \pi^{-1}(\phi)$
is considered as the {\it tangent space} in $\phi \in \Diff(M)$. 
We have natural left and right
actions of $\Diff(M)$ on $T(\Diff(M))$ by 
$$ \phi.X = T(\phi) \circ X \quad \hbox{ and } \quad X.\phi := X \circ \phi. $$
The action 
$$ \Ad \: \Diff(M) \times {\cal V}(M) \to {\cal V}(M), \quad 
(\phi, X) \mapsto \phi_*X := \Ad(\phi).X := T(\phi) \circ X \circ \phi^{-1} $$
is called the {\it adjoint action of $\Diff(M)$ on ${\cal V}(M)$}. 

Smooth curves $\phi \: J \subeq \R\to \Diff(M)$ have (left) logarithmic
derivatives 
$$ \delta(\phi)\: J \to {\cal V}(M), \quad 
\delta(\phi)_t := \phi(t)^{-1}.\phi'(t)$$ 
which are smooth curves
in the Lie algebra ${\cal V}(M)$ of smooth vector fields on $M$, i.e., 
time-dependent vector fields. For general $N$, the logarithmic derivatives 
are ${\cal V}(M)$-valued $1$-forms on~$N$, defined as follows: 

If $\phi \: N \to \Diff(M)$ is smooth and 
$\hat\phi_1 \: N \times M \to M, (n,x) \mapsto \phi(n)(x)$, 
then we have a smooth 
tangent map 
$$ T(\hat\phi_1) \: T(N \times M) \cong T(N) \times T(M) \to T(M), $$
and for each $v \in T_p(N)$ the partial map 
$$ T_p(\phi)v \: M \to T(M), \quad 
m \mapsto T_{(p,m)}(\hat\phi_1)(v,0) $$
is an element of $T_{\phi(p)}(\Diff(M))$. We thus obtain a {\it tangent map} 
$$ T(\phi) \: T(N) \to T(\Diff(M)), \quad v \in T_p(N) \mapsto T_p(\phi)v. $$

\Definition II.3.2. 
We define the {\it (left) logarithmic derivative of $\phi$ in $p$} by 
$$ \delta(\phi)_p \: T_p(N) \to {\cal V}(M), \quad 
v \mapsto \phi(p)^{-1}.T_p(\phi)(v) = 
T(\phi(p)^{-1}) \circ T_p(\phi)(v). $$
It can be shown that $\delta(\phi)$ is a smooth 
${\cal V}(M)$-valued $1$-form on $N$ (see [GN06] for details), 
but recall that ${\cal V}(M)$ need not be a topological Lie algebra if 
$M$ is not finite-dimensional (Remark~I.5.3). 
\qed

For calculations, it is convenient to observe the 
Product- and Quotient Rule, both easy consequences of the Chain Rule:

\Lemma II.3.3.  
For two smooth maps $f, g \:  N \to \Diff(M)$, define 
$(fg)(n) := f(n)\circ g(n)$ and $(fg^{-1})(n) := f(n)\circ g(n)^{-1}$. 
Then we have the 
\litem{(1)} Product Rule: $\delta(fg) = \delta(g) + \Ad(g^{-1}).\delta(f)$, and the 
\litem{(2)} Quotient Rule: $\delta(fg^{-1}) = \Ad(g).(\delta(f) - \delta(g))$, 
\par\nin where we write 
$(\Ad(f).\alpha)_n := \Ad(f(n)).\alpha_n$ for a ${\cal V}(M)$-valued 
$1$-form $\alpha$ on $N$.  
\qed

\Remark II.3.4. Although we shall only use the left logarithmic derivative, 
we note that one can also define the {\it right logarithmic derivative} of a smooth 
map $\phi \: N \to \Diff(M)$ by 
$$ \delta^r(\phi)_p(v) = \big(T_p(\phi)v\big) \circ \phi(p)^{-1}, $$
which also defines an element of $\Omega^1(N,{\cal V}(M))$, satisfying 
$\delta^r(\phi) = \Ad(\phi).\delta(\phi) = - \delta(\phi^{-1}).$

We then have for two smooth maps $f,g \: N \to \Diff(M)$ the 
\litem{(1)} Product Rule: $\delta^r(fg) = \delta^r(f) + \Ad(f).\delta^r(g)$, 
and the 
\litem{(2)} Quotient Rule: $\delta^r(fg^{-1}) = \delta^r(f) - \Ad(fg).\delta^r(g).$
\qed

The following lemma generalizes Lemma~7.4 in [Mil84] which deals with Lie group-valued 
curves. 

\Lemma II.3.5. {\rm(Uniqueness Lemma)} 
Suppose that $N$ is connected. For two smooth maps $f, g \:  N \to \Diff(M)$, the 
relation $\delta(f)= \delta(g)$ is equivalent to the
existence of $\phi \in \Diff(M)$ with 
$g(p) = \phi \circ f(p)$ for all $p \in N$. 
In particular, $g(p_0) = f(p_0)$ for some $p_0 \in N$ implies $f = g$. 

\Proof. If $g(p) = \phi \circ f(p)$ for each $p \in N$, then 
$T_p(g)  = \phi(p).T_p(f)$, and therefore 
$\delta(g) = \delta(f)$. 

If, conversely, $\delta(g) = \delta(f)$ and $\gamma := gf^{-1}$, 
then the Quotient Rule implies 
$\delta(\gamma) = \delta(gf^{-1}) = 0$, which in turn implies that 
 for each $x \in M$ the map $p \mapsto \gamma(p)(x)$ has vanishing derivative, hence is
locally constant (Proposition~I.3.4). 
Since $N$ is connected, $\gamma$ is constant. 
We conclude that $g = \phi \circ f$ for some $\phi \in \Diff(M)$. 
\qed

The Uniqueness Lemma is a key tool which implies in particular 
that solutions to certain initial value problems
are unique whenever they exist (which need not be the case). 
In this generality, this is quite 
remarkable because there are ordinary linear differential 
equations with constant coefficients on Fr\'echet spaces $E$ for which
solutions are not unique (cf.\ Example II.3.11 below). 
Nevertheless, the Uniqueness Lemma implies 
that solutions of the corresponding operator-valued 
initial value problems on the group $\GL(E) \subeq \Diff(E)$ are unique whenever they exist. 

\Remark II.3.6. Smooth maps with values in $\Diff(M)$ can be 
specialized in several ways: 

(a) Let $E$ be a locally convex space and $\GL(E)$ the group of linear topological 
automorphisms of $E$. Then $\GL(E)$ consists of all diffeomorphisms of $E$ commuting 
with the scalar multiplications $\mu_t(v) = tv$, $t \in \K^\times$, and 
$\gl(E) = ({\cal L}(E),[\cdot,\cdot])$ can be identified with the 
Lie subalgebra  of ${\cal V}(E)$ consisting of linear vector fields which can be characterized 
in a similar way. 
This observation implies that the logarithmic derivative of a smooth map 
$\phi \: N \to \GL(E)$ is a $\gl(E)$-valued $1$-form on $N$ and that the 
Uniqueness Lemma applies to $\GL(E)$-valued smooth maps. 

(b) If $K$ is a Lie group with Lie algebra $\k$, 
then we consider the group $C^\infty(M,K)$ of smooth maps, 
endowed with the pointwise bracket, as a subgroup of $\Diff(M \times K)$, by 
letting $f \in C^\infty(M,K)$ act on $M \times K$ by 
$\tilde f(m,k) := (m,f(m)k)$. 
The corresponding Lie algebra of vector fields on 
$M \times K$ is $C^\infty(M,\k)$, where 
$\xi \in C^\infty(M,\k)$ corresponds to the vector field given by 
$$\tilde \xi(m,k) = T_\1(\rho_k)\xi(m) \in T_k(K) \subeq T_{(m,k)}(M \times K).$$ 

A map $\phi \: N \to C^\infty(M,K)$ is smooth as a map into $\Diff(M \times K)$ 
if and only if the map 
$$  N \times M \times K \to (M \times K)^2, \quad 
(n,m,k) \mapsto \big((m,\phi(n)(m)k), (m,\phi(n)(m)^{-1}k)\big) $$
is smooth, which in turn means that the map  
$\hat \phi \: N \times M \to K, (n,m) \mapsto \phi(n)(m)$
is smooth. Hence the Uniqueness Lemma also applies 
to  functions $\phi \: N \to C^\infty(M,K)$ which are smooth in the sense that 
$\hat\phi$ is smooth. 
Their logarithmic derivatives $\delta(\phi)$ can be viewed as 
$C^\infty(M,\k)$-valued $1$-forms on~$N$. 

(c) If $G$ is a Lie group, then $G$ itself can be identified with 
the subgroup $\{ \lambda_g \: g \in G \}$ of $\Diff(G)$, consisting of all 
left translations. On the Lie algebra level, this 
corresponds to the embedding $\L(G) \into {\cal V}(G)$ as the right invariant 
vector fields. 
Then a map $\phi \: N \to G \subeq \Diff(G)$ is smooth if and only if it is 
smooth as a $G$-valued map, and we thus obtain a 
Uniqueness Lemma for $G$-valued smooth maps and $\L(G)$-valued $1$-forms. 
\qed
 
\Remark II.3.7. (a) The Uniqueness Lemma implies in particular that a smooth left action of 
a connected Lie group $G$ on a smooth manifold $M$, given by a homomorphism  
$\sigma \: G \to \Diff(M)$, is uniquely determined by 
the corresponding homomorphism of Lie algebras 
$$\dot\sigma := -\delta(\sigma)_\1 \: \L(G) \to {\cal V}(M)$$ 
because $\delta(\sigma)$ is a left invariant ${\cal V}(M)$-valued 
$1$-form on $G$, hence determined by its value in $\1$. 

(b) It likewise follows that any smooth representation 
$\pi \: G \to \GL(E)$ of a connected Lie group $G$ on some locally convex space $E$ 
is uniquely determined by its derived representation 
$$\L(\pi) := \delta(\pi)_\1 \: \L(G) \to \gl(E) \subeq {\cal V}(E).
\qeddis 

\Remark II.3.8. (Complete vector fields) (a) Another 
consequence of the Uniqueness Lemma is that we may define a 
{\it complete vector field $X$ on $M$} as a vector field for which 
there exists a 
smooth one-parameter group $\gamma_X \: \R \to \Diff(M)$
with $\gamma_X'(0) = X$. 
In this sense, we consider the complete vector fields as the 
domain of the exponential
function $\exp(X) := \gamma_X(1)$ 
of $\Diff(M)$. 

(b) Likewise, 
the domain of the exponential function of $\GL(E)$, $E$ a locally convex space, is
the set of all continuous linear operators $D$ on $E$ for which the corresponding
linear vector field $X_D(v) = Dv$ is complete, i.e., there exists a 
smooth representation $\alpha \: \R \to \GL(E)$ with 
$\alpha'(0) = D$. We call these operators $D$ {\it integrable}. 

(c) We may further define for each
Lie group $G$ the domain of the exponential function of $G$ as those
elements $x \in \L(G)$ for which the corresponding left invariant
vector field $x_l$ is complete. 
\qed

\Example II.3.9. (The adjoint representation of a Lie group)  
The adjoint action encodes a good deal of structural information of a Lie group $G$. 
It provides a linearized picture of the non-commutativity of $G$. 

For each $g \in G$, the map 
$c_g \: G \to G, x \mapsto gxg^{-1},$ 
is a smooth automorphism of $G$, hence induces a continuous linear
automorphism 
$$ \Ad(g) := \L(c_g) \: \L(G) \to \L(G). $$
We thus obtain a smooth action 
$ G \times \L(G) \to \L(G), (g, x) \mapsto \Ad(g).x, $ 
called the {\it adjoint action} of $G$ on $\L(G)$.
By considering the Taylor expansion of the map 
$(g,h) \mapsto ghg^{-1}$, one shows that the derived representation 
of $\L(G)$ on $\L(G)$ satisfies 
$$ \L(\Ad) = \ad, \quad \hbox{ i.e.,} \quad \L(\Ad)(x)(y) = [x,y] \quad \hbox{ 
for } \quad x,y \in \L(G). \leqno(2.3.1) $$

If $\L(G)' := {\cal L}(\L(G),\K)$ denotes the {\it topological dual of $\L(G)$}, then we also obtain 
a representation of $G$ 
on $\L(G)'$ by $\Ad^*(g).f := f \circ \Ad(g)^{-1}$, called the 
{\it coadjoint action}. Since we do not endow $\L(G)'$ with a topology, we will not specify 
any smoothness or continuity properties of this action. 
\qed 

The following lemma shows that, whenever there is a smooth curve 
$\gamma \: J \to \Diff(M)$ satisfying the initial value problem 
$$ \gamma(0) = \id_M \quad \hbox{ and } \quad 
\gamma'(t) = X_t \circ \gamma(t) \leqno(2.3.2) $$
for a time-dependent vector field $X \: J \to {\cal V}(M)$, then all integral 
curves of $X$ on $M$ are of the form 
$$\eta(t) = \gamma(t)(m), \leqno(2.3.3) $$
hence unique. It follows in particular, that the existence 
of multiple integral curves of $X$ implies that (2.3.2) has no solution. 
Below we shall see examples where this situation arises, even for linear 
differential equations. 

\Lemma II.3.10. Let $J \subeq \R$ be an interval containing $0$ 
and $\gamma \: J \to \Diff(M)$ 
be a smooth curve with $\gamma(0) = \id_M$. 
Let $X_t := \delta^r(\gamma)_t$ be the corresponding time-dependent 
vector field on $M$ with $X_t \circ \gamma(t) = \gamma'(t)$, 
$m_0 \in M$, and assume that $\eta \: J \to M$ 
is a solution of the initial value problem: 
$$ \eta(0)= m_0 \quad \hbox{ and } \quad 
\eta'(t) = X_t(\eta(t)) \quad \hbox{ for } \quad t \in J. $$
Then $\eta(t) = \gamma(t)(m_0)$ holds for all $t \in J$. 

\Proof. The smooth curve $\alpha \: J \to M, t \mapsto \gamma(t)^{-1}(\eta(t))$ 
satisfies $\alpha(0) = m_0$ and 
$$ \eqalign{ \alpha'(t) 
&= (\gamma^{-1})'(\eta(t)) + T(\gamma(t)^{-1})(\eta'(t)) 
= T(\gamma(t)^{-1})\big(\delta(\gamma^{-1})_t(\eta(t)) + \eta'(t)\big) \cr
&= T(\gamma(t)^{-1})\big(-\delta^r(\gamma)_t(\eta(t)) + \eta'(t)\big) 
= T(\gamma(t)^{-1})\big(-X_t(\eta(t)) + \eta'(t)\big) =0. \cr} $$
Hence $\alpha$ is constant $m_0$, and the assertion follows. 
\qed

In [OMYK82], one finds the particular version of the preceding lemma 
dealing with solutions of the initial value problem 
$$ \eta'(t) = [\eta(t),\xi(t)] + \eta(t), \quad \eta(0) = x $$
in the Lie algebra of a regular Lie group (see also [KYMO85, 2.5/2/6]). 

\Example II.3.11. (A linear ODE with multiple solutions) (cf.\ [Ham82, 5.6.1], [Mil84]) 
We give an example of a linear ODE for which solutions 
to initial value problems exist, but are 
not unique. We consider the Fr\'echet space $E := C^\infty([0,1],\R)$ of smooth 
functions on the closed unit interval, and 
the continuous 
linear operator $Df := f'$ on~$E$. We are asking for solutions of the 
initial value problem 
$$ \dot\gamma(t) = D\gamma(t), \quad \gamma(0) = v_0, \quad 
\gamma \: I \subeq \R \to E. \leqno(2.3.4) $$
As a consequence of E.~Borel's Theorem that each power series is the Taylor 
series of a smooth function, each $v_0 \in E$ 
has an extension to a  smooth function on~$\R$. 
Let $h$ be such a function and consider the curve 
$$ \gamma \: \R \to E, \quad \gamma(t)(x) := h(t + x). $$
Then $\gamma(0) = h\res_{[0,1]} = v_0$ and 
$\dot\gamma(t)(x) = h'(t + x) = \gamma(t)'(x) = (D \gamma(t))(x)$. It is clear that these 
solutions of (2.3.4) depend on the choice of the extension $h$ of $v_0$. 

Lemma II.3.10 and the discussions preceding it now imply that $D$ is not integrable. 
In fact, for any smooth homomorphism 
$\alpha \: \R \to \GL(E)$ with $\alpha'(0) = D$, we would have 
$\delta^r(\alpha) = D$, so that any solution of (2.3.4) is of the form 
$\gamma(t) = \alpha(t).v_0$, contradicting the existence of multiple solutions. 
\qed

\Example II.3.12. (A linear ODE without solutions; [Mil84]) 
We identify $E := C^\infty(\SS^1,\C)$ 
with the space of $2\pi$-periodic smooth functions on the real line. 
We consider the linear operator $Df := -f''$ and the equation (2.3.4), 
which in this case 
is the heat equation with reversed time. 
If $\gamma$ is a solution of (2.3.4) and 
$\gamma(t)(x) = \sum_{n \in \Z} a_n(t) e^{inx}$
its Fourier expansion, then 
$a_n'(t)= n^2 a_n(t)$ for each $n \in \Z$ leads to 
$a_n(t) = a_n(0)e^{tn^2}$. 
If the Fourier coefficients $a_n(0)$ of $\gamma_0$ do not satisfy 
$\sum_n |a_n(0)| e^{\eps n^2} < \infty$
for any $\eps > 0$ (which need not be the case for a smooth function 
$\gamma_0$), then (2.3.4) does not have a solution on $[0,\eps]$. 

As a consequence, the operator $\exp(tD)$ is not defined in $\GL(E)$ for any $t > 0$. 
Nevertheless, we may use the Fourier series expansion to see that 
$\beta(t) := (1 + it^2)\1 + tD$ defines a smooth curve 
$\beta \: \R \to \GL(E)$. 
We further have $\beta'(0) = D$, so that $D$ arises as the 
tangent vector of a smooth curve in $\GL(E)$, but not of any 
smooth one-parameter group. 
\qed

The following example is of some interest for the integrability of 
Lie algebras of formal vector fields (Example~VI.2.8). 

\Example II.3.13. We consider the space 
$E := \R[[x]]$ of formal power series $\sum_{n = 0}^\infty a_n x^n$ in one variable.
We endow it with the Fr\'echet topology for which the map 
$\R^{\N_0} \to \R[[x]], (a_n) \mapsto \sum_n a_n x^n$ is a topological isomorphism. 
Then $Df := f'$ with 
$f'(x) := \sum_{n=1}^\infty a_n n x^{n-1} = \sum_{n = 0}^\infty a_{n+1} (n+1) x^n$ 
for $f(x) = \sum_{n = 0}^\infty a_n x^n$ defines a continuous linear 
operator on $E$. We claim that this operator is not integrable. 

We argue by contradiction, and assume that $\alpha \: \R \to \GL(E)$ 
is a smooth $\R$-action of $E$ with $\alpha'(0) = D$. For each 
$n \in \N$, the curve $\gamma \: \R \to E, \gamma(t) := (x+t)^n$, 
satisfies 
$\dot\gamma(t) = n (x+t)^{n-1} = D\gamma(t)$, so that Lemma~II.3.10 
implies that $\alpha(t)x^n = (x+t)^n$ for all $t \in \R$. 
Then we obtain $\alpha(1)x^n = 1 + n x+ \ldots$. 
In view of $\lim_{n \to \infty} x^n \to 0$ in $E$, this contradicts the 
continuity of the operator $\alpha(1)$. Therefore $D$ is not integrable. 
\qed

\Example II.3.14. Let $M$ be a compact manifold and 
$\g = {\cal V}(M)$, the Lie algebra of smooth vector fields on $M$. 
We now sketch how the group 
$G := \Diff(M)$ can be turned into a Lie group, modeled on~${\cal V}(M)$, 
endowed with its natural Fr\'echet topology (Definition~I.5.2) ([Les67]). 

If $\Fl^X \: \R \times M \to M, (t,m) \mapsto \Fl^X_t(m)$ 
denotes the flow of the vector field $X$, then 
the exponential function of the group $\Diff(M)$ should be given by the time-$1$-map 
of the flow of a vector field: 
$$ \exp_{\Diff(M)} \: {\cal V}(M) \to \Diff(M), \quad 
X \mapsto \Fl^X_1. $$
For the Lie group structure described below,  
this is indeed the case. Unfortunately, 
it is not a local diffeomorphism 
of a $0$-neighborhood in ${\cal V}(M)$ onto any identity neighborhood in $\Diff(M)$. 
Therefore we cannot 
use it to define a chart around $\1 = \id_M$ 
(cf.\ [Grab88], [Pali68/74], and also [Fre68], which deals with local smooth diffeomorphisms 
in two dimensions). 

Fortunately, 
there is an easy way around this problem. Let $g$ be a Riemannian metric on $M$ and 
$\Exp \: TM \to M$
be its exponential function, which assigns to $v \in T_m(M)$ the point 
$\gamma_v(1)$, where $\gamma_v \: [0,1] \to M$ is the geodesic segment with 
$\gamma_v(0) = m$ and $\gamma_v'(0) = v$. We then obtain a smooth map 
$$ \Phi \: TM \to M \times M, \quad v \mapsto (m, \Exp v), \quad v \in T_m(M). $$
There exists an open neighborhood $U \subeq TM$ of the zero section such that 
$\Phi$ maps $U$ diffeomorphically onto an open neighborhood of the diagonal in $M \times M$.
Now 
$$ U_\g := \{ X \in {\cal V}(M) \: X(M) \subeq U\} $$
is an open subset of the Fr\'echet space ${\cal V}(M)$, and we define a map 
$$ \phi \: U_\g \to C^\infty(M,M), \quad \phi(X)(m) := \Exp(X(m)). $$
It is clear that $\phi(0) = \id_M$. 
One can show that after shrinking $U_\g$ to a sufficiently small 
$0$-neighborhood in the compact open $C^1$-topology on ${\cal V}(M)$, 
we achieve that $\phi(U_\g) \subeq \Diff(M)$. 
To see that $\Diff(M)$ carries a Lie group 
structure for which $\phi$ is a chart, one has to verify that the 
group operations are smooth in a $0$-neighborhood 
when transferred to $U_\g$ via $\phi$, so that Theorem~II.2.1 applies. 
We thus obtain a Lie group structure on $\Diff(M)$ (cf.\ [Omo70], [GN06]). 
    
From the smoothness of the map 
$U_\g \times M \to M, (X,m) \mapsto \phi(X)(m) = \Exp(X(m))$ it follows 
that the canonical left action $\sigma \: \Diff(M) \times M \to M, 
(\phi,m) \mapsto \phi(m)$ is smooth in an identity neighborhood of 
$\Diff(M)$, and hence smooth, because it is an action by smooth maps. 
The corresponding homomorphism of Lie algebras 
$\dot\sigma \: \L(\Diff(M)) \to {\cal V}(M)$ (Remark~II.3.7(a))  
is given by 
$$ \dot\sigma(X)(m) = -T\sigma(X,0_m)  = -X(m), $$
i.e., $\dot\sigma = -\id_{{\cal V}(M)}$, which leads to 
$$ \L(\Diff(M)) = ({\cal V}(M),[\cdot,\cdot])^{\rm op}, $$
where $\g^{\rm op}$ is the {\it opposite} of the  Lie algebra $\g$ with the bracket 
$[x,y]_{\rm op} := [y,x]$. 

This ``wrong'' sign is caused by the fact that we consider 
$\Diff(M)$ as a group acting on $M$ from the left. If we consider 
it as a group acting on the right, we obtain the opposite multiplication 
$\phi * \psi := \psi \circ \phi$ 
and 
$\L(\Diff(M)^{\rm op}) \cong  ({\cal V}(M),[\cdot,\cdot])$. 
Here we write $G^{\rm op}$ for the {\it opposite group} with the order of 
multiplication reversed. 
\qed

\subheadline{II.4. Applications to Lie group-valued smooth maps} 

In this subsection, we describe some applications of the Uniqueness Lemma 
to Lie group-valued smooth maps (cf.\ Remark~II.3.6(c)).  

Let $G$ be a Lie group with Lie algebra $\g = \L(G)$. 
The {\it Maurer--Cartan form} $\kappa_G \in \Omega^1(G,\g)$ is the unique left invariant 
$1$-form on $G$ with $\kappa_{G,\1} = \id_\g$, i.e., 
$\kappa_G(v) = g^{-1}.v$ for $v \in T_g(G)$. In various disguises, 
this form plays a central role in the approach to (local) (Banach--)Lie groups 
via partial differential equations ([Mau88], [CaE01], 
[Lie95], [Bir38], [MicA48], [Lau56]). 

Identifying $G$ with the subgroup of left translations in $\Diff(G)$, 
the concepts of the preceding subsection apply to any smooth map $f \: M \to G$ 
(Remark~II.3.6(c)). The logarithmic derivative of $f$ can be described as a 
pull-back of the Maurer--Cartan form: 
$$ \delta(f)  = f^*\kappa_G \in \Omega^1(M,\g). $$

\Proposition II.4.1. Let $G$ and $H$ be Lie groups. 
\litem{(1)} If $\phi \: G \to H$ is a morphism of Lie groups, then 
$\delta(\phi) = \L(\phi) \circ \kappa_G.$ 
\litem{(2)} If $G$ is connected and $\phi_1, \phi_2 : G \to H$ are
morphisms of Lie groups with $\L(\phi_1) = \L(\phi_2)$, then 
$\phi_1 = \phi_2$. 
\litem{(3)} For a smooth function 
$f \: G \to H$ with $f(\1) = \1$, the following are equivalent: 
\litemitem{(a)} $\delta(f)$ is a left invariant $1$-form. 
\litemitem{(b)} $f$ is a group homomorphism.

\Proof. (1) is a simple computation, and (2) follows with (1) and the 
Uniqueness Lemma (cf.\ Remark~II.3.6(c)). 

The proof of (3) follows a similar pattern, applying the Uniqueness Lemma 
to the relations $\lambda_g^*\delta(f) = \delta(f)$. 
\qed

Applying (2) to the conjugation automorphisms $c_g \in \Aut(G)$, we obtain: 
\Corollary II.4.2. If $G$ is a connected Lie group, then 
$\ker \Ad = Z(G)$. 
\qed

It follows in particular, that the adjoint action of a connected Lie group 
$G$ is trivial if and only if $G$ is abelian. In view of Remark~II.3.6(b), 
this is equivalent to the triviality of the corresponding derived action, 
which is the adjoint action of $\L(G)$ (Example II.3.9). We thus obtain the following 
affirmative answer to a question of {\smc J.~Milnor} ([Mil84]): 

\Proposition II.4.3. A connected Lie group is abelian if and only 
if its Lie algebra is abelian. 
\qed

This argument can be refined by investigating the structure of 
logarithmic derivatives of iterated commutators of smooth curves in 
a Lie group $G$. A systematic use of the Uniqueness Lemma then leads to the following result 
(see [GN06]): 

\Theorem II.4.4. A connected Lie group $G$ is nilpotent, resp., solvable, if 
and only if its Lie algebra $\L(G)$ is nilpotent, resp., solvable. 
\qed

\subheadline{II.5. The exponential function and regularity} 

In the Lie theory of finite-dimensional and Banach--Lie groups, the 
exponential function is a central tool used to pass information from 
the group to the Lie algebra and vice versa. Unfortunately, the exponential 
function is less powerful in the context of locally convex Lie groups. 
Here we take a closer look at its basic properties, and in Section IV below 
we study the class of locally exponential Lie groups for which the exponential 
function behaves well in the sense that it is a local diffeomorphism in~$0$.

\Definition II.5.1. For a Lie group $G$ with Lie algebra $\g = \L(G)$, we call a 
smooth function $\exp_G \: \g \to G$ an 
{\it exponential function for $G$} if for each $x \in \g$ the curve 
$\gamma_x(t) := \exp_G(tx)$ is a one-parameter group 
with $\gamma_x'(0) = x$. 
\qed

It is easy to see that any such curve is a solution of the 
initial value problem (IVP)  
$$ \gamma(0) = \1, \quad \delta(\gamma) = x,$$
so that the Uniqueness Lemma implies that solutions are unique whenever 
they exist. Hence a Lie group $G$ has at most one 
exponential function. 

The question for the existence of an exponential function 
leads to the more general question when for a smooth curve 
$\xi \in C^\infty(I,\g)$ ($I = [0,1]$), the initial value problem (IVP) 
$$ \gamma(0) = \1, \quad \delta(\gamma) = \xi, \leqno(2.5.1) $$
has a solution. If this is the case for constant functions $\xi(t) = x$, 
the corresponding solutions are the curves $\gamma_x$ required to obtain an 
exponential function. The solutions of (2.5.1) are unique by the 
Uniqueness Lemma (Remark~II.3.6(c)). 

\Definition II.5.2. A Lie group $G$ is called {\it regular} if for each
$\xi \in C^\infty(I,\g)$, the initial value problem (2.5.1) 
has a solution $\gamma_\xi \in C^\infty(I,G)$, and the
{\it evolution map} 
$$ \evol_G \: C^\infty(I,\g) \to G, \quad \xi \mapsto \gamma_\xi(1) $$
is smooth. 
\qed

\Remark II.5.3. (a) If $G$ is regular, then $G$ has a smooth exponential 
function, given by 
$$ \exp_G(x) := \evol_G(\xi_x), $$
where $\xi_x(t) = x$ for $t \in I$. 

(b) For any Lie group $G$, the logarithmic derivative 
$$\delta \: C^\infty_*([0,1],G) \to C^\infty(I,\L(G)) \cong \Omega^1(I,\L(G))$$ 
is a smooth map with 
$T_\1(\delta)\xi = \xi'$. If $G$ is regular, this fact can be used to show that 
$T_\1(\delta)$ is surjective, hence that $\L(G)$ is Mackey complete 
(cf.\ [GN06]). 
\qed

\Remark II.5.4. As a direct consequence of the existence of solutions to 
ordinary differential equations on open domains of Banach spaces and 
their smooth dependence on parameters (cf.\ [La99]), every Banach--Lie group is regular. 
\qed

All Lie groups known to the author which are modeled on Mackey complete spaces are regular. 
In concrete situations, it is sometimes hard to verify regularity, and in some case 
it is not known if the Lie groups under consideration are regular.  
We shall take a closer look at criteria for regularity in 
Section III below. In particular, we shall see that essentially all groups 
belonging to the major classes discussed in the introduction are in fact regular. 

\Example II.5.5. If the model space is no longer assumed to be Mackey complete, 
one can construct non-regular Lie groups as follows (cf.~[Gl02b, Sect.~7]): 
Let $A \subeq C([0,1],\R)$ denote the unital subalgebra of all rational functions, 
i.e., of all quotients $p(x)/q(x)$, where $q(x)$ is a polynomial without zeros in 
$[0,1]$. We endow $A$ with the induced norm $\|f\| := \sup_{0 \leq t \leq 1} |f(t)|$.  
If an element $f \in A$ is invertible in $C([0,1],\R)$, then it has no zero 
in $[0,1]$, which implies that it is also invertible in $A$, i.e., 
$$ A^\times = C([0,1],\R)^\times \cap A. $$
This shows that $A^\times$ is open in $A$, and since the Banach algebra 
$C([0,1],\R)$ is a CIA, the continuity of the inversion  is inherited by 
$A$, so that $A$ is a CIA. In particular, $A^\times$ is a Lie group (Example II.1.4). 

Let $f \in A$ and assume that there exists a smooth homomorphism 
$\gamma_f \: \R \to A^\times$ 
with $\gamma_f'(0) =f$. Then 
Proposition II.4.1, applied to $\gamma_f$ as a map $\R \to C([0,1],\R)^\times$, 
leads to $\gamma_f(t) = e^{tf}$ for each $t \in \R$. 
Since $e^f$ is not rational if $f$ is not constant, we conclude that 
$f$ is constant. Therefore the Lie group $A^\times$ does not have an 
exponential function and in particular it is not regular. 
\qed

The following proposition illustrates the relation between regularity and Mackey 
completeness. 

\Proposition II.5.6. The additive Lie group $(E,+)$ of a locally convex space 
$E$ is regular if and only if $E$ is Mackey complete. 

\Proof. For a smooth curve $\xi \: I \to E$, any solution 
$\gamma_\xi \: I \to E$ of (2.5.1) satisfies $\gamma_\xi' = \xi$ and vice versa. 
Therefore regularity implies that $E$ is Mackey complete (Definition~I.1.4). 
Conversely, Mackey completeness of $E$ implies that 
$\evol_G(\xi) := \int_0^1 \xi(s)\, ds$
defines a continuous linear map 
$\evol_G \: C^\infty(I,E) \to E$, so that it is in particular smooth. 
\qed

\Proposition II.5.7. Suppose that the Lie group $G$ has a
smooth exponential function $\exp_G \: \g \to G$. Then its 
logarithmic derivative is given by 
$$ \delta(\exp_G)(x) = \int_0^1 \Ad(\exp_G(-tx))\, dt, \leqno(2.5.2) $$
where the operator-valued integral is defined pointwise, i.e., 
$$ \delta(\exp_G)(x)y = \int_0^1 \Ad(\exp_G(-tx))y\, dt 
\quad \hbox{ for each } y \in \g. $$

\Proof. ([Grab93]) For $t,s \in \R$, we consider the three smooth functions 
$f, f_t,f_s  \: \g \to G,$ given by 
$$ f(x) := \exp_G((t+s)x), \quad 
f_t(x) := \exp_G(tx) \quad \hbox{and} \quad f_s(x):= \exp_G(sx),  $$
satisfying $f = f_t f_s$ pointwise on $\g$. The 
Product Rule (Lemma~II.3.3) implies that 
$$ \delta(f) = \delta(f_s) + \Ad(f_s)^{-1} \delta(f_t). $$

For the smooth curve $\psi \: \R \to \g, \psi(t) := \delta(\exp_G)_{tx}(ty)$,  
we therefore obtain 
$$ \psi(t+s) = \delta(f)_x(y) 
= \delta(f_s)_x(y) + \Ad(f_s)^{-1}.\delta(f_t)_x(y) 
= \psi(s) + \Ad(\exp_G(-sx)).\psi(t). \leqno(2.5.3) $$
We have $\psi(0) = 0$ and 
$\psi'(0) = \lim_{t \to 0} \delta(\exp_G)_{tx}(y) = \delta(\exp_G)_0(y) = y, $ 
so that taking derivatives with respect to $t$ in $0$, (2.5.3) leads to 
$\psi'(s) = \Ad(\exp_G(-sx)).y.$
Now the assertion follows by integration from 
$\delta(\exp_G)_x(y) = \psi(1) = \int_0^1 \psi'(s)\, ds$. 
\qed

If $\g$ is integrable to a group with exponential function, 
then the one-parameter groups 
$\Ad(\exp_G(tx))$ have the infinitesimal generator 
$\ad x$ (Remark~II.3.7(b), Example~II.3.9), so that we may also write 
$$ \Ad(\exp_G(tx)) = e^{t\ad x}. \leqno(2.5.4)  $$
If, in addition, $\g$ is Mackey complete, then the operator-valued integral 
$$ \kappa_\g(x) := \int_0^1 e^{-t\ad x}\, dt \leqno(2.5.5) $$
exists pointwise because the curves $t \mapsto e^{-t\ad x}y$ are smooth, 
and the preceding theorem states that  for each $x \in \g$: 
$$\delta(\exp_G)_x = \kappa_\g(x). \leqno(2.5.6) $$ 
The 
advantage of $\kappa_\g(x)$ is that it is expressed completely in 
Lie algebraic terms. 

\Remark II.5.8. If $\g$ is a Banach--Lie algebra, then $\kappa_\g(x)$ can be represented by 
a convergent power series 
$$ \kappa_\g(x) 
= \int_0^1 e^{-t\ad x}\, dt 
= {\1 - e^{-\ad x} \over \ad x} 
= \sum_{k = 0}^\infty {(-1)^k\over (k+1)!} (\ad x)^k. $$
This means that $\kappa_\g(x) = f(\ad x)$ holds for the entire function 
$$ f(z) :=  \sum_{k = 0}^\infty {(-1)^k\over (k+1)!} z^k = {1 - e^{-z} \over z}. $$
As $ f^{-1}(0) = 2\pi i \Z \setminus \{0\}$, 
and $\Spec(\kappa_\g(x)) = f(\Spec(\ad x))$
by the Spectral Mapping Theorem, 
we see that $\kappa_\g(x)$ is invertible if and only if 
$\Spec(\ad x) \cap 2 \pi i \Z \subeq \{0\}.$

Part of this observation can be saved in the general case. 
If $\g$ is Mackey complete one can show that 
$\kappa_\g(x)$ is not injective if and only if there exists 
some $n \in \N$ with 
$$  \ker((\ad x)^2 + 4 \pi^2 n^2\1) \not= \{0\}.$$ 
If $\g$ is a complex Lie algebra, this means that some  
$2\pi i n \in 2 \pi i \Z \setminus \{0\}$ is an eigenvalue of $\ad x$ 
(cf.\ [GN06] for details). 
\qed

\Examples II.5.9. Let $\alpha \: \R \to \GL(E)$ be a smooth 
representation of $\R$ on the Mackey complete locally convex space $E$ with the 
infinitesimal generator $D = \alpha'(0)$. 
Then the semi-direct product group 
$$ G := E \rtimes_\alpha \R, \quad (v,t) (v',t') = (v + \alpha(t)v', t + t') $$
is a Lie group with Lie algebra $\g = E \rtimes_D \R$ and exponential function 
$$ \exp_G(v,t) 
= \big( \beta(t)v,  t\big) \quad \hbox{ with } \quad 
\beta(t) =  \int_0^1 \alpha(st)\, ds = \cases{ 
\id_E & for $t = 0$ \cr 
{1 \over t} \int_0^t \alpha(s)\, ds & for $t \not=0$. \cr} $$
From this formula it is clear that $(w,t) \in \im(\exp_G)$ is equivalent to 
$w \in \im(\beta(t))$. We conclude that 
$\exp_G$ is injective on some $0$-neighborhood 
if and only if $\beta(t)$ is injective for $t$ close to $0$, 
and it is surjective onto some $\1$-neighborhood in $G$ if 
and only if $\beta(t)$ is surjective for $t$ close to $0$ (cf.\ Problem IV.4 below). 

Note that the eigenvector equation $Dv = \lambda v$ for $t\lambda \not=0$ implies that 
$$ \beta(t)v = \int_0^1 e^{st\lambda}v\, ds  = {e^{t\lambda} - 1\over t\lambda}v, $$
so that $\beta(t)v = 0$ is equivalent to $t\lambda \in 2 \pi i \Z \setminus \{0\}$. 

(a) For the Fr\'echet space $E = \C^\N$ 
and the diagonal operator $D$ given by 
$D(z_n) = (2\pi in z_n)$, we see that 
$\beta({1\over n}) e_n = 0$ holds for $e_n = (\delta_{mn})_{m \in \N}$, 
and $e_n \not\in \im\big(\beta({1\over n})\big)$. 
We conclude that $(e_n, {1\over n})$ is not contained in the image of $\exp_G$, and since 
$(e_n, {1\over n})\to (0,0)$, the identity of $G$, 
$\im(\exp_G)$ does not contain any identity neighborhood of $G$. 
Hence the exponential function of the Fr\'echet--Lie group 
$G = E \rtimes_\alpha \R$ is neither locally injective nor locally surjective in $0$. 

(b) For the Fr\'echet space $E = \R^\N$ 
and the diagonal operator $D$ given by 
$D(z_n) = (n z_n)$, it is easy to see that all operators $\beta(t)$ are invertible 
and that $\beta \: \R \to \GL(E)$ is a smooth map. This implies that 
$\exp_G \: \g \to G$ is a diffeomorphism. 
\qed

\Remark II.5.10. If $G$ is the unit group $A^\times$ of a Mackey complete CIA, then 
we identify $T(G) \subeq T(A) \cong A \times A$ with $A^\times \times A$ and note that 
$\ad x = \lambda_x - \rho_x$ and $e^{\ad x}y = e^x ye^{-x}.$
Therefore (2.5.2) can be written as 
$$ T_x(\exp_G)y 
= e^x \int_0^1 e^{-tx} y e^{tx}\, dt 
= \int_0^1 e^{(1-t)x} y e^{tx}\, dt $$
(cf.\ [MicA45] for the case of Banach algebras). 
\qed

A closer investigation of (2.5.5) leads to the 
following results on the behavior of the exponential function 
(cf.\ [GN06]; and [LaTi66] for the finite-dimensional case): 

\Proposition II.5.11. Let $G$ be a Lie group with Lie algebra $\g$ 
and a smooth exponential function. Then the following assertions hold for 
$x,y \in \g$: 
\litem{(1)} If $\kappa_\g(x)y = 0$, 
then 
$$ \exp_G(e^{t\ad y}.x) = \exp_G(x) \quad \hbox{ for all } \quad t \in \R.$$
\litem{(2)} If $\kappa_\g(x)$ is not injective and $\g$ is Mackey complete, 
then $\exp_G$ is not injective in any neighborhood of $x$. 
\litem{(3)} If $\kappa_\g(x)$ is injective, then 
\litemitem{(a)} $\exp_G(y) = \exp_G(x)$ implies $[x,y]=0$ and $\exp_G(x-y) = \1$. 
\litemitem{(b)} $\exp_G(x) \in Z(G)$ implies $x \in \z(\g)$ and equivalence holds 
if $G$ is connected. 
\litemitem{(c)} $\exp_G(x) = \1$ implies $x \in \z(\g)$. 
\litem{(4)} Suppose that $0$ is isolated in $\exp_G^{-1}(\1)$. 
Then $x$ is isolated in $\exp_G^{-1}(\exp_G(x))$ 
if and only if $\kappa_\g(x)$ is injective. 
\litem{(5)} If $\a \subeq \g$ is an abelian subalgebra, then 
$\exp_\a := \exp_G \res_\a \: \a \to G$ is a morphism of Lie groups. 
Its kernel $\Gamma_\a := \ker(\exp_\a)$ is a closed subgroup of $\a$ 
in which all $C^1$-curves are constant. It intersects each finite-dimensional 
subspace of $\a$ in a discrete subgroup. 
\qed

\Remark II.5.12. (a) Let $U \subeq \g$ be a $0$-neighborhood with the property that 
$\kappa_\g(z)$ is injective for each $z \in U-U$. 
Then the preceding proposition implies for $x,y \in U$ with 
$\exp_G x = \exp_G y$ that $\exp_G(x-y) = \1$, $[x,y] = 0$, and since $x - y \in U$, 
it further follows that $x - y \in \z(\g)$. 
If we assume, in addition, that the closed subgroup 
$\Gamma_{\z(\g)} := \ker(\exp_{\z(\g)})$ intersects 
$U - U$ only in~$\{0\}$, $\exp_G$ is injective on $U$. 

(b) If $G$ is a Banach--Lie group and 
$\g = \L(G)$ carries a norm with $\|[x,y]\| \leq \|x\|\cdot\|y\|$, 
then $\| \ad x\| \leq \|x\|$. Therefore $\|x\| < 2\pi$ implies 
that $\kappa_\g(x)$ is invertible (Remark~II.5.8). If $\exp_G\res_{\z(\g)}$ is injective, 
i.e., $Z(G)$ is simply connected, the preceding remark implies 
that $\exp_G$ is injective on the open ball $B_\pi := \{ x \in \g \: \|x\| < \pi\}$ 
(cf.\ [LaTi66]). In general, we may put 
$$ \delta_G := \inf \{ \|x\|\: 0 \not= x \in \Gamma_{\z(\g)}\} $$
to see that $\exp_G$ is injective on the ball of radius 
$r := \min \{ \pi, {\textstyle{\delta_G\over 2}}\}$
(cf.\ [GN03, [Bel04, Rem. 2.3]). 
\qed

\Example II.5.13. In [Omo70], [Ham82] and [Mil84], it is shown that for the group 
$G := \Diff(\SS^1)$ of diffeomorphisms of the circle, the image of the exponential 
function is not a neighborhood of $\1$ (cf.\ also [KM97, Ex. 43.2], [PS86, p.~28]). 
Small perturbations of  
rigid rotations of order $n$ lead to a sequence of diffeomorphisms converging 
to $\id_{\SS^1}$ which do not lie on any one-parameter group.  

More generally, for any compact manifold $M$, 
the image of the exponential function of $\Diff(M)$ does not contain any 
identity neighborhood (cf.\ [Grab88], [Pali68/74], and [Fre68] for some $2$-dimensional cases). 

Identifying the Lie algebra $\g := {\cal V}(\SS^1)$ of $\Diff(\SS^1)$ 
with smooth $2\pi$-periodic functions on~$\R$, the Lie bracket corresponds to 
$$ [f,g] = fg'-f'g. $$
For the constant function $f_0 = 1$ and 
$c_n(t) := \cos(nt)$ and $s_n(t) = \sin(nt)$, this leads to 
$$ [f_0, s_n] = n c_n \quad \hbox{ and }  \quad [f_0, c_n] = -n s_n, $$
so that $\span\{f_0, s_n, c_n\} \subeq {\cal V}(\SS^1)$ is a Lie subalgebra 
isomorphic to $\sL_2(\R)$. It further follows that 
$((\ad f_0)^2 + n^2\1) s_n = 0$, so that 
$\kappa_\g({2\pi\over n}f_0)s_n = 0$ implies that $\exp_G$ is 
not injective in any neighborhood of ${2\pi\over n}f_0$ (Proposition~II.5.11(1)) 
(cf.\ [Mil82, Ex. 6.6]). 
Therefore $\exp_G$ is neither locally surjective nor injective. 
\qed

\Remark II.5.14. (Surjectivity of $\exp_G$) The global behavior of the exponential function and in 
particular the question of its surjectivity is a quite complicated issue, depending 
very much on specific properties of the groups under consideration  
(cf.\ [W\"u03/05]). 

(a) For finite-dimensional Lie groups, the most basic general result is that 
if $G$ is a connected Lie group with compact Lie algebra $\g$, then 
$\exp_G$ is surjective. Since the compactness of $\g$ is equivalent to the 
existence of an $\Ad(G)$-invariant scalar product, which in turn leads to a 
biinvariant Riemannian metric on $G$, the surjectivity of $\exp_G$ can be 
derived from the Hopf--Rinow Theorem in Riemannian geometry. 

(b) A natural generalization of the notion of a compact Lie algebra to 
the Banach context is to say that a real Banach--Lie algebra $(\g,\|\cdot\|)$ is 
{\it elliptic} if the norm on $\g$ is 
invariant under the group $\Inn(\g) := \la e^{\ad \g}\ra \subeq \Aut(\g)$ 
of {\it inner automorphisms} (cf.\ [Ne02c, Def.~IV.3]). 
A finite-dimensional Lie algebra $\g$ is elliptic with respect to some
norm if and only if it is compact. 
In this case, the requirement of an invariant scalar product
leads to the same class of Lie algebras, but in the
infinite-dimensional context this is different. Here the existence of 
an invariant scalar product turning $\g$ into a real Hilbert space
leads to the structure of an $L^*$-algebra. Simple $L^*$-algebras can be
classified, and each $L^*$-algebra is a Hilbert space direct sum of
simple ideals and its center (cf.\ [Sc60/61], 
[dlH72], [CGM90], [Neh93], [St99]). In particular, the
classification shows that every $L^*$-algebra can be
realized as a closed subalgebra of the $L^*$-algebra $B_2(H)$ of
Hilbert--Schmidt operators on a complex Hilbert space $H$. Therefore 
the requirement of an invariant scalar product on $\g$ leads to the
embeddability into the Lie algebra $\uu_2(H)$ of skew-hermitian
Hilbert--Schmidt operators on a Hilbert space $H$. 

The class of elliptic Lie algebras is much bigger. It contains
the algebra $\uu(A)$ of skew-hermitian elements of any $C^*$-algebra
$A$ and in particular the Lie algebra $\uu(H)$ of the full unitary
group $\UU(H)$ of a Hilbert space $H$. 

Although finite-dimensional connected 
Lie groups with compact Lie algebra have a surjective exponential
function, this is no longer true for 
connected Banach--Lie groups with elliptic Lie
algebra. This is a quite remarkable phenomenon discovered  by {\smc Putnam} and
{\smc Winter} in [PW52]: the orthogonal group $\OO(H)$ of a real 
infinite-dimensional Hilbert 
space is a connected Banach--Lie group with elliptic Lie algebra, but
its exponential function is {\sl not} surjective. This contrasts the fact that 
the exponential function of the unitary group $\UU(H)$ of a complex Hilbert space 
is always surjective, as follows from the spectral theory of unitary operators.
\qed

\subheadline{II.6. Initial Lie subgroups} 

It is one of the fundamental problems of Lie theory (FP5) to understand to which extent 
subgroups of Lie groups carry natural Lie group structures. In this subsection,  
we briefly discuss the rather weak concept of an initial Lie subgroup. As a consequence 
of the universal property built into its definition, such a structure is unique 
whenever it exists. 
As the discussion in Remark II.6.5 and further results in Section IV 
below show, it is hard to prove that a subgroup does 
not carry any initial Lie group structure (cf.\ Problem~II.6). 

\Definition II.6.1. An injective morphism $\iota \: H \to G$ of Lie groups 
is called {\it an initial Lie subgroup} 
if $\L(\iota) \: \L(H) \to \L(G)$ is injective, and for each 
$C^k$-map  
$f \: M \to G$ ($k \in \N \cup \{ \infty\}$) from a $C^k$-manifold $M$ 
to $G$ with $\im(f) \subeq H$, the corresponding map 
$\iota^{-1} \circ f \: M \to H$ is $C^k$. 
\qed

The following lemma shows that the existence of an initial Lie group structure 
only depends on the subgroup $H$, considered as a subset of $G$. 

\Lemma II.6.2. Any subgroup $H$ of a Lie group $G$ carries at most one structure 
of an initial Lie subgroup.

\Proof. If $\iota' \: H' \into G$ is another initial Lie subgroup with the same range 
as $\iota \: H \to G$, then 
$\iota^{-1} \circ \iota' \: H' \to H$ and 
$\iota'^{-1} \circ \iota \: H \to H'$ are smooth 
morphisms of Lie groups, so that $H$ and $H'$ are isomorphic. 
\qed

A priori, any subgroup $H$ 
of a Lie group $G$ can be an initial Lie subgroup. 
A first step to a better understanding of initial subgroups is 
to find a natural candidate for the Lie algebra 
of such a subgroup. In the following, we write $C^1_*(I,G)$ for the set of 
all $C^1$-curves $\gamma \: I = [0,1] \to G$ with $\gamma(0) = \1$. 
Then the following definition works well for all subgroups 
(cf.\ [Lau56]; see also [vN29; pp.\ 18/19]):

\Proposition II.6.3. Let $H \subeq G$ be a subgroup of the Lie group $G$. Then the 
differential tangent set 
$$ \L^d(H) := \{ \alpha'(0) \in \L(G) = T_\1(G) \: \alpha \in C^1_*([0,1],G),\ \im(\alpha) \subeq H\} $$
is a Lie subalgebra of $\L(G)$. If, in addition, $H$ carries the structure 
$\iota_H \: H \to G$ of an initial Lie subgroup, then $\L^d(H) = \im(\L(\iota_H))$. 

\Proof. If $\alpha,\beta \in C^1_*(I,G)$, then 
$(\alpha\beta)'(0) = \alpha'(0) + \beta'(0)$, $(\alpha^{-1})'(0) = -\alpha(0),$
and for $0 \leq \lambda \leq 1$ the curve $\alpha_\lambda(t) := \alpha(\lambda t)$ 
satisfies $\alpha_\lambda'(0) = \lambda \alpha'(0)$. This implies that 
$\L^d(H)$ is a real linear subspace of $\L(G)$. 

Next we recall that $[x,y]$ is the lowest order term in the Taylor expansion of 
the commutator map $(x,y) \mapsto xyx^{-1}y^{-1}$ 
in any local chart around $\1$ (Remark~II.1.8). This implies that the curve 
$$\gamma(t) := \alpha(\sqrt t)\beta(\sqrt t)\alpha(\sqrt t)^{-1}\beta(\sqrt t)^{-1} $$
with $\gamma(0) = \1$ is $C^1$ with $\gamma'(0) 
= [\alpha'(0), \beta'(0)]$.\footnote{$^1$}{\eightrm Note that in general this curve is not 
twice differentiable.} 
We conclude that $\L^d(H)$ is a Lie subalgebra of $\L(G)$. 

If, in addition, $H$ is initial, then 
$C^1_*([0,1],H) = \{ \alpha \in C^1_*([0,1],G) \: \im(\alpha) \subeq H\}$ 
implies that $\L(\iota_H)(\L(H)) = \L^d(H)$. 
\qed

We put the superscript $d$ (for differentiable) 
to distinguish $\L^d(H)$ from the Lie algebra of a Lie group. Later, we shall encounter 
another approach to the Lie algebra of a subgroup which works well for closed 
subgroups of locally exponential Lie groups. 

\Remark II.6.4. Let $\alpha \in C^1_*([0,1],G)$ and $H \subeq G$ be a subgroup. 
If $\im(\alpha) \subeq H$, then the image of the continuous curve 
$\delta(\alpha) \in C([0,1],\L(G))$ is contained in $\L^d(H)$. 
If, conversely, $\im(\delta(\alpha)) \subeq \L^d(H)$, then it is not clear why 
this should imply that $\im(\alpha) \subeq H$. We shall see below that 
the concept of regularity helps to deal with this problem. 
\qed

\Remark II.6.5. The following facts demonstrate that it is not easy to find 
subgroups with no initial Lie subgroup structure ([Ne05, Lemma~I.7]): 

(a) Let $H \subeq G$ be a subgroup such that all $C^1$-arcs in $H$ are constant. Then 
the discrete topology defines on $H$ an initial Lie subgroup structure. 

(b) If $\dim G < \infty$, then any subgroup $H \subeq G$ carries an initial Lie group 
structure: According to Yamabe's 
Theorem ([Go69]), the arc-component $H_a$ of $G$ 
is of the form $\la \exp \h \ra$ for some Lie 
subalgebra $\h \subeq \L(G)$, which can be identified with $\L^d(H)$. 
To obtain the initial Lie group structure on $H$, we endow $H_a$ with its intrinsic 
Lie group structure and extend it with Corollary II.2.3 to all of $H$. 

(c) If the connected Lie group $G$ has a smooth exponential function, 
the center $\z(\g)$ of $\g = \L(G)$ is Mackey complete, and 
the subgroup $\Gamma_{\z(\g)} := \exp_G^{-1}(\1) \cap \z(\g)$ 
is discrete, then $Z(G)$ carries an initial Lie group structure with 
Lie algebra $\z(\g).$

We endow $\exp_G(\z(\g)) \cong \z(\g)/\Gamma_{\z(\g)}$ with the quotient 
Lie group structure (Corollary~II.2.4) and use Corollary~II.2.3 
to extend it to all of $Z(G)$.  
\qed

\Remark II.6.6. If $H \subeq G$ is a Lie subgroup in the sense of Remark~II.2.5(b), 
then some identity neighborhood of $H$ is a submanifold of $G$ and its intrinsic 
Lie group structure turns $H$ into an initial Lie subgroup of $G$. 
\qed

\subheadline{Open Problems for Section II}

\Problem II.1. Show that every Lie group $G$ modeled on a Mackey complete 
locally convex space has a smooth exponential function, or find a counterexample 
(cf.\ Example~II.5.5). 
\qed

The following assertion is even stronger: 
\Problem II.2. ([Mil84]) Show that every Lie group $G$ modeled on a Mackey complete 
locally convex space is regular, or find a counterexample. 

The assumption of Mackey completeness of $\L(G)$ is necessary because for any 
regular Lie group the 
differential of the evolution map $\evol_G \: C^\infty([0,1],\g) \to G$ 
is given by 
$$ T_0(\evol_G) \xi = \int_0^1 \xi(t)\, dt. $$
Therefore the regularity of $G$ implies the Mackey completeness of $\L(G)$ 
(cf.\ Proposition~II.5.6). 
\qed

\Problem II.3. ([Mil84]) Show that two $1$-connected Lie groups $G$ with isomorphic 
Lie algebras are isomorphic. For groups with Mackey complete Lie algebras, this 
would follow from Theorem~III.1.5 and a positive solution to Problem~II.2. 
\qed

\Problem II.4. Prove or disprove the following claims 
for all Lie groups $G$ with a smooth 
exponential function $\exp_G \: \g = \L(G) \to G$: 
\litem{(1)} $0$ is isolated in $\exp_G^{-1}(\1)$. 
\litem{(2)} $0$ is isolated in $\Gamma_{\z(\g)} := \exp_G^{-1}(\1) \cap \z(\g)$, where 
$\z(\g)$ denotes the center of $\g$. 

In view of Remark~II.6.5(c), a solution of (2) would be of particular interest 
to classify classes of extensions of Lie groups by non-abelian Lie 
groups (cf.\ Theorem~V.1.5). Note that (2) is equivalent to the discreteness of the group 
$\Gamma_{\z(\g)}$. We know that all $C^1$-curves in this closed subgroup 
of $\z(\g)$ are constant and that all intersections with finite-dimensional 
subspaces are discrete (Proposition~II.5.11(5)). 
\qed

\Problem II.5. (Small Torsion Subgroup Problem; (FP8)) 
Show that for any Lie group $G$ there exists an 
identity neighborhood $U$ such that $\1$ is the only element of finite order 
generating a subgroup lying in $U$. 

If the answer to Problem II.4(1) is negative for some Lie group 
$G$, then each identity neighborhood contains the range of a 
homomorphism $\T \cong \R/\Z \to G$ obtained by $\exp_G(\R x)$ for 
$x \in \exp_G^{-1}(\1)$ sufficiently close to $0$. This implies in 
particular that each identity neighborhood contains non-trivial torsion subgroups.

It is a classical result that Banach--Lie groups do not contain 
small subgroups, i.e., there exists a $\1$-neighborhood $U$ 
for which $\{\1\}$ is the only subgroup  contained in $U$. 
This is no longer true for locally convex vector groups, such 
as $G = \R^\N$, with the product topology. Then each $0$-neighborhood 
contains non-zero vector subspaces, so that $G$ has small subgroups. 
However, $G$ is torsion free. 

For a locally convex space $E$, the non-existence of small subgroups 
is equivalent to the existence of a continuous norm on $E$. Every 
locally exponential Lie group 
$G$ for which $\L(G)$ has a continuous norm has no small subgroups 
(cf.\ Section IV). 
Since any real vector space is torsion free, this implies that 
no locally exponential Lie group contains small torsion subgroups. 
For strong ILB--Lie groups, it is also known that they do not contain 
small subgroups (cf.\ Theorem~III.2.3), and this implies in particular 
that for each compact manifold $M$ the group $\Diff(M)$ does not contain small subgroups. 
We also know that direct limits of finite-dimensional Lie groups do not contain 
small subgroups (Theorem~VII.1.3). 
\qed

\Problem II.6. (Initial Subgroup Problem) 
Give an example of a subgroup $H$ of some infinite-dimensional 
Lie group which does not possess any initial Lie subgroup structure. 

We think that such examples exist, but in view of Remark II.6.5(b), 
there is no such example in any finite-dimensional Lie group. Moreover, 
Theorem~IV.4.17 below implies that all closed subgroups of Banach--Lie groups 
carry initial Lie subgroup structures. Therefore the most natural candidates of groups
to consider are (non-closed) subgroups of Banach spaces which are connected 
by smooth arcs. For $E := C([0,1],\R)$, the subgroup $H \subeq E$ generated by the smooth 
curve $\gamma \: [0,1] \to E, \gamma(t)(x) := e^{tx}-1$ is a natural candidate. 
Since the values $\gamma(t)$ for $t > 0$ 
are linearly independent, $H = \sum_{t \in ]0,1]} \Z \gamma(t)$ 
is a free abelian group. 
\qed

\Problem II.7. (Canonical factorization for Lie groups) Let 
$\phi \: G \to H$ be a morphism of Lie groups. 
Does the quotient group $G/\ker \phi \cong \phi(G) \subeq H$ carry a natural Lie group 
structure for which the induced map $G/\ker(\phi) \to H$ is smooth and  
each other morphism $\psi \: G \to H'$ with $\ker \phi \subeq \ker \psi$ factors 
through $G/\ker(\phi)$? 
Does $\phi(G)$ carry the structure of an initial Lie subgroup of $H$? 
Maybe it helps to assume that $G$ is a regular Lie group (cf.\ Section III below). 
\qed

\Problem II.8. (Locally Compact Subgroup Problem; (FP9)) 
Show that any locally compact subgroup of a Lie group $G$ is a (finite-dimensional) 
Lie group. Since locally compact subgroups are Lie groups if and only if they 
have no small subgroups, this is closely related to Problem II.5. 
We shall see below that this problem has a positive solution 
for most classes of concrete groups 
(cf.\ Theorem IV.3.15 for locally exponential Lie groups; 
[MZ55, Th.~5.2.2, p.~208] for the Lie group $\Diff(M)$ of diffeomorphisms of a compact 
manifold,  
and Theorem VII.1.3 for direct limits of finite-dimensional Lie groups). 
\qed

\Problem II.9. (Completeness of Lie groups) Suppose that the Lie algebra 
$\L(G)$ of the Lie group $G$ is a complete locally convex space. Does this 
imply that the group $G$ is complete with respect to the left, resp., right uniform 
structure? 
\qed

\Problem II.10. (Large tori in Lie groups) Suppose that $G$ is a Lie group 
with a smooth exponential function and that 
$\a \subeq \L(G)$ is a closed abelian subalgebra for which the closed subgroup 
$\Gamma_\a := \exp_G^{-1}(\1) \cap \a$ spans a dense subspace of $\a$. 
Then the exponential 
function $\Exp_\a := \exp_G\res_\a \: \a \to G$ factors through a continuous map 
$\a/\Gamma_\a \to G$. Characterize the groups $A :=\a/\Gamma_\a$ for which this may happen. 

If $\a$ is finite-dimensional, then $A$ is a torus (Proposition II.5.11(5)), so that 
we may think of these groups $A$ as {\it generalized tori}. If 
$\Gamma_\a$ is discrete, then $A$ is a Lie group. If, in addition, 
$\a$ is separable, then $\Gamma_\a$ is a free group 
([Ne02a, Rem. 9.5(c)]). If Problem II.4 has a positive 
solution, then $\Gamma_\a$ is always discrete. 

An interesting example in this context is $E = \R^\N$ with the closed subgroup $\Gamma_E := \Z^\N$. 
In this case, the quotient $E/\Gamma_E \cong \T^\N$ is the compact torus which is 
not a Lie group because it is not locally contractible. Do pairs 
$(\a,\Gamma_\a) \cong (\R^\N, \Z^\N)$ occur? 
For the free vector space $E = \R^{(\N)}$ over $\N$, the subgroup 
$\Gamma_E := \Z^{(\N)}$ is discrete and $E/\Gamma_E$ is a Lie group, 
a direct limit of finite-dimensional tori (cf.\ Theorem VII.1.1). 
\qed

\Problem II.11. Does the adjoint group $\Ad(G) \subeq \Aut(\L(G))$ of a Lie group $G$ 
always carry a natural Lie group structure for which the adjoint representation 
$\Ad \: G \to \Ad(G)$ is a quotient morphism of topological groups? 
Since, in general, the group $\Aut(\L(G))$ is not a Lie group if $\L(G)$ is not Banach, this 
does not follow from a positive solution of Problem~II.7. 
Closely related is the question if $\exp_G x \in Z(G)$ for $x$ small enough implies  
$x \in \z(\g)$. 
\qed

\Problem II.12. Let $G$ be a connected Lie group with a smooth exponential function 
and $\a \subeq \L(G)$ a Mackey complete abelian subalgebra for which the group 
$\Gamma_\a := \exp_G^{-1}(\1) \cap \a$ is discrete. Then 
$\exp_G \res_\a$ factors through an injective smooth map 
$A := \a/\Gamma_\a \into G$ and $A$ carries a natural Lie group structure 
(Corollary~II.2.4). Is this Lie group always initial? According to 
Remark~II.6.5(c), this is the case for $\a = \z(\g)$. 
\qed

\sectionheadline{III. Regularity} 

\nin In this section, we discuss regularity of Lie groups in some more detail. 
In particular, we shall see how regularity of a Lie group can be used to 
obtain a Fundamental Theorem of Calculus for Lie group-valued smooth functions. 
This implies solutions to many integrability questions. 
For example, for each homomorphism $\psi \:  \L(G) \to \L(H)$ from the Lie algebra of 
a $1$-connected Lie group $G$ into the Lie algebra of a regular 
Lie group $H$, there exists a unique morphism of Lie groups $\phi$ with 
$\L(\phi) = \psi$. In Section III.2, we turn to the concepts of strong ILB--Lie groups 
and $\mu$-regularity and their relation to our context. In the remaining 
two subsections III.3 and III.4, we discuss some applications to 
groups of diffeomorphisms and groups of smooth maps, resp., gauge groups.  

\subheadline{III.1. The Fundamental Theorem for Lie group-valued functions} 

\Definition III.1.1. Let $G$ be a Lie group with Lie algebra $\g = \L(G)$. 
We call a $\g$-valued $1$-form 
$\alpha \in \Omega^1(M,\g)$ {\it integrable} if there exists a smooth 
function $f \: M \to G$ with $\delta(f) = \alpha$. 
The $1$-form $\alpha$ is said to be 
{\it locally integrable} if each point $m \in M$ has an open 
neighborhood $U$ such that $\alpha\res_U$ is integrable. 
\qed

We recall from Definition~I.4.1(b) the brackets 
$\Omega^p(M,\g) \times \Omega^q(M,\g) \to \Omega^{p+q}(M,\g)$. 
If $f$ is a solution of the equation $\delta(f) = f^*\kappa_G = \alpha \in \Omega^1(M,\g)$, 
then the fact that $\kappa_G$ satisfies the Maurer--Cartan equation 
$d\kappa_G + \frac{1}{2} [\kappa_G, \kappa_G] =0$
implies that so does $\alpha$: 
$$d\alpha + \frac{1}{2} [\alpha,\alpha] =0. \leqno(MC) $$

The following theorem is a version of the Fundamental Theorem of Calculus 
for functions with values in regular Lie groups ([GN06]). 

\Theorem III.1.2. {\rm(Fundamental Theorem for Lie group-valued functions)} 
Let $M$ be a smooth manifold, $G$ a Lie group and 
$\alpha \in \Omega^1(M,\L(G))$. Then the following assertions hold: 
\litem{(1)} If $G$ is regular and $\alpha$ satisfies the Maurer--Cartan equation, 
then $\alpha$ is locally integrable. 
\litem{(2)} If $M$ is $1$-connected and $\alpha$ is locally integrable, then it is 
integrable. 
\litem{(3)} If $M$ is connected, $m_0 \in M$, and $\alpha$ is locally integrable, 
then there exists a homomorphism 
$$ \per_\alpha \: \pi_1(M,m_0) \to G $$
that vanishes if and only if $\alpha$ is integrable. 
For a piecewise smooth representative $\sigma \: [0,1] \to M$ of a loop 
in $M$, the element $\per_\alpha([\sigma])$ is given by 
$\gamma(1)$ for $\gamma \: [0,1] \to G$ satisfying 
$\delta(\gamma) = \sigma^*\alpha$. 
\qed

\Remark III.1.3. If $M$ is one-dimensional, 
then each $\g$-valued $2$-form on $M$ vanishes, so that 
$[\alpha,\beta] = 0 = d\alpha$ 
for $\alpha, \beta \in \Omega^1(M,\g)$. Therefore 
all $1$-forms trivially satisfy the Maurer--Cartan equation. 
\qed 

This remark applies in particular to the manifold with boundary 
$M = I = [0,1]$. The requirement that for each smooth curve $\xi \in C^\infty(I,\g) 
\cong \Omega^1(I,\g)$,  
the IVP 
$$ \gamma(0) = \1, \quad \gamma'(t) = \gamma(t).\xi(t)\quad \hbox{ for } \quad t \in I, $$
has a solution depending smoothly on $\xi$ leads to the concept 
of a regular Lie group.

\Remark III.1.4. (a) If $M$ is a complex manifold, $G$ is a complex Lie group and 
$\alpha \in \Omega^1(M,\g)$ is a holomorphic $1$-form, then 
for any smooth function $f \: M \to G$ with $\delta(f) = \alpha$,  
the differential of $f$ is complex linear in each point, so that 
$f$ is holomorphic. Conversely, the logarithmic derivative of any holomorphic 
function $f$ is a holomorphic $1$-form. 

If, in addition, $M$ is a one-dimensional complex manifold, 
then for each holomorphic $1$-form 
$\alpha \in \Omega^1(M,\g)$ the $2$-forms 
$d\alpha$ and $[\alpha,\alpha]$ are holomorphic, which implies that they 
vanish. Therefore the Maurer--Cartan equation is automatically satisfied by all holomorphic 
$1$-forms. 
\qed

The following theorem is one of the main motivations for introducing the notion of 
regularity. It was  
proved in [OMYK82] under the stronger assumption of $\mu$-regularity (cf.\ Subsection III.2 below) and by {\smc Milnor} (who attributed it to {\smc Thurston}) in the following form ([Mil82/84]):  

\Theorem III.1.5.  If $H$ is a regular Lie group, $G$ is a $1$-connected Lie group,
and $\phi \: \L(G)\to \L(H)$ is a continuous homomorphism of Lie algebras,
then there exists a unique Lie group homomorphism $f \: G \to H$
with $\L(f) = \phi$. 

\Proof. This is Theorem 8.1 in [Mil84] (see also [KM97, Th.~40.3]). 
The uniqueness assertion follows from Proposition~II.4.1 
and does not require the regularity of $H$. 

On $G$, we consider the smooth $\L(H)$-valued $1$-form 
$\alpha := \phi \circ \kappa_G$ and it is easily verified that 
$\alpha$ satisfies the MC equation. 
Therefore the Fundamental Theorem implies the existence of a unique 
smooth function $f \: G \to H$ with $\delta(f) = \alpha$ and $f(\1_G) = \1_H$. 
In view of Proposition~II.4.1(3),  the function $f$ is a homomorphism of 
Lie groups with $\L(f) = \alpha_\1 = \phi$. 
\qed

\Corollary III.1.6. If $G_1$ and $G_2$ are regular $1$-connected 
Lie groups with isomorphic Lie algebras, then $G_1$ and $G_2$ are isomorphic. 
\qed

\Corollary III.1.7. Let $G$ be a connected Lie group with Lie algebra 
$\g$ and $\n \trile \g$ a closed ideal which is not $\Ad(G)$-invariant. 
Then the quotient Lie algebra $\g/\n$ is not integrable to a regular Lie group. 

\Proof. If $Q$ is a regular Lie group with Lie algebra $\q := \g/\n$, 
then the quotient map $q \: \g \to \n$ integrates to a morphism 
of Lie groups $\phi \: \tilde G \to Q$ with 
$\L(\phi) = q$ (Theorem~III.1.5), so that $\n = \ker(\L(\phi))$, contradicting its 
non-invariance under $\Ad(\tilde G) = \Ad(G)$. 
\qed

\Remark III.1.8. Let $G$ be a regular Lie group and 
$\h \subeq \L(G)$ a closed Lie subalgebra. Let  
$\iota \: H \to G$ be a regular connected initial Lie subgroup of $G$ with 
$\L^d(H) = \h$. 
Then for each smooth curve $\gamma \: I \to H$ the curve 
$\delta(\gamma)$ has values in $\h$, and, conversely, for any smooth curve 
$\xi \: I \to \h$, the regularity of $H$ and the Uniqueness Lemma imply that the 
corresponding curve $\gamma_\xi$ has values in $H$. Hence $H$ coincides 
with the set of endpoints of all curves $\gamma_\xi$, $\xi \in C^\infty(I,\h)$. 
In particular, $H$ is uniquely determined by the Lie algebra $\h$. 
\qed

For the groups of smooth maps on a compact manifold, it is quite easy to 
find charts of the corresponding mapping groups, such as $C^\infty(M,K)$, 
by composing with charts of $K$ (Theorem~II.2.8). This does no longer work for 
non-compact manifolds, as the discussion in Remark~II.2.9(a) shows. 
The Fundamental Theorem implies that for any 
regular Lie group $K$ with Lie algebra $\k$, any $1$-connected manifold $M$, $m_0 \in M$ and 
$$C^\infty_*(M,K) := \{ f \in C^\infty(M,K) \: f(m_0) = \1\},$$ 
the map 
$$ \delta \: C^\infty_*(M,K) \to \{ \alpha \in \Omega^1(M,\k) \: 
d\alpha + {\textstyle{1\over 2}}[\alpha,\alpha]= 0\} $$
is a bijection, which can be shown to be a homeomorphism. 
If the solution set of the MC equation carries a natural manifold 
structure, we thus obtain a manifold structure on the group $C^\infty_*(M,K)$ and 
hence on $C^\infty(M,K)$. This is the case if $K$ is abelian, $M$ is one-dimensional 
(all $2$-forms vanish), and for holomorphic $1$-forms on complex one-dimensional 
manifolds (cf.\ Remark~III.1.4). 
Following this strategy and using Gl\"ockner's Implicit Function Theorem 
to take care of the period conditions if $M$ is not simply  connected, 
we get the following result ([NeWa06b]). To formulate the real and complex case 
in one statement, let $\K \in \{\R,\C\}$, $K$ be a $\K$-Lie group, and  
$C^\infty_\K(M,K)$ be the group of $\K$-smooth $K$-valued maps. For $\K = \C$, these 
are the holomorphic maps, and in this case the smooth $C^\infty$-topology on 
$C^\infty_\C(M,K) = {\cal O}(M,K)$ coincides with the compact open topology. 

\Theorem III.1.9. Let $K$ be a regular $\K$-Lie group and $M$ 
a finite-dimensional connected $\sigma$-compact $\K$-manifold. 
We endow the group $C^\infty_\K(M,K)$ with the compact open $C^\infty$-topology, 
turning it into a topological group. This topology is compatible with a Lie group 
structure if 
\litem{(1)} $\dim_\K M = 1$, $\pi_1(M)$ is finitely generated and $K$ is a Banach--Lie group. 
\litem{(2)} $H^1_{\rm sing}(M,\Z)$ is finitely generated and $K$ is abelian. 
\litem{(3)} $H^1_{\rm sing}(M,\Z)$ is finitely generated and 
$K$ is finite-dimensional and solvable. 
\litem{(4)} $K$ is diffeomorphic to a locally convex space. 
\qed

\subheadline{III.2. Strong ILB--Lie groups and $\mu$-regularity} 

An important criterion for regularity of a Lie group rests on the concept of a 
(strong) ILB--Lie group, a concept developed by {\smc Omori} by abstracting a common feature 
from groups of smooth maps and diffeomorphism groups ([Omo74]). 
The bridge from ILB--Lie groups to Milnor's regularity concept was built in 
[OMYK82], where even a stronger regularity concept, called $\mu$-regularity below, is used.  
In this subsection, we explain some of the key results concerning $\mu$-regularity 
and how they apply to diffeomorphism groups. In his book [Omo97], 
{\smc Omori} works with a slight variant of the axiomatics of $\mu$-regularity, 
as defined below, but since it is quite close to the original concept, 
we shall not go into details on this point. 

\Definition III.2.1. (ILB--Lie groups; [Omo74, p.2])

(a) An {\it ILB chain} is a sequence $(E_n)_{n \geq d}$, $d \in \N$, of Banach spaces 
with continuous dense inclusions $\eta_n \: E_n \into E_{n+1}$. The projective 
limit $E := \prolim E_n$ of this system is a Fr\'echet space. 

Realizing $E$ as $\{ (x_n)_{n \geq d} \in \prod_{n \geq d} E_n \: (\forall n)\ 
\eta_n(x_n) = x_{n+1}\}$, we see that for each $k \geq d$ the projection map 
$$ q_k \: E \to E_k, \quad (x_n) \mapsto x_k $$
is injective. We may therefore think of $E$ and all spaces $E_n$ as subspaces 
of $E_d$, which leads to the identification of $E$ with the intersection $\bigcap_{n \geq d} E_n$.  

(b) A topological group $G$ is called an {\it ILB--Lie group  
modeled on the ILB chain $(E_n)_{n \geq d}$}  
if there exists a sequence of topological groups $G_n$, $n \geq d$, 
satisfying the conditions (G1)--(G7) below. If, in addition, (G8) holds, then 
$G$ is called a {\it strong ILB--Lie group}. 
\litemindent=0.8cm
\litem{(G1)} $G_n$ is a smooth Banach manifold modeled on $E_n$. 
\litem{(G2)} $G_{n+1}$ is a dense subgroup of $G_n$ and the inclusion map $G_{n+1} \into G_n$ 
is smooth. 
\litem{(G3)} $G = \prolim G_n$ as topological groups, so that we may identify 
$G$ with $\bigcap_{n \geq d} G_n \subeq G_d$. 
\litem{(G4)} The group multiplication of $G$ extends to a $C^{\ell}$-map 
$\mu_G^{n,\ell} \: G_{n+\ell} \times G_{n} \to G_n$. 
\litem{(G5)} The inversion map of $G$ extends to a $C^\ell$-map $G_{n+\ell} \to G_n$. 
\litem{(G6)} The right translations in the groups $G_n$ are smooth. 
\litem{(G7)} The tangent map $T(\mu_G^{n,\ell})$ induces a $C^\ell$-map 
$T_\1(G_{n+\ell}) \times G_n \to T(G_n)$. 
\litem{(G8)} There exists a chart $(\phi_d, U_d)$ of $G_d$ with $\1 \in U_d$ 
and $\phi_d(\1) = 0$ such that $U_n := U_d \cap G_n$ and 
$\phi_n := \phi_d \res_{U_n}$ define an $E_n$-chart $(\phi_n, U_n)$ of $G_n$. 
\litemindent=0.7cm

If all spaces $E_n$ are Hilbert, we call $(E_n)_{n \geq d}$ an {\it ILH chain} and 
$G$ an {\it ILH--Lie group}. 
\qed

\Remark III.2.2. (Omori) Every strong ILB--Lie group $G$ 
carries a natural Fr\'echet--Lie group structure with 
$\L(G) \cong \bigcap_{n\geq d} E_n = E$. A chart $(\phi, U)$ in the identity is obtained 
by $U := G \cap U_d$ and $\phi := \phi_d\res_U$ (notation as in (G8)). 
\qed

A complete solution to (FP8) for a large class of Lie groups is provided by: 
\Theorem III.2.3. {\rm([Omo74, Th.1.4.2])} Strong ILB-Lie groups have no small subgroups, 
i.e., there exists an identity neighborhood containing no non-trivial subgroups. 
\qed

We now give the slightly involved definition of the regularity concept 
introduced in \break [OMYK82/83a].

\Definition III.2.4. Let $G$ be a Lie group with Lie algebra $\g$ and 
$\Delta := \{t_0, \ldots, t_m\}$ a division of the real 
interval $J := [a,b]$ with 
$a = t_0$ and $b = t_m$. We write 
$$|\Delta| := \max\{t_{j+1} -t_j \; j = 0,\ldots, m-1\}.$$ 
For $|\Delta| \leq \eps$, a pair $(h, \Delta)$ is called a {\it step function} on $[0,\eps] \times J$ 
if $h \: [0,\eps] \times J \to G$ is a map satisfying 
\litem{(1)} $h(0,t) = \1$ for all $t \in J$ and all maps 
$h^t(s) := h(s,t)$ are $C^1$. 
\litem{(2)} $h(s,t) = h(s,t_j)$ for $t_j \leq t < t_{j+1}$. 

\nin For a step function $(h,\Delta)$, we define the {\it product integral} 
$\prod_a^t (h,\Delta) \in G$ by 
$$ \prod_a^t (h,\Delta) := h(t-t_k, t_k) h(t_k-t_{k-1}, t_{k-1})\cdots h(t_1 - t_0, t_0)
\quad \hbox{ for } \quad t_k \leq t < t_{k+1}. $$ 

Now let $(h_n,\Delta_n)$ be a sequence of step functions with $|\Delta_n| \to 0$ 
for which the sequence $(h_n, {\partial h_n \over \partial s})$ converges 
uniformly to a pair $(h, {\partial h \over \partial s})$ for a 
function $h \: [0,\eps] \times J \to G$. 
Then the limit function $h$ is a {\it $C^1$-hair in $\1$}, 
i.e., it is continuous, differentiable with respect to $s$, and 
${\partial h \over \partial s}$ is continuous on $[0,\eps] \times J$. 

The Lie group $G$ is called {\it $\mu$-regular}\footnote{$^1$}{\eightrm $\mu$ stands for 
``multiplicative''} (called ``regular'' in 
[OMYK82/83a]) if the product integrals $\prod_a^t (h_n,\Delta_n)$ converges  
uniformly on $J$ for each sequence $(h_n, \Delta_n)$ converging in the sense explained 
above to some $C^1$-hair in $\1$. 
Then the limit is denoted $\prod_a^t (h, d\tau)$ and called the {\it product integral of $h$}.
\qed

\Remark III.2.5. The First Fundamental Theorem in [OMYK82] asserts that the product integral 
$\prod_a^t (h, d\tau)$ is $C^1$ with respect to $t$ and satisfies 
$$ {d\over dt} \prod_a^t (h, d\tau) = u(t)\cdot \prod_a^t (h, d\tau)
\quad \hbox{ for } \quad 
u(t) = {\partial h \over \partial s}(0,t), $$
where $u \in C(J,\g)$ is a continuous curve. Hence the product integral is 
the unique $C^1$-curve $\gamma_u \: J \to G$ with $\gamma_u(a) = \1$ and 
$\delta^r(\gamma_u) = u$. The Second Fundamental Theorem in [OMYK82] is that the right 
logarithmic derivative 
$$ \delta^r \: C^1_*(J,G) \to C(J,\g) $$
is a $C^\infty$-diffeomorphism, where 
$C^1_*(J,G)$ is the group 
of $C^1$-paths $\gamma \: J \to G$ with $\gamma(a) = \1$, 
endowed with the compact open $C^1$-topology (cf.\ Theorem~II.2.8). 

Since the inclusion map 
$C^\infty([0,1],\g) \to C^0([0,1],\g)$
is continuous and the evaluation map $\ev_1 \: C^1_*([0,1],G) \to G, \gamma \mapsto \gamma(1)$ 
is smooth, it follows in particular that each $\mu$-regular Lie group is regular. 
\qed

\Theorem III.2.6. {\rm([OMYK82, Th.~6.9])} Strong ILB-Lie groups are $\mu$-regular, hence 
in particular regular. 
\qed

\Lemma III.2.7. {\rm([OMYK83a, Lemma 1.1])} In 
each $\mu$-regular Fr\'echet--Lie group $G$, we have for each $C^1$-curve 
$\gamma \: [0,1] \to G$ with $\gamma(0) = \1$ the relation 
$$ \lim_{n \to \infty} \gamma\big({\textstyle{t\over n}}\big)^n = \exp_G(t \gamma'(0)) \quad 
\hbox{ for } \quad 0 \leq t \leq 1. 
\qeddis 

\Theorem III.2.8. {\rm([OMYK83a, Th.~4.2])} Let $G$ be a $\mu$-regular Fr\'echet--Lie group. 
For each closed finite-codimensional subalgebra $\h \subeq \L(G)$,  
there exists a connected Lie group $H$ with $\L(H) = \h$ and an injective 
morphism of Lie groups $\eta_\h \: H \to G$ for which $\L(\eta_\h) \: \L(H) \to \L(G)$ 
is the inclusion of $\h$. 
\qed

\Theorem III.2.9. {\rm([OMYK83a, Prop.~6.6])} Let $M$ be a compact manifold, 
$G$ a locally exponential $\mu$-regular Fr\'echet--Lie group, $r \in \N_0 \cup \{\infty\}$, and 
$q \: \G \to M$ a smooth fiber bundle whose fibers are groups isomorphic to $G$, for  
which the transition functions are group automorphisms. Then the group 
$C^r(M,\G)$ of $C^r$-sections of this bundle is a group with respect to pointwise 
multiplication, and it carries a natural Lie group structure, turning it into a 
$\mu$-regular Fr\'echet--Lie group. 
\qed

A slightly weaker version of the preceding theorem can already be found in [Les68]. 
Note that it applies in particular to gauge groups of $G$-bundles over $M$. 
We have added the assumption that $G$ is locally exponential because 
this is needed for the standard constructions of charts of the group 
$C^r(M,\G)$ (cf.\ Theorem IV.1.12 below for gauge groups). 

\Theorem III.2.10. {\rm([OMYK83a, Prop.~2.4])} 
Let $G$ be a $\mu$-regular Fr\'echet--Lie group 
and $H \subeq G$ a subgroup for which there exists an identity neighborhood $U^H$ 
whose smooth arc-component of $\1$ is a submanifold of $G$. Then $H$ carries the structure 
of an initial Lie subgroup {\rm(Remark~II.6.6)} which is $\mu$-regular. 
\qed

\subheadline{III.3. Groups of diffeomorphisms} 

As an important consequence of Theorem III.2.6, several classes of 
groups of diffeomorphisms are regular: 

\Theorem III.3.1. Let $M$ be a compact smooth manifold. Then the following 
groups carry natural structures of strong ILH--Lie groups, and hence are $\mu$-regular: 
\litem{(1)} $\Diff(M)$. 
\litem{(2)} $\Diff(M,\omega) := \{ \phi \in \Diff(M)\: \phi^*\omega = \omega\}$, where 
$\omega$ is a symplectic $2$-form on $M$. 
\litem{(3)} $\Diff(M,\mu)$, where $\mu$ is a volume form on $M$.
\litem{(4)} $\Diff(M,\alpha)$, where $\alpha$ is a contact form on $M$. 

\nin The corresponding Lie algebras are ${\cal V}(M)$, 
${\cal V}(M,\omega) := \{ X \in {\cal V}(M) \: {\cal L}_X\omega = 0\}$, 
${\cal V}(M,\mu)$, resp., ${\cal V}(M,\alpha)$. 
\qed

The preceding results on the Lie group structure of groups of diffeomorphisms 
have a long history. 
The Lie group structure on $\Diff(M)$ for a compact manifold $M$ has first 
been constructed by {\smc J.~Leslie} in [Les67], and {\smc Omori} proved in 
[Omo70] that $\Diff(M)$ can be given the structure of a strong 
ILH--Lie group (cf.\ also [Eb68] for the ILH structure). {\smc Ebin} and {\smc Marsden} 
extended the ILH results to compact manifolds with boundary ([EM69/70]). 
Later expositions of this result can be found in [Gu77], [Mi80] and [Ham82]. 
The regularity of $\Diff(M)$ is proved in [Mil84/82] with direct arguments, 
not using ILB techniques.

In [Arn66], {\smc Arnold} studies the group $\Diff(M,\mu)$, 
where $\mu$ is a volume form on the compact manifold $M$,  
as the configuration space of a perfect fluid. Arguing by analogy with 
finite-dimensional groups, he showed that, for a 
suitable right invariant Riemannian metric on this group, the Euler 
equation of a perfect fluid 
corresponds to the geodesic equation for a left invariant 
Riemannian metric on $\Diff(M,\mu)$. This was made rigorous 
by Marden and Abraham in [MA70]. 
For a more recent survey on this circle of ideas, we refer to [EMi99]. 

Using Hodge theory, {\smc Ebin} and {\smc Marsden} show in [EM70] that 
if $\omega$ either is a volume form or a symplectic form on a compact 
manifold $M$, then 
$\Diff(M,\omega)$ carries the structure of an ILH--Lie group (see also [Wei69]). 
They further show that 
the group $\Diff_+(M)$ of orientation preserving diffeomorphisms of $M$ is 
diffeomorphic to the direct product 
$\Diff(M,\mu) \times \Vol_1(M)$, where $\Vol_1(M)$ denotes the convex 
set of volume forms of total mass $1$ on 
$M$. This refines a result of {\smc Omori} on the topological 
level  (cf.\ [KM97, Th.~43.7]). 
For the symplectic case, a more direct proof of the regularity 
assertion can be found in [KM97, Th.~43.12], where it is also shown that 
$\Diff(M,\omega)$ is a submanifold of $\Diff(M)$. 

In [EM70], one also finds that the following groups are ILH--Lie groups: 

\litem{(1)} $\Diff(M,N) := \{ \phi \in \Diff(M) \: \phi(N) = N\}$ 
and 
$\Diff_N(M) := \{ \phi \in \Diff(M) \: \phi\res_N = \id_N\},$  
where $N \subeq M$ is a closed submanifold and $M$ compact without boundary. 
\litem{(2)} 
$\Diff_{\partial M}(M) := \{ \phi \in \Diff(M) \: (\forall x \in \partial M)
\ \phi(x) = x\}$, if $M$ has a boundary. 
\litem{(3)} $\Diff(M,\mu)$ and $\Diff_{\partial M}(M,\mu) := 
\Diff(M,\mu) \cap \Diff_{\partial M}(M)$
for any volume form $\mu$ on $M$. 
\litem{(4)} If, in addition, $\omega = d\theta$ 
is an exact symplectic $2$-form on $M$, then 
$\Diff_{\partial M}(M,\omega)$ and 
$$ \Ham(M,\omega) := \{ \phi \in  \Diff_{\partial M}(M) \: \phi^*\theta -\theta 
\in B^1_{\rm dR}(M,\R)\} $$
are ILH--Lie groups (see also Remark~V.2.14(c)). 
\msk 

The following result of {\smc Michor} ([Mi91]) 
concerns the Lie group structure of a gauge group in a setting where the gauge group 
is a Lie subgroup of a diffeomorphism group of a compact manifold. 

\Theorem III.3.2. If $q \: B \to M$ is a locally trivial fiber bundle over the 
compact manifold $M$ with compact fiber $F$, then the gauge group 
$\Gau(B)$ is a split submanifold of the regular Fr\'echet--Lie group $\Diff(B)$. 
\qed

\subheadline{III.4. Groups of compactly supported smooth maps and diffeomorphisms} 

In the preceding subsection, we discussed diffeomorphisms of compact manifolds. 
We now briefly take a look at the corresponding 
picture for compactly supported maps on $\sigma$-compact manifolds.  

It is interesting that if $M$ is a $\sigma$-compact finite-dimensional manifold, 
then for each locally convex space $E$, the space $C^\infty_c(M,E)$ has two natural topologies. 
The first one is the locally convex direct limit structure 
$$ C^\infty_c(M,E) = \indlim C^\infty_{M_n}(M,E), $$
where $(M_n)_{n \in \N}$ is an exhaustion of $M$, which for 
the case that $E$ is Fr\'echet, defines an LF space structure on $C^\infty_c(M,E)$ 
(cf.\ Examples I.1.3 and Theorem II.2.8). 
The other locally convex topology is obtained by endowing for each $r \in \N$ 
the space $C^r_c(M,E)$ 
with the direct limit structure $\indlim C^r_{M_n}(M,E)$ and 
then topologize $C^\infty_c(M,E)$ as the projective limit 
$\prolim C^r_c(M,E)$ of these spaces. 
These two topologies do not coincide (cf.\ [Gl02b], see also [Gl06a]). 

Similar phenomena occur for the space $C^\infty_c(M,\E)$ of smooth compactly supported 
sections of a vector bundle $\E \to M$ whose fibers are locally convex spaces. 
In the context of Lie algebras, this problem affects the 
model spaces $C^\infty_c(M,\k)$ of the Lie groups 
$C^\infty_c(M,K)$ and the space ${\cal V}_c(M)$ of compactly supported vector fields 
on $M$. For the natural LF space structure on ${\cal V}_c(M)$, 
the corresponding Lie group structure on $\Diff_c(M)$ has been constructed 
by {\smc Michor} in [Mi80,pp.~39, 197], 
where he even endows $\Diff(M)$ with the Lie group structure 
for which $\Diff_c(M)$ is an open subgroup (Corollary~II.2.3). 

The following theorem complements Theorem~II.2.8 in a natural way (cf.~[GN06], based on 
[Gl02d]). 

\Theorem III.4.1. Let $M$ be a $\sigma$-compact finite-dimensional smooth manifold 
and $K$ a regular Lie group. 
Both natural topologies turn $C^\infty_c(M,\L(K))$ into a topological Lie algebra. 
Accordingly, the group $C^\infty_c(M,K)$ 
carries two regular Lie group structures for which the Lie algebra is 
$C^\infty_c(M,\L(K))$, endowed with these two topologies. 
If $M$ is non-compact, these two regular Lie groups are {\sl not} isomorphic. 
\qed

The corresponding result for diffeomorphism groups is proved by {\smc Gl\"ockner} in [Gl02b] 
(for corresponding statements without proof see also [Mil82]). 

\Theorem III.4.2. Let $M$ be a $\sigma$-compact finite-dimensional manifold. 
Both natural topologies turn ${\cal V}_c(M)$ into a topological Lie algebra, 
and the group $\Diff_c(M)^{\rm op}$ 
carries two corresponding regular Lie group structure turning ${\cal V}_c(M)$ into 
its Lie algebra. For $M$ non-compact, these two regular Lie groups are {\sl not} isomorphic. 
\qed

\subheadline{Open Problems for Section III} 

\Problem III.1. Show that every abelian Lie group $G$ modeled on a Mackey complete 
locally convex space $\g$ is regular. 

We may w.l.o.g.\ assume that $G$ is $1$-connected (cf.\ Theorem~V.1.8 below). 
Then the regularity of the additive group of 
$\g = \L(G)$ (Proposition~II.5.6) implies that $\id_{\g}$ integrates to a smooth 
homomorphism 
$\Log_G \: G \to \g$ (Theorem~III.1.5), 
so that the assumption implies the existence of a logarithm function, but 
it is not clear how to get an exponential function (cf.\ [Mil82, p.36]). 
One would like to show that $\Log_G$ is an isomorphism 
of Lie groups, but also weaker information would be of interest: 
Is $\Log_G$ surjective or injective? 

If $H := \im(\Log_G)$ were a proper subgroup of $\g$, it would 
be a strange object:  
Since $\L(\Log_G) = \id_{\g}$, we have $\L^d(H) = \g$ (Remark~II.6.4). 
For any $\alpha \in C^1_*([0,1], H)$, the relation 
$\alpha'(0) = \lim_{n \to \infty} n \alpha({\textstyle{1\over n}})$
implies that $H$ is dense in $\g$. Is $H$ a vector space? 
Let 
$$P := \{ \xi \in C^\infty([0,1],\g) \: (\exists \gamma \in C^\infty([0,1],G)\, \delta(\gamma) 
= \xi\}.$$ 
For $\gamma(0) = \1$ and $\delta(\gamma) = \xi$, we then have 
$\Log_G(\gamma(1)) = \int_0^1 \xi(t)\, dt.$
Therefore $H$ is the image of the additive group $P$ under the integration map. 
Is $P$ a vector subspace of $C^\infty([0,1],\g)$? 
\qed

\Problem III.2. Let $G$ be a regular Lie group (not necessarily Fr\'echet or $\mu$-regular). 
Show that for each closed finite-codimensional subalgebra $\h \subeq \L(G)$ 
there exists a connected Lie group $H$ with $\L(H) = \h$ and an injective 
morphism of Lie groups $\eta_\h \: H \to G$ for which $\L(\eta_\h) \: \L(H) \to \L(G)$ 
is the inclusion of $\h$ (cf.\ Theorem III.2.8). 
\qed

\Problem III.3. Let $\K \in \{\R,\C\}$, $M$ be a $\sigma$-compact finite-dimensiomal 
$\K$-manifold, and $K$ a (finite-dimensional) $\K$-Lie group. 
We endow the group $C^\infty_\K(M,K)$ of $K$-valued $\K$-smooth 
functions $M \to K$ with the compact open $C^\infty$-topology, turning it into a topological 
group. For $\K = \C$, this is the group of holomorphic functions  
and the compact open $C^\infty$-topology coincides with 
the compact open topology (cf.\ [NeWa06b]). When is this topology on the group 
$C^\infty_\K(M,K)$ compatible with a Lie group structure? See Theorem III.1.9, for partial 
results in this direction.  
\qed

\Problem III.4. Consider a topological group 
$G = \prolim G_j$ which is a projective limit of the Banach--Lie groups $G_j$ 
(or more general locally exponential groups). 
\litem{(1)} Characterize the situations where $G$ is locally exponential 
in the sense that it carries a 
(locally exponential) Lie group structure (cf.\ Remark~IV.1.22). 
For the case where all $G_j$ are finite-dimensional, this is done in 
[HoNe06] (cf.\ Theorem~X.1.9). 
\litem{(2)} Can we say more in the special case $G = \prod_{j \in J} G_j$? 
\litem{(3)} Suppose that $G$ carries a compatible Fr\'echet--Lie group 
structure. Does this imply that $G$ is regular? 
(cf.\ [Ga97] for some results in this direction). 
\qed

\Problem III.5. Let $\iota_j \: H_j \to G$, $j = 1,2$, 
be two initial Lie subgroups of the Lie group $G$ with 
$\L^d(H_1) = \L^d(H_2)$. 
Which additional assumptions are necessary to conclude that $H_1 = H_2$ as subgroups of 
$G$, hence that $H_1$ and $H_2$ are isomorphic as Lie groups? 
Note that Remark III.1.8 implies that this is the case if $H_1$ and $H_2$ are $\mu$-regular 
or at least if the maps $\delta \: C^1_*([0,1],H_j) \to C^0([0,1],\L^d(H_j))$ are surjective. 
\qed

\sectionheadline{IV. Locally exponential Lie groups} 

\nin In this section, we turn to Lie groups with an exponential function 
$\exp_G \: \L(G) \to G$ which is 
well-behaved in the sense that it maps a $0$-neighborhood in $\L(G)$ 
diffeomorphically onto a $\1$-neighborhood in~$G$. We call 
such Lie groups {\it locally exponential}. 

This class of Lie groups has been introduced by {\smc Milnor} in 
[Mil84]\footnote{$^1$}{\eightrm In [Rob96/97],  
these groups are called ``of the first kind.''}, where 
one finds some of the basic results explained below. In [GN06], we devote a long chapter 
to this important class of infinite-dimensional Lie groups, which properly contains 
the class of BCH--Lie groups as those for which the BCH--series defines an analytic 
local multiplication on a $0$-neighborhood in $\L(G)$ (cf.\ [Gl02c] for basic 
results in the BCH context). In particular, it contains all Banach--Lie groups, 
but also many other interesting 
types of groups such as unit groups of Mackey complete CIAs, 
groups of the form $C^\infty_c(M,K)$, where $M$ is $\sigma$-compact 
and $K$ is locally exponential, 
and moreover, all projective limits of nilpotent Lie groups. It therefore includes 
many classes of ``formal'' Lie groups. The appeal of this 
class is due to its large scope and the strength of the general Lie theoretic 
results that can be obtained for these groups. Up to certain refinements 
of assumptions, a substantial part of the 
theory of Banach--Lie groups carries over to locally exponential groups. 

One of the most important structural consequences of local exponentiality 
is that it provides canonical local coordinates given by the exponential function. 
This in turn permits us to develop a good theory of subgroups and there even is a characterization 
of those subgroups for which we may form Lie group quotients. Moreover, we shall 
see in Section VI below that integrability of a locally exponential 
Lie algebra (to be defined below) can be characterized similarly as for Banach algebras. 

Not all regular Lie groups are locally exponential. 
The simplest examples can be found among groups of the form 
$G = E \rtimes_\alpha \R$ for a smooth $\R$-action on $E$ (Example~II.5.9). 
Another prominent example of a regular Lie group which is not locally exponential 
is the group $\Diff(\SS^1)$ of diffeomorphisms of the circle (Example~II.5.13). 

\subheadline{IV.1. Locally exponential Lie groups and BCH--Lie groups} 

\Definition IV.1.1. We call a Lie group $G$ {\it locally exponential} if 
it has a smooth exponential function $\exp_G \: \L(G) \to G$ 
which is a local diffeomorphism in $0$, i.e., 
there exists an open $0$-neighborhood $U \subeq \L(G)$ 
mapped diffeomorphically onto an open $\1$-neighborhood of $G$. 

A Lie group is called {\it exponential} if, in addition, $\exp_G$ is a global 
diffeomorphism. 
\qed

If $\exp_G \: \L(G) \to G$ is an exponential function, then 
$ T_0(\exp_G) = \id_{\L(G)}$
by definition. This observation is particularly useful 
in the finite-dimensional or Banach context, where it 
follows from the Inverse Function Theorem that 
$\exp_G$ is a local diffeomorphism in $0$, 
so that we can use the exponential function to obtain charts around $\1$:  

\Proposition IV.1.2. Banach--Lie groups are locally exponential. 
\qed

We shall see below that a similar conclusion does not work for general 
Fr\'echet--Lie groups, because in this context there is no general Inverse Function 
Theorem. From that it follows that to integrate a Lie algebra 
homomorphism $\phi \: \L(G) \to \L(H)$ 
to a group homomorphisms, it is in general not enough 
to start with the prescription $\exp_G x \mapsto \exp_H \phi(x)$ to obtain a 
local homomorphism, because $\exp_G(\L(G))$ need not be a $\1$-neighborhood in $G$ 
(cf.\ Example~II.5.9). 

For Banach--Lie groups, the existence of ``canonical'' coordinates provided by the 
exponential map leads to a description of the local multiplication in a canonical 
form, given by the BCH series:

\Definition IV.1.3. For two elements $x,y$ in a Lie algebra $\g$, we define 
$$ H_1(x,y) := x + y, \quad H_2(x,y) := {1\over 2}[x,y], $$
and for $n \geq 3$: 
$$ H_n(x,y) := \sum_{k,m\ge 0\atop p_i+q_i>0}{(-1)^k\over
(k+1)(q_1+\ldots+q_k+1)}{(\ad x)^{p_1}(\ad y)^{q_1}\ldots
(\ad x)^{p_k}(\ad y)^{q_k}(\ad x)^m\over p_1!q_1!\ldots p_k!q_k!m!}y, $$
where the sum is extended over all summands with 
$p_1 + q_1 + \ldots + p_k + q_k + m + 1 = n.$
The formal series 
$\sum_{n = 1}^\infty H_n(x,y)$
is called the {\it Baker--Campbell--Hausdorff series}. 
\qed

There are many different looking ways to write the polynomials $H_n(x,y)$. We have chosen the one 
obtained from the integral formula 
$$ x * y = x + \int_0^1 \psi(e^{\ad x}e^{t\ad y})y\, dt, \leqno(4.1.1) $$
where $\psi$ denotes the analytic function $\psi(z) := {z \over z-1} \log z$, defined 
in a neighborhood of $1$. Formula (4.1.1) is valid for sufficiently 
small elements $x$ and $y$ in a Banach--Lie algebra, because 
we may use functional calculus in Banach algebras to make sense 
of $\psi(e^{\ad x}e^{\ad y})$ for $x,y$ close to $0$. Then the explicit expansion 
of the BCH series is obtained from the series expansion of $\psi$ and the exponential series 
of $e^{\ad x}$ and $e^{\ad ty}$. 

\Remark IV.1.4. (History of the BCH series) 
In [SchF90a], {\smc F.~Schur} derived recursion formulas for the summands of 
the series describing the multiplication of a Lie group in canonical 
coordinates (i.e., in an exponential chart). He also proved the local convergence 
of the series given by this recursion relations, which can in turn be used to obtain 
the associativity of the BCH multiplication (cf.\ [BCR81, p.\ 93], [Va84, Sect.~2.15]). 
His approach is quite 
close to our treatment of locally exponential Lie algebras in the sense that 
he derived the series from the Maurer--Cartan form by integration of 
a partial differential equation of the form 
$f^*\kappa_\g = \kappa_\g$ with $f(0) = x$, whose unique solution 
is the left multiplication $f = \lambda_x$ in the local group. 

The BCH series was made more explicit by {\smc Campbell} in [Cam97/98], 
and in [Hau06] {\smc Hausdorff} approached the BCH series on a formal level, 
showing that the formal expansion of $\log(e^x e^y)$ 
can be expressed in terms of Lie polynomials. Part of his results 
had been obtained earlier by {\smc Baker} ([Bak01/05]). See [Ei68] for a more 
recent short argument that all terms in the BCH series are Lie brackets. 
\qed

\Definition IV.1.5. A topological Lie algebra $\g$ is called 
{\it BCH--Lie algebra} if there exists an open $0$-neighborhood $U \subeq \g$ 
such that for $x,y \in U$ the BCH series 
$$\sum_{n = 1}^\infty H_n(x,y) $$
converges and defines an analytic function $U \times U \to \g, (x,y) \mapsto x * y$ 
(cf.\ Definition~I.2.1). In view of [Gl02a, 2.9], the analyticity 
of the product $x*y$ is automatic if $\g$ is a Fr\'echet space. 
\qed

\Example IV.1.6. (a) If $\g$ is a nilpotent locally convex Lie algebra of nilpotency class 
$m$, then the BCH series defines a polynomial multiplication 
$$ x * y = x + y +  {1\over 2}[x,y] + \sum_{n \leq m} H_n(x,y) $$
on $\g$. From the structure of the series it follows immediately that 
for $t,s \in \R$ and $x \in \g$ we have 
$$ tx * sx = (t+s)x, $$
so that $(\g,*)$ is an exponential nilpotent Lie group. 

(b) If $\g$ is a Banach--Lie algebra whose norm is submultiplicative in the sense that 
$\|[x,y]\| \leq \|x\|\cdot\|y\|$ for $x, y \in \g$, then the BCH series 
$x*y = \sum_{n = 1}^\infty H_n(x,y)$ 
converges for  $\|x\|,\|y\| < {1\over 3} \log({3\over 2})$ ([Bir38]). 
\qed

The following result is quite useful to show that certain Lie algebras are not BCH: 
\Theorem IV.1.7. {\rm(Robart's Criterion; [Rob04])} If $\g$ is a 
sequentially complete BCH--Lie algebra, then there exists a 
$0$-neighborhood $U \subeq \g$ such that 
$f(x,y) := \sum_{n = 0}^\infty (\ad x)^n y$
converges and defines an analytic function on $U \times \g$. 
\qed

On the global level we have the following result whose proof requires 
the uniqueness assertion from Theorem~IV.2.8 below:  

\Theorem IV.1.8. For a Lie group $G$ the following are equivalent: 
\litem{(1)} $G$ is analytic with an analytic 
exponential function which is a local analytic 
diffeomorphism in $0$. 
\litem{(2)} $G$ is locally exponential and $\L(G)$ is BCH. 
\qed

In Examples IV.1.14(b) and IV.1.16 
below, we describe an analytic Lie group with an analytic exponential 
function which is a smooth diffeomorphism, but such that $\L(G)$ is not BCH. This is a 
negative answer to a question raised in [Mil84, p.31].

\Definition IV.1.9. A 
group satisfying the equivalent conditions of the preceding theorem is 
called a {\it BCH--Lie group}. 
\qed

Our introductory discussion now can be stated as: 
\Corollary IV.1.10. Each Banach--Lie group is BCH. 
\qed

The Lie group concept used in [BCR81] is stronger than our concept of a BCH--Lie group 
because additional properties of the Lie algebra 
are required, namely that it is a so-called AE--Lie algebra, a property 
which encodes the existence of certain seminorms, compatible with the Lie bracket. 

The following two theorems show that many interesting classes of Lie groups 
are in fact BCH. 

\Theorem IV.1.11. If $A$ is a Mackey complete CIA, then its unit group $A^\times$ is 
BCH. If, in addition, $A$ is sequentially complete, then $A^\times$ is regular. 

\Proof. (Sketch) If $A$ is a Mackey 
complete complex CIA, then the fact that $A^\times$ is open implies that for each 
$a \in A$ the spectrum $\Spec(a)$ is a compact subset of $\C$, 
and the holomorphic functional calculus works as for Banach algebras 
(cf.\ [Wae54a/b]\footnote{$^1$}{\eightrm Waelbroeck even introduces a functional 
calculus in several variables for tuples in complete locally convex 
algebras which 
are not necessarily CIAs, but where spectra and resolvents satisfy certain 
regularity conditions.}, [Al65], [Gl02b]). This provides an analytic exponential function, and 
on the open star-like subset 
$$ U := \{ a \in A \: \Spec(a) \cap ]{-\infty},0] = \eset\} \subeq A^\times $$
we likewise obtain an analytic logarithm function $\log\: U \to A$. 
From that, local exponentiality for complex CIAs follows easily. 

That the multiplication 
on  $A^\times$ is analytic follows from its bilinearity on $A$, 
and the analyticity of the inversion is obtained from functional calculus, 
which in turn leads to the expansion by the Neumann series 
$(\1 - x)^{-1} = \sum_{n = 0}^\infty x^n.$
We conclude that $A^\times$ is a BCH--Lie group. 

The real case can be reduced to the complex case, because for each real 
CIA $A$ its complexification $A_\C$ is a complex CIA ([Gl02b]). 

Now assume that $A$ is sequentially complete. For $u \in C([0,1],A)$ we want to solve the 
linear initial value problem 
$$ \gamma(0) = \1,\quad \gamma'(t) = \gamma(t)u(t). \leqno(4.1.2) $$
According to an idea of {\smc T. Robart} ([Rob04]), 
the BCH property of $A^\times$ implies that 
this can be done by Picard iteration: 
$$ \gamma_0(t) := \1, \quad 
\gamma_{n+1}(t) := \1 + \int_0^t \gamma_n(\tau)u(\tau)\, d\tau, $$
which leads to 
$$ \gamma_n(t) =  \1 + \sum_{k = 1}^n \int_0^t \int_0^{\tau_n} \cdots \int_0^{\tau_2} 
u(\tau_1)u(\tau_2)\cdots u(\tau_n)\ d\tau_1\, d\tau_2 \cdots d\tau_n. $$
Now one argues that the analyticity of the function 
$(\1 - x)^{-1} = \sum_{n = 0}^\infty x^n$
implies that all sums of the form 
$\sum_{n = 0}^\infty x_{n1}\cdots x_{nn}$ converge for $x_{ij}$ in some sufficiently small 
$0$-neighborhood. A closer inspection of the limiting process implies 
that the limit curve $\gamma := \lim_{n \to \infty} \gamma_n$ is $C^1$, 
solves the initial value problem (4.1.2), and depends analytically on $u$. This implies 
the regularity of $A^\times$. 
\qed

\Theorem IV.1.12. If $K$ is a locally exponential Lie group and 
$q \: P \to M$ a smooth $K$-principal bundle over the $\sigma$-compact finite-dimensional 
manifold $M$, then the group $\Gau_c(P)$ of compactly supported gauge transformations 
is a locally exponential Lie group. 
In particular, the Lie group $C^\infty_c(M,K)$ is locally exponential. 

If, in addition, $K$ is regular, then $\Gau_c(P)$ is regular and 
if $K$ is BCH, then so is $\Gau_c(P)$. 

\Proof. (Sketch; cf. [GN06] and Theorem II.2.8) 
Let $\exp_K \: \L(K) \to K$ be the exponential function of $K$ 
and realize $\Gau(P)$ as the subgroup 
$C^\infty(P,K)^K$ of $K$-fixed points in $C^\infty(P,K)$ with respect to the 
$K$-action given by $(k.f)(p) := k f(p.k)k^{-1}$. 
Then we put 
$$ \gau(P) := C^\infty(P,\L(K))^K = \{ \xi \in C^\infty(M,\L(K)) \: 
(\forall p \in P)(\forall k \in K)\, \Ad(k).\xi(p.k) = \xi(p) \}, $$
and observe that for the group $G := \Gau_c(P)$ the map 
$$ \exp_G \:\g := \gau_c(P) \to G, \quad \xi \mapsto \exp_K \circ \xi $$
is a local homeomorphism in $0$. Using Theorem~II.2.1, 
this can be used to define a Lie group structure on $G$. 
Then $\exp_G$ is an exponential function of $G$, and, by construction, 
it is a local diffeomorphism in $0$. 
\qed

Various special cases of the preceding theorem can be found in the literature: 
[OMYK82], [Sch04] (for $M$ compact, $K$ finite-dimensional), [Mil84] (without proofs), 
[KM97, 42.21] (in the convenient setting), and 
[Gl02c], [Wo05a]). That $\Gau(P)$ is $\mu$-regular if $K$ is $\mu$-regular follows from 
Theorem~III.2.9. 

\Example IV.1.13. (Pro-nilpotent Lie groups) If $\g = \prolim \g_j$ is a projective limit of a family of 
nilpotent Lie algebras $(\g_j)_{j \in J}$ (a so-called {\it pro-nilpotent Lie algebra}), 
then the corresponding connecting homomorphisms 
of Lie algebras are also morphisms for the corresponding group structures (Example IV.1.6(a)), 
so that 
$(\g,*) := \prolim (\g_j,*)$ defines on $\g$ a Lie group structure 
with $\L(\g,*) = \g$. We thus obtain an exponential Lie group 
$G := (\g,*)$ with $\exp_G = \id_\g$. This group is {\it pro-nilpotent} in the sense 
that it is a projective limit of nilpotent Lie groups. 

\Example IV.1.14. (Formal diffeomorphisms) (a) Important examples of pro-nilpotent Lie groups arise 
as certain groups of formal diffeomorphisms. 
We write $\Gf_n(\K)$ for the group of formal diffeomorphisms of 
$\K^n$ fixing $0$, where $\K \in \{\R,\C\}$. 
The elements of this group are represented by formal power series of the form 
$$ \phi(x) = gx + \sum_{|{\bf m}| > 1} c_{\bf m} x^{\bf m}, $$
where $g \in \GL_n(\K)$, 
$$ {\bf m} = (m_1,\ldots, m_n) \in \N_0^n, \quad 
|{\bf m}| := m_1 + \ldots + m_n, \quad 
x^{\bf m} := x_1^{m_1} \cdots x_n^{m_n}, \quad c_{{\bf m}} \in \K^n, $$
and the group operation is given by composition of power series. 
We call $\phi$ {\it pro-unipotent} if $g = \1$. 
It is easy to see that the pro-unipotent formal diffeomorphisms form a pro-nilpotent 
Lie group $\Gf_n(\K)_1 = \prolim G_k$, where 
$G_k$ is the finite-dimensional nilpotent group obtained by 
composing polynomials of the form 
$$ \phi(x) = x + \sum_{1 < |{\bf m}| \leq k} c_{\bf m} x^{\bf m} $$
modulo terms of order $> k$. The group $\Gf_n(\K)$ of all formal diffeomorphisms 
of $\K^n$ fixing $0$ is a semidirect product 
$$ \Gf_n(\K) \cong \Gf_n(\K)_1 \rtimes \GL_n(\K), \leqno(4.1.3) $$
where the group $\GL_n(\K)$ of linear automorphisms acts by conjugation. 
As this action is smooth, $\Gf_n(\K)$ is a Fr\'echet--Lie group. 

These groups are $\mu$-regular Lie groups (cf.\ [Omo80]): 
In view of the semidirect decomposition and the fact that 
$\mu$-regularity is an extension property (Theorem V.1.8), 
it suffices to observe that pro-nilpotent Lie groups are $\mu$-regular, which follows 
by an easy projective limit argument. 

The  group $\Gf_n(\K)$ has been studied by {\smc Sternberg} in [St61], where he shows in particular 
that for $\K = \C$ and $n = 1$ the elements 
$$ \phi_m(x) = e^{2 \pi i \over m} x + p x^{m+1}, \quad m \in \N\setminus \{1\}, p \in \C^\times, $$
are not contained in the image of the exponential function. This is of particular 
interest because $\phi_m \to \1$ in the Lie group $\Gf_n(\C)$, 
so that the image of the exponential function in this group 
is not an identity neighborhood. A detailed analysis of the exponential function 
of this group can also be found in {\smc Lewis}'s paper [Lew39].  

To see that $\phi_m$ is not in the image of the exponential function 
of $\Gf_n(\C)$, 
it suffices to verify this in the finite-dimensional solvable quotient group 
$G_{m+1} \rtimes \C^\times$, i.e., modulo terms of order $m+2$. 
The subgroup $\C x^{m+1} \rtimes \C^\times$ is isomorphic to 
$\C \rtimes \C^\times$ with the multiplication 
$$ (z,w)(z',w') = (z + w^mz', ww') $$
and the exponential function 
$$ \exp(z,w) = \Big( {e^{wm} - 1\over wm}z, e^w\Big)= \Big( {(e^w)^m - 1\over wm}z, 
e^w\Big), $$
showing that $\phi_m$ is not contained in the exponential image 
of this subgroup. However, one can use Proposition~II.5.11(3) to 
see that any element $\xi$ with $\exp\xi = \phi_m$ must be contained in 
the plane $\span \{x, x^{m+1}\}$. This completes the proof. 

(b) For $\K = \R$ the identity component $\Gf_1(\R)_0$ is exponential and analytic, but not BCH. 
For $n \geq 2$ the group $\Gf_n(\R)$ is analytic, but not locally exponential. 
If a subgroup $H \subeq \GL_n(\R)$ consists of matrices with real eigenvalues, 
then the subgroup $\Gf_n(\R)_1 \rtimes H \subeq \Gf_n(\R)$ 
is locally exponential ([Rob02, Ths.~6/7]).
\qed

\Example IV.1.15. Let $F= \K[x_1,\ldots, x_n]$ be the free associative algebra in $n$ generators 
$S := \{x_1,\ldots, x_n\}$. 
Then $F$ has a natural filtration 
$$F_k := \span\{ s_1\cdots s_m \: s_i \in S, m \geq k \}. $$
Each quotient $F/F_m$ is a finite-dimensional unital algebra, hence a CIA. 
Therefore the algebra 
$\hat F := \prolim F/F_n$, which can be identified with the algebra 
of non-commutative formal power series in 
the generators $x_1, \ldots, x_n$, is a complete CIA ([GN06]). 

We conclude that the unit group $\hat F^\times$ is a BCH--Lie group (Theorem IV.1.11).  
Let $\eps \: \hat F \to \K$ denote the homomorphism sending each $x_i$ to $0$. 
Then the normal subgroup $U := \1 + \ker \eps$ is a pro-nilpotent Lie group 
and $\hat F^\times \cong U \rtimes \K^\times$. In particular, the exponential function 
$$ \Exp \: \ker \eps \to U = \1 + \ker \eps, \quad x \mapsto \sum_{k = 0}^\infty {x^k \over k!} $$
is an analytic diffeomorphism with the analytic inverse 
$\Log(x) := \sum_{k = 1}^\infty {(-1)^{k-1}\over k} (x - \1)^k.$
We thus obtain on $\ker \eps$ a global analytic multiplication 
$$ x * y := \Log(\Exp x \Exp y) $$
given by the BCH series, so that its values lie in the 
completion $\hat L$ of the free Lie algebra $L$ generated by $x_1,\ldots, x_n$, 
which is a closed Lie subalgebra of $\hat F$. 
\qed

In Section XI below, we shall discuss more aspects of projective limits of finite- and 
infinite-dimensional Lie groups.

\Example IV.1.16. We recall the group $G = E \rtimes_\alpha \R$ from Example~II.5.9(b), 
where $E = \R^\N$ and $\alpha(t) = e^{tD}$ with the diagonal operator 
$D(z_n) = (n z_n)$. Then the Lie group structure on $G$ is analytic and 
the explicit formula shows that 
$$\exp_G \: \L(G) \to G, \quad (v,t) \mapsto (\beta(t)v,t) $$ 
is analytic. Further, it is a smooth diffeomorphism whose inverse 
$$\log_G \: G \to \L(G), \quad (v,t) \mapsto (\beta(t)^{-1}v,t) $$ 
is smooth but not analytic. In fact, 
$\beta(t)^{-1}e_n = {tn \over e^{nt} - 1}e_n,$
and for $n$ fixed, the radius of convergence of the Taylor series of this function in $0$ 
is ${2\pi \over n}$. A similar argument shows that the corresponding 
global multiplication $x*y := \log_G(\exp_G(x)\exp_G(y))$ 
on $\L(G)$ is smooth but not analytic. With Robart's Criterion 
(Theorem IV.1.7), this follows from the fact that the power series 
$$ \sum_{k = 0}^\infty t^k D^k e_n = (1-tn)^{-1} e_n $$
is not convergent for $|t| > {1\over n}$. 
\qed

For details concerning the following results, we refer to [GN06] (see also 
[Mil82, 4.3] for some of the statements). 
Minor modifications of the corresponding argument 
for finite-dimensional, resp., Banach--Lie groups lead to the following lemma, 
which in turn is the key to the following theorem: 

\Lemma IV.1.17. Let $G$ be a locally exponential Lie group. 
For $x, y \in \L(G)$, we have the {\it Trotter Product
Formula} 
$$ \exp_G(x + y) = \lim_{n \to \infty} \Big(\exp_G\big(\frac{x}{n}\big)
\exp_G\big(\frac{y}{n}\big)\Big)^n$$
and the {\it Commutator Formula} 
$$ \exp_G([x,y]) = \lim_{n \to \infty} \Big(
\exp_G\big(\frac{x}{n}\big)\exp_G\big(\frac{y}{n}\big)
\exp_G\big(-\frac{x}{n}\big)\exp_G\big(-\frac{y}{n}\big)\Big)^{n^2}. 
\qeddis 

\Theorem IV.1.18. {\rm(Automatic Smoothness Theorem)} Each continuous homomorphism 
$\phi \: G \to H$ of locally exponential (BCH) Lie groups is smooth (analytic). 
\qed

For (local) Banach--Lie groups, Theorem IV.1.18 can already be found in [Bir38], 
and for BCH--Lie groups in [Gl02c] (see also [Mil84] for the statement without proof). 
The special case of one-parameter groups $\R \to A^\times$, where $A$ is a Banach algebra 
is due to {\smc Nagumo} ([Nag36]) and {\smc Nathan} ([Nat35]). 
Mostly such ``automatic smoothness'' theorems concern continuous 
homomorphisms $\phi \: G\to H$ of Lie groups, where $H$ is a Lie group 
with an exponential function for which each continuous one-parameter group 
is of the form $\gamma_x(t) = \exp_H(tx)$ and $G$ is locally exponential. 
We then obtain a map $\L(\phi) \: \L(G) \to \L(H)$ by 
$\phi(\exp_G(tx)) = \exp_H(t \L(\phi)x)$ for $t \in \R$ and $x \in \L(G)$, 
and then it remains to show that $\L(\phi)$ is continuous and linear.

The following theorem is due to {\smc Maissen} for Banach--Lie groups ([Mais62, Satz~10.3]): 

\Theorem IV.1.19. Let 
$G$ and $H$ be Lie groups and $\psi \: \L(G) \to \L(H)$ a continuous 
homomorphism of Lie algebras. 
Assume that $G$ is locally exponential and $1$-connected 
and that $H$ has a smooth exponential function. 
Then there exists a unique morphism of Lie groups $\phi \: G \to H$ with 
$\L(\phi) = \psi$. 

\Proof. (Idea) Let $U_\g \subeq \g = \L(G)$ be a convex balanced $0$-neighborhood 
mapped diffeomorphically by the exponential function to an open subset $U_G$ of $G$. 

First one observes that $\psi^*\kappa_{\L(H)} = \psi \circ \kappa_{\L(G)}$ (cf.\ (2.5.5)). 
For the map $$f \: U_G \to H, \quad \exp_G(x) \mapsto \exp_H(\psi(x)),$$ 
this leads to 
$\delta(f) = f^*\kappa_H = \psi \circ \kappa_G,$
showing that the ${\L(H)}$-valued $1$-form $\psi \circ \kappa_G$ is locally integrable. 
Since this form on $G$ is left invariant and $G$ is $1$-connected, it 
is globally integrable to a function 
$\phi  \: G \to H$ with $\phi(\1) = \1$ and $\delta(\phi) = \psi \circ \kappa_G$ 
(Theorem~III.1.2). Now Proposition~II.4.1 implies that $\phi$ is a group 
homomorphism, and by construction $\L(\phi) =\alpha_\1 =  \psi$. 
\qed

Since we do not know if all Lie groups with an exponential function are regular, 
the preceding theorem is not a consequence of Theorem~III.1.5. 

\Corollary IV.1.20. If $G_1$ and $G_2$ are locally exponential $1$-connected 
Lie groups with isomorphic Lie algebras, then $G_1$ and $G_2$ are isomorphic. 
\qed

\Remark IV.1.21. It is instructive to compare Corollary~IV.1.20 with 
the corresponding statement for regular Lie groups (Corollary~III.1.6). 
They imply that there exists for each locally convex Lie algebra $\g$ at most one 
$1$-connected Lie group $G$ which is regular and at most one $1$-connected 
locally exponential Lie group $H$ with $\L(G) = \L(H) = \g$. 
The regularity of $G$ implies that 
$\id_\g$ integrates to a smooth homomorphism 
$\phi \: H \to G$, but we do not know if there is a morphism 
$\psi \: G \to H$ with $\L(\psi) = \id_\g$ (cf.\ Problem~III.1). 

Presently we do not know if all locally exponential Lie groups (modeled on Mackey complete 
spaces) are regular, therefore it is still conceivable that there might be locally 
exponential Lie algebras which are the Lie algebra of a $1$-connected 
regular Lie group and a non-isomorphic $1$-connected locally exponential Lie group 
which is not regular. 
\qed

\Remark IV.1.22. Theorem IV.1.18 implies in particular that being a locally exponential 
Lie group is a {\sl topological property}: Any topological group $G$ 
carries at most one structure of a locally exponential Lie group. We thus 
adjust our terminology in the sense that we call a topological group {\it locally 
exponential} if it carries a locally exponential Lie group structure compatible 
with the topology. 

Forgetting the differentiable structure on $G$, it becomes an interesting issue how 
to recover it. In view of Theorem~IV.1.18, we recover the Lie algebra 
$\L(G)$, as a set, by identifying $x \in \L(G) \cong T_\1(G)$ with the 
corresponding one-parameter group $\gamma_x(t) = \exp_G(tx)$. 
Starting from $G$, as a topological group, we may then put 
${\frak L}(G) := \Hom_c(\R,G)$, the set of continuous homomorphisms $\R \to G$. 
The scalar multiplication of ${\frak L}(G)$ can be written as 
$$ (\lambda\alpha)(t) := \alpha(\lambda t), \quad \lambda \in \R, \alpha \in \Hom_c(\R,G), 
\leqno(4.1.4) $$ 
and, in view of Lemma IV.1.17, addition and Lie bracket may be written on the  level 
of one-parameter groups by 
$$ (\alpha + \beta)(t) := \lim_{n \to \infty} 
\Big(\alpha({t \over n})\beta({t\over n})\Big)^n\leqno(4.1.5) $$
and 
$$ [\alpha,\beta](t^2) := \lim_{n \to \infty} \Big(
\alpha({t\over n})\beta({t\over n})
\alpha(-{t\over n})\beta(-{t\over n})\Big)^{n^2}. \leqno(4.1.6) $$
We can also recover the topology on $\L(G)$ as the compact open 
topology on $\fL(G)$, and the exponential function as the evaluation map 
$$ \exp_G \: \Hom_c(\R,G) \to G, \quad \gamma \mapsto \gamma(1). \leqno(4.1.7) $$
\qed

In [HoMo05/06], {\smc Hofmann} and {\smc Morris} use (4.1.4-7) as the starting 
point in the investigation of a remarkable class of topological groups: 

\Definition IV.1.23. Let $G$ be a topological group 
and ${\frak L}(G) := \Hom_c(\R,G)$ the set of one-parameter groups, 
endowed with the compact open topology. 
Then $G$ is said to be 
{\it a topological group with Lie algebra} if 
the limits in (4.1.5/6) exist for $\alpha,\beta \in \Hom_c(\R,G)$ and 
define elements of ${\frak L}(G)$, addition and bracket are continuous maps 
${\frak L}(G) \times {\frak L}(G) \to {\frak L}(G),$
and ${\frak L}(G)$ is a real Lie algebra
with respect to the scalar multiplication (4.1.4), the addition (4.1.5), 
and the bracket (4.1.6). This implies 
that $\fL(G)$ is a topological Lie algebra. 
The exponential function of $G$ is defined by (4.1.7). 
\qed

In [BCR81], {\smc Boseck}, {\smc Czichowski} and {\smc Rudolph} define smooth 
functions on a topological group in terms of restrictions to  
one-parameter groups, which leads them to (4.1.4-7), together 
with the assumption that ${\frak L}(G)$ can be identified with the set of 
derivations of the algebra of germs of smooth functions in $\1$ 
([BCR81, Sect.~1.5]). 

We have just seen that any locally exponential Lie group 
is a topological group with Lie algebra. Since $\R$ is connected, a topological 
group $G$ has a Lie algebra if and only if its identity component $G_0$ does. 
In [HoMo05, Th.~2.3], it is also 
observed that any abelian topological group is a group with Lie algebra,
where the addition on ${\frak L}(G)$ is pointwise multiplication and the bracket 
is trivial (cf.\ Problem IV.7). 

\Theorem IV.1.24. Each $2$-step nilpotent topological group has a Lie 
algebra. 

\Proof. (Sketch) The commutator map $c \: G \times G \to Z(G)$ is an 
alternating bihomomorphism. Then direct calculations lead to the formulas 
$$ (\alpha + \beta)(t) = \alpha(t) \beta(t) c(\alpha(t), \beta(-{\textstyle{t\over 2}})) 
\quad \hbox{ and } \quad  [\alpha,\beta](t) = c(\alpha(1), \beta(t)), $$
which can be used to verify all requirements. 
\qed

We shall return to topological groups with Lie algebras 
in our discussion of projective limits in 
Section~X.

\Remark IV.1.25. We have seen in  Corollary~IV.1.20 that 
a $1$-connected locally exponential Lie group 
$G$ is completely determined up to isomorphism (as a topological group) 
by its Lie algebra. 

If $G$ is connected but not simply connected, then we have 
a universal covering morphism $q_G \: \tilde G \to G$ 
and $\ker q_G \cong \pi_1(G)$ is a discrete central subgroup of $\tilde G$ 
with $G \cong \tilde G/\ker q_G$. 
It is easy to see that two discrete central subgroups 
$\Gamma_1, \Gamma_2 \subeq Z(\tilde G)$ lead to isomorphic quotient groups 
$\tilde G/\Gamma_1$ 
and $\tilde G/\Gamma_2$ if and only if there exists an automorphism
 $\phi \in \Aut(\tilde G) \cong \Aut(\L(G))$ with 
$\phi(\Gamma_1) = \Gamma_2$. Therefore the isomorphism classes of 
connected Lie groups $G$ with a given Lie algebra $\g \cong \L(G)$ are parametrized 
by the orbits of $\Aut(\g) \cong \Aut(\tilde G)$ in the set of discrete central 
subgroup of $\tilde G$. 

If $G_0$ and a discrete group $\Gamma$ are given, then the determination of all Lie groups 
$G$ with identity component $G_0$ and component group $\pi_0(G)\cong \Gamma$ 
corresponds to the classification of all Lie group extensions 
$$ \1 \to G_0 \into G \onto \Gamma \to \1, $$
i.e., to a description of the set $\Ext(\Gamma,G_0)$. Extension problems of this type 
are discussed in Section V.1 below. 
\qed

\subheadline{IV.2. Locally exponential Lie algebras} 

We now turn to the Lie algebras which are candidates for Lie algebras 
of locally exponential Lie groups. We call these Lie algebras ``locally exponential''. 
They are defined by the requirement that some $0$-neighborhood carries a local 
group structure in ``canonical''  coordinates, i.e., the additive one-parameter groups 
$t \mapsto tx$, should also be one-parameter groups for the local group structure 
(cf.\ [Bir38], [Lau55]). 

\Definition IV.2.1. A locally convex Lie algebra $\g$ 
is called {\it locally exponential} if there exists a circular 
convex open $0$-neighborhood $U \subeq \g$ and an open subset 
$D \subeq U \times U$ on which we have a smooth map 
$$ m_U \: D \to U, \quad (x,y) \mapsto x * y $$
such that $(U,D,m_U,0)$ is a local Lie group (Definition~II.1.10) satisfying:
\litem{(E1)} For $x \in U$ and $|t|,|s|, |t+s|\leq 1$, we have 
$(tx,sx) \in D$ with 
$tx * sx = (t+s)x.$ 
\item{(E2)} The second order term in the Taylor expansion of $m_U$ is 
$b(x,y) = \frac{1}{2}[x,y].$ 

\nin The Lie algebra $\g$  is called {\it exponential} if $U = \g$ and $D = \g \times \g$.
\qed

Since any local Lie group on an open subset of a locally convex space 
$V$ leads to a Lie algebra structure on $V$ (Definition~II.1.10), 
condition (E2) only ensures that $\g$ is the Lie algebra of the local group 
(cf.\ Remark~II.1.8). 

Using exponential coordinates, we directly get: 

\Lemma IV.2.2. The Lie algebra $\L(G)$ of a locally exponential Lie group $G$ 
is locally exponential. 
\qed

\Definition  IV.2.3. We call a locally exponential Lie algebra $\g$ 
{\it enlargeable} if it is integrable to a locally exponential Lie group $G$. 
As we shall see in Remark~IV.2.5 below, this 
is equivalent to the enlargeability of some associated local group in $\g$. 
\qed

\Examples IV.2.4. (a) All BCH--Lie algebras, hence in particular all 
Banach--Lie algebras and therefore all finite-dimensional 
Lie algebras are locally exponential (Example~IV.1.6). 

A different existence proof for the local multiplication on a Banach--Lie 
algebra $\g$ is given by {\smc Laugwitz} ([Lau56]): As a first step, 
we observe that 
$\kappa_\g(x) := {\1 - e^{-\ad x} \over \ad x}$ defines a smooth map 
$\kappa_\g \: \g \to {\cal L}(\g)$ with $\kappa_\g(0) = \id_\g$, 
so that its values are  invertible on some $0$-neighborhood. 
We consider $\kappa_\g$ as a $\g$-valued $1$-form on $\g$. 
Then one verifies that $\kappa_\g$ satisfies the Maurer--Cartan equation, 
which implies the existence of an open $0$-neighborhood $U$ such that 
for each $x \in U$ the (partial differential) equation 
$$ f^*\kappa_\g = \kappa_\g, \quad f(0) = x $$ 
has a unique solution $f_x$ on $U$. For $x,y$ close to $0$, the composition 
$f_x \circ f_y$ is then defined on some $0$-neighborhood and satisfies 
$f_x\circ f_y(0) = f_x(y) = f_{f_x(y)}(0)$ as well as 
$(f_x \circ f_y)^*\kappa_\g = f_y^*f_x^*\kappa_\g = \kappa_\g$, which implies 
$f_x \circ f_y = f_{f_x(y)}$ on some $0$-neighborhood. 
For $x * y := f_x(y)$, this leads to the associativity condition 
$$ x * (y * z) = (x * y) * z $$
on some $0$-neighborhood in $\g$, hence to a local group structure.
As $\kappa_\g$ satisfies $\kappa_\g(x)x = x$ for each $x \in \g$, the curves 
$t \mapsto tx$ are local one-parameter groups. This corresponds to 
condition (E1). 

(b) If $\g$ is locally exponential and $M$ a compact manifold, then 
$C^\infty(M,\g)$ is also locally exponential with respect to 
$(x*y)(m) := x(m) * y(m)$ for all $m \in M$ and $x,y$ close to $0$.
\qed

\Remark IV.2.5. A similar reasoning as in the proof of Theorem~IV.1.19 implies that 
any morphism $f \: \g \to \h$ of locally exponential Lie algebras 
satisfies  $f(x*y) = f(x)*f(y)$ for $x,y$ close to $0$. 
Applying this to $f = \id_\g$ shows in particular that the 
Lie algebra $\g$ determines the germ of the local multiplication 
$x*y$ (cf.\ [Lau56] for the Banach case). 
We know that this multiplication need not be analytic, 
not even if it is defined on all of $\g \times \g$ (Example~IV.1.16). 
\qed

Suppose that $\g$ is an exponential 
Lie algebra for which the 
group $(\g,*)$ is regular. Then $(\g,*)$ is the unique $1$-connected regular 
Lie group with Lie algebra $\g$. If $G$ is any $1$-connected Lie group (regular or not) 
with $\L(G) = \g$, and $G$ has an exponential function $\exp_G \: \L(G) =\g\to G$, 
then $\exp_G$ is a group homomorphism $(\g,*) \to G$ 
(cf.\ Propositions~II.4.1 and II.5.7). 
The regularity of $(\g,*)$ implies the existence of a 
unique homomorphism $\Log_G \: G \to (\g,*)$ with $\L(\Log_G) = \id_\g$,  
and the uniqueness assertion of Proposition~II.4.1 yields 
$\Log_G \circ \exp_G = \id_\g$ and $\exp_G \circ \Log_G = \id_G$. 
Since on any Mackey complete nilpotent Lie algebra $\g$, the BCH multiplication 
defines a regular Lie group structure 
([GN06]), these arguments lead to the following theorem:

\Theorem IV.2.6. If $G$ is a connected nilpotent Lie group with a smooth exponential 
function and $\L(G)$ is Mackey complete, then the exponential 
function 
$$ \exp_G \: (\L(G),*) \to G $$
is a covering morphism of Lie groups. In particular,  
$G \cong (\L(G),*)/\Gamma$ for a discrete subgroup $\Gamma \subeq \z(\g)$, isomorphic 
to $\pi_1(G)$. Moreover, $G$ is regular and locally exponential. 
\qed

This generalizes a result of {\smc Michor} and {\smc Teichmann} who showed in 
[MT99] that any connected regular abelian Lie group $G$ is of the form 
$\L(G)/\Gamma$ for a discrete subgroup $\Gamma \cong \pi_1(G)$ of $\L(G)$. 
Related results can be found in [Ga96], where locally exponential abelian 
Fr\'echet--Lie groups are studied as projective limits of Banach--Lie groups. 

Without any completeness assumption we obtain the following 
very natural intrinsic characterization of the BCH series as the only Lie 
series which leads on nilpotent Lie algebras to a group multiplication 
satisfying (E1). 

\Proposition IV.2.7. If $G$ is a $1$-connected exponential nilpotent Lie group, 
then $G \cong (\L(G),*)$, where $*$ denotes the (polynomial) BCH multiplication on 
$\L(G)$. 
\qed

We have already seen that the Lie bracket on a locally exponential Lie algebra $\g$ 
determines the germ of the corresponding local multiplication (Remark~IV.2.5), 
hence in particular its Taylor series in $(0,0)$. 
The preceding proposition is the key step to the following theorem, 
identifying this series as the BCH series. In the Banach context, 
the corresponding result is due to Birkhoff ([Bir38]). 
Its statement can be found in [Mil82] as Lemma~4.4, 
with the hint that it can be proved with the methods used in [HS68] 
in the finite-dimensional case, which is based on formula (4.1.1). 
Since the spectra of the operators $\ad x$ and $\ad y$ are possibly  unbounded, 
formula (4.1.1) makes no sense for general locally exponential Lie algebras. 
The situation is much better if $\g$ is nilpotent. In this 
case, the operators $e^{\ad x}$ are unipotent, so that 
$\psi(e^{\ad x}e^{\ad y})$ is a polynomial in $x$ and $y$. 
The reduction to this case is a key point in the proof of the following theorem. 

\Theorem IV.2.8. {\rm(Universality Theorem)} If $\g$ is locally exponential, then the 
Taylor series of the local multiplication $x*y$ in $(0,0)$ is the BCH series. 

\Proof. (Sketch) A central idea is the following. For each Lie algebra 
we obtain by extension of scalars from $\R$ to the two-dimensional algebra $\R[\eps]$ 
of dual numbers ($\eps^2 = 0$), the Lie algebra $T(\g) := \g \otimes_\R \R[\eps]$. 
One can show that $T(\g)$ is also locally exponential. The local multiplication 
$m_{T(\g)}$ is the tangent map of the local multiplication $m_{\g}$ of $\g$ and 
$U_{T(\g)} = T(U_\g) = U_\g \times \g$ is the tangent bundle of $U_\g \subeq \g$. 

Iterating this procedure, we obtain a sequence of locally exponential Lie algebras 
$$ T^n(\g) := \g\otimes \R[\eps_1,\ldots, \eps_n] \quad \hbox{ with } 
\quad \eps_i \eps_j = \eps_j \eps_i, \ \eps_i^2 = 0, $$
whose local multiplication $T^n(m_\g)$ induces a global multiplication on the nilpotent ideal 
$J \trile T^n(\g)$ which is the kernel of the augmentation map $T^n(\g) \to \g$. 
Applying Proposition~IV.2.7 to $(J, T^n(m_\g))$ now shows that the 
$n$-th order Taylor polynomial of $m_\g$ in $(0,0)$ is given by the BCH series. 
\qed
 
For the discussion of quotients of locally exponential groups below, the 
following theorem is crucial: 

\Theorem IV.2.9. {\rm(Quotient Theorem for locally exponential Lie algebras)} 
Let $\g$ be a locally exponential Mackey complete Lie algebra and $\n \trile \g$ a closed ideal. 
Then $\g/\n$ is locally exponential if and only if
\litem{(1)} $\n$ is stable, i.e., $e^{\ad x}(\n) = \n$ for each $x \in \g$, and 
\litem{(2)} $\kappa_\g(x)\n = \n$ for each $x$ in some $0$-neighborhood of $\g$.

If $\n$ is the kernel of a morphism $\phi \: \g \to \h$ of locally exponential Lie algebras, 
then $\n$ is locally exponential and 
both conditions are satisfied, so that $\phi$ factors through the quotient map 
$q \: \g \to \g/\n$ to an injective morphism 
$\oline\phi \: \g/\n \to \h$ of locally exponential Lie algebras. 
\qed

The preceding result is nicely complemented by the following observation 
on extensions: 

\Theorem IV.2.10. {\rm([GN07])} Let $\g$ be a locally exponential Lie algebra and 
$q \: \hat\g \to \g$ be a central extension, i.e., a quotient morphism with 
central kernel $\z$. If $\z$ is Mackey complete, then $\hat\g$ is locally exponential. 
\qed

\Remark IV.2.11. In [Hof72/75], {\smc K.~H.~Hofmann} advocates an approach to Banach--Lie groups 
by defining a Banach--Lie group as a topological group possessing 
an identity neighborhood isomorphic (as a topological local group) 
to the local group defined by the BCH multiplication in a 
$0$-neighborhood of a Banach--Lie algebra. A key point of this 
perspective is that Banach--Lie groups form a full sub-category of the category 
of topological groups (Theorem~IV.1.18, Remark~IV.1.22). 
Due to the analyticity of the BCH multiplication, this approach works 
quite well for Banach--Lie groups, and also for the larger class of 
BCH Lie groups which behave in almost all respects like Banach--Lie groups. 
Although one may think that one can adopt a similar point of 
view for locally exponential 
Lie groups, a closer analysis of the arguments used in this theory to pass from 
infinitesimal to local information  
shows that the behavior of locally exponential groups is far from being controlled 
by topology. Actually the arguments we use are much closer to the 
original approach to Lie theory via differential equations (cf.\ Examples~IV.2.4). 

Giving up the analyticity requirement of the local 
multiplication in an identity neighborhood implies that we have 
to work in a smooth category to prove uniqueness assertions. 
The Maurer--Cartan form and the Uniqueness Lemma are the fundamental tools. 
In the analytic context, one can often argue quite directly by analytic continuation. 
\qed

\subheadline{IV.3. Locally exponential Lie subgroups} 

It is a well-known result in finite-dimensional Lie theory 
that each closed subgroup $H$ of a Lie group $G$ carries a natural Lie group 
structure turning it into a submanifold of $G$ (see [vN29] for 
closed subgroups of $\GL_n(\R)$). This becomes already 
false for closed subgroups of infinite-dimensional Hilbert spaces, which contain 
contractible subgroups not containing any smooth arc. 
Therefore additional assumptions on closed subgroups are needed to 
make them accessible by Lie theoretic methods. Since we already know that each 
topological group carries at most one locally exponential Lie group 
structure, it is clear that a closed subgroup deserves to be called a 
{\it Lie subgroup} if it is a locally exponential Lie group with respect to the 
induced topology. For Banach--Lie groups, this is precisely {\smc Hofmann}'s approach, 
and for several of the results described below, Banach versions can be found in 
[Hof75]. 

\subheadline{Lie subgroups and factor groups}

\Lemma IV.3.1. For every 
closed subgroup $H$ of the locally exponential Lie group $G$, we have 
$$ \L^d(H) = \L^e(H) := \{ x \in \L(G) \: \exp_G(\R x) \subeq H \},  $$
and this is a closed Lie subalgebra of $\L(G)$. 

\Proof. The equality $\L^e(H) = \L^d(H)$ follows from 
$\lim_{n \to \infty} \gamma\big({\textstyle{t\over n}}\big)^n = \exp_G(t \gamma'(0))$
for each curve $\gamma \: [0,1] \to H$ with $\gamma(0) = \1$ which is differentiable in $0$ 
because we can write it on some interval $[0,\eps]$ as $\gamma = \exp_G \circ \eta$ 
with some $C^1$-curve $\eta$ in $\L(G)$ with $\eta(0) = 0$. 
The closedness follows from the obvious closedness of $\L^e(H)$. 
\qed

In the following we shall keep the notation 
$\L^e(H) = \{ x \in \L(G) \: \exp_G(\R x) \subeq H \}$ for a subgroup 
$H$ of a Lie group $G$ with an exponential function, because if $H$ is not closed or 
$G$ is not locally exponential, it is not clear that this set coincides with 
$\L^d(H)$

\Definition IV.3.2. A closed subgroup $H$ of a locally exponential Lie group $G$ 
is called a {\it locally exponential Lie subgroup}, or simply a {\it Lie subgroup},  if 
$H$ is a locally exponential Lie group with respect to the induced topology 
(cf.\ Remark~IV.1.22). 
\qed

A Banach version of the following theorem is Proposition~3.4 in [Hof75]. 
\Theorem IV.3.3. For a closed subgroup $H$ of the locally exponential Lie group $G$ 
the following are equivalent: 
\litem{(1)} $H$ is a locally exponential Lie group. 
\litem{(2)} There exists an open $0$-neighborhood $V \subeq \L(G)$ such that 
$\exp_G\res_V$ is a diffeomorphism onto an open $\1$-neighborhood in $G$ and 
$\exp_G(V \cap \L^e(H)) = \exp_G(V) \cap H.$

In particular, each locally exponential Lie subgroup is a submanifold of $G$.
\qed

\Proposition IV.3.4. If $\phi \: G' \to G$ is a morphism of
locally exponential Lie groups and $H \subeq G$ is a locally exponential Lie subgroup, 
then $H' := \phi^{-1}(H)$ is a locally exponential Lie subgroup. In particular,  
$\ker \phi$ is a locally exponential Lie subgroup of $G'$. 
\qed

The preceding proposition implies in particular that if a quotient 
$G/N$ by a closed normal subgroup $N$ is locally exponential, then $N$ is a 
locally exponential Lie subgroup. But the converse is more subtle: 

\Theorem IV.3.5. {\rm(Quotient Theorem for locally exponential groups)} 
For a closed normal subgroup  $N \trile G$ the following are equivalent: 
\litem{(1)} $G/N$ is a locally exponential Lie group. 
\litem{(2)} $N$ is a locally exponential Lie subgroup and $\L(G)/\L^e(N)$ is a locally exponential Lie algebra. 
\litem{(3)} $N$ is a locally exponential Lie subgroup and 
$\kappa_{\L(G)}(x)(\L^e(N)) = \L^e(N)$ for $x \in \L(G)$ sufficiently close to $0$. 

\msk 
If $N$ is the kernel of a morphism $\phi \: G \to H$ of locally exponential Lie groups, 
then $G/\ker \phi$ is a Lie group, so that $\phi$ factors through a quotient map 
$G \to G/\ker \phi$ and an injective morphism 
$\oline\phi \: G/\ker \phi$ of locally exponential Lie groups. 
\qed

Since quotients of BCH--Lie algebras are BCH--Lie algebras, no matter whether they are complete or 
not ([Gl02c, Th.~2.20]), we get the following corollary, whose Banach version is also contained 
in [Hof75, Prop.~3.6] and [GN03]. 

\Corollary IV.3.6. {\rm(Quotient Theorem for BCH--Lie groups)} 
A closed normal subgroup $N$ of a BCH--Lie group $G$ is a BCH--Lie group if and only if 
the quotient $G/N$ is a BCH--Lie group. 
\qed

If $\phi \: G \to H$ is an injective morphism of locally exponential Lie groups, 
then the preceding theorem provides no additional information. In Section IV.4, we 
shall encounter this situation for integral subgroups of $G$. 
The following example provides a bijective morphism $\phi$ for which $\L(\phi)$ is not surjective.  
The only way to avoid this pathology is to assume that 
$\L(G)$ is separable 
(cf.\ Theorems IV.4.14/15 below). That not all surjective morphisms 
of locally exponential Lie groups are quotient morphisms can already be seen 
for surjective continuous linear maps between non-Fr\'echet spaces. 

\Example IV.3.7. We give an example of a proper closed subalgebra $\h$ of the Lie algebra 
$\L(G)$ of some Banach--Lie group $G$ for which $\la \exp \h \ra = G$ ([HoMo98, p.157]). 

We consider the abelian Lie group $\g := \ell^1(\R,\R) \times \R$, where 
the group structure is given by addition. We write $(e_r)_{r
\in \R}$ for the canonical topological basis elements of
$\ell^1(\R,\R)$. Then the subgroup $D$ generated by the pairs 
$(e_r, -r)$, $r \in \R$, is closed and discrete, so that 
$G := \g/D$ is an abelian Lie group. Now we consider the closed
subalgebra $\h := \ell^1(\R,\R)$ of $\g$. As $\h + D = \g$, we have 
$H := \exp_G \h = G$, and therefore $(0,1) \in \L^e(H) \setminus \h.$ 

The map $\phi := \exp_G \res_\h \: (\h,+) \to G$ is a surjective morphism of Lie groups 
for which $\L(\phi)$ is the inclusion of the proper subalgebra $\h$. 
\qed

That for connected Banach--Lie groups $G$ the center $Z(G) = \ker \Ad$ is a 
locally exponential Lie subgroup follows immediately from Proposition~IV.3.4 (cf.\ [Lau55]). 
For non-Banach--Lie algebras $\g$, 
$\Aut(\g)$ carries no natural Lie group 
structure, so that Proposition~IV.3.4 does not apply. This makes the 
following theorem quite remarkable. The crucial point in its proof is to show that 
for  the exponential function 
$$ \Exp \: \g/\z(\g) \to \Aut(\g), \quad x \mapsto e^{\ad x}, \leqno(4.3.1) $$
the point $0$ is isolated in $\Exp^{-1}(\id_\g)$ (cf.\ Problems II.4 and IX.1). 

\Theorem IV.3.8. Let $\g$ be a locally exponential Lie algebra. 
Then the adjoint group 
$G_{\rm ad} := \la e^{\ad \g} \ra \subeq \Aut(\g)$ carries the structure of a 
locally exponential Lie group whose Lie algebra is the quotient $\g_{\rm ad} := \g/\z(\g)$ 
and {\rm(4.3.1)} its exponential function. 
\qed

Combining the preceding theorem with Proposition IV.3.4, we get: 
\Corollary IV.3.9. If $G$ is a connected locally exponential Lie group, then its 
center $Z(G) = \ker \Ad$ is a locally exponential Lie subgroup. 
\qed

\subheadline{Algebraic subgroups} 

The concept of an algebraic subgroup of a Banach--Lie algebra, introduced by 
{\smc Harris} and {\smc Kaup} ([HK77]),  
provides very convenient criteria which in many concrete
cases can be used to verify that a closed subgroup $H$ of a Banach--Lie group
is a Banach--Lie subgroup. 

\Definition IV.3.10. Let $A$ be a unital Banach algebra. 
A subgroup $G \subeq A^\times$ is called {\it algebraic} if there exists a
$d \in \N_0$ and a set 
${\cal F}$ of Banach space-valued polynomial functions on $A\times A$ 
of degree $\leq d$ such that 
$$ G = \{ g \in A^\times \: (\forall f \in {\cal F})\ f(g,g^{-1}) = 0\}. 
\qeddis 

\Theorem IV.3.11. {\rm([HK77], [Ne04b, Prop.~IV.14])} 
Every algebraic subgroup $G
\subeq A^\times$ of the unit 
group $A^\times$ of a Banach algebra $A$ is a Banach--Lie subgroup. 
\qed

\Corollary IV.3.12. Let $E$ be a Banach space and $F \subeq E$ a closed
subspace. Then 
$$ \GL(E,F) := \{ g \in \GL(E) \: g(F) = F\} $$ 
is  a Banach--Lie subgroup of $\GL(E)$. 
\qed

\Corollary IV.3.13. Let $E$ be a Banach space and $v\in E$. Then 
$$ \GL(E)_v := \{ g \in \GL(E) \: g(v) = v\} $$
is  a Banach--Lie subgroup of $\GL(E)$. 
\qed

\Corollary IV.3.14. For each continuous 
bilinear map $\beta \: E \times E \to E$ on a Banach space $E$, the group 
$$ \Aut(E,\beta) := \{ g \in \GL(E)\: \beta \circ (g \times g) = g \circ \beta\} $$
is a Banach--Lie subgroup of $\GL(E)$ with Lie algebra 
$$ \der(E,\beta) := \{ D \in \gl(E) \: (\forall v,w \in E)\ 
D.\beta(v,w) = \beta(D.v,w) + \beta(v,D.w)\}. 
\qeddis 

\Corollary IV.3.15. For each bilinear map $\beta \: E \times E \to \K$, the group 
$$ \OO(E,\beta) := \{ g \in \GL(E)\: \beta \circ (g \times g) = \beta\} $$
is a Banach--Lie subgroup of $\GL(E)$ with Lie algebra 
$$ \oo(E,\beta) := \{ D \in \gl(E) \: (\forall v,w \in E)\ 
\beta(D.v,w) + \beta(v,D.w)=0\}. 
\qeddis

\subheadline{Closed subgroups versus Lie subgroups} 

For finite-dimensional Lie groups, closed subgroups are
Lie subgroups (cf.\ [vN29]), but for Banach--Lie groups this is no longer true. What
remains true is that locally compact subgroups (which are closed in particular) are Lie subgroups. 
For subgroups of Banach algebras the following theorem is due to 
{\smc Yosida} ([Yo36]) and for general Banach--Lie groups to  {\smc Laugwitz} 
([Lau55], [Les66]).
Although the arguments in the Banach case do not immediately carry over 
because unit spheres for seminorms are no longer bounded, one can 
use Gl\"ockner's Implicit Function Theorem ([Gl03a]) to get: 

\Theorem IV.3.16. {\rm([GN06])} Each locally compact subgroup of a locally exponential 
Lie group is a finite-dimensional Lie subgroup. 
\qed

\Remark IV.3.17. How bad closed subgroups may behave is
illustrated by the following example ([Hof75, Ex.~3.3(i)]): 
We consider the real Hilbert space $G := L^2([0,1],\R)$ as a
Banach--Lie group. Then the subgroup $H := L^2([0,1],\Z)$ of all those
functions which almost everywhere take values in $\Z$ is a closed
subgroup. Since the one-parameter subgroups of $G$ are of the form 
$\R f$, $f \in G$, we have $\L^e(H) = \{0\}$. 
On the other hand, the group $H$ is arc-wise connected. It is 
contractible, because the map $F \: [0,1] \times H \to H$ given by 
$$ F(t,f)(x) := \cases{ 
f(x) & for $0 \leq x \leq t$ \cr 
0& for $t < x \leq 1$ \cr} $$
is continuous with $F(1,f) = f$ and $F(0,f) = 0$. 
\qed

\subheadline{IV.4. Integral subgroups} 

It is a well-known result in finite-dimensional Lie theory 
that for each subalgebra $\h$ of the Lie algebra $\g = \L(G)$ of a finite-dimensional 
Lie group $G$, there exists a Lie group $H$ with Lie algebra $\h$ together with 
an injective morphism of Lie groups $\iota \: H \to G$ for which 
$\L(\iota) \: \h \to \g$ is the inclusion map. As a group,  
$H$ coincides with $\la \exp \h \ra$, the analytic subgroup corresponding to $\h$, 
and $\h$ can be recovered from this subgroup as the set 
$\L^e(H) = \{ x \in \L(G) \: \exp(\R x) \subeq H \}.$
This nice and simple theory of analytic subgroups and integration of 
Lie algebra inclusions $\h \into \L(G)$ becomes much more subtle 
for infinite-dimensional Lie groups. 
Even for Banach--Lie groups some pathologies arise. Here any inclusion 
$\h \into \L(G)$ of Banach--Lie algebras integrates to an ``integral'' subgroup 
$H \into G$, but if the Banach--Lie algebra $\h$ is not separable, then 
it may happen that $\h$ cannot be recovered from the abstract subgroup $H$ of $G$. 
In Example~IV.3.7, it even occurs that $\h \not= \L(G)$ and $H = G$.

\Definition IV.4.1. Let $G$ be a Lie group with an exponential function, so that 
we obtain for each $x \in \L(G)$ an automorphism 
$e^{\ad x} := \Ad(\exp_G x) \in \Aut(\L(G))$. 
A  subalgebra $\h \subeq \L(G)$ is called {\it stable} if 
$$ e^{\ad x}(\h) = \Ad(\exp_G x)(\h) = \h \quad \hbox{ for all } x \in \h. $$
An ideal $\n \trile \L(G)$ is called {\it a stable ideal} if 
$e^{\ad x}(\n) = \n$ for all $x \in \L(G)$. 
\qed 

The following lemma shows that stability of kernel and range 
is a necessary requirement for the integrability of a homomorphism 
of Lie algebras. 

\Lemma IV.4.2. If $\phi \: G \to H$ is a morphism of Lie groups with an 
exponential function, then 
$\im(\L(\phi))$ is a stable subalgebra of $\L(H)$, and 
$\ker(\L(\phi))$ is a stable ideal of $\L(G)$. 

\Proof. We have $\phi \circ \exp_G = \exp_H \circ \L(\phi),$
which leads to 
$$ \eqalign{ \L(\phi) \circ e^{\ad x} 
&= \L(\phi) \circ \Ad(\exp_G x) 
= \L(\phi \circ c_{\exp_G x}) 
= \L(c_{\phi(\exp_G x)} \circ \phi) \cr
&= \Ad(\exp_H \L(\phi)x) \circ \L(\phi) 
= e^{\ad \L(\phi)(x)} \circ \L(\phi).  \cr} $$
We conclude in particular that $\im(\L(\phi))$ is a stable subalgebra and 
that $\ker(\L(\phi))$ is a stable ideal. 
\qed
 
\Lemma IV.4.3. {\rm(a)} Each closed subalgebra which is 
finite-dimensional or finite-codimensional is stable. 

{\rm(b)} Let $\g$ be a BCH Lie algebra. 
Then each closed subalgebra $\h \subeq \g$ and each closed ideal $\n \trile \g$ 
is stable. 

{\rm(c)} If $\h \into \g$ is a continuous inclusion of locally convex 
Lie algebras such that 
for each $x \in \h$ the operators $\ad_\h x$ and $\ad_\g x$ are integrable on $\h$, resp., 
$\g$, then $\h$ is a stable subalgebra of $\g$. 

\Proof. (a) (cf.\ [Omo97, Lemma~III.4.8]) If $\h$ is finite-dimensional, then the Uniqueness Lemma 
implies for $x \in \h$ the relation $e^{\ad x}\res_\h = e^{\ad x\,\res\,\,\h}$. 
If $\h$ is finite-codimensional and $q \: \g \to \g/\h$ the projection map, 
then the curve $\gamma(t) := q(e^{t \ad x}y)$ satisfies the  linear ODE 
$\gamma'(t) = \ad_{\g/\h}(x) \gamma(t)$, hence vanishes for $y \in \h$. 

(b) Since $\g$ is BCH, the map 
$x \mapsto e^{\ad x}y = x * y * (-x)$ is analytic on some open $0$-neighborhood, 
hence given by the power series $\sum_{n = 0}^\infty {1\over n!} (\ad x)^n y$. 
Therefore the closedness of $\h$ implies that for $x$ close to $0$ we have 
$e^{\ad x}(\h) \subeq \h$. This implies stability. 
A similar argument yields the stability of closed ideals. 

(c) We apply Lemma~II.3.10 to see that $e^{\ad_\g x}y= e^{\ad_\h x}y \in \h$ 
holds for $x,y \in \h$. 
\qed

In view of the preceding lemma, stability causes no problems for BCH--Lie 
algebras, but the condition becomes crucial in the non-analytic context. 

\Example IV.4.4. The first example of a closed Lie subalgebra $\h$ 
of some $\L(G)$ which does not integrate to any group homomorphism is due 
to {\smc H.~Omori} (cf.\ [Mil84, 8.5]). 

We consider the group $G := \Diff(\T^2)$ of diffeomorphisms of the $2$-dimensional 
torus and use coordinates $(x,y)\in [0,1]^2$ corresponding to the identification 
$\T^2 \cong \R^2/\Z^2$. Then 
$$ \h := \Big\{ f {\partial \over \partial x} + g {\partial \over \partial y} \: 
{1\over 2} \leq x \leq 1 \Rarrow g(x,y) = 0 \Big\} $$
is easily seen to be a closed Lie subalgebra of $\g = {\cal V}(\T^2)$. 
The vector field $X := {\partial \over \partial x}$ generates the 
smooth action $\alpha \: \T \to \Diff(\T^2)$ of $\T$ on 
$\T^2$ given by $[z].([x],[y]) = ([x+z],[y])$. 
This vector field is contained in $\h$, but 
$$ e^{{1\over 2} \ad X}\h = \Ad(\alpha({\textstyle{1\over 2}}))\h 
= \Big\{ f {\partial \over \partial x} + g {\partial \over \partial y} \: 
0 \leq x \leq {1\over 2}  \Rarrow g(x,y) = 0 \Big\} \not=\h.$$
This shows that $\h$ is not stable and hence that it does not integrate to any 
subgroup of $\Diff(\T^2)$ with an exponential function. 
\qed

\Example IV.4.5. Let $E := C^\infty(\R,\R)$ and 
consider the one-parameter group $\alpha \: \R \to \GL(E)$ given by 
$\alpha(t)(f)(x) = f(x+t)$. Then $\R$ acts smoothly on $E$, so that we can 
form the corresponding semi-direct product group 
$G := E \rtimes_\alpha \R.$
This is a Lie group with a smooth exponential function given by 
$$ \exp_G(v,t) 
= \Big( \int_0^1 \alpha(st).v\, ds, t\Big), 
\quad \hbox{ where } \quad 
 \Big(\int_0^1 \alpha(st).v\, ds\Big)(x) = \int_0^1 v(x + st)\, ds. $$
The Lie algebra $\g = \L(G)$ has the corresponding semi-direct product structure 
$\g = E \rtimes_D \R$ with $Dv = v'$, i.e., 
$$ [(f,t), (g,s)] = (tg'-sf', 0). $$

In $\g$, we now consider the subalgebra $\h := E_{[0,1]} \rtimes \R$, 
where 
$$ E_{[0,1]} := \{ f \in E \: \supp(f) \subeq [0,1]\}. $$
Then $\h$ is a closed subalgebra of $\g$. It is not stable because 
$\alpha(-t) E_{[0,1]} = E_{[t,t+1]}.$
The subgroup of $G$ generated by $\exp_G \h$ contains $\{0\} \rtimes \R$, 
$E_{[0,1]}$, and hence all subspaces $E_{[t,t+1]}$, which implies that 
$\la \exp_G \h \ra = C^\infty_c(\R) \rtimes \R$. 

Lemma~IV.4.2 implies that the inclusion $\h \into \g$ does not 
integrate to a homomorphism $\phi \: H \to G$ for any Lie group $H$ with an exponential 
function and $\L(H) = \h$. 
\qed 

\Example IV.4.6. Let 
$E \subeq C^\infty(\R,\C)$ be the closed subspace of $1$-periodic functions, 
$\mu \in \R^\times$, and 
consider the homomorphism $\alpha \: \R \to \GL(E)$ given by 
$$ (\alpha(t)f)(x) := e^{\mu t} f(x + t). $$
That the corresponding $\R$-action on $E$ is smooth follows from the smoothness of 
the translation action and one can show that the group 
$G := E \rtimes_\alpha \R$ is exponential with 
Lie algebra $\g = E \rtimes_D \R$ and 
$Df = \mu f + f'$ ([GN06]). In particular, the product 
$x * y := \exp_G^{-1}(\exp_G(x)\exp_G(y))$ is globally defined on $\g$. 

Let $M \subeq [0,1]$ be an open subset which is not dense and put 
$$E_M := \{ f \in E \: f\res_M = 0\}.$$ 
Then $E_M$ is a closed subspace of 
$E$ with $DE_M \subeq E_M$ but 
$\alpha(t)(E_M) = E_{M-t}\not\subeq E_M$ for some $t \in \R$. 
Therefore $\h_M := E_M \rtimes_D \R \subeq \g = E \rtimes_D \R$ is a closed subalgebra of 
the exponential Fr\'echet--Lie algebra $\g$ which is not stable. 
Since $e^{\ad x}y = x * y * (-x)$ for all $x,y \in \g$, this implies in particular  
that $\h$ is not closed under the $*$-multiplication. 
\qed

\Definition IV.4.7. 
Let $G$ be a Lie group. An {\it integral subgroup} 
is an injective morphism $\iota \: H \to G$ of 
Lie groups such that $H$ is connected and the differential $\L(\iota) \: \L(H) \to \L(G)$ 
is injective. 
\qed

\Remark IV.4.8. Let $\iota \: H \to G$ be an integral subgroup and assume that $H$ and $G$ 
have exponential functions. Then the relation 
$$\exp_G \circ \L(\iota)= \iota \circ \exp_H  \leqno(4.4.1) $$
 implies that 
$\ker(\L(\iota)) = \L(\ker \iota) = \{0\},$
so that 
$\L(\iota) \: \L(H) \to \L(G)$ is an injective morphism of topological 
Lie algebras, which implies in particular that $\h := \im(\L(\iota))$ is a stable 
subalgebra of $\L(G)$ (Lemma IV.4.2). 
Moreover, (4.4.1) shows that the subgroup $\iota(H)$ of $G$ coincides, 
as a set, with the subgroup 
$\la \exp_G \h\ra$ of $G$ generated by $\exp_G\h$. 
Therefore a locally exponential integral subgroup can be viewed as  
a locally exponential Lie group structure on the subgroup of $G$ generated by
$\exp_G \h$. 
\qed

\Theorem IV.4.9. {\rm(Integral Subgroup Theorem)} Let 
$G$ be a Lie group with a smooth exponential function 
and $\alpha \: \h \to \L(G)$ an injective morphism of topological 
Lie algebras, where $\h$ is locally exponential. 
We assume that the closed subgroup 
$$\Gamma := \{x \in \z(\h) \: \exp_G(\alpha(x))= \1\} $$
is discrete. Then 
there exists a locally exponential integral subgroup $\iota \: H \to G$ with 
$\L(H) = \h$ and $\L(\iota) = \alpha$. 
In particular, $\h$ is integrable to a locally 
exponential Lie group. 
\qed

The discreteness of $\Gamma$ is automatic in the following two special 
cases (cf.\ Problem~II.4). 

\Corollary IV.4.10. Let $G$ be a Lie group with a smooth exponential 
function and $\h$ locally exponential. 
Then any injective morphism 
$\alpha \: \h \to \L(G)$ integrates to a locally exponential integral subgroup if 
one of the following two conditions is satisfied: 
\litem{(1)} $\z(\h)$ is finite-dimensional.
\litem{(2)} $G$ is locally exponential. 
\qed

Corollary IV.4.10(1) is a substantial generalization of the main result 
of [Pe95b] which assumes that $G$ is regular and $\h$ is Banach. 
Other special cases can be found in many places in the literature, such as 
[MR95, Th.~2]. The versions given in [RK97, Th.~2] and [Rob97, Cor.~2] 
contradict the existence 
of unstable closed subalgebras in locally exponential Lie algebras (Example IV.4.6). 
For Banach--Lie groups it is contained in [EK64], and for BCH--Lie groups 
in [Rob97]. 

\Remark IV.4.11. (a) In [Rob97], {\smc Robart} gives a criterion for the existence of 
integral subgroups of a locally exponential Lie group $G$ 
for a prescribed injective morphism $\alpha \: \h \to \g = \L(G)$: 
The Lie algebra morphism $\alpha$ can be integrated to an integral subgroup if 
and only if $\h/\z(\h)$ is the Lie algebra of a locally exponential Lie group 
isomorphic to $H_{\rm ad} :=\la e^{\ad \h} \ra \subeq \Aut(\h)$ with 
exponential function as in (4.3.1).  In view of Theorems IV.2.10 and IV.3.8, for Mackey complete 
Lie algebras, this condition is equivalent to $\h$ being locally exponential. 
This argument shows in particular, that Robart's concept of a {\it Lie algebra of the 
first kind} coincides with our concept of a locally exponential Lie algebra. 
In the light of this remark, Theorem~5 in [Rob97] can be read as a version of 
our Corollary~IV.4.10(2), whereas 
Corollary~IV.4.10(1) corresponds to his Theorem~8. We do not understand 
the precise meaning of his remark concerning a generalization 
to the case where $\z(\h)$ is infinite-dimensional by simply refining the topology. 

(b) Even for a closed subalgebra $\h \subeq \g := \L(G)$, the condition that it 
is locally exponential is quite subtle. 
It means that for $x,y$ sufficiently close to $0$ in $\h$, 
we have $x * y \in \h$. If $\g$ is BCH and $\h$ is closed, this is clearly satisfied, 
but if $\g$ is not BCH, not every closed subalgebra satisfies this condition  
because it implies stability (Example~IV.4.5). 
\qed

To verify this condition, one would like 
to show that the integral curve $\gamma(t) := x * ty$ of the left 
invariant vector field $y_l$ through $x$ does not leave the closed subspace 
$\h$ of $\g$. This leads to the necessary condition 
$T_0(\lambda_x)(\h) \subeq \h$, which, under the assumption that $\h$ is stable, 
means that the operator $\kappa_\g(x) = \int_0^1 e^{-t\ad x}\, dt$ 
satisfies $\kappa_\g(x)(\h) = \h$ for $x \in \h$ sufficiently close to $0$ 
(cf.\ Theorem~IV.3.8 and Problem IV.5 below). 

\Example IV.4.12. (a) Applying 
Corollary IV.4.10 to the CIA $\hat F$ obtained by completing the free associative 
algebra in $n$ generators $x_1,\ldots, x_n$ (Example IV.1.15), it follows that the 
closed Lie subalgebra generated by $x_1,\ldots, x_n$, i.e., the completion 
of the free Lie algebra, integrates to a subgroup. As $\hat F$ is topologically 
isomorphic to $\R^\N$, each closed subspace is complemented 
([HoMo98, Th.~7.30(iv)]), so that the existence of the corresponding integral subgroup 
could also be obtained by the methods developed in [Les92, Sect.\ 4] which require complicated 
assumptions on groups and Lie algebras. 

(b) If $K$ is a Banach--Lie group with Lie algebra $\k$ and $M$ a compact manifold, 
then the group $C^\infty(M,K)$ is BCH (Theorem IV.1.12), so that 
the Integral Subgroup Theorem also applies to each closed subalgebra 
$\h \subeq C^\infty(M,\k)$ (cf.\ [Les92, Sect.\ 4]). 
\qed

\Remark IV.4.13. (a) In [La99], S.\ Lang calls a subgroup $H$ of a Banach--Lie 
group $G$ a ``Lie subgroup'' if $H$ carries a Banach--Lie group structure for
which there exists an immersion $\eta \: H \to G$. This 
requires the Lie algebra $\L(H)$
of $H$ to be a closed subalgebra of $\L(G)$ which is complemented in the
sense that there exists a closed vector space complement (cf.\ Remark~I.2.7). 
From that, it follows that his Lie subgroups coincide with 
the integral subgroups with closed complemented Lie algebra (cf.\ Corollary~IV.4.10). 

The advantage of Lang's more restrictive concept is that 
for a closed complemented Lie subalgebra $\h \subeq \L(G)$ one obtains 
the existence of corresponding integral subgroups from the 
Frobenius Theorem for Banach manifolds ([La99, Th.~VI.5.4]). But 
it excludes in particular closed non-complemented 
subspaces of Banach spaces. 

 (b) The most restrictive 
concept of a Lie subgroup is the one used in [Bou89, Ch.~3], 
where a ``Lie subgroup'' of a Banach--Lie group $G$ 
is a Banach--Lie subgroup $H$ with the additional 
property that $\L(H)$ is complemented, i.e., 
$H$ is required to be a split submanifold of $G$. 
This concept has the advantage that it implies that the 
quotient space $G/H$ carries a natural manifold structure for which the
quotient map $q \: G \to G/H$ is a submersion 
([Bou89, Ch.~3, \S 1.6, Prop.~11]). However, the 
condition that $\L(H)$ is complemented is very hard to check in concrete situations 
and, as the Quotient Theorem and the Integral Subgroup Theorem show, 
not necessary. 

(c) For closed subalgebras which are not necessarily complemented, the Integral 
Subgroup Theorem can already be found in [Mais62] who also shows that 
kernels are Banach--Lie subgroups and that $G/N$ is a Lie group if $N$ 
is a Banach--Lie subgroup with complemented Lie algebra as in (b). 
This case is also dealt with in [Hof68], [Hof75, Th.~4.1], and a local version can be 
found in [Lau56]. 
\qed

The following theorem 
generalizes [Hof75, Prop.~4.3] from 
Banach--Lie groups to locally exponential ones, which is quite straightforward ([GN06]). 
The necessity of the separability assumption follows from Example~IV.3.7. 

\Theorem IV.4.14. {\rm(Initiality Theorem for integral subgroups)} 
Let $G$ be a Lie group with a smooth exponential 
function $\exp_G \: \L(G) \to G$ which is injective on some $0$-neighborhood. 
Further let $\iota_H \: H \into G$ be a locally exponential integral 
subgroup whose Lie algebra $\L(H)$ is separable. 
Then the subgroup $\iota_H(H)$ of $G$ satisfies 
$$  \L^e(\iota_H(H)) = \{ x \in \L(G) \: \exp_G(\R x) \subeq \iota_H(H)\} = \im(\L(\iota_H)). $$
In particular, the surjectivity of $\iota_H$ implies the surjectivity of $\L(\iota_H)$. 

If, in addition, $G$ is locally exponential and $\L(H)$ is a closed subalgebra of 
$\L(G)$, then $\iota_H \: H \to G$ is an initial Lie subgroup of $G$. 
\qed

\Theorem IV.4.15. {\rm([Hof75, Prop.~4.6])} Let $G$ be a separable Banach--Lie 
group and assume that ${\iota_H \: H \to G}$ is an integral subgroup with closed range. Then 
$\iota_H$ is an embedding. In particular, $H$ is a Banach--Lie subgroup of $G$. 
\qed

We conclude this section with a discussion of initial Lie subgroup structures on 
closed subgroups of Banach--Lie groups and locally convex spaces. 

\Theorem IV.4.16. {\rm(Initiality Theorem for closed subgroups of Banach--Lie 
groups)} Let $G$ be a Banach--Lie group and $H \subeq G$ a closed 
subgroup. Then $H$ carries the structure of an initial Lie subgroup with 
Lie algebra $\L^e(H) = \L^d(H)$. Its identity component is an 
integral Lie subgroup for the closed Lie subalgebra $\L^e(H)$ of $\L(G)$. 

\Proof. (Idea) We know from Lemma~IV.3.1 that $\L^d(H) = \L^e(H)$ is a closed Lie subalgebra 
of $\L(G)$. 
Let $\iota \: H_0 \to G$ be the corresponding integral Lie subgroup. 
Since each smooth curve $\gamma \: [0,1] \to \iota(H_0)\subeq H$ satisfies 
$\delta(\gamma)(t) \in \L^d(H)$ for each $t$ and $H_0$ is regular (Remark~II.5.4), 
$\gamma$ is smooth as a curve to $H_0$, and this further 
permits us to conclude that $H_0$ is initial and coincides with the 
smooth arc--component of $H$. Now one uses Corollary~II.2.3 to extend the 
Lie group structure to all of $H$. 
\qed

\Theorem IV.4.17. {\rm(Initiality Theorem for closed subgroups of locally convex 
spaces)} Let $E$ be a locally convex space and $H \subeq (E,+)$ a closed subgroup. 
Then $H$ carries an initial Lie group structure, for which 
$H_0 = \L^d(H) = \L^e(H)$ is the largest vector subspace contained in~$H$. 

\Proof. (Idea) For each curve $\alpha \in C^1_*(I,E)$ with $\im(\alpha) \subeq H$ we have 
$tx = \lim_{n \to \infty} n\alpha({t\over n}) \in H,$
which leads to $\L^d(H) = \L^e(H)$, a closed subspace of $E$. 
For each $C^1$-curve $\gamma \: [0,1] \to E$ with range in $H$,  
all tangent vectors lie in $\L^d(H)$. This implies that $\gamma$ lies in a 
coset of $\L^d(H)$. Defining the Lie group structure in such a way 
that $\L^d(H)$ becomes an open subgroup of $H$, it follows easily 
that $H$ is initial (Corollary~II.2.3). 
\qed

\Remark IV.4.18. (Stability and distributions) (a) For 
a subset $D \subeq {\cal V}(M)$, we call the subset 
$\Delta_D \subeq T(M)$ defined by $\Delta_D(m) := \span \{ X(m) \: X \in D\}$ 
the corresponding {\it smooth distribution}. Conversely, 
we associate to a (smooth) distribution $\Delta \subeq T(M)$ the subspace 
$D_\Delta := \{X \in {\cal V}(M) \: (\forall m \in M)\ X(m) \in \Delta\}$. 
The distribution $\Delta$ is said to be {\it involutive} if $D_\Delta$ is a 
Lie subalgebra of ${\cal V}(M)$. A smooth distribution $\Delta_D$ is 
called {\it $D$-invariant} 
if it is preserved by the local flows generated by elements of $D$,  
and {\it integrable} it possesses (maximal) integral submanifolds through each 
point of~$M$.  

{\smc Sussman}'s Theorem asserts that 
$\Delta_D$ is integrable if and only if it is $D$-invariant ([Sus73, Th.~4.2]). 
As a special case, where all subspaces $\Delta_D(m)$, $m \in M$, are 
of the same dimension, we obtain {\smc Frobenius}' Theorem. 
The invariance condition on $\Delta_D$ implies that it is involutive, 
but the converse does not hold. E.g., consider on $M = \R^2$ the 
set $D$, consisting of two vector fields 
$$ {\partial \over \partial x_1}, \quad f(x_1) {\partial \over \partial x_2} 
\quad \hbox{ with } \quad f^{-1}(0) = ]-\infty,0] $$
(see also Example~IV.4.4). 

If $M$ is analytic and $D$ consists of analytic vector fields, 
then {\smc Nagano} shows in [Naga66] that the involutivity of the 
corresponding distribution is sufficient for the existence of 
integral submanifolds. 

(b) The invariance condition for a distribution is quite analogous 
to the stability condition for a Lie subalgebra $\h \subeq \L(G)$.
If, furthermore, $G$ is analytic with an analytic exponential function, then each 
closed subalgebra is stable, which is analogous to Nagano's result 
(cf.\ Lemma~IV.4.3(b)). 

To relate this to stability of Lie algebras of vector fields, assume that $M$ is compact. 
If $\h \subeq {\cal V}(M)$ is a stable subalgebra (not necessarily  closed), 
then the corresponding 
distribution $\Delta_\h$ is stable, and Sussman's Theorem implies that its maximal 
integral submanifolds are the orbits of the subgroup 
$H := \la \exp_{\Diff(M)}(\h)\ra$ of $\Diff(M)$ generated by the flows of 
elements of $\h$ ([KYMO85, Sect.\ 3.1]). 

If, conversely, we start with a smooth distribution $\Delta$, 
then the closed space $\h_\Delta$ of all vector fields with values in $\Delta$ 
is a Lie subalgebra if and only if $\Delta$ is involutive. 
Furthermore, it is not hard to see that Sussman's Theorem implies that $\Delta$ 
is involutive if and only if $\h_\Delta$ is stable. 
If this is the case, the corresponding subgroup $H_\Delta$ of $\Diff(M)$ satisfies 
$\h_\Delta \subeq \L^e(H_\Delta) \subeq \L^d(H_\Delta) \subeq \h_\Delta.$
Now is is a natural question whether 
$H_\Delta$ carries the structure of a Lie group. 
For more details and related examples, we refer to Section~3.1 in [KYMO85]. 
\qed

\subheadline{Open Problems for Section IV} 

\Problem IV.1. Show that for each subgroup $H$ 
of a locally exponential Lie group $G$, the set 
$\L^e(H) = \{ x \in \L(G)  \: \exp_G(\R x) \subeq H\}$
is a Lie subalgebra of $G$ or find a counterexample. 

If $G$ is finite-dimensional, then Yamabe's Theorem implies that 
the arc-component of $H$ is an integral subgroup (Remark II.6.5) which proves 
the assertion in this case. If $H$ is closed, then $\L^e(H)$ is a Lie subalgebra 
by Lemma~IV.3.1. 

We further observe that $\L^e(H)$ is invariant under all operators $\Ad(\exp_G(x)) 
= e^{\ad x}$ for $x \in \L^e(H)$. If $\L^e(H)$ is closed (which is the case for 
each closed subgroup) and closed under addition, then it is a closed 
vector subspace of $\L(G)$, and for $x,y \in \h$ it 
contains the derivative of the curves $t \mapsto e^{t\ad x}y$ in $0$. 
This implies that it is a Lie subalgebra. 
\qed

\Problem IV.2. Show that for each closed subgroup $H$ of a locally exponential Lie group 
$G$ the closed Lie subalgebra $\L^e(H) \subeq \L(G)$ is locally exponential. 
\qed

\Problem IV.3. Find an example of a locally exponential normal Lie subgroup $N$ of a 
locally exponential Lie group $G$ for which $\L(G)/\L(N)$ is not locally exponential 
or prove that it always is. In view of Theorem~IV.3.5, this would imply that 
$G/N$ is locally exponential. 
\qed 

\Problem IV.4. (One-parameter groups and local exponentiality) 
Let $\alpha \: \R \to \GL(E)$ define a smooth action of 
$\R$ on the Mackey complete 
locally convex space and $D := \alpha'(0)$ be its infinitesimal 
generator. We then obtain a $2$-step solvable Lie group 
$G := E \rtimes_\alpha \R$ with 
the product 
$$ (v,t) (v',t') = (v + \alpha(t).v', t + t') $$
and the Lie algebra 
$\g = E \rtimes_D \R$. Characterize local exponentiality of 
$G$ in terms of the infinitesimal generator $D$. 

Writing the exponential function as 
$\exp_G(v,t) = (\beta(t).v, t)$ with
$\beta(t) = \int_0^1 \alpha(st)\, ds,$
we obtain the curve $\beta \: \R \to {\cal L}(E)$. 
We are looking for a characterization of those operators 
$D$ for which there exists some $T > 0$ such that 
\litem{(1)} $\beta(]{-T},T[) \subeq \GL(E)$, and 
\litem{(2)} $\tilde\beta \: ]{-T},T[\times E \to E, (t,v) 
\mapsto \beta(t)^{-1}v$ is smooth.  

\nin Note that $(t,v) \mapsto \beta(t)v$ is always smooth. 
If $E$ is a Banach space, then $G$ is a Banach--Lie group, 
hence locally exponential. In this case, $D$ is a bounded operator 
and we have for each $t \not=0$: 
$$ \beta(t) 
= {1\over t} \int_0^t e^{sD}\, ds 
= {1\over t} {e^{tD} - \1 \over D} 
= {e^{tD} - \1 \over tD} 
= \sum_{k = 0}^\infty {t^k \over (k+1)!} D^k. $$
As $\beta \: \R\to {\cal L}(E)$ is analytic w.r.t.\ to the operator norm on 
${\cal L}(E)$, $\beta(0) = \1$, and $\GL(E)$ is open, 
conditions (1) and (2) follow immediately. 
Moreover, the Spectral Theorem implies that 
$$ \Spec(\beta(t)) = 
\Big\{ {e^{t\lambda} - 1 \over t\lambda} \: \lambda \in \Spec(D)\Big\}, $$
which means that $\beta(t)$ is invertible for 
$|t|  < {2 \pi \over \sup \{ \Im(\lambda) \: \lambda \in \Spec(D)\}}.$ 
\qed

\Problem IV.5. (Invariant subspaces) 
Let $\alpha \: \R \to \GL(E)$ be a smooth action of 
$\R$ on the Mackey complete locally convex space $E$ and $\beta(t)$ as in 
Problem IV.4. 

Suppose that $F \subeq E$ is a closed invariant subspace. Then we also 
have $\beta(t)(F) \subeq F$ for each $t \in \R$. 
Assume that for some $\eps > 0$ the operator $\beta(t)$ is invertible for $|t| \leq \eps$. 
Show that $\beta(t)^{-1}(F) \subeq F$ for $|t| \leq \eps$ 
or find a counterexample. Note that this is trivially the case if 
$F$ is of finite dimension or codimension. 
\qed

\Problem IV.6. Show that BCH--Lie groups are regular. 
In [Rob04], {\smc Robart} has obtained substantial results in this direction, 
including that for each BCH--Lie group $G$ with Lie algebra $\g$ 
and each smooth path 
$\xi \in C^\infty(I,\g)$, the initial value problem 
$$ \eta(0)= \id_{\g}, \quad \eta'(t) = [\eta(t),\xi(t)] $$
has a solution in ${\cal L}(\g)$. Unfortunately, it is not clear whether 
these solutions define curves in $\GL(\g)$. 
\qed

\Problem IV.7. Show that each nilpotent topological group is a topological 
group with Lie algebra in the sense of Definition IV.1.23 (cf.\ Theorem IV.1.24). 
\qed

\Problem IV.8. In Theorem IV.3.3, we have seen that a locally exponential Lie subgroup 
$H$ of a locally exponential Lie group $G$ is a submanifold, where the submanifold 
chart in the identity can be obtained from the exponential function of $G$. 

It is an interesting question whether every Lie group $H$ which is a submanifold 
of a locally exponential Lie group $G$ is in fact a locally 
exponential Lie group. This is true if $G$ is a Banach--Lie group because 
this property is inherited by every subgroup which is a submanifold. 

This point also concerns the use of the term ``Lie subgroup,'' which would 
also be natural for subgroups which are submanifolds. 
\qed

\Problem IV.9. Develop a theory of algebraic subgroups for CIAs in the context of 
locally exponential, resp., BCH--Lie groups. A typical question such a theory should 
answer is: For which linear actions of a locally exponential Lie group $G$ on a 
locally convex space $E$ are the stabilizers $G_v$, $v \in E$, locally exponential 
Lie subgroups? 
\qed

\Problem IV.10. Show that Theorems~IV.4.15/16 
remain valid for locally exponential Lie groups. 
\qed

\Problem IV.11. Let $G$ be a regular Lie group and $H \subeq G$ a closed subgroup. 
Then $\L^e(H)$ is a closed subset of $\L(G)$, stable under scalar multiplication. 
On the other hand, $\L^d(H)$ is a Lie subalgebra containing $\L^e(H)$. 
Do these two sets always coincide? 
If, in addition, $G$ is $\mu$-regular, this follows from Lemma~III.2.7.
\qed

\Problem IV.12. (a) Let $G$ be a Lie group with a smooth exponential function. 
Find examples where the Trotter Formula and/or the Commutator Formula 
(Lemma~IV.1.17) do not hold. 
For which classes of groups (beyond the locally exponential ones) 
are these formulas, or the more general Lemma~III.2.7, valid? What about the group $\Diff_c(M)$? 

(b) What can be said about the sequence of power maps $p_n(x) = x^n$ in a local Lie group? 
In local coordinates with $0$ as neutral element, it is interesting 
to consider for an element $x$ the sequence 
$\big({x\over n}\big)^n$. Does it converge (to $x$)?
\qed

\Problem IV.13. Let $\Delta$ be an integrable distribution on the compact manifold~$M$ 
and $H_\Delta$ the corresponding subgroup of $\Diff(M)$ preserving the maximal 
integral submanifolds of $\Delta$ (Remark~IV.4.18(b)). Show that $H_\Delta$ is a 
regular Lie group. 
\qed

\Problem IV.14. (cf.\ [Rob97, p.837, Prop.~3]) Let 
$G$ be a $\mu$-regular Lie group and $\h \subeq \L(G)$ a closed 
stable subalgebra. Does $\h$ integrate to an integral Lie subgroup?  
Since product integrals converge in $G$, for two smooth 
curves $\alpha, \beta \: [0,1]\to \h$, the curve 
$$ (\alpha * \beta)(t) := \alpha(t) + \Ad(\gamma_\alpha(t)).\beta(t) $$
has values in $\h$, which leads to a Lie group structure on $C^\infty(I,\h)$ 
with Lie algebra $C^\infty([0,1],\h)$, where the bracket is given by 
$$ [\xi,\eta](t) := \Big[\xi(t), \int_0^t \eta(\tau)\, d\tau\Big] + 
\Big[\int_0^t \xi(\tau)\, d\tau, \eta(t)\Big] $$
([Rob97, Th.~9]). 
The map $E \: \alpha \mapsto \gamma_\alpha(1)$ is a group homomorphism, so that 
the problem is to see that the quotient group 
$(C^\infty(I,\h),*)/\ker E \cong \im(E)$ carries a natural Lie group structure. 
\qed

\sectionheadline{V. Extensions of Lie groups}  

\nin In this section, we turn to some general results on extensions of infinite-dimensional 
Lie groups. In Section V.1, we explain how an extension of $G$ by $N$ is described in 
terms of data associated to $G$ and $N$. This description is easily adapted from the abstract group 
theoretic setting (cf.\ [ML68]). 
In Section V.2, we describe the appropriate 
cohomological setup for Lie theory and explain criteria for the integrability of 
Lie algebra cocycles to group cocycles. This is applied in Section V.3 to 
integrate abelian extensions of Lie algebras to corresponding group 
extensions. 

\subheadline{V.1. General extensions} 

\Definition V.1.1. An {\it extension of Lie groups} is a short exact sequence 
$$ \1 \to N \sssmapright{\iota} \hat G \sssmapright{q} G \to \1 $$
of Lie group morphisms, for which 
$\hat G$ is a smooth (locally trivial) principal 
$N$-bundle over $G$ with respect to the right action of $N$ given by 
$(\hat g,n) \mapsto \hat gn$. In the following, we identify 
$N$ with the subgroup $\iota(N) \trile \hat G$. 

We call two extensions 
$N \into \hat G_1 \onto G$ and 
$N \into \hat G_2 \onto G$ of the Lie group $G$ by the 
Lie group $N$ {\it equivalent} if there exists a Lie group morphism 
$\phi \: \hat G_1 \to \hat G_2$ such that the following diagram commutes: 
$$ \matrix{
 N & \into& \hat G_1 & \onto & G \cr 
\mapdown{\id_N} & & \mapdown{\phi} & & \mapdown{\id_G} \cr 
 N & \into& \hat G_2 & \onto & G. \cr } $$
It is easy to see that any such $\phi$ is an isomorphism of
 Lie groups and that we actually obtain an equivalence relation. 
We write $\Ext(G,N)$ for the set of equivalence classes of
Lie group extensions of $G$ by~$N$.\footnote{$^1$}{\eightrm From the description of Lie 
group extensions as in Theorem~V.1.4 below, one obtains cardinality 
estimates showing that 
the equivalence classes actually form a set.} 

We call an extension $q \: \hat G \to G$ with $\ker q = N$ 
{\it split} if there exists a Lie group morphism 
$\sigma \: G \to \hat G$ with $q \circ \sigma = \id_G$. This implies 
that $\hat G \cong N \rtimes_S G$ for 
$S(g)(n) := \sigma(g)n\sigma(g)^{-1}$. 
\qed 

\Remark V.1.2. A Lie group extension  
$N \into \hat G \onto G$ can also be described in terms of data associated to $G$ and $N$ as 
follows: 
Let $q\: \hat G \to G$ be a Lie group extension of $G$ by $N$.
By assumption, the map $q$ has a smooth local section. 
Hence there exists a global section 
$\sigma \: G \to \hat G$ smooth in an identity neighborhood and {\it normalized} by 
$\sigma(\1) = \1$. Then the map 
$$ \Phi \: N \times G \to \hat G, \quad (n,g) \mapsto n \sigma(g)$$ is
a bijection which restricts to a local diffeomorphism on an identity neighborhood. 
In general, $\Phi$ is not continuous, but we may nevertheless use it to identify 
$\hat G$ with the product set $N \times G$, endowed with the multiplication 
$$ (n,g) (n',g') = (n S(g)(n') \omega(g,g'), gg'), \leqno(5.1.1) $$
where 
$$ S := C_N \circ \sigma \: G \to \Aut(N) \quad \hbox{ for } \quad 
C_N \: \hat G \to \Aut(N),\ C_N(g) = gng^{-1}, \leqno(5.1.2) $$
and 
$$ \omega \: G \times G \to N, \quad 
(g,g') \mapsto \sigma(g) \sigma(g')\sigma(gg')^{-1}. \leqno(5.1.3) 
$$
Note that $\omega$ is smooth in an identity neighborhood and that the map 
$\hat S \: G \times N \to N, (g,n) \mapsto S(g)(n)$ is smooth in a set of the form 
$U_G \times N$, where $U_G$ is an identity neighborhood of $G$. 
The maps $S$ and $\omega$ satisfy the relations 
\litem{(C1)}$\sigma(g) \sigma(g') = \omega(g,g')\sigma(gg')$, 
\litem{(C2)} $S(g)S(g') = C_N(\omega(g,g')) S(gg')$, 
\litem{(C3)} $\omega(g,g') \omega(gg',g'') = S(g)\big(\omega(g',g'')\big) \omega(g, g'g'').$ 
\qed

\Definition V.1.3. Let $G$ and $N$ be Lie groups. A {\it smooth outer action of 
$G$ on $N$} is a map $S \: G \to \Aut(N)$ with $S(\1) = \id_N$ for which 
$$ \hat S \: G \times N \to N, \quad (g,n) \mapsto S(g)(n) $$
is smooth on a set of the form $U_G \times N$, where $U_G \subeq G$ is an open 
identity neighborhood, and for which there exists a map 
$\omega \: G \times G \to N$ with $\omega(\1,\1) = \1$, smooth in 
an identity neighborhood, such that (C2) holds. We call 
$(S,\omega)$ a {\it locally smooth non-abelian $2$-cocycle}. 

We define an equivalence relation on the set of all smooth outer actions of 
$G$ on $N$ by $S' \sim S$ if $S' = (C_N \circ \alpha) \cdot S$ for some 
map $\alpha \: G \to N$ with $\alpha(\1) = \1$, smooth in an identity 
neighborhood. We write $[S]$ for the equivalence class of $S$. 
\qed

Remark V.1.2 implies that for each extension $q \: \hat G \to G$ of Lie groups and 
any section $\sigma \: G \to \hat G$ which is smooth in an identity neighborhood with 
$\sigma(\1) = \1$, (5.1.2) defines a smooth outer action of $G$ on $N$. 
Different choices of such sections lead to equivalent 
outer actions. 

\Theorem V.1.4. Let $G$ be a connected Lie group, $N$ a Lie group and 
$(S,\omega)$ a smooth outer action of $G$ on $N$. Then {\rm(5.1.1)} defines a 
group structure on $N \times G$ if and only if {\rm(C3)} holds. If this is the case, 
then this group carries a unique Lie group structure, 
denoted $N \times_{(S,\omega)} G$, 
for which the identity $N \times G \to N \times_{(S,\omega)} G$ is smooth in an identity 
neighborhood and 
$$q \: N \times_{(S,\omega)} G \to G, \quad (n,g) \mapsto g $$
defines a Lie group extension of $G$ by $N$. 
\qed

All Lie group extensions of $G$ by $N$ arise in this way, so that we obtain a partition 
$$ \Ext(G,N) = \bigcup_{[S]} \Ext(G,N)_{[S]}, $$
where $\Ext(G,N)_{[S]}$ denote the set of equivalence classes of extensions 
corresponding to the equivalence class $[S]$. 

If $N$ is abelian, then each class $[S]$ contains a unique representative $S$, which is a 
smooth action of $G$ on $N$. Fixing $S$, the set $\Ext(G,N)_S$ carries a natural 
abelian group structure, where the addition is given by the {\it Baer sum}: 
For two extensions $q_1 \: \hat G_1 \to G$, $q_2 \: \hat G_2 \to G$ of $G$ by $N$, the 
Baer sum is defined by 
$$ \hat G := \{ (\hat g_1, \hat g_2) \in \hat G_1 \times \hat G_2 \: q_1(\hat g_1) = q_2(\hat g_2)\}
/\Delta_N', \quad \Delta_N' := \{(n,n^{-1}) \:n \in N\}, $$
and the projection map $q(\hat g_1, \hat g_2) := q_1(\hat g_1)$. 
This defines an abelian group structure on the set 
$\Ext(G,N)_S$ whose neutral element is 
the class of the split extension $\hat G = N \rtimes_S G$ (cf.\ [ML63, Sect.~IV.4]). 
In Theorem V.2.8 below, we shall recover this group structure in terms of group 
cohomology. 

\Theorem V.1.5. {\rm([Ne05])} Assume that $Z(N)$ carries an initial Lie subgroup structure 
{\rm(Remark~II.6.5)}. Then 
each class $[S]$ determines a smooth $G$-action on $Z(N)$ by 
$g.z := S(g)(z)$ and the abelian group $\Ext(G,Z(N))_S$ acts simply transitively 
on $\Ext(G,N)_{[S]}$ by 
$$ [H] * [\hat G] := [(\alpha^*\hat G)/\Delta_{Z(N)}'], $$
where $\alpha \: H \to G$ is a Lie group extension of $G$ by the $G$-module $Z(N)$, 
$$\alpha^*\hat G = \{ (\hat g_1, \hat g_2) \in H \times \hat G \: \alpha(\hat g_1) = q(\hat g_2)\} 
\quad \hbox{ and } \quad \Delta_{Z(N)}' := \{(n,n^{-1}) \:n \in Z(N)\}.
\qeddis 

\Examples V.1.6. Interesting classes of extensions of Lie groups arise as follows. 

(a) Projective unitary representations: Let $H$ be a complex Hilbert space, 
$\UU(H)$ its unitary group with center $Z(\UU(H)) = \T \1$, and $\PU(H) := \UU(H)/\T \1$ the 
projective unitary group (all these groups are Banach--Lie groups). If 
$H$ is a complex Hilbert space and $\pi \: G \to \PU(H)$ a projective representation of 
the Lie group $G$ with at least one smooth orbit in the projective space 
$\P(H)$, then the pull-back diagram 
$$ \matrix{ \T = \R/\Z & \into& \UU(H) & \onto & \PU(H)  \cr 
\mapup{=} & & \mapup{} & & \lmapup{\pi}\cr 
\T & \into& \hat G & \onto & G  \cr} $$
defines a central Lie group extension of $G$ by the circle group (cf.\ [Lar99]). 
This leads to a 
partition of the set of equivalence classes of projective unitary representations according to the 
set $\Ext(G,\T)$ of central extensions of $G$ by $\T$ (cf. [Ne02a]). 

(b) Hamiltonian actions: Let $G$ be a connected Lie group, $M$ a locally convex 
manifold and $\omega \in \Omega^2(M,\R)$ a closed $2$-form which is the curvature 
of a pre-quantum line bundle $p \: P \to M$ with connection $1$-form $\theta \in \Omega^1(P,\R)$. 
Assume further that $\sigma \: G \to \Diff(M)$ defines a smooth action of $G$ on $M$ for which 
all associated vector fields $\dot\sigma(x) \in {\cal V}(M)$ are Hamiltonian, i.e., 
the closed $1$-forms $i_{\dot\sigma(x)}\omega$ are exact. Then the range of the map 
$\sigma \: G \to \Diff(M)$ lies in the group $\Ham(M,\omega)$ of Hamiltonian automorphisms of 
$(M,\omega)$, and the diagram 
$$ \matrix{ 
 \T & \into& \Aut(P,\theta) & \onto & \Ham(M,\omega) \cr 
\mapup{=} & & \mapup{} & & \lmapup{\sigma}\cr 
\T & \into& \hat G & \onto & G  \cr} $$
defines a central Lie group extension of $G$ by $\T$ 
(cf.\ [Kos70], [RS81], [NV03]). 

(c) (Extensions by gauge groups) Let 
$q \: P \to M$ be a smooth $K$-principal bundle, where $M$ is compact and 
$K$ is locally exponential (cf.\ Theorem~IV.1.12). 
Further, let $\sigma \: G \to \Diff(M)$ define a smooth action of $G$ on $M$, whose range lies in 
the subgroup $\Diff(M)_{[P]}$ of diffeomorphisms $\phi$ with $\phi^*P \sim P$, i.e., 
fixing the equivalence class $[P]$. Then the diagram 
$$ \matrix{ 
 \Gau(P) & \into& \Aut(P) = \Diff(P)^K & \onto & \Diff(M)_{[P]}  \cr 
\mapup{=} & & \mapup{} & & \lmapup{\sigma}\cr 
\Gau(P) & \into& \hat G & \onto & G  \cr}  $$
defines an extension of $G$ by the gauge group $\Gau(P)$ of this bundle. 

(d) If $q \: \V \to M$ is a finite-dimensional vector bundle and $K = \GL(V)$, the preceding remark 
applies to the corresponding frame bundle, and leads to the diagram 
$$ \matrix{ 
 \Gau(\V) & \into& \Aut(\V) & \onto & \Diff(M)_{[\V]}  \cr 
\mapup{=} & & \mapup{} & & \lmapup{\sigma}\cr 
\Gau(\V) & \into& \hat G & \onto & G . \cr} $$
In view of [KYMO85, p.89], the Lie group $\Aut(\V)$ is $\mu$-regular if $M$ is compact. 

(e) (Non-commutative generalizations; cf.\ [GrNe06], [KYMO85, Sect.~3.2])  
Let $A$ be a CIA and $E$ a finitely generated projective right $A$-module. 
Then the group 
$\GL_A(E)$ of $A$-linear endomorphisms of $E$ is a Lie group, 
which is a subgroup of the larger group 
$$ \GaL(E) := \{ \phi \in \GL(E) \: (\exists \phi_A \in \Aut(A))
(\forall s \in E)(\forall a \in A)\ \phi(s.a) = \phi(s).\phi_A(a)\} $$
of semilinear automorphisms of $E$. These are the linear automorphisms 
$\phi \in \GL(E)$ for which there exists an automorphism $\phi_A$ of $A$ with 
$\phi(s.a) = \phi(s).\phi_A(a)$ for $s \in E, a \in A$. 
Then for each homomorphism $\sigma \: G \to \Aut(A)$, $G$ a connected Lie group, 
whose range lies in the set $\Aut(A)_E$ of those automorphisms of $A$ 
preserving $E$ by pull-backs, the diagram 
$$ \matrix{ 
 \GL_A(E) & \into& \GaL(E) & \onto & \Aut(A)_E  \cr 
\mapup{=} & & \mapup{} & & \lmapup{\sigma}\cr 
\GL_A(E) & \into& \hat G & \onto & G  \cr} $$
defines an extension of $G$ by the linear Lie group $\GL_A(E)$. 

For the special case $A = C^\infty(M,\R)$, finitely generated projective modules 
correspond to vector bundles over $M$ (cf.\ [Ros94], [Swa62]), so that (e) specializes to (d). 
\qed

\Remark V.1.7. Non-abelian extensions of Lie groups also play a crucial 
role in the structural analysis of the group of 
invertible Fourier integral operators of order zero on a compact manifold~$M$ 
([OMYK81], [ARS86a/b]), which is an extension of a group of symplectomorphisms of the complement 
of the zero section in the cotangent bundle $T^*(M)$ whose Lie algebra 
corresponds to smooth functions homogeneous of degree $1$ on the fibers. 
\qed

The following theorem is an important tool to verify that 
given Lie groups are regular (cf.\ [KM97], [OMYK83a, Th.~5.4] in the context 
of $\mu$-regularity, and [Rob04]). A variant of this result for ILB--Lie groups 
is Theorem~3.4 in [ARS86b]. 

\Theorem V.1.8. Let $\hat G$ be a Lie group extension of the Lie group 
$G$ by $N$. Then $\hat G$ is regular if and only if the groups 
$G$ and $N$ are regular. 
\qed

\Remark V.1.9. A typical class of examples illustrating the difference between abelian 
and central extensions of Lie groups arises from 
abelian principal bundles. If $q \: P \to M$ is a smooth principal bundle 
with the abelian structure group $Z$ over the compact connected manifold $M$, 
then the group $\Aut(P) = \Diff(P)^Z$ of all diffeomorphisms of $P$ commuting with 
$Z$ (the automorphism group of the bundle) is an extension of 
the  open subgroup $\Diff(M)_{[P]}$ of $\Diff(M)$ by the gauge group 
$\Gau(P) \cong C^\infty(M,Z)$ of the bundle (Example~V.1.6(c)). Here the conjugation action of 
$\Diff(M)$ on $\Gau(P)$ is given by composing functions with diffeomorphisms. 
Central extensions corresponding to the bundle $q \: P \to M$ are obtained by 
choosing a principal connection $1$-form 
$\theta \in \Omega^1(P,\z)$. Let $\omega \in \Omega^2(M,\z)$ denote the corresponding 
curvature form. Then the subgroup 
$\Aut(P,\theta)$ of those elements of $\Aut(P)$ preserving $\theta$ 
is a central extension of an open subgroup of $\Diff(M,\omega)$, which is substantially 
smaller that $\Diff(M)$. This example shows that the passage from central 
extensions to abelian extensions is similar to the passage from 
symplectomorphism groups to diffeomorphism groups (see also [Ne06a]). 

As the examples of principal bundles over compact manifolds show, abelian 
extensions of Lie groups occur naturally in geometric contexts and in particular 
in symplectic geometry, where the pre-quantization problem is to find for a 
symplectic manifold $(M,\omega)$ a $\T$-principal bundle with curvature $\omega$, 
which leads to an abelian extension of $\Diff(M)_0$ by the group $C^\infty(M,\T)$. 
Conversely, every abelian extension $q \: \hat G \to G$ 
of a Lie group $G$ by an abelian Lie group $A$ is in particular an $A$-principal 
bundle over $G$. This leads to an interesting interplay between abelian extensions 
of Lie groups and abelian principal bundles over (finite-dimensional) 
manifolds. 

A shift from central to abelian extensions occurs naturally as follows: 
Suppose that a connected Lie group $G$ acts on a 
smooth manifold $M$ which is 
endowed with a $Z$-principal bundle $q \: P \to M$ ($Z$ an abelian Lie group) 
and that each element of $G$ lifts 
to an automorphism of the bundle. If all elements of $G$ lift to elements of 
the group $\Aut(P,\theta)$ for some principal connection 
$1$-form $\theta$, then we obtain a central extension as in 
Example~V.1.6(b). But if there is no such connection $1$-form $\theta$, 
then we are forced to consider the much larger abelian extension of 
$G$ by the group $\Gau(P) \cong C^\infty(M,Z)$ or at least some subgroup containing 
non-constant functions. 
The case where $M$ is a restricted 
Gra\3mannian of a polarized Hilbert space and the groups are 
restricted operator groups of Schatten class $p > 2$, resp., mapping groups 
$C^\infty(M,K)$, where $K$ is finite-dimensional and $M$ is a compact manifold 
of dimension $\geq 2$, is discussed in detail in [Mick89] 
(see also [PS86] for a discussion of related topics). 
\qed

\subheadline{V.2. Cohomology of Lie groups and Lie algebras} 

Any good setting for a cohomology theory on Lie groups 
should be fine enough to take the smooth structure into account 
and flexible enough 
to parameterize equivalence classes of group extensions. All these 
criteria are met by the locally smooth cohomology we describe in this 
subsection (cf.\ [Ne02a], [Ne04a]). The traditional approach in finite dimensions 
uses globally smooth cochains ([Ho51]), which is too restrictive in infinite dimensions.

\subheadline{From Lie group cohomology to Lie algebra cohomology} 

\Definition V.2.1. (a) Let 
$\g$ be a topological Lie algebra and $E$ a locally convex space. 
We call $E$ a {\it topological $\g$-module} if $E$ is a $\g$-module for which the action map 
$\g \times E \to E$ is continuous. 

(b) Let $G$ be a Lie group and $A$ an abelian Lie group. 
We call $A$ a {\it smooth $G$-module} if it is endowed with a
$G$-module structure defined by a smooth action map 
$G \times A \to A$. 
\qed

\Definition V.2.2. Let $\g$ be a topological Lie algebra and 
$E$ a topological $\g$-module. 
For $p \in \N_0$, let $C^p_c(\g,E)$ denote the space of continuous 
alternating maps $\g^p \to E$, i.e., 
the {\it Lie algebra $p$-cochains with values in the module $E$}. 
We identify $C^0_c(\g,E)$ with $E$ and put 
$C_c^\bullet(\g,E) := \bigoplus_{p= 0}^\infty C_c^p(\g,E)$. 
We then obtain a cochain complex with the {\it Lie algebra differential}
$d_\g \: C^p_c(\g,E) \to C^{p+1}_c(\g,E)$
given on $f \in C^p_c(\g,E)$ by 
$$ \eqalign{ 
(d_\g f)(x_0, \ldots, x_p) 
&:= \sum_{j = 0}^p (-1)^j x_j.f(x_0, \ldots, \hat x_j, \ldots, 
x_p) \cr
& + \sum_{i < j} (-1)^{i + j} f([x_i, x_j], x_0, \ldots, \hat
x_i, \ldots, \hat x_j, \ldots, x_p),\cr} $$
where $\hat x_j$ indicates omission of $x_j$ ([ChE48]). 
In view of $d_\g^2= 0$, the space 
$Z^p_c(\g,E) := \break \ker(d_\g\res_{C^p_c(\g,E)})$
of {\it $p$-cocycles} contains the space 
$B^p_c(\g,E) :=  d_\g(C^{p-1}_c(\g,E))$
of {\it $p$-coboundaries}. The quotient 
$$ H^p_c(\g,E) := Z^p_c(\g,E)/B^p_c(\g,E) $$
is the {\it $p$-th continuous cohomology space of $\g$ with values in
the $\g$-module $E$}. We write $[f] := f + B^p_c(\g,E)$ for the 
cohomology class of the cocycle $f$. 
\qed

\Definition V.2.3. Let $G$ be a Lie group and 
$E$ a {\it smooth locally convex $G$-module}, i.e. a smooth $G$-module which is a 
locally convex space. We write 
$$ \rho_E \: G \times E \to E, \quad (g,v) \mapsto \rho_E(g,v) =: \rho_E(g)(v) =: g.v$$ 
for the action map. 
We call a $p$-form $\alpha \in \Omega^p(G,E)$ {\it equivariant} if we have 
for each $g \in G$ the relation 
$$ \lambda_g^*\alpha = \rho_E(g) \circ \alpha. $$
If $E$ is a trivial module, then an equivariant form is 
a left invariant $E$-valued form on $G$. 

We write $\Omega^p(G,E)^G$ for the subspace of equivariant $p$-forms in 
$\Omega^p(G,E)$ and note that this is the space of $G$-fixed elements with 
respect to the action given by $g.\alpha := \rho_E(g)\circ (\lambda_{g^{-1}})^*\alpha$. 
The subcomplex $(\Omega^\bullet(G,E)^G,d)$ of 
equivariant differential forms in the $E$-valued de Rham complex on $G$ 
has been introduced in the 
finite-dimensional setting by Chevalley and Eilenberg in [ChE48].

Let $\g := \L(G) \cong T_\1(G)$. 
An equivariant $p$-form $\alpha$ is uniquely determined by the corresponding element 
$\alpha_\1 \in C^p_c(\g,E)$: 
$$ 
\alpha_g(g.x_1, \ldots, g.x_p) 
= \rho_E(g) \circ \alpha_\1(x_1, \ldots, x_p) \quad \hbox{ 
for } \quad g \in G, x_i \in \g.\leqno(5.2.1)$$
Conversely, (5.2.1) provides for each $\omega \in
C^p_c(\g,E)$ a unique 
equivariant $p$-form $\omega^{\rm eq}$ on $G$ with $\omega^{\rm eq}_\1 = \omega$.
\qed

The following observation is due to Chevalley/Eilenberg ([ChE48, Th.~10.1]). For an 
adaptation to the infinite-dimensional setting, we refer to [Ne04a]. 

\Proposition V.2.4. The evaluation maps 
$$ \ev_\1 \: \Omega^p(G,E)^G \to C^p_c(\g,E), \quad 
\omega \mapsto \omega_\1 $$
define an isomorphism from the cochain complex $(\Omega^\bullet(G,E)^G,d)$ 
of equivariant
$E$-valued differential forms on $G$ to the continuous $E$-valued Lie algebra
complex $(C^\bullet_c(\g,E),d_\g)$. 
\qed

\Definition V.2.5. Let $A$ be a smooth $G$-module and 
$C^n_s(G,A)$ denote the space of all functions $f \: G^n \to A$
which are smooth in an identity neighborhood and normalized in the
sense that $f(g_1,\ldots, g_n)$ vanishes if $g_j = \1$ holds for some
$j$. We call these functions {\it normalized locally smooth group
cochains}. 
The differential $d_G \: C^n_s(G,A) \to C^{n+1}_s(G,A)$, defined by  
$$ \eqalign{ (d_G f)(g_0, \ldots, g_n) 
&:= g_0.f(g_1,\ldots, g_n) \cr
&\ \ \ \ + \sum_{j=1}^n (-1)^j f(g_0, \ldots, g_{j-1} g_j,\ldots, g_n) 
+ (-1)^{n+1} f(g_0, \ldots, g_{n-1}). \cr} $$
turns 
$(C^\bullet_s(G,A), d_G)$ into a differential complex. 
We write $Z^n_s(G,A) := \ker(d_G\res_{C^n_s(G,A)})$ for the corresponding group of cocycles, 
$B^n_s(G,A) := d_G(C^{n-1}_s(G,A))$ for the subgroup of coboundaries, and 
$$ H^n_s(G,A) := Z^n_s(G,A)/B^n_s(G,A) $$
for the {\it $n$-th locally smooth cohomology group with values in the smooth module $A$}. 
\qed

Let $M_1, \ldots, M_k$ be smooth manifolds, $A$ an abelian Lie group 
and $f \: M_1 \times \ldots \times M_k \to A$ a smooth function. 
For $v_k \in T_{m_k}(M_k)$ we obtain a smooth function 
$$ \partial_k(v_k)f \: M_1 \times \ldots \times M_{k-1} \to \a := \L(A), \quad 
(m_1,\ldots, m_{k-1}) \mapsto \delta(f)_{(m_1,\ldots, m_k)}(0, \ldots, 0, v_k). $$
Iterating this process, we obtain for each tuple 
$(m_1,\ldots, m_k) \in M_1 \times \ldots \times M_k$ a continuous $k$-linear map 
$$ T_{m_1}(M_1) \times \ldots \times T_{m_k}(M_k) \to \a, \quad 
(v_1,\ldots, v_k) \mapsto \big(\partial_1(v_1)\cdots \partial_k(v_k) f\big)(m_1,\ldots, m_k). $$

The following theorem describes the natural map from Lie group to Lie algebra cohomology 
([Ne04a, App.~B]; see also [EK64]): 
\Theorem V.2.6. For $f \in C^n_s(G,A)$, $n \geq 1$, 
and $x_1,\ldots, x_n \in \g \cong T_\1(G)$ we put
$$ (D_n f)(x_1,\ldots, x_n) 
:= \sum_{\sigma \in S_n} \sgn(\sigma) 
\big(\partial_1(x_{\sigma(1)})\cdots \partial_n(x_{\sigma(n)})f\big)(\1,\ldots, \1). $$
Then $D_n(f) \in C^n_c(\g,\a)$, and these maps induce 
a morphism of cochain complexes 
$$ D \: (C^n_s(G,A), d_G)_{n \geq 1} \to (C^n_c(\g,\a), d_\g)_{n \geq 1} $$
and in particular homomorphisms
$D_n \: H^n_s(G,A) \to H^n_c(\g,\a)$ for $n \geq 2.$

For $A = \a$ these assertions hold for all $n \in \N_0$, and if 
$A \cong \a/\Gamma_A$ holds for a discrete subgroup $\Gamma_A$ of
$\a$, then $D_1$ also induces a homomorphism 
$D_1 \: H^1_s(G,A) \to H^1_c(\g,\a), [f] \mapsto [df(\1)].$
\qed

\subheadline{Integrability of Lie algebra cocycles} 

We have seen above that for $n \geq 2$ there is a natural {\it derivation map} 
$$ D_n \: H^n_s(G,A) \to H^n_c(\g,\a) $$
from locally smooth Lie group cohomology to continuous Lie algebra cohomology. 
Since 
the Lie algebra cohomology spaces $H^n_c(\g,\a)$ are
much better accessible by algebraic means than those of $G$, it is important to understand 
the amount of information lost by the map $D_n$. Thus one
is interested in kernel and cokernel of $D_n$. A determination of
the cokernel can be considered as describing integrability conditions
on cohomology classes $[\omega] \in H^n_c(\g,\a)$ which have to be
satisfied to ensure the existence of $f \in Z^n_s(G,A)$ with $D_n f =
\omega$.

Before we turn to the complete solution for $n = 2$ ([Ne04a, Sect. 7]), 
we take a closer look at the much simpler case $n = 1$. 

\Remark V.2.7. For $n = 1$ we consider the more general setting of a  
Lie group action on a non-abelian group: Let $G$ and $N$ be Lie groups with 
Lie algebras $\g$ and $\n$ and 
$\sigma \: G \times N \to N$ a smooth action of $G$ on $N$ by automorphisms. 
A {\it crossed homomorphism}, or a {\it $1$-cocycle}, is a smooth map $f \: G \to N$ with 
$$ f(gh) = f(g) \cdot g.f(h) \quad \hbox{ for } \quad g,h \in G, $$
which is equivalent to $(f, \id_G) \: G \to N \rtimes G$ being a group homomorphism. 
We note that for a $1$-cocycle smoothness in an identity neighborhood implies smoothness 
and write $Z^1_s(G,N)$ for the set of smooth $1$-cocycles $G \to N$. 

It is easy to see that for each crossed homomorphism $f \: G \to N$,  
the logarithmic derivative $\delta(f) \in \Omega^1(G,\n)$ is 
an equivariant $1$-form with values in the smooth $G$-module $\n$, hence uniquely 
determined by $D_1(f) := T_\1(f) \: \g \to \n$. Conversely, an easy application of the 
Uniqueness Lemma shows that if $G$ is connected, then a smooth 
function $f \: G \to N$ is a crossed homomorphism if and only if $f(\1) = \1$ and 
$\delta(f)$ is an equivariant $1$-form. 

To see the infinitesimal picture, we call a 
continuous linear map $\alpha \: \g \to \n$ a {\it crossed homomorphism},  
or a {\it $1$-cocycle}, if 
$$ \alpha([x,y]) = x.\alpha(y) - y.\alpha(x) + [\alpha(x), \alpha(y)], \leqno(5.2.2) $$
which is equivalent to $(\alpha, \id_\g) \: \g \to \n \rtimes \g$ being a homomorphism 
of Lie algebras. We write $Z^1_c(\g,\n)$ for the set of continuous $1$-cocycles $\g \to \n$. 
If $G$ is connected, we obtain an injective map 
$$ D_1 \: Z^1_s(G,N) \to Z^1_c(\g,\n). $$
The cocycle condition (5.2.2) for $\alpha$ holds if and only if 
$\alpha^{\rm eq} \in \Omega^1(G,\n)$ satisfies the Maurer--Cartan equation.  
Therefore the Fundamental Theorem (Theorem~III.1.2) shows that if $G$ is connected and 
$N$ is regular, then a $1$-cocycle $\alpha \in Z^1_c(\g,\n)$ is integrable to some group $1$-cocycle 
if and only if the period homomorphism 
$$ \per_\alpha := \per_{\alpha^{\rm eq}} \: \pi_1(G) \to N $$
vanishes. This can also be expressed by the exactness of the sequence
$$ \0 \to Z^1_s(G,N) \smapright{D_1} Z^1_c(\g,\n) \smapright{\per} \Hom(\pi_1(G),N) $$
which already gives 
an idea of what kind of obstructions to expect for $2$-cocycles. 
\qed

The special importance of the group $H^2_s(G,A)$ stems from the following theorem, 
which can be derived easily from the construction in Section V.1. 

\Theorem V.2.8. If $G$ is a connected Lie group and $S \: G \to \Aut(A)$ a smooth action of $G$  
on the abelian Lie group $A$, then we obtain an isomorphism of abelian groups 
$$ H^2_s(G,A) \to \Ext(G,A)_S, \quad 
[\omega] \mapsto A \times_{(S,\omega)} G. 
\qeddis 

\Remark V.2.9. (a) If the group $G$ is not connected, then condition (L3) in Theorem~II.2.1 
requires an additional smoothness condition on cocycles, which is equivalent to 
smoothness of the functions  
$$ f_g \: G \to A, \quad f_g(g') := f(g,g') - f(gg'g^{-1},g) $$
 on an identity neighborhood for each $g \in G$. 
For $g \in G_0$ this is automatically the case for each $f \in Z^2_s(G,A)$. 
We write $Z^2_{ss}(G,A) \subeq Z^2_s(G,A)$ for the set of all cocycles satisfying 
this additional condition. Then $B^2_s(G,A) \subeq Z^2_{ss}(G,A)$, and we put 
$H^2_{ss}(G,A) := Z^2_{ss}(G,A)/B^2_s(G,A)$. 
Theorem V.2.8 remains valid for general $G$ with $H^2_{ss}(G,A)$ instead of $H^2_s(G,A)$. 

(b) The second cohomology groups do not only classify abelian extensions
of $G$. In view of Theorem V.1.5, the sets $\Ext(G,N)_{[S]}$ are principal 
homogeneous spaces of the groups $\Ext(G,Z(N))_S \cong H^2_{ss}(G,Z(N))$, provided 
$Z(N)$ is an initial Lie subgroup of the Lie group $N$ (Remark~II.6.5(c)). 
Therefore the knowledge of the  
second cohomology groups is also crucial for an 
understanding of non-abelian extension classes. 
\qed

On the Lie algebra level, we similarly have for topologically split extensions of 
Lie algebras (cf.\ Remark I.2.7): 

\Proposition V.2.10. Let $(\a,S)$ be a topological $\g$-module, where 
$S \: \g \to \End(\a)$ denotes the module structure, and
write $\Ext(\g,\a)_S$ for the set of all equivalence classes of
topologically split 
$\a$-extensions $\hat\g$ of $\g$ for which the adjoint action of
$\hat\g$ on $\a$ induces the given $\g$-module structure on $\a$. 
Then the map 
$$ Z^2_c(\g,\a) \to \Ext(\g,\a)_S, \quad 
\omega \mapsto [\a \oplus_\omega \g], $$ 
where $\a \oplus_\omega \g$ denotes $\a \times \g$, endowed with the bracket 
$$ [(a,x), (a',x')] := (x.a' - x'.a + \omega(x,x'),[x,x']), $$
factors through a bijection 
$H^2_c(\g,\a) \to \Ext(\g,\a)_S, [\omega] \mapsto [\a \oplus_\omega \g].$ 
\qed

We now turn to the description of the obstruction for the integrability of 
Lie algebra $2$-cocycles. 

\Theorem V.2.11. {\rm(Approximation Theorem; [Ne02a; Th.~A.3.7])} Let 
$M$ be a compact manifold. Then the inclusion map 
$C^\infty(M,G) \into C(M,G)$
is a morphism of Lie groups which is a weak homotopy equivalence, i.e., it induces 
isomorphisms of homotopy groups 
$$ \pi_k(C^\infty(M,G)) \to \pi_k(C(M,G)) $$
for each $k \in \N_0$. In particular, we have 
$$ [M,G] \cong \pi_0(C^\infty(M,G))$$
for the group $[M,G]$ of homotopy classes of maps $M \to G$.
\qed

Below, $\a$ denotes a smooth Mackey complete $G$-module. 

\Definition V.2.12. (a) If 
$M$ is a compact oriented manifold of dimension $k$ and 
$\Omega \in \Omega^k(G,\a)$ a closed $\a$-valued $k$-form, then the map 
$$ \tilde \per_\Omega \: C^\infty(M,G) \to \a, \quad 
\sigma \mapsto \int_\sigma \Omega := \int_M \sigma^*\Omega $$
is locally constant. If, in addition, $\Omega$ is equivariant, then its values lie in 
the closed subspace $\a^\g$ of $\g$-fixed elements in $\a$,  hence defines a period map 
$[M,G] \to \a^\g$ ([Ne02a, Lemma 5.7]). 
If $M = \SS^k$ is a sphere, so that $\pi_k(G) \subeq [\SS^k,G]$ corresponds to base point 
preserving maps, then restriction to $\pi_k(G)$ defines a group homomorphism 
$$ \per_\Omega \: \pi_k(G) \to \a^\g, $$
called the {\it period homomorphism} defined by $\Omega$. 

(b) For $k = 2$ and $\omega \in Z^2_c(\g,\a)$, we obtain a 
Lie algebra $1$-cocycle 
$$f_\omega \: \g \to C^1_c(\g,\a)/d_\g \a, \quad  x \mapsto [i_x \omega], $$
and it is shown in [Ne04a, Lemma 6.2] that this $1$-cocycle gives rise to a well-defined period 
homomorphism, called the {\it flux homomorphism},   
$$ F_\omega \: \pi_1(G) \to H^1_c(\g,\a) $$
as follows. For each piecewise smooth loop $\gamma \: \SS^1 \to G$, we define a $1$-cocycle 
$$ I_\gamma \: \g \to \a,\quad I_\gamma(x) := \int_\gamma i_{x_r} \omega^{\rm eq}, $$
where $x_r$ is the right invariant vector field on $x$ with $x_r(\1) = x$, 
and put $F_\omega([\gamma]) := [I_\gamma]$. 
\qed

Now we turn to the main result of this section ([Ne04a, Th.~7.2]): 

\Theorem V.2.13. Let $G$ be a connected Lie group, $A$ a smooth $G$-module 
of the form $A \cong \a/\Gamma_A$, where $\Gamma_A \subeq \a$ is a
discrete subgroup of the Mackey complete space $\a$ 
and $q_A \: \a \to A$ the quotient map. Then the map 
$$ \tilde P \: Z_c^2(\g,\a) \ \to  \Hom\big(\pi_2(G),A\big)
\times \Hom\big(\pi_1(G), H^1_c(\g,\a)\big), \quad 
\tilde P(\omega) = (q_A \circ \per_\omega, F_\omega) $$
factors through a homomorphism 
$$ P \: H_c^2(\g,\a) \ \to  \Hom\big(\pi_2(G),A\big)
\times \Hom\big(\pi_1(G), H^1_c(\g,\a)\big), \quad 
P([\omega]) = (q_A \circ \per_\omega, F_\omega), $$
and the following sequence is exact: 
$$ \eqalign{
\0 &\to  H^1_s(G,A) \sssmapright{I} H^1_s(\tilde  G,A) \sssmapright{R}  
H^1\big(\pi_1(G),A\big)^G \cong \Hom\big(\pi_1(G),A^G\big) 
\ssmapright{\delta} \cr 
&\ \ \ \ \ssmapright{\delta} H^2_s(G,A) 
\ssmapright{D_2} H_c^2(\g,\a) \ssmapright{P} \Hom\big(\pi_2(G),A\big)
\times \Hom\big(\pi_1(G), H^1_c(\g,\a)\big).  \cr}$$
Here the map $\delta$ assigns to a group homomorphism $\gamma \: \pi_1(G)
\to A^G$ the quotient of the semi-direct product 
$A \rtimes \tilde G$ by the graph $\{ (\gamma(d),d) \: d \in
\pi_1(G)\}$ of $\gamma$ which is a discrete central subgroup, 
$I$ denotes the inflation map and 
$R$ the restriction map to the subgroup $\ker q_G \cong \pi_1(G)$ of 
$\tilde G$. 

If, in particular,  $\pi_1(G)$ and $\pi_2(G)$ vanish, we obtain an isomorphism 
$$ D_2 \: H^2_s(G,A) \to H^2_c(\g,\a). 
\qeddis 

\Remark V.2.14. (a) If $G$ is $1$-connected, things become much simpler and the 
criterion for the integrability of a Lie algebra cocycle 
$\omega$ to a group cocycle is that $\im(\per_\omega) \subeq \Gamma_A$. 
Similar conditions arise in the theory of abelian principal bundles 
on  smoothly paracompact presymplectic manifolds $(M,\Omega)$ ($\Omega$ is a 
closed $2$-form on $M$). Here the integrality of the cohomology class $[\Omega]$ 
is equivalent to the existence of a 
pre-quantum bundle, i.e., a $\T$-principal bundle $\T \into
\hat M \onto M$ whose curvature $2$-form is $\Omega$ (cf.\ [Bry93], [KYMO85]). 

(b) For finite-dimensional Lie groups the integrability 
criteria simplify significantly because $\pi_2(G)$ vanishes 
([CaE36]). This has been used by {\smc \'E.~Cartan} 
to construct central extensions and thus 
to prove that each finite-dimensional Lie algebra belongs to a 
global Lie group, which is known as Lie's third theorem (cf.\ [CaE30/52], [Est88]).  

(c) Let $(M,\omega)$ denote a compact symplectic manifold and 
$\tilde\Diff(M,\omega)$ the universal covering group of 
the identity component $\Diff(M,\omega)_0$ of the Fr\'echet--Lie group $\Diff(M,\omega)$ 
(Theorem~III.3.1). 
Then the Lie algebra homomorphism 
$$ f_\omega \: {\cal V}(M,\omega) := \{ X \in {\cal V}(M) \: {\cal L}_X\omega = 0\} 
\to H^1_{\rm dR}(M,\R), \quad 
X \mapsto [i_X\omega], \leqno(5.2.3) $$
where $H^1_{\rm dR}(M,\R)$ is considered as an abelian Lie algebra, 
integrates to a Lie group homomorphism
$$ {\cal F} \: \tilde\Diff(M,\omega) \to H^1_{\rm dR}(M,\R), $$
whose restriction $\per_{f_\omega}$ to the discrete subgroup $\pi_1(\Diff(M,\omega))$ 
is called the {\it flux homomorphism}. 
Let $\ham(M,\omega) := \ker f_\omega$ denote the Lie subalgebra of Hamiltonian 
vector fields. In [KYMO85, 2.2], it is shown that 
$$ \tilde\Ham(M,\omega) := \ker {\cal F} $$
is a $\mu$-regular Lie group. 

Recently, there has been quite some activity concerning the flux homomorphism for symplectic 
manifolds and generalizations thereof ([Ban97, Ch.~3], [KKM05], [Ne06a]), including a proof of the 
flux conjecture ([Ono04]), formulated by Calabi ([Cal70]). It asserts that  
the image of the flux homomorphism $\per_{f_\omega}$ 
is discrete for each compact symplectic manifold (cf.\ [MD05], [LMP98] for a survey). 

(d) In [RS81], {\smc Ratiu} and {\smc Schmid} address the existence problem of ILH--Lie group 
structures for the following three classes of groups:  
Under the assumption that the image of the flux homomorphism 
is discrete, which is always the case ([Ono04]), 
they show that the group $\Ham(M,\omega)$ of Hamiltonian diffeomorphisms carries 
an ILH--Lie group structure. If $q \: P \to M$ is a pre-quantum $\T$-bundle 
with curvature $\omega$ and connection $1$-form $\theta$, 
they further show that the quantomorphism group $\Aut(P,\theta)$, 
a central $\T$-extension of $\Ham(M,\omega)$ (cf.\ Example~IV.1.6(b)), 
is an ILH--Lie group, 
and they obtain an ILH--Lie group $G$ for the Lie algebra of 
real-valued smooth functions on $T^*(M)$ which are homogeneous of degree $1$  
with respect to the Poisson bracket. The latter group is of particular 
interest for the Lie group structure on the group of invertible 
Fourier--integral operators of order zero, which is a Lie group 
extension of $G$ ([ARS84,86a/b]). 

For a discussion of the relation between quantomorphisms and Hamiltonian 
diffeomorphisms, extending some of these structures, such as {\smc Kostant}'s Theorem 
([Kos70]) to infinite dimensional manifolds, we refer to [NV03]. 

(e) The period and the flux homomorphism annihilate the torsion subgroups of 
$\pi_2(G)$, resp., $\pi_1(G)$. Hence they factor through the rational homotopy groups 
$\pi_2(G) \otimes \Q$, resp., $\pi_1(G) \otimes \Q$. 

(f) If $\a$ is a trivial module and $\omega \in Z^2_c(\g,\a)$, 
then $\hat \g := \a \oplus_\omega \g$ is a central extension of $\g$, 
and 
$x.(a,y) := (\omega(x,y), [x,y])$
turns $\hat\g$ into a topological $\g$-module. A $1$-cocycle 
$f \: \g \to \a$ is the same as a Lie algebra homomorphism, and $B^1_c(\g,\a)= \{0\}$, 
so that $H^1_c(\g,\a) = \Hom_{\rm Lie}(\g,\a) \cong {\cal L}(\g/\oline{[\g,\g]},\a)$. 
In this case, the flux homomorphism 
$$ F_\omega \: \pi_1(G) \to \Hom_{\rm Lie}(\g,\a) $$
vanishes if and only if the action of $\g$ on $\hat\g$ integrates to a smooth 
action of the group $G$ on $\hat\g$ ([Ne02a, Prop.~7.6]). 
\qed

We emphasize that Theorem V.2.13 holds for Lie groups which are not necessarily 
smoothly paracompact. On these groups de Rham's Theorem is 
not available, so that one has to get along without it and to use more direct methods. 
This is important because 
many interesting Banach--Lie groups are not smoothly paracompact, because 
their model spaces do not permit smooth bump functions (cf.\ Remark~I.4.5). 

\Remark V.2.15. Let $G$ be a connected smoothly paracompact Lie group and $A$ a smooth
$G$-module of the form $\a/\Gamma_A$, where $\Gamma_A$ is a discrete subgroup of $\a$. 
Let $Z^2_{gs}(G,A)$ denote the
group of smooth $2$-cocycles $G \times G \to A$ 
and $B^2_{gs}(G,A) \subeq Z^2_{gs}(G,A)$ the cocycles of the form $d_G
h$, where $h \in C^\infty(G,A)$ is a smooth function with $h(\1) =
0$. Then we have an injection 
$$ H^2_{gs}(G,A) := Z^2_{gs}(G,A)/B^2_{gs}(G,A) \into H^2_s(G,A), $$
and the space $H^2_{gs}(G,A)$ classifies those $A$-extensions of $G$ with
a smooth global section. We further have an exact sequence 
$$ \Hom(\pi_1(G),\a^G) \sssmapright{\delta} H^2_{gs}(G,A) 
\sssmapright{D} H^2_c(\g,\a) \ssmapright{P} H^2_{\rm dR}(G,\a) 
\times \Hom\big(\pi_1(G), H^1_c(\g,\a)\big), $$
where $P([\omega]) = ([\omega^{\rm eq}], F_\omega)$ 
(cf.\ Section 8 in [Ne02a] and Remark~8.5 in [Ne04a]). 
\qed

\Remark V.2.16. Let $G$ be a Lie group with Lie algebra 
$\g$ and $\a := C^\infty(G,\R)$, endowed with the compact open $C^\infty$-topology. 
Note that $G$ acts on $\a$ by $(g.f)(x) := f(xg)$, and that the corresponding 
action of $\g$ corresponds to the embedding $\g \to {\cal V}(G), x \mapsto x_l$. 
Using the left trivialization of $T(G)$, 
we see that $\R$-valued $p$-forms are in one-to-one correspondence 
with those smooth functions $G \times \g^p \to \R$ which are $p$-linear 
and alternating in the last $p$ arguments. This implies in particular, that each 
$p$-form $\omega \in \Omega^p(G,\R)$ defines an element of $C^p_c(\g,\a)$, 
and it is easy to see that this leads to an injection of cochain complexes 
$$ \eta \: (\Omega^\bullet(G,\R),d) \into (C^\bullet(\g,\a), d_\g).$$ 

If $G$ is a Fr\'echet--Lie group, then the cartesian closedness 
argument from the convenient calculus (cf.\ [KM97, p.30])  
implies that $\eta$ is bijective, which leads to an isomorphism 
$$ H^p_{\rm dR}(G,\R) \cong H^p_c(\g,\a). $$
If, in addition, $G$ is smoothly paracompact, we thus obtain a description of 
real-valued singular cohomology of $G$ in terms of Lie algebra cohomology 
(cf.\ [Mi87], [Ne04a, Ex.~7.6]). 

In [Mi87], {\smc Michor} applies this construction in particular 
to the group $\Diff(M)$ for a compact manifold $M$. 
For more detailed information on the de Rham cohomology of 
groups like $\Diff(M)$ or $C^\infty(M,K)$, where $M$ is compact and $K$ finite-dimensional, 
we refer to [Beg87]. 
\qed

We have seen above that period homomorphisms arise naturally in the integration theory 
of Lie algebra extensions to group extensions. 
Below we describe some interesting classes of 
Lie algebra $2$-cocycles which have some independent topological interpretation. 

Let $G = C^\infty_c(M,K)$, where $K$ is a Lie group with Lie algebra $\k$ and 
$M$ is a $\sigma$-compact finite-dimensional manifold, so that 
$\g = \L(G) \cong C^\infty_c(M,\k)$, endowed with the natural locally convex 
direct limit structure (Theorem~II.2.8). 
For detailed proofs of the results below we refer to [MN03] for the compact case and 
to [Ne04c] for the non-compact case. 

The Lie algebra cocycles we are interested in are those of 
{\it product type}, constructed as follows. 
Let $E$ be a Mackey complete space and 
$\kappa\: \k \times \k \to E$ an invariant continuous symmetric bilinear form. 
Then the quotient space $\z := \Omega^1_c(M,E)/dC^\infty_c(M,E)$ carries a natural 
locally convex topology because the space of exact forms is closed 
with respect to the natural direct limit topology. 
We then obtain a continuous Lie algebra cocycle 
$$\omega \in Z^2_c(\g,\z) \quad \hbox{ by } \quad 
\omega(\xi,\eta) := [\kappa(\xi,d\eta)]. $$
Of particular interest is the case $E = V(\k)$, where 
$\kappa \: \k \times \k \to V(\k)$ is the universal invariant symmetric bilinear form 
and the case $E = \R$, where $\kappa$ is the Cartan--Killing form of a finite-dimensional 
Lie algebra. We write $\Pi^M_\kappa \subeq \z$ for the corresponding period group. 
Other types of cocycles, which are not of product type, 
are described in [NeWa06a]. If $\k$ is a compact simple Lie algebra and $M = \SS^1$, 
then $H^2_c(\g,\R)$ is one-dimensional and generated by the cocycle 
defined by the Cartan--Killing form~$\kappa$. The associated central extensions 
and their integrability to Lie groups are discussed in some detail in Section~4.2 in~[PS86]. 

\Theorem V.2.17. The following assertions hold:
\litem{(1)} For $M = \SS^1$ we have $\z \cong E$, and the period group 
$\Pi^{\SS^1}_\kappa$ is the image of a homomorphism 
$$ \per_\kappa \: \pi_3(K) \to E, $$
obtained by identifying the subgroup $\pi_2(C_*(\SS^1,K))$ of $\pi_2(G)$ with $\pi_3(K)$. 
\litem{(2)} $\Pi_\kappa^M$ is contained in $H^1_{\rm dR,c}(M,E)$ 
and coincides with the set of all those cohomology classes 
$[\alpha]$ for which integration over circles and properly embedded 
copies of $\R$, we obtain elements of $\Pi^{\SS^1}_\kappa$. 
\litem{(3)} $\Pi^M_\kappa$ is discrete if and only if $\Pi^{\SS^1}_\kappa$ is discrete. 
\litem{(4)} If $\dim K < \infty$ and $\kappa$ is universal, 
then $\Pi^{\SS^1}_\kappa \subeq V(\k)$ is discrete.
\litem{(5)} If $K$ is compact and simple, then the Cartan--Killing form $\kappa$ is universal,  
for a suitable normalization of $\kappa$ we have $\Pi^{\SS^1}_\kappa = \Z$, and 
$\Pi^M_\kappa \subeq H^1_{\rm dR,c}(M,\R)$ is the subgroup of all cohomology 
classes with integral periods in the sense of {\rm(2)}. 
\qed

\Remark V.2.18. A 
particularly interesting class of corresponding central extensions has been studied 
by {\smc Etingof} and {\smc Frenkel} in [EF94]. They investigate the situation 
where $M$ is a compact complex manifold, $K$ is a simple complex Lie group, 
$\kappa$ is the Cartan--Killing form, and by projecting onto the subspace of 
$H^1_{\rm dR}(M,\C) \subeq \Omega^1(M,\C)/dC^\infty(M,\C)$ generated by the holomorphic 
$1$-forms, they obtain a central extension of the complex Lie group 
$C^\infty(M,K)$ by a compact complex Lie group, which in some cases is an elliptic curve 
or an abelian variety. 
\qed

\subheadline{V.3. Abelian extensions of Lie groups} 

In this section, we use the results of the preceding section to 
integrate abelian extensions of Lie algebras to Lie group extensions. 

If $S \: G \to \Aut(A)$ defines on $A$ the structure of a smooth $G$-module, 
$G$ is connected and $A \cong \a/\Gamma_A$ with $\Gamma_A \subeq \a$ discrete,  
then $H^2_s(G,A) \cong \Ext(G,A)_S$ (Theorem~V.2.8), so that Theorem~V.2.13 provides 
in particular necessary and sufficient conditions for a Lie algebra 
cocycle $\omega \in Z^2_c(\g,\a)$ to correspond to a global Lie group 
extension ([Ne04a, Th.~6.7]): 

\Theorem V.3.1. {\rm(Integrability Criterion)} Let $G$ be a connected
Lie group and $A$ a smooth $G$-module with 
$A_0  \cong \a/\Gamma_A$, 
where $\Gamma_A$ is a discrete subgroup of the Mackey complete space
$\a$. For each $\omega \in Z^2_c(\g,\a)$, the abelian Lie algebra extension 
$\a \into \hat\g := \a \oplus_\omega \g \onto \g$ 
integrates to a Lie group extension $A \into \hat G \onto G$ with a 
connected Lie group $\hat G$ if and only if 
\litem{(1)} $\Pi_\omega := \im(\per_\omega) \subeq \Gamma_A$, and 
\litem{(2)} there exists a surjective homomorphism $\gamma \: \pi_1(G) \to \pi_0(A)$ 
such that the flux homomorphism $F_\omega \: \pi_1(G) \to H^1_c(\g,\a)$ 
is related to the characteristic homomorphism 
$$\oline\theta_A \: \pi_0(A) \to H^1_c(\g,\a), \quad [a] \mapsto [D_1(d_G(a))] 
\quad \hbox{ by } \quad  F_\omega = \oline\theta_A \circ \gamma. $$ 

\nin If $A$ is connected, then {\rm(2)} is equivalent to $F_\omega = 0$. 
\qed

\Corollary V.3.2. Let $G$ be a connected Lie group, $\a$ a smooth Mackey complete 
$G$-module and $\omega \in Z^2_c(\g,\a)$. 
Then there exists a smooth $G$-module $A$ with Lie algebra $\a$ such that 
the abelian Lie algebra extension $\a \into \hat\g := \a \oplus_\omega \g \onto \g$ 
integrates to a Lie group extension $A \into \hat G \onto G$ with a 
connected Lie group $\hat G$ if and only if 
$\Pi_\omega$ is a discrete subgroup of $\a^G$. 

\Proof. The necessity is immediate from Theorem V.3.1. For the converse, we first use this 
theorem to find an extension 
$q_0 \: G^\sharp \to \tilde G$ of the universal covering group $\tilde G$ of $G$ 
by the smooth $G$-module $A_0 := \a/\Pi_\omega$. 
Then $A := q_0^{-1}(\pi_1(G)) \subeq G^\sharp$ is a Lie group with identity 
component $A_0$, so that $G^\sharp$ is an $A$-extension of $G$. 
\qed

Note that it may happen that the group $A$ constructed in the preceding proof is not 
abelian. Since $A_0$ and $\pi_1(G)$ are abelian, it is at most $2$-step nilpotent. 

\Remark V.3.3. (a) Suppose that only (1) in Theorem V.3.1 is satisfied, 
and that $A$ is connected. Consider 
the corresponding extension $q^\sharp \: G^\sharp \to \tilde G$ of 
$\tilde G$ by $A \cong \a/\Gamma_A$. 
Then $G \cong G^\sharp/{\hat\pi_1(G)}$, where 
${\hat\pi_1(G)} := (q^\sharp)^{-1}(\pi_1(G))$ 
is a central $A$-extension of $\pi_1(G)$, hence
$2$-step nilpotent. This group is abelian if and only if 
the induced commutator map 
$$ C \: \pi_1(G) \times \pi_1(G) \to A $$
vanishes. 
It is shown in [Ne04a, Rem.~6.8] that, up to sign, this map is given by 
$$ C([\gamma],[\eta]) = \int_{\gamma * \eta} \omega^{\rm eq}, 
\quad \hbox{ where } \quad 
\gamma * \eta \: \T^2 \to G,\quad (t,s) \mapsto \gamma(t)\eta(s). $$
  
(b) According to a result of {\smc H.~Hopf} ([Hop42]), we have for each arcwise connected 
topological space $X$ an exact sequence 
$$ \0 \to H^2(\pi_1(X), A) \to 
H^2_{\rm sing}(X,A) \cong \Hom(H_2(X), A) \to \Hom(\pi_2(X), A) \to \0 $$
(cf.\ [ML78, p.5]). If $G$ is smoothly paracompact, then 
the closed $2$-form $\omega^{\rm eq}$ defines a singular cohomology 
class in $H^2_{\rm sing}(G,\a) \cong \Hom(H_2(M), \a)$ and after composition 
with the quotient 
map $q_A \: \a \to A$, a singular cohomology class $c_\omega \in H^2_{\rm sing}(G,A)$. 
The inclusion $\Pi_\omega \subeq \Gamma_A$ means that this class vanishes on 
the spherical cycles, i.e., the image of $\pi_2(G)$ in $H_2(G)$. Hence it determines 
a central extension of $\pi_1(G)$ by $A$, and if $A$ is divisible, this central extension 
is determined by the commutator map $C \: \pi_1(G) \times \pi_1(G) \to A$. 
If this map vanishes, then $c_\omega = 0$,
but Example~V.3.5(b) below shows that this does not imply the existence of a corresponding 
global group cocycle. If $G$ is $1$-connected, then $c_\omega$ vanishes if and 
only if $\omega$ integrates to a group cocycle (cf.\ [EK64]), 
but in general this simple criterion fails. 

(c) If $F_\omega'([\gamma]) \in H^1_{\rm dR}(G,\a)$ denotes the de Rham class 
obtained as in Proposition V.2.4, then we have for each 
piecewise smooth loop $\eta \: \SS^1  \to G$ the formula 
$\int_\eta F_\omega'(\gamma)= \int_{\gamma*\eta} \omega^{\rm eq}$. 
\qed

The following proposition displays another facet of Hopf's result mentioned 
under (b) above for the special case of topological groups 
(cf.\ [Ne04a, Prop.~6.11]). In the context of rational homotopy theory,  
it can be extended to the Cartan--Serre Theorem, that the rational homology 
algebra of an arcwise connected topological group is generated 
by the homology classes defined by maps $\SS^k \to G$, $k \in \N$ 
(cf.\ [BuGi02, Th.~3.17]). 

\Proposition V.3.4. Let $G$ be a topological group, 
$S_2(G) \subeq H_2(G)$ the subgroup of {\it spherical $2$-cycles}, i.e., 
the image of $\pi_2(G)$ under the Hurewicz homomorphism $\pi_2(G) \to H_2(G)$,  
and $\Lambda_2(G) := H_2(G)/S_2(G)$ the quotient group. Then 
$\Lambda_2(G)$ is generated by the images of cycles defined by maps of the form 
$$ \alpha * \beta \: \T^2\to G, \quad (t,s) \mapsto \alpha(t) \beta(s), $$
where $\alpha, \beta \: \T \to G$ are loops in $G$. 
\qed

\Example V.3.5. (a) Let $G := \Diff(M)_0^{\rm op}$ be the opposite group of the 
identity component of $\Diff(M)$ for a connected compact manifold $M$. 
Recall that its Lie algebra is $\g := {\cal V}(M)$ (Example~II.3.14). For each 
Fr\'echet space $E$, the abelian Lie group $\a = C^\infty(M,E)$ 
is a smooth $G$-module with respect to $\phi.f := f \circ \phi$. 
Each closed $E$-valued $2$-form $\omega_M$ defines a 
continuous Lie algebra $2$-cocycle by $\omega(X,Y) := \omega_M(X,Y)$.
In this case, the period 
map and the flux cocycle can be described in geometrical terms. 
In [Ne04a, Sect. 9], it is shown that the period map 
$$ \per_{\omega} \: \pi_2(\Diff(M)) \to \a^\g = C^\infty(M,E)^{{\cal V}(M)} = E $$
factors for each $m_0 \in M$ through the evaluation map $\ev_{m_0} 
\: \Diff(M) \to M, \phi \mapsto \phi(m_0),$  to the map 
$$ \per_{\omega_M} \: \pi_2(M,m_0) \to E, \quad [\sigma] \mapsto \int_\sigma \omega_M. $$
Likewise, the flux homomorphism can be interpreted as a map 
$$ F_\omega \: \pi_1(\Diff(M)) \to H^1_{\rm dR}(M,E) \cong \Hom(\pi_1(M),E), $$
that vanishes if and only if all integrals of the $2$-form $\omega_M$ over smooth cycles of the 
form $H \: \T^2 \to M, (s,t) \mapsto \alpha(s).\beta(t)$ with loops 
$\alpha$ in $\Diff(M)$ and $\beta$ in $M$ vanish. 

This easily leads to the sufficient condition for the integrability of 
$\omega$ that the period group $\Gamma_E$ of the $2$-form $\omega_M$ should be discrete in $E$. 
This in turn implies the existence of a $Z$-principal bundle for $Z := E/\Gamma_E$ with 
curvature $\omega_M$ over $M$, and the identity component of the group 
$\Aut(P) = \Diff(P)^Z$ is a 
Lie group extension of $G$ by $\Gau(P) \cong C^\infty(M,Z)$, integrating $\omega$ 
(Example~V.1.6(c)). 

It would be very interesting to understand to which extent the discreteness of 
the periods of $\omega_M$ is necessary for the discreteness of the period group 
of $\omega$ (see also the discussion in [KYMO85, p.86] and Problem~V.4).

(b) We consider the special case $M = \T^2$, realized as the unit torus in $\C^2$ and let 
$\omega_M$ be an invariant $2$-form on $M$ with $\int_M \omega_M = 1$. 

Since $\pi_2(M,m_0)$ is trivial, $\per_{\omega}$ vanishes. 
By $\alpha(z)(w_1, w_2) = (z w_1, w_2)$, we obtain a loop 
$\alpha$ in $\Diff(M)$, and the loop 
$\beta(z) := (1,z)$ in $M$ satisfies 
$\alpha(z_1).\beta(z_2) = (z_1, z_2)$, so that 
$$ \int_{\alpha*\beta} \omega_M = 1. $$
We conclude that $F_\omega\not=0$. Hence the Lie algebra cocycle 
$\omega$ on ${\cal V}(M)$ does not integrate to a group cocycle with values 
in the connected group $\a= C^\infty(\T^2,\R)$. 

Since $\omega_M$ is integral, it is the curvature of a natural $\T$-bundle 
$q \: P \to M$, which leads to an abelian extension 
$$ \1 \to A := \Gau(P) \cong C^\infty(M,\T) \into \hat \Diff(M)_0 \onto \Diff(M)_0 \to \1 $$
whose Lie algebra cocycle coincides with $\omega$. Note that 
$\pi_0(A) \cong [\T^2, \T] \cong \Z^2$ is non-trivial. 

(c) The same phenomenon occurs already for the subgroup $T := \T^2$, acting on itself 
by translations, and accordingly on $\a$. By restriction, we obtain an abelian extension 
$$ \1 \to A =C^\infty(\T^2,\T) \into \hat \T^2 \onto \T^2 \to \1 $$
whose flux homomorphism 
$F_\omega \: \pi_1(\T^2) \to H^1(\R^2, \a) \cong H^1_{\rm dR}(\T^2, \R) \cong \R^2$ 
is injective. In this case, there is a reduction of the 
extension of $T$ to an extension by the subgroup 
$$ B := \T \times \Hom(T,\T) \cong \T \times \Z^2  \subeq A =C^\infty(\T^2,\T), $$
generated by the constant maps and the characters of $T$. The corresponding extension 
$\tilde T$ of $T$ by $B$ is isomorphic to the Heisenberg group
$$ H :=\Bigg\{ \pmatrix{1 & a & c \cr  0 & 1 & b \cr 0 & 0 & 1 \cr} \: a,b,c \in \R \Bigg\} 
= \pmatrix{ 1 & \R & \R \cr 0 & 1 & \R \cr 0 & 0 & 1\cr} 
\quad \hbox{ modulo } \quad 
\pmatrix{ 1 & 0& \Z \cr 0 & 1 & 0 \cr 0 & 0 & 1\cr}. 
\qeddis

\Example V.3.6. Let $G := \SL_2(\R)$. From the natural action of $G$ on 
$\P_1(\R) \cong \SS^1$, we derive an action on the space 
$\a := \Omega^1(\SS^1,\R)$. In Section 10 of [Ne04a], it is shown that 
there exists a non-trivial $\omega \in Z^2(\sL_2(\R), \a)$ which 
integrates to an abelian extension 
$$\hat \SL_2(\R) = \Omega^1(\SS^1,\R) \rtimes_f \SL_2(\R), $$
so that we obtain a non-trivial infinite-dimensional abelian extension of 
$\SL_2(\R)$ which is a Fr\'echet--Lie group. 

Since all finite-dimensional 
Lie group extensions of $\SL_2(\R)$ by vector spaces 
split on the Lie algebra level, this example illustrates the difference between 
the finite- and infinite-dimensional theory. 
\qed

For more references dealing specifically with central extensions, we refer to 
[Ne02a]. In particular, [CVLL98] is a nice survey on central $\T$-extensions of 
Lie groups and their role in quantum physics (see also [Rog95]). It also contains 
a description of the universal central extension for finite-dimensional groups. 
For infinite-dimensional groups, universal central extensions are constructed 
in [Ne02d], and for root graded Lie algebras in [Ne03] 
(cf.\ Subsection~VI.1). 

\Example V.3.8. (The Virasoro group) Let $G := \Diff_+(\T)$ be the group of orientation
preserving diffeomorphisms of the circle $\T$. Then the inclusion 
$\T \into G$ of the rigid rotations is a homotopy equivalence, so that $\pi_2(G)$ vanishes 
and $\pi_1(G) \cong \Z$ (cf.\ [Fu86, p.~302]).

Furthermore, $H_c^2(\g,\R) = \R [\omega]$ is one-dimensional ([PS86]), 
and the corresponding flux homomorphism $F_\omega$ vanishes ([Ne02a, Ex.~9.3]), 
so that Theorem~V.3.1 implies the existence of a corresponding 
central $\R$-extension of $G$, called the {\it Bott--Virasoro group} $\Vir$. 
Remark~V.2.15 implies that this extension has a smooth global section, hence can 
be described by a smooth global cocycle. 
Such a cocycle, and other related ones, are described explicitly 
by {\smc Bott} in [Bo77]. A more direct construction of this and related 
cocycles has been described recently by {\smc Billig} ([BiY03]). 

In [Se81], {\smc G.~Segal} studies projective unitary representations of
$\Diff(\SS^1)$ via representations of loop groups, which implicitly 
define unitary representations of the Bott--Virasoro group.
In [GW84/85], {\smc Goodman and Wallach} give 
an analytic construction of the unitary highest weight
representations of $\Vir$ by directly
integrating the corresponding Lie algebra representation on the 
representations of loop groups,  using scales of Banach spaces. 

The Bott--Virasoro group is also a very interesting geometric object. 
One aspect of its rich geometric structure is that, although it is only a smooth 
real Lie group which is not analytic (Remark~VI.2.3 below), 
it carries the structure of a complex Fr\'echet manifold, which is 
obtained by  identifying it with the complement 
of the zero section in the holomorphic line bundles over $\Diff_+(\SS^1)/\T$ ([Lem95], [KY87]).

\subheadline{Open Problems for Section V} 

\Problem V.1. Generalize Theorem V.2.13 in an appropriate way to non-connected 
Lie groups $G$ and $A$. 

The generalization to non-connected Lie groups $G$ means to 
derive accessible criteria for
the extendibility of a $2$-cocycle on the identity component $G_0$ to the whole group $G$. From the short exact sequence 
$G_0 \into G \onto \pi_0(G)$, we obtain maps 
$$ H^2(\pi_0(G), A) \sssmapright{I} H^2(G, A) \sssmapright{R} H^2(G_0, A)^{G}, $$
but it is not clear how to describe the image of the restriction map $R$ from $G$ to $G_0$. 

If $A$ is a trivial module, one possible approach is to introduce additional structures on 
a central extension $\hat G$ of $G_0$ by $A$, so that the map 
$q \: \hat G \to G$ describes a crossed module, which requires an extension of the 
natural $G_0$-action of $G$ on $\hat G$ to an action of $G$ (cf.\ [Ne05]). 

To deal with non-connected groups $A$ seems to be tractable if we assume that 
$A_0 \cong \a/\Gamma_A$ as in Theorem V.2.13. 
Under the assumption that $G$ is connected, the crucial information is 
contained in an exact sequence 
$$ \0 \to H^2_s(G,A_0) \to H^2_s(G,A) \sssmapright{\gamma} \Hom(\pi_1(G), \pi_0(A)) \to
H^3_s(G, A_0), $$
where $\gamma$ assigns to an extension of $G$ by $A$ the corresponding connecting 
homomorphism $\pi_1(G) \to \pi_0(A)$ in the long exact homotopy sequence
(cf.\ [Ne04a, App.~E]). To determine $H^2_s(G,A)$ in terms of $H^s(G,A_0)$ and
known data, one has to determine the image of $H^2_s(G,A)$ in 
$\Hom(\pi_1(G), \pi_0(A))$. 
\qed

\Problem V.2. Do the spaces $Z^2_s(G,A)$ and $Z^2_{ss}(G,A)$ (Remark~V.2.9) coincide for 
each non-connected Lie group $G$ and each smooth $G$-module $A$? 

This is true if $G$ is connected ([Ne04a, Prop.~2.6]), 
but in general we do not know if $Z^2_{ss}(G,A)$ 
is a proper subgroup of $Z^2_s(G,A)$, which is equivalent to 
$H^2_{ss}(G,A)$ being a proper subgroup of $H^2_s(G,A)$. 
In terms of abelian extensions, this means that there exists an abelian extension  
$\hat G$ of $G$ by the $G$-module $A$ for which the restriction $\hat G_0$ to the identity 
component $G_0$ is a Lie group extension, but $\hat G$ cannot be turned into a 
Lie group because for certain elements $\hat g \in \hat G$ the conjugation action 
on $\hat G_0$ is not an action by smooth group automorphisms 
(cf.\ condition (L3) in Theorem~II.2.1). 
\qed

\Problem V.3. Give an explicit description of kernel and cokernel of the derivation maps 
$$ D_n \: H^n_s(G,A) \to H^n_c(\g,\a) \quad \hbox{ for } \quad n \geq 3. $$
For $A \cong \a/\Gamma_A$ for some discrete subgroup $\Gamma_A \subeq \a$, 
the first necessary condition for $[\omega] \in H^n_c(\g,\a)$ to lie in the image of $D$, 
one obtains quite easily 
is that the range of the period homomorphism 
$$ \per_\omega \: \pi_n(G) \to \a $$
must be contained in $\Gamma_A \cong \pi_1(A)$ (cf.\ [GN07]). 
\qed

\Problem V.4. An interesting special case of the preceding problem 
arises for $G = \Diff(M)_0^{\rm op}$, $M$ a compact manifold, 
$\a = C^\infty(M,\R)$, where $G$ acts by $(\phi.f)(m) := f(\phi(m))$, 
and $\omega \in \Omega^2(M,\R)$ is a closed $2$-form. Then 
$\omega$ defines a Lie algebra cocycle in $Z^2_c({\cal V}(M), \a)$, 
and it is an interesting question when this cocycle integrates to a group 
cocycle on $G$. We know that this is the case if the period group 
$\la [\omega], H_2(M)\ra \subeq \R$ is discrete, but this is not necessary 
(cf.\ [KYMO85, p,86]). The approach described in Example~V.3.5 may be useful 
to analyze this problem. The crucial point is to understand the range of the 
homomorphism $\pi_2(\Diff(M)) \to \pi_2(M,m_0)$ and of the natural map  
$\pi_1(\Diff(M)) \times \pi_1(M,m_0) \to [\T^2, M] \to H_2(M)$ (Example~V.3.5) 
(see [Ban97, Ch.~3] for more details on such maps).  
\qed

\Problem V.5. Give a characterization of those principal $K$-bundles 
$q \: P \to M$ for which the extension 
$\Aut(P)$ of the subgroup $\Diff(M)_{[P]}$ by the gauge 
group $\Gau(P)$ splits on the group level (cf.\ Example V.1.6). 
On the Lie algebra level, such conditions are given by {\smc Lecomte} in 
[Lec85]. Note that this is obviously the case if the bundle is trivial, which implies 
$\Aut(P) \cong C^\infty(M,K) \rtimes \Diff(M)$. It is also the case 
for natural bundles to which the action of $\Diff(M)$ lifts, such as the frame 
bundle and other natural bundles. 
\qed

\sectionheadline{VI. Integrability of locally convex Lie algebras} 

\nin In this section, we take a systematic look at the integrability 
problem for locally convex Lie algebras with an emphasis on locally exponential ones, 
because they permit a quite satisfying general theory. For Lie algebras which are not 
locally exponential only isolated results are available. 

\subheadline{VI.1. Enlargeability of locally exponential Lie algebras} 

\Definition VI.1.1. A locally convex Lie algebra $\g$ is said to be 
{\it integrable} if there exists a Lie group $G$ with 
$\L(G) \cong \g$. It is called 
{\it locally integrable} if there exists a local Lie group 
$(G,D,m_G,\1)$ with Lie algebra $\L(G) \cong \g$. 
A  locally exponential Lie algebra 
is called {\it enlargeable} if it is integrable to a locally exponential Lie group, 
i.e., if some of the corresponding local 
groups are enlargeable (cf.\ Definition~IV.2.3). 
\qed

Although every finite-dimensional Lie algebra is integrable, 
integrability of infinite-di\-men\-sio\-nal Lie algebras turns out to be a very subtle 
property. 

\Examples VI.1.2. (a) If $\g$ is a finite-dimensional Lie algebra, endowed with its 
unique locally convex topology, then $\g$ is integrable. 
This is Lie's Third Theorem. One possibility to prove this is first to 
use Ado's Theorem to find an embedding $\g \into \gl_n(\R)$ and then to 
endow the integral subgroup $G := \la \exp \g \ra \subeq \GL_n(\R)$ with a Lie group 
structure such that $\L(G) = \g$ (cf.\ Corollary~IV.4.10). 

(b) If $\g$ is locally exponential, then it is locally integrable by definition. 
In particular, every Banach--Lie algebra is locally integrable (Examples~IV.2.4). 
\qed

\subheadline{Enlargeability and generalized central extensions}

The criteria described in Section V.3 provide good tools 
to understand the difference between the group and Lie algebra picture 
for abelian extensions. However, not all quotient maps 
$q \: \hat\g \to \g$ of Lie algebras are topologically split in the sense 
that there is a continuous linear section, therefore they are not extensions 
of the type just discussed. An important example 
is the map $\ad \: \g \to \g/\z(\g)$, where $\z(\g)$ is the center. 
The fact that for each locally exponential Lie algebra 
$\g$, the Lie algebra $\g_{\rm ad} := \g/\z(\g)$ is always integrable  
(Theorem IV.3.8) shows that the question of the integrability of central extensions 
has to be addressed even for those which are not topologically split. 
Fortunately, there is a method to circumvent the problems caused by 
this topological difficulty by reducing all assertions to topologically split 
central extensions. The key concept is that of a generalized central extension 
(cf.\ [Ne03], [GN07]).

\Definition VI.1.3. A morphism $q \: \hat\g\to\g$ of 
locally convex Lie algebras is 
called a {\it generalized central extension} if it has dense range and 
there exists a continuous bilinear map $b \: \g \times \g \to \hat\g$ for which 
$b \circ (q \times q)$ is the Lie bracket on $\hat\g$. 
It is called a {\it central extension} if, in addition, $q$ is a quotient map. 
\qed

The subtlety of generalized central extensions 
is that $q$ need not be surjective and if it is surjective, it need not be 
a quotient map. Fortunately,  
these difficulties are compensated by the following nice fact. 
Let us call a locally convex Lie algebra $\g$ {\it topologically perfect} 
if its commutator algebra is dense. We call a generalized central extension 
$q_\g \: \tilde\g \to \g$ {\it universal} if for any generalized central extension 
$q \: \hat\g \to \g$ there exists a unique morphism of locally convex 
Lie algebras $\alpha \: \tilde\g \to \hat\g$ with $q \circ \alpha = q_\g$. 
Then one can show that each topologically perfect locally convex Lie algebra 
$\g$ has a universal generalized central extension (unique up to isomorphism). 
For the basic results on generalized central extensions we refer to [Ne03, Sect.~III], where 
one also finds descriptions of the universal generalized central extensions 
of several classes of Lie algebras. 

\Remark VI.1.4. If $q \: \hat\g \to \g$ is a central extension, 
then $q \times q \: \hat\g \times \hat\g \to \g \times \g$ also is a 
quotient map. Therefore the Lie bracket of $\hat\g$ factors through 
a continuous bilinear map $b \: \g\times \g \to \hat\g$ 
with $b(q(x), q(y)) = [x,y]$ for $x,y \in \hat\g$, 
showing that $q$ is a generalized central extension of $\g$. 
\qed

\Proposition VI.1.5. {\rm([Ne03, Lemma~III.4])\footnote{$^1$}{\eightrm For the case of 
central extensions of Banach--Lie algebras, part of the assertions below can be found in a 
footnote in [ES73, p.58].}} For a generalized central extension
$q \: \hat\g\to\g$ the following assertions hold: 
\litem{(1)} The corresponding map $b$ is a Lie algebra cocycle in $Z^2_c(\g,|\hat\g|)$, where 
$|\hat\g|$ denotes $\hat\g$, considered as a trivial $\g$-module. 
\litem{(2)} If $|\g|$ denotes the space $\g$, endowed with the trivial Lie bracket, 
then the maps 
$$\psi \: \hat\g \to |\hat\g| \oplus_b \g,  \quad x \mapsto (x,q(x)) \quad \hbox{ and } 
\quad \eta \: |\hat\g| \oplus_b \g \to |\g|, \quad (x,y) \mapsto y - q(x) $$
are homomorphisms of Lie algebras, $\psi$ is a topological embedding, 
$\eta$ is a quotient map, and the sequence 
$$ \0 \to \hat\g \sssmapright{\psi} |\hat\g| \oplus_b \g \sssmapright{\eta} |\g| \to \0 $$
is exact. 
\qed

For the following theorem from [GN07], we recall that central extensions of locally exponential 
Lie algebras by Mackey complete spaces are locally exponential (Theorem IV.2.10).

\Theorem VI.1.6. {\rm(Enlargeability criterion for generalized central extensions)} Let 
$G$ be a connected locally exponential Lie group with Lie algebra $\g$ and 
$q \: \hat\g\to \g$ a generalized central extension for which $\hat\g$ is Mackey complete. 
Let $\omega \in Z^2_c(\g,|\hat\g|)$ be the associated 
Lie algebra cocycle and $\per_\omega \: \pi_2(G) \to |\hat\g|$ the corresponding period homomorphism. 
Then the following assertions hold: 
\litem{(1)} $\Pi_\omega := \im(\per_\omega)$ is contained in $\z$. 
\litem{(2)} $\hat\g$ is enlargeable if $\Pi_\omega$ is discrete. 
\litem{(3)} If $q$ is a central extension, then $\hat\g$ is enlargeable if and only if $\Pi_\omega$ is discrete. 

\Proof. (1) follows from the fact that the cocycle $q \circ \omega = - d_\g \id_\g$ is trivial. 
It is the Lie bracket of $\g$. 

(2) Corollary V.3.2 implies that $\tilde\g := |\hat\g| \oplus_b \g$ 
is enlargeable if and only if $\Pi_\omega$ is discrete. 
If this is the case, then the closed ideal $\hat\g$ of $\tilde\g$ is also enlargeable 
because $\hat\g \cong \ker \eta$ implies that it is locally exponential (Theorem~IV.2.9), 
so that Corollary~IV.4.10 applies.

(3) Suppose that $q$ is a quotient map, i.e., a central extension, and that $\hat\g$ is enlargeable. 
Since the cocycle $\tilde b := q^*b$ coincides with the Lie bracket on $\hat\g$, 
the corresponding central extension 
$\hat\g^\sharp := |\hat\g| \oplus_{\tilde b} \hat\g$ is split by the section 
$\sigma(x) := (x,x)$, hence is enlargeable. 
Furthermore, 
$$\tilde\g = |\hat\g| \oplus_b \g \cong \hat\g^\sharp/(\{0\} \times \z) $$
is locally exponential by Theorem IV.2.9, which applies in particular to all quotients by 
central ideals. In view of Theorem VI.1.10 below, it now 
suffices to show that the integral subgroup $Z$ generated by $\z$ 
is a locally exponential Lie subgroup. But this follows from the fact that the projection 
onto $|\hat\g| \times \{0\}$ along $\im(\sigma)$ restricts to a homeomorphism 
on $\z$. Hence the corresponding subgroup is a locally exponential Lie subgroup, 
and this completes the proof. 
\qed

The preceding theorem applies in particular to central extensions 
$\z \into \hat\g\onto \g = \L(G)$ of Banach--Lie algebras, for which 
it characterizes integrability in terms of the discreteness of $\Pi_b$. 
In this case, a similar criterion is given by {\smc van Est} and {\smc Korthagen} 
in [EK64]. On the surface, their criterion has the same formulation, but their 
period homomorphism arises as an element of 
$H^2_{\rm sing}(G,\z) \cong \Hom(H_2(G),\z)$ obtained from the enlargeability 
theory of local groups ([Est62]). Under their assumption 
that $G$ is $1$-connected, the Hurewicz homomorphism $\pi_2(G) \to H_2(G)$ 
is an isomorphism, so that their period homomorphism also is a 
homomorphism $\pi_2(G) \to \z$, and one can even show that both coincide up to sign. 
We think that the definition of the period homomorphism 
in terms of integration of differential forms makes it much more accessible 
than the implicit construction in [EK64]. 

\Definition VI.1.7. Let $\g$ be a locally exponential Lie algebra and consider the 
central extension 
$$ \0 \to \z(\g) \to \g \to \g_{\rm ad} := \g/\z \to \0. $$
Let $G_{\rm ad} \subeq \Aut(\g)$ be endowed with its locally exponential group 
structure with Lie algebra $\g_{\rm ad}$ (Theorem IV.3.8) and 
$$ \per_\g \: \pi_2(G_{\rm ad}) \to \z(\g) $$
the corresponding period homomorphism (Theorem VI.1.6(1)). 
We write $\Pi(\g) := \im(\per_\g)$ for its image and call it the {\it period group of $\g$}. 
\qed

The following theorem generalizes the enlargeability criterion of [EK64] for Banach algebras. 
It follows immediately from Theorem IV.3.8 on the integrability of $\g/\z(\g)$ 
and Theorem~VI.1.6. 

\Theorem VI.1.8. {\rm(Enlargeability Criterion for locally exponential Lie algebras)} 
A Mackey complete locally exponential Lie algebra $\g$ is enlargeable if and only if 
its period group $\Pi(\g)$ is discrete. 
\qed

\Proposition VI.1.9. If $\g$ is a separable locally exponential Lie algebra, then 
$\Pi(\g)$ is countable. If, in addition, $\g$ is Fr\'echet, then $\Pi(\g)$ is 
closed if and only if it is discrete. 

\Proof. If $\g$ is separable, then the same holds for the connected group 
$G_{\rm ad}$ and hence for the identity component $C_*(\SS^1,G_{\rm ad})_0$ of the loop 
group. Its universal covering group is also separable, so that its fundamental 
group, which is isomorphic to $\pi_2(G_{\rm ad})$, is countable. 
This implies that $\Pi(\g)$ is countable. 

If $\Pi(\g)$ is closed and $\g$ is Fr\'echet, it is a countable complete metric space, 
hence discrete. 
\qed

For the second part of the preceding proposition, the Fr\'echet assumption 
on $\g$ is crucial: the space $\R^\R$ contains a non-discrete closed 
subgroup isomorphic to $\Z^{(\N)}$ ([HMP04, Cor.~3.2(i)]). 

Combining the fact that kernels of morphisms are locally exponential Lie subgroups 
(Proposition~IV.3.4) and Theorem IV.1.19 on the integration of morphisms of 
Lie algebras, one obtains the equivalence of (1) and (2) in the 
following integrability criterion for quotient algebras ([GN06]): 

\Theorem VI.1.10. {\rm(Enlargeability Criterion for quotients)} Let 
$G$ be a $1$-connected locally exponential Lie group 
and $\n \trile \g$ a closed ideal for which 
the quotient Lie algebra $\q := \g/\n$ is locally 
exponential. Let 
$$Z(G,\n) := \{ g \in G \: (\Ad(g) - \1)(\g) \subeq \n\}. $$ 
Then $Z(G,\n) \trile G$ is a normal locally exponential Lie subgroup with 
Lie algebra 
$$\z(\g,\n) := \{ x \in \g \: [x,\g] \subeq \n\},$$ 
and the Lie algebra homomorphism 
$q \: \z(\g,\n) \to \z(\q)$ defines a period homomorphism 
$$ \per_q  \: \pi_1(Z(G,\n)) \to \z(\q), \quad 
\per_q([\gamma]) = \int_0^1 q(\delta(\gamma)_t)\, dt, $$
where $\gamma \:[0,1] \to Z(G,\n)$ is a piecewise smooth loop. 
The following assertions are equivalent: 
\litem{(1)} The locally exponential Lie algebra $\q = \g/\n$ is enlargeable. 
\litem{(2)} The normal integral 
subgroup $N := \la \exp_G\n\ra \trile G$ is a locally exponential Lie subgroup.   
\litem{(3)} The image of $\per_q$ is a discrete subgroup of $\z(\q)$. 
\qed

\Remark VI.1.11. In addition to the assumptions of 
the preceding theorem, suppose that $G$ is a Fr\'echet--Lie group. 
We then have 
$Q_{\rm ad} \cong G/Z(G,\n)$, and since Michael's Selection Theorem  ([MicE59]) 
applies to the quotient map $\g\to \q$, 
this leads to a surjective homomorphism 
$$ \delta \: \pi_2(Q_{\rm ad}) \to \pi_1(Z(G,\n)). $$
The surjectivity of $\delta$ follows from the $1$-connectedness of $G$ 
and the exactness of the long exact homotopy sequence of the bundle 
$G \to Q_{\rm ad}$. Then it is not hard to see that 
$$ \per_q \circ \delta  = \per_{\q} \: \pi_2(Q_{\rm ad}) \to \z(\q), $$
which shows that (3) in Theorem~VI.1.10 is equivalent to the discreteness of the $\Pi(\q)$ 
(Theorem~VI.1.8). 
\qed

\Proposition VI.1.12. {\rm(Functoriality of the period group)} 
Let $\phi \: \g\to \h$ be a morphism of Mackey complete locally exponential 
Lie algebras with $\phi(\z(\g)) \subeq \z(\h)$ 
and $\phi_{\rm ad} \: \g_{\rm ad} \to \h_{\rm ad}$ the induced homomorphism. 
Then $\phi_{\rm ad}$ integrates to a group homomorphism 
$\tilde\phi_{\rm ad} \: \tilde G_{\rm ad} \to \tilde H_{\rm ad}$, 
$\phi(\Pi(\g)) \subeq \Pi(\h)$, and the following diagram commutes 
$$ \matrix{ 
\pi_2(G_{\rm ad}) &
\smapright{\pi_2(\tilde\phi_{\rm ad})} & \pi_2(H_{\rm ad}) \cr 
\mapdown{\per_\g}  & {}  & \mapdown{\per_\h} \cr 
\z(\g)  & \smapright{\phi} & \z(\h)}.  $$
\qed

\Corollary VI.1.13. If $\g_1, \g_2$ are Mackey complete locally exponential Lie algebras, then 
$$ \Pi(\g_1 \times \g_2) = \Pi(\g_1) \times \Pi(\g_2). 
\qeddis 

\Remark VI.1.14. (Constructing non-enlargeable Lie algebras) 
Suppose that $\g$ is a locally exponential Lie algebra with 
$\Pi(\g) \cong \Z$. Let $\theta \in \R \setminus \Q$. 
Then $\z := \{ (x,\theta x) \: x \in \z(\g)\}$ is a central ideal of 
$\g \times \g$, so that $\h := (\g \times \g)/\z$ is locally exponential. 
Corollary VI.1.13 and Proposition VI.1.12 imply that, writing $\Pi(\g) = \Z d$, 
we get 
$$ \Pi(\h) \cong \Z[(d,0)] + \Z[(0,d)] 
= (\Z + \Z \theta) [(d,0)], $$ 
which is not discrete. Hence $\h$ is not integrable. 
\qed

Using the construction of the group $G$ via Theorem~VI.1.6 and the long exact homotopy sequence,  
one can identify the period group of enlargeable Fr\'echet--Lie 
algebras in terms of the center: 

\Proposition VI.1.15. If $G$ is a locally exponential $1$-connected Fr\'echet--Lie group and 
$\g = \L(G)$ its Lie algebra, then 
$$ \Pi(\g) = \ker(\exp_G \res_{\z(\g)}) \cong \pi_1(Z(G)). 
\qeddis

\Example VI.1.16. The first example of a non-enlargeable Banach--Lie algebra 
was given by van Est and Korthagen with the method described in Remark~VI.1.14 ([EK64]). 
It is the central extension 
$\g$ of the Banach--Lie algebra $C^1(\SS^1, \su_2(\C))$ by $\R$, defined by the cocycle 
$$ \omega(f,g) := \int_0^1 \tr(f(t)g'(t))\, dt, $$
where we identify functions on $\SS^1 \cong \R/\Z$ with $1$-periodic functions on $\R$. 
Then $\g_{\rm ad} \cong C^1(\SS^1, \su_2(\C))$ and 
$G_{\rm ad} \cong C^1(\SS^1, \SU_2(\C))$ leads to 
$\pi_2(G_{\rm ad}) \cong \pi_3(\SU_2(\C)) \cong \pi_3(\SS^3) \cong \Z$.  
Now one shows that $\per_\g = \per_\omega$ is non-trivial to verify that 
$\Pi(\g) \cong \Z$. 

Using Kuiper's Theorem ([Ku65]),  
{\smc Douady} and {\smc Lazard} gave a simpler example ([DL66]): 
by observing that the $1$-connectedness of the unitary group $\UU(H)$ of an 
infinite-dimensional complex Hilbert space $H$ implies that its Lie algebra 
$\uu(H) := \{ X \in {\cal L}(H) \: X^* = - X\}$ satisfies 
$$ \Pi(\uu(H)) \cong \pi_1(Z(\UU(H))) = \pi_1(\T) \cong \Z $$
(Proposition~VI.1.15). 

Based on the fact that $\UU(H)$ is $1$-connected, one can also give the following 
direct argument. For any irrational $\theta \in \R \setminus \Q$ the line 
$\n := \R i (\1,\theta \1)$ generates a dense subgroup of the center 
$Z(\UU(H) \times \UU(H)) \cong \T^2$ of the $1$-connected group $\UU(H) \times \UU(H)$, 
so that Theorem~VI.1.10 implies that the quotient Lie algebra 
$(\uu(H) \times \uu(H))/\n$ is not enlargeable. 
\qed

\subheadline{Enlargeability of quotients} 

One may take Theorem VI.1.10 as a starting point of a theory 
of certain topological groups which are more general than Lie groups, namely 
quotients of Lie groups. This leads to the concept of a {\it scheme of Lie groups}, or 
{\it $S$-Lie group} (cf.\ [Ser65], [DL66] and [Est84]). The strength of this concept 
for Banach--Lie algebras and, more generally, locally exponential Lie algebras, 
follows from the fact that each such Lie algebra is a quotient of an enlargeable one:  

\Theorem VI.1.17. {\rm([Swi71] for the Banach case)} For each locally exponential 
Fr\'echet--Lie algebra $\g$, the Lie algebra 
$$ \Lambda(\g) := C_*([0,1],\g) := \{ \gamma \in C([0,1],\g) \: \gamma(0) = 0\} $$
is enlargeable. 

\Proof. Clearly, $\z(\Lambda(\g)) = \Lambda(\z(\g))$, so that 
$\Lambda(\g)_{\rm ad} \cong \Lambda(\g_{\rm ad})$ follows from Michael's Theorem 
([MicE59]). The corresponding 
group $C_*([0,1],G_{\rm ad})$ is contractible, and this leads to 
$\Pi(\Lambda(\g)) = \{0\}$, which implies enlargeability. 
\qed

A central point of the preceding theorem is that it implies that each 
locally exponential Fr\'echet--Lie algebra $\g$ 
is a quotient of an enlargeable Fr\'echet--Lie algebra (cf.\ [Rob02, Th.~5]):  
the evaluation map $\ev_1 \: \Lambda(\g) \to \g, \gamma \mapsto \gamma(1)$ 
is a quotient map. Now one can address the enlargeability problem along the 
lines of Theorem VI.1.10.

\Remark VI.1.18. (a) In [Woj06], {\smc Wojty\'nski} describes a variant of this approach for 
Banach--Lie algebras. 
Instead of considering the Banach--Lie algebra $\Lambda(\g)$, he considers 
analytic paths $\gamma(t) := \sum_{n = 1}^\infty a_n t^n$, for which 
$\|\gamma\|_1 := \sum_{n = 1}^\infty \|a_n\|$ is finite. Identifying these curves 
with their coefficient sequences, we denote this space by 
$\ell^1(\g) := \ell^1(\N,\g)$. The Lie bracket on this sequence space is 
given by 
$$ [(a_n), (b_n)] = (c_n)\quad \hbox{ with} \quad c_n = \sum_{j = 1}^{n-1} [a_j, b_{n-j}]. 
\leqno(6.1.1) $$

With the same Lie bracket, we also turn the full sequence space 
$\g^\N$ into a pro-nilpotent Fr\'echet--Lie algebra, which is exponential for trivial 
reasons. Since the Banach--Lie algebra $\ell^1(\g)$ injects into the 
exponential Lie algebra $\g^\N$, it is enlargeable 
by Corollary IV.4.10. Again, we have an evaluation map 
$$ q \: \ell^1(\g) \to \g, \quad (a_n) \mapsto \sum_{n = 1}^\infty a_n,  $$
which is a quotient morphism of Lie algebras and since 
the subgroup 
$\la \exp \ell^1(\g)\ra$ is contractible (cf.\ [Woj06]), one may proceed with 
Theorem VI.1.10 as for $\Lambda(\g)$. 

(b) In [Pe93a/95a], {\smc Pestov} shows that if $E$ is a Banach space of 
$\dim E > 1$, then the free Banach--Lie 
algebra over $E$ has trivial center. As a consequence, every Banach--Lie algebra 
$\g$ of dimension $> 1$ is a quotient of a centerless Banach--Lie algebra 
$F(\g)$, the free Banach--Lie algebra over the Banach space $\g$, which is enlargeable 
because its center is trivial (Theorem~IV.3.8). Again, we can proceed 
with Theorem~VI.1.10 to obtain enlargeability criteria. 
\qed

The following enlargeability criterion of {\smc Swierczkowski} for 
extensions by not necessarily abelian ideals is a powerful tool. 
It would be very interesting to see if it can be extended to the locally exponential 
setting. It applies in particular to all situations where 
$\q$ is finite-dimensional or abelian (cf.\ [Swi65]; Remark~V.2.14(b)). 

\Theorem VI.1.19. {\rm([Swi67, Th., Sect. 12])} Suppose that $\g$ is a Banach--Lie algebra 
and $\n \trile \g$ a closed enlargeable ideal for which $\q := \g/\n$ is enlargeable 
to some group $Q$ with vanishing $\pi_2(Q)$, then $\g$ is enlargeable. 
\qed

\Definition VI.1.20. A 
Banach--Lie algebra is said to be {\it lower solvable} if there exists 
an ordinal number $\alpha$ and an ascending chain of closed subalgebras 
$$ \{0\} = \g_0 \subeq \g_1 \subeq \g_2 \subeq \ldots \subeq \g_\beta \subeq 
\g_{\beta+1} \subeq \ldots \subeq \g_\alpha = \g $$
such that 
\litem{(a)} If $\beta \leq \alpha$ is not a limit ordinal, then 
$X_{\beta-1}$ is an ideal of $X_\beta$ containing all commutators. 
\litem{(b)} If $\beta \leq \alpha$ is a limit ordinal, then 
$X_{\beta}$ is the closure of $\bigcup_{\gamma < \beta} X_\gamma$. 
\qed

The following theorem is an immediate consequence of Theorem VI.1.19, 
applied to the situation where $\q$ is abelian: 

\Theorem VI.1.21. {\rm([Swi65, Th.~2])} Each lower solvable 
Banach--Lie algebra is enlargeable.
\qed

Some of the methods used above for Banach--Lie groups have some potential to work 
in greater generality. Here are some ideas:

\Remark VI.1.22. If $G$ is a Lie group with Lie algebra $\g$, then 
$$ P(G) := C^\infty_*(I,G) := \{ \gamma \in C^\infty(I,G) \: \gamma(0) = \1 \} $$
also is a Lie group with Lie algebra $P(\g) := C^\infty_*(I,\g)$, endowed with the pointwise bracket. 
The logarithmic derivative 
$\delta \: P(G) \to C^\infty(I,\g)$ 
is a smooth map satisfying 
$\delta(\alpha\beta) = \delta(\beta) + \Ad(\beta)^{-1}.\delta(\alpha)$ and 
$T_\1(\delta)(\xi) = \xi'.$ (Lemma~II.3.3). 
As $[\xi,\eta]' = [\xi', \eta] + [\xi,\eta']$, it follows that 
$T_1(\delta) \: P(\g) \to C^\infty(I,\g)$ becomes a topological isomorphism of Lie algebras if 
$C^\infty(I,\g)$ is endowed with the bracket 
$$ [\xi,\eta](t) := \Big[\xi(t), \int_0^t \eta(\tau)\, d\tau\Big] + 
\Big[\int_0^t \xi(\tau)\, d\tau, \eta(t)\Big]. \leqno(6.1.2) $$
The evaluation map $\ev_1 \: P(\g) \to \g$ corresponds to the quotient map 
$$ C^\infty(I,\g) \to \g, \quad \xi \mapsto \int_0^1 \xi(\tau)\, d\tau. $$

If, in addition, $G$ is regular, then $\delta$ is a diffeomorphism, 
and it follows that $C^\infty(I,\g)$, endowed with the bracket (6.1.2), is 
integrable. Since this property is clearly necessary for the regular integrability of 
$\g$, Lie algebras with this property are called {\it pre-integrable} 
in [RK97] (see also [Les93]). 
\qed

If $G$ is a real BCH--Lie group, then a morphism 
$\kappa_G \: G \to G_\C$ to a complex BCH--Lie group $G_\C$ is called a 
{\it universal complexification} if for each other morphism 
$\alpha \: G \to H$ to a complex BCH--Lie group $H$, there exists a unique morphism 
$\beta \: G_\C \to H$ with $\alpha = \beta \circ \kappa_G$. It is well known that if 
$G$ is finite-dimensional, then a universal complexification always exists  
(cf.\ [Ho65] , [Ne99, Th.~XIII.5.6]), but it need not be locally injective, 
so that it may occur that $\dim_\C G_\C < \dim_\R G$. The following theorem shows 
that, due to the existence of non-enlargeable Lie algebras,  the situation 
becomes more complicated in infinite dimensions. 

\Theorem VI.1.23. {\rm(Existence of universal complexifications; [GN03], [Gl02c])} 
Given a real BCH--Lie group~$G$, 
let $N_G$ be the intersection of all kernels
of smooth homomorphisms from~$G$ to complex BCH--Lie
groups. Then~$G$ has a universal complexification 
 if and only if $N_G$ is a BCH--Lie subgroup of~$G$ and the complexification of
$\L(G)/\L^e(N_G)$ is enlargeable. 
\qed

Note that Theorem~VI.1.10 implies that if $G$ is $1$-connected, the 
existence of a universal complexification is equivalent to the enlargeability 
of $\L(G)/\L^e(N_G)$. 
In [GN03], one finds an example of a Banach-Lie group for which
$N_G$ fails to be a Lie subgroup ([GN03, Sect.~V]) 
and also examples where $N_G=\{\1\}$ but $\L(G)_\C$ is not enlargeable. 
The setting of BCH--Lie groups is the natural one for complexifications 
because if $\g$ is a locally exponential Lie algebra for which $\g_\C$ 
is locally exponential as a complex Lie algebra, then the local multiplication 
is complex analytic. This implies that $\g_\C$ is BCH which in turns entails that 
$\g$ is BCH. 

\subheadline{Localizing enlargeability} 

We call a norm $\|\cdot\|$ on a Lie algebra $\g$ {\it submultiplicative} if 
$\|[x,y]\| \leq \|x\|\|y\|$ for all $x,y \in \g$. A Banach--Lie algebra 
$(\g,\|\cdot\|)$ is called {\it contractive} if its norm is submultiplicative. 
For any contractive Lie algebra, we define 
$$\delta_\g := \inf \{ \|x\|\: 0 \not= x \in \Pi(\g)\} \in [0,\infty] $$
and note that $\g$ is enlargeable if and only if $\delta_\g > 0$, which is equivalent 
to the discreteness of the period group $\Pi(\g)$ 
(Theorem~VI.1.8). The following theorem is originally due to 
{\smc Pestov} who proved it with non-standard methods. A ``standard'' proof 
has been given in [Bel04] by {\smc Beltita}. 

\Theorem VI.1.24. {\rm(Pestov's Local Theorem on Enlargeability)} 
A contractive Banach--Lie algebra $\g$ is enlargeable if and only if 
there exists a directed family ${\cal H}$ of closed subalgebras of $\g$ 
for which $\bigcup {\cal H}$ is dense in $\g$ and 
$\inf \{ \delta_\h \: \h \in {\cal H}\} > 0$. 
\qed

Since for each finite-dimensional Lie algebra $\g$ the period group is 
trivial, we have $\delta_\g = \infty$, and the preceding theorem, applied to the 
directed family of finite-dimensional subalgebras of $\g$ leads to:

\Corollary VI.1.25. {\rm([Pe92], [Bel04])} If $\g$ is a Banach--Lie algebra containing 
a locally finite-dimensional dense subalgebra, then $\g$ is enlargeable. 
\qed

\Corollary VI.1.26. {\rm([Pe93b, Th.~7])} A Banach--Lie algebra $\g$ is 
enlargeable if and only if all its separable closed subalgebras are. 
\qed

\subheadline{Period groups for continuous inverse algebras} 

Another interesting class of cocycles arises for complete CIAs~$A$ ([Ne06c]). 
A continuous alternating bilinear map $\alpha \: A \times A \to E$, 
$E$ a locally convex space, is said to be a {\it cyclic $1$-cocycle} 
if 
$$ \alpha(ab,c) + \alpha(bc,a) + \alpha(ca,b) = 0 \quad \hbox{ for } \quad 
a,b,c \in A. $$
We write $ZC^1(A,E)$ for the set of all cyclic $1$-cocycles with values in $E$. 
Let $A_L = \gl_1(A)$ denote the Lie algebra $(A,[\cdot,\cdot])$ obtained by endowing 
$A$ with the commutator bracket. Then each cyclic cocycle defines a Lie algebra 
cocycle $\alpha \in Z^2_c(A_L,E)$ with respect to the trivial module structure on~$E$. 
To describe the universal cyclic cocycle, we endow 
$A \otimes A$ with the projective tensor topology and 
define $\la A, A \ra$ as the completion of the quotient space 
$$ (A \otimes A)/\oline{\span\{a \otimes a, a b \otimes c + bc \otimes a + ca \otimes
b;a,b, c \in A\}}. $$
We write $\alpha_u(a,b) := \la a,b \ra$ for the image of $a \otimes b$ in 
$\la A, A \ra$. 
Then the universal property of the projective tensor product implies that 
$$ {\cal L}(\la A,A \ra, E) \to ZC^1(A,E), \quad f \mapsto f \circ \alpha_u  $$
is a bijection for each complete locally convex space $E$, so that $\alpha_u$ is a 
universal cyclic $1$-cocycle. Of particular interest is the map 
$b_A \: \la A, A \ra \to A, \la a,b \ra \to [a,b]$ 
defined by the commutator bracket. Its kernel 
$$ HC_1(A) := \ker b_A \subeq \la A, A \ra $$
is the {\it first cyclic homology space of $A$} (cf.\ [Lo98]). 
We write $\omega_u$ for the universal cyclic $1$-cocycle, interpreted 
as a Lie algebra $2$-cocycle. Then the corresponding period map 
$$\per_{\omega_u} \: \pi_2(A^\times) \to \la A, A \ra $$ 
actually has values in the subspace $HC_1(A)$, which leads to a homomorphism 
$$ \per_{\omega_u} \: \pi_2(A^\times) \to HC_1(A). $$

It is a remarkable fact that this structure behaves nicely if we replace $A$ 
by a matrix algebra $M_n(A)$. 
Let $\eta_n \: A  \to M_n(A), a \mapsto a E_{11}$ denote the natural inclusion map 
and observe that it induces  maps  
$\la A, A \ra \to \la M_n(A), M_n(A)\ra$, taking 
$HC_1(A)$ into $HC_1(M_n(A))$. In the other direction, we have maps 
$$ \tr^{(2)} \: \la M_n(A), M_n(A)\ra \to \la A, A \ra, \quad 
\la (a_{ij}), (b_{ij}) \ra \mapsto \sum_{i,j=1}^n \la a_{ij}, b_{ji}\ra, $$
and the topological version of the Morita invariance of
cyclic homology ([Lo98, Th.~2.2.9]) asserts that these maps restrict to isomorphisms 
$HC_1(M_n(A)) \to HC_1(A)$. This leads to extensions of the universal cocycle 
to a cocycle $\omega_u^n \in Z^2_c(\gl_n(A), \la A, A \ra)$ with 
$\eta_n^*\omega_u^n = \omega_u$ for each $n \in \N$. In terms of the 
tensor product structure $\gl_n(A) \cong A \otimes \gl_n(\K)$, it is given by 
$$ \omega_u^n(a \otimes x, b \otimes y) = \tr(xy)\la a,b\ra. $$

To explain the corresponding compatibility on the level of period homomorphisms, 
we define the {\it topological $K$-groups} of $A$ by 
$$ K_{i+1}(A) := \indlim \pi_i(\GL_n(A)) \quad \hbox{ for } \quad i \in \N_0, $$
where the direct limit on the right hand side corresponds to the embeddings 
$$\GL_n(A) \to \GL_{n+1}(A), \quad a \mapsto \pmatrix{a & 0 \cr 0 & \1\cr}, $$
induced by the corresponding embeddings $M_n(A) \into M_{n+1}(A)$. The group $K_0(A)$ is defined 
as the Grothendieck group of the abelian monoid $\indlim \pi_0(\Idem(M_n(A)))$, endowed 
with the addition $[e]+ [f] := \Big[\pmatrix{e & 0 \cr 0 & f\cr}\Big]$ (cf.\ [Bl98]). 

The naturality of the universal cocycles now implies 
that the period maps 
$$ \per_{\omega_u^n} \: \pi_2(\GL_n(A)) \to HC_1(A)  $$
combine to a group homomorphism 
$$ \per_A^1 \: K_3(A) = \indlim \pi_2(\GL_n(A)) \to HC_1(A), $$
which is a natural transformation from  the functor $K_3$ with values in abelian 
groups to the functor $HC_1$ with values in complete locally convex spaces. 

It is of some interest to know whether the group 
$$ \Pi^1_A := \im(\per_A^1) \subeq HC_1(A) $$
is discrete. 
If this is the case, then each period homomorphism 
$\per_{\omega_u^n}$ has discrete image, which implies that the 
corresponding central extension $\hat\gl_n(A)$ of the Lie algebra $\gl_n(A)$ by 
$\la M_n(A), M_n(A)\ra$ is enlargeable. 

This central extension is of particular interest when restricted 
to the subalgebra $\sL_n(A) := \oline{[\gl_n(A), \gl_n(A)]}$. We 
define a Lie bracket on $\la M_n(A), M_n(A) \ra$ 
by 
$$ [\la  a,b\ra, \la a',b'\ra] := \la [a,b], [a',b']\ra, $$
turning it into a locally convex Lie algebra. 
Now the bracket map of $M_n(A)$ induces a generalized central extension 
$$ q \: \hat\sL_n(A) := \la M_n(A), M_n(A) \ra \to \sL_n(A), \quad 
\la a, b \ra \mapsto [a,b] $$
with $\ker q = HC_1(M_n(A)) \cong HC_1(A)$,  
which is a universal generalized central extension, called the 
{\it topological Steinberg--Lie algebra} ([Ne03, Ex.~4.10]). 
The enlargeability criterion in Theorem~VI.1.6 immediately leads to: 

\Theorem VI.1.27. If the subgroup $\Pi^1_A$ of $HC_1(A)$ is discrete, then all 
Steinberg--Lie algebras $\hat\sL_n(A)$ are enlargeable. 
\qed

\Remark VI.1.28. Let $\omega \in ZC^1(A,E)$ be a cyclic cocycle, considered as a 
Lie algebra cocycle on $(A,[\cdot,\cdot])$. Then 
the adjoint action of $A$ on the Lie algebra $\hat A_L := E \oplus_\omega A_L$ 
integrates to an action of $A^\times$ by 
$$ g.(z,a) = (z - \omega(ag^{-1},g),gag^{-1}). $$
In view of Remark~V.2.14(f), this implies the triviality of the flux homomorphism 
$$F_\omega \: \pi_1(A^\times) \to H^1_c(A_L,E) \subeq {\cal L}(A,E). 
\qeddis

According to [Bos90], we have for each complex complete CIA $A$ natural isomorphisms 
$$ \beta_A^i \: K_i(A) \to K_{i+2}(A),\quad i \in \N_0. $$ 
This is an abstract version of Bott periodicity. In particular, the range of 
$$ P_A := \per^1_A \circ \beta_A^1 \: K_1(A) \to HC_1(A) $$
coincides with $\Pi_1^A$. The main advantage of this picture is that 
natural transformations from $K_1$ to $HC_1$ are unique, which leads to 
the explicit formula 
$$ P_A([g]) = \sum_{i,j} \la (g^{-1})_{ij},g_{ji} \ra \quad \hbox{ for } \quad 
[g] \in K_1(A), g \in \GL_n(A) $$
(cf.\ [Ne06c]). 
If, in addition, $A$ is commutative, then 
$HC_1(A)$ is the completion of the quotient $\Omega^1(A)/d_A(A)$, where 
$\Omega^1(A)$ is the (topological) universal differential module of $A$. 
In these terms, we then have 
$$ P_A(g) = \la \det(g)^{-1}, \det(g)\ra = [\det(g)^{-1} d_A(\det(g))], $$
which leads to 
$$ \im(P_A) = P_A([A^\times]) = \{ [a^{-1} d_A(a)] \: a \in A^\times \}. $$

\Examples VI.1.29. (1) For $A=C^\infty_c(M,\C)$, $M$ a $\sigma$-compact finite-dimensional 
manifold, we have 
$$ HC_1(A) \cong \Omega^1_c(M,\C)/dC_c^\infty(M,\C)$$ 
([Co94], [Mai02]). 
Moreover, 
$ M_n(A) \cong C_c^\infty(M, M_n(\C))$ 
and 
$$\omega_u^n(f,g) = [\tr(f \cdot dg)] $$
is a cocycle of product type, which implies that its period group coincides with 
the group 
$$ \im(P_A) = \delta(C^\infty_c(M,\C^\times))/dC_c^\infty(M,\C)$$ 
of integral cohomology classes in $H^1_{\rm dR,c}(M,\C)$, which is 
discrete (Theorem~V.2.17). 

(2) For $A = C(X,\C)$, where $X$ is a compact space, {\smc Johnson}'s Theorem 
entails that $\Omega^1(A)$, and hence $HC_1(A) \cong \Omega^1(A)/d_A(A)$, 
vanish ([BD73, Th.\ VI.12]). This further implies that for 
each $C^*$-Algebra $A$ the homomorphism $P_A$ vanishes. 
\qed

A particularly interesting class of Fr\'echet CIAs are the {\it $d$-dimensional smooth 
quantum tori}. 
These algebras are parametrized by 
skew-symmetric matrices $\Theta \in \Skew_d(\R)$, as follows. 
They are topologically generated by 
$d$ invertible elements $u_1,\ldots, u_d$, together with their inverses, 
satisfying the commutation relations 
$$u_p u_q = e^{2\pi i \Theta_{pq}} u_q u_p \quad \hbox{for} \quad 1 \leq p,q\leq d. $$
Moreover, 
$$ A_\Theta = \Big\{ 
\sum_{I \in \Z^d} \alpha_I u^I \: 
(\forall k \in \N) \sum_I |I|^k |\alpha_I| < \infty\Big\}, $$
where $|I| = i_1 + \ldots + i_d$ and $u^I := u_1^{i_1}\cdots u_d^{i_d}$, so that, 
as a Fr\'echet space, $A_\Theta$ is isomorphic to the space of smooth functions on the 
$d$-dimensional torus. 
In particular, we have the commutative case $A_0 \cong C^\infty(\T^d,\C)$. 
The following theorem characterizes those for which the image of $P_A$ is discrete ([Ne06c]): 

\Theorem VI.1.30. For the $d$-dimensional smooth quantum torus $A_\Theta$,  
the group $\im(P_{A_\Theta})$ is discrete if and only if 
$d \leq 2$ or the matrix $\Theta$ has rational entries. 
\qed

An interesting consequence of the preceding theorem is that 
there exists a CIA $A$ for which $\im(P_A)$ is not discrete. The smallest examples 
are of the form $A := C^\infty(\T,A_\Theta)$, where 
$\Theta 
= \pmatrix{ 0 & \theta \cr -\theta & 0\cr}, \theta \in\R\setminus \Q$, so that 
$A_\Theta$ is a so-called {\it irrational rotation algebra}.

\subheadline{VI.2. Integrability of non-locally exponential Lie algebras} 

After the discussion of the enlargeability of locally exponential and Banach--Lie algebras 
in the preceding subsection, we now turn to more general classes of Lie algebras. 
Unfortunately, there is no general theory beyond the locally exponential class, 
so that all positive and negative results are quite particular. 

We start with a discussion of some obstructions to the integrability to an analytic 
Lie group, then turn to complexifications of Lie algebras of vector fields, and 
finally to Lie algebras of formal vector fields, resp., Lie algebras of germs. 

\Proposition VI.2.1. {\rm([Mil84, Lemma 9.1])} Let $G$ be a connected analytic Lie group.
Then each closed ideal $\n \trile \L(G)$ is invariant under $\Ad(G)$. 
\qed

\Corollary VI.2.2. If $\g$ is a Lie algebra containing 
a closed ideal which is not stable, then $\g$ is not integrable to an analytic 
Lie group with an analytic exponential function. 
\qed

\Remark VI.2.3. Proposition~VI.2.1 implies that 
the Lie group $\Diff(M)$ of all diffeomorphisms of a compact manifold $M$ 
does not possess an analytic Lie group structure for which its Lie algebra 
is ${\cal V}(M)$. Indeed,  
for each non-dense open subset $K \subeq M$, the subspace 
$$ {\cal V}(M)_K := \{ X \in {\cal V}(M) \: X\res_K = 0\} $$
is a closed ideal of ${\cal V}(M)$ not invariant under $\Ad(\Diff(M))$ because 
$\Ad(\phi).{\cal V}(M)_K = {\cal V}(M)_{\phi(K)}$  for 
$\phi \in \Diff(M)$. 
\qed

The situation improves if we restrict our attention to analytic diffeomorphisms: 

\Theorem VI.2.4. {\rm([Les82/83])} Let $M$ be a compact analytic manifold and 
${\cal V}^\omega(M)$ the Lie algebra of analytic vector fields on $M$. 
Then ${\cal V}^\omega(M)$ carries a natural Silva space structure, turning it into a 
topological Lie algebra, and the group $\Diff^\omega(M)$ of 
analytic diffeomorphisms 
carries a smooth Lie group structure 
for which ${\cal V}^\omega(M)^{\rm op}$ is its Lie algebra. 
\qed

It is shown by {\smc Tognoli} in [Ta68] that the group $\Diff^\omega(M)$, 
$M$ a compact analytic manifold, carries no analytic Lie group structure 
(cf.\ [Mil82, Ex.~6.12]). That there is no analytic Lie group with an 
analytic exponential function and Lie algebra ${\cal V}^\omega(M)$ 
can be seen by verifying that the map $(X,Y) \mapsto \Ad(\Fl^X_1).Y$ 
is not analytic on a $0$-neighborhood in ${\cal V}^\omega(M) \times 
{\cal V}^\omega(M)$ (cf.\ [Mil82, Ex.~6.17]).

The following non-integrability result is quite strong because it does not 
assume the existence of an exponential function. Its outcome is that complexifications 
of Lie algebras of vector fields are rarely integrable. For complexifications 
of Lie algebras of ILB--Lie groups, similar results 
are described by {\smc Omori} in [Omo97, Cor.~4.4].

\Theorem VI.2.5. {\rm([Lem97])} Let $M$ be a compact manifold of positive dimension. 
Then the complexifications $\g_\C$ of the following Lie algebras $\g$ are not integrable: 
\litem{(1)} The Lie algebra ${\cal V}(M)$ of smooth vector fields on $M$. 
\litem{(2)} If $M$ is analytic, the Lie algebra ${\cal V}^\omega(M)$ 
of analytic vector fields on $M$. 
\litem{(3)} If $\Omega$ is a symplectic $2$-form on $M$, the Lie algebra 
${\cal V}(M,\Omega) := \{ X \in {\cal V}(M) \: {\cal L}_X\Omega = 0\}$ 
of symplectic vector fields on $M$. 
\litem{(4)} If $M$ is analytic and $\Omega$ is an analytic 
symplectic $2$-form on $M$, the Lie algebra 
${\cal V}^\omega(M,\Omega)$ of analytic symplectic vector fields on $M$. 

\Proof. (Idea) Lempert's 
proof is based on the following result, which is obtained by PDE methods: 
If $\xi \: \R \to \g_\C$ is a smooth curve such that for each 
$x \in \g_\C$ the IVP 
$$ \gamma(0) = x, \quad \dot\gamma(t) = [\xi(t),\gamma(t)] $$
has a smooth solution, then $\xi(0) \in \g$. 

For (1) he gives another argument, based on the fact that 
$$ \Aut({\cal V}(M)_\C) 
\cong \Aut({\cal V}(M)) \rtimes \{\1,\sigma\} 
\cong \Diff(M) \rtimes \{\1,\sigma\}, \leqno(6.2.1) $$
where $\sigma$ denotes the complex conjugation on ${\cal V}(M)_\C$. 
The first isomorphism is obtained in [Lem97], using [Ame75], and 
the second is an older result of {\smc Pursell} and {\smc Shanks} ([PuSh54]; 
cf.\ Theorem~IX.2.1). 

Clearly (6.2.1) implies that for any 
connected Lie group $G$ with Lie algebra $\L(G) = {\cal V}(M)_\C$, the group 
$\Ad(G) \subeq \Aut({\cal V}(M)_\C)$ preserves the real 
subspace ${\cal V}(M)$. Taking derivatives of orbit maps, this 
leads to the contradiction 
$[{\cal V}(M)_\C, {\cal V}(M)] \subeq {\cal V}(M)$.
\qed

\Theorem VI.2.6. {\rm([Omo81])} For any non-compact $\sigma$-compact smooth 
manifold $M$ of positive dimension, 
the Lie algebra ${\cal V}(M)$ is not integrable to any Lie group 
with an exponential function. 

\Proof. (Sketch) If $G$ is a Lie group with Lie algebra $\L(G) = {\cal V}(M)$ 
and an exponential function, then for each 
$X \in {\cal V}(M)$ we obtain a smooth $1$-parameter group 
$t \mapsto \Ad(\exp_G(tX))$ of automorphisms of ${\cal V}(M)$ 
with generator $\ad X$. By [Ame75, Thm.~2], 
$\Aut({\cal V}(M)) \cong \Diff(M)$, so that we obtain a 
one-parameter group $\gamma_X$ of $\Diff(M)$ which then is shown to coincide 
with the flow generated by $X$ (cf.\ Lemma~II.3.10; [KYMO85, Sect.~3.4]). 
This contradicts the existence of
non-complete vector fields on $M$. 
\qed

Since the BCH series can be used to defined a Lie group structure on 
any nilpotent locally convex Lie algebra, all these Lie algebras are integrable. 
The following theorem shows that the integrability problem 
for solvable locally convex Lie algebras contains the integrability problem for 
continuous linear operators on locally convex spaces, which is highly 
non-trivial (Problem VI.1). E.g., if $M$ is a finite-dimensional $\sigma$-compact manifold 
and $X \in {\cal V}(M)$ a vector field, then the corresponding derivation of the Fr\'echet 
algebra $C^\infty(M,\R)$ is integrable if and only if the vector field 
$X$ is complete. 

\Theorem VI.2.7. Let $E$ be a locally convex space and $D \in \gl(E)$. 
Then the solvable Lie algebra $\g := E \rtimes_D \R$ 
with the bracket 
$[(v,t), (v',t')] := (tDv' - t'Dv, 0)$ is integrable 
if and only if $D$ is integrable to a smooth $\R$-action on $E$. 

\Proof. If $D$ is integrable to a smooth representation 
$\alpha \: \R \to \GL(E)$ with $\alpha'(0) = D$, 
then  the semi-direct product $G := E \rtimes_\alpha \R$ is a Lie group with 
the Lie algebra $\g$. 

Suppose, conversely, that $G$ is a connected Lie group with Lie algebra $\g$. 
Replacing $G$ by its universal covering group, we may 
assume that $G$ is $1$-connected. Then the regularity of the additive 
group $(\R,+)$ implies the existence of a smooth homomorphism 
$\chi \: G \to \R$ with $\L(\chi) = q$, where 
$q(v,t) = t$ is the projection $\g = E \rtimes_D \R \to \R$ (Theorem~III.1.5). 

Using Gl\"ockner's Implicit Function Theorem ([Gl03a]), 
it follows that $\ker \chi$ is a submanifold of $G$ and there exists a smooth 
curve $\gamma \: \R \to G$ with $\gamma(0) = \1$ and $\chi \circ \gamma = \id_\R$. 

Next we observe that $[\g,\g] \subeq E$ implies that 
$E$ is $\Ad(G)$-invariant, so that $\Ad_E(g) := \Ad(g)\res_E$ defines 
a smooth action of $G$ on $E$ whose derived representation is given by $\ad_E(x,t) = tD$. 
We now put 
$\alpha(t) := \Ad_E(\gamma(t))$ and observe that 
$$ \delta(\alpha)(t) = \ad_E(\delta(\gamma)(t)) 
= q(\delta(\gamma)(t))\cdot D 
= \delta(\chi \circ \gamma)(t)\cdot D = D. $$
Hence $D$ is integrable. 
\qed

\Example VI.2.8. Let $\gf_n(\R)_{-1} := \R^n[[x_1,\ldots, x_n]]$ denote 
the space of all $\R^n$-valued formal power series in $n$ variables, 
considered as the Lie algebra of formal vector fields, endowed with the bracket 
$$ [f,g](x) := dg(x)f(x) - df(x)g(x), $$
which makes sense on the formal level because if $f$ is homogeneous of degree $p$ 
and $g$ is homogeneous of degree $q$, then $[f,g]$ is of degree $p+q-1$. 

We have already seen in Example IV.1.14 that the subalgebra 
$\gf_n(\R)$ of all elements with vanishing constant term is the Lie algebra 
of the Fr\'echet--Lie group $\Gf_n(\R)$ of formal diffeomorphisms of $\R^n$ 
fixing $0$. We obviously have the split short exact sequence 
$$ \0 \to \gf_n(\R) \into \gf_n(\R)_{-1} \onto \R^n \to \0, $$
where $\R^n$ is considered as an abelian Lie algebra, corresponding to the constant 
vector fields. 

We claim that the Lie algebra $\gf_n(\R)_{-1}$ is not integrable to any Lie group 
with an exponential function. This strengthens a statement in [KYMO85, p.80], 
that it is not integrable to a $\mu$-regular Fr\'echet--Lie group. 
Let us assume that $G$ is a Lie group with 
Lie algebra $\gf_n(\R)_{-1}$. With a similar argument as in the proof of 
Theorem~VI.2.7, one can show that for each constant function $x$ the operator $\ad x$ 
on $\gf_n(\R)_{-1}$ is integrable. 
We consider the constant function $e_1$. Then 
$[e_1, g] = {\partial g \over \partial x_1},$
and we can now justify as in Example II.3.13 that $\ad e_1$ is not integrable, 
hence that $\gf_n(\R)_{-1}$ cannot be integrable to any Lie group with an 
exponential function. 
\qed

The preceding example shows that the constant terms create problems in integrating 
Lie algebras of formal vector fields, which is very natural because the formal completion 
distinguishes the point $0 \in \R^n$. 
A similar phenomenon arises in the context of groups of germs of local diffeomorphisms. 
For germs of functions in $0$, the non-integrability of vector fields with 
non-zero constant term follows from the fact that all automorphisms preserve the 
unique maximal ideal of functions vanishing in $0$ (cf.\ [GN06]). 

Let $\gs_n(\R)_{-1}$ denote the space of germs of smooth maps $\R^n \to \R^n$ in $0$, 
identified with germs of vector fields in $0$. 
According to [RK97, Sect.~5.2], this space carries a natural Silva structure, turning it into a 
locally convex Lie algebra. Let $\gs_n(\R)$ denote the subspace of all germs vanishing in $0$ 
and $\gs_n(\R)_1$ the set of germs vanishing of second order in $0$. 

\Theorem VI.2.9. {\rm([RK97, Th.~3])} The 
group $\Gs_n(\R)$ of germs of diffeomorphism of $\R^n$ in $0$ fixing $0$ 
carries a Lie group structure for which the Lie algebra is the space $\gs_n(\R)$ of 
germs of vector fields vanishing in $0$. 

We have a semidirect product decomposition $\Gs_n(\R) \cong \Gs_n(\R)_1 \rtimes\GL_n(\R)$, 
where $\Gs_n(\R)_1$ is the normal subgroup of those germs $[\phi]$ for which 
$\phi - \id_{\R^n}$ vanishes of order $2$. The map 
$$ \Phi \: \gs_n(\R)_{1} \to \Gs_n(\R)_{1}, \quad \xi \mapsto \id + \xi $$
is a global diffeomorphism. 
\qed

In view of the preceding theorem, it is a natural problem 
to integrate Lie algebras of germs of vector fields 
vanishing in the base point to Lie groups of germs of diffeomorphisms. 
This program is carried out by {\smc Kamran} and {\smc Robart} 
in several papers (cf.\ [RK97], [KaRo01/04], [Rob02]). It results in several interesting classes 
of Silva--Lie groups of germs of smooth and also analytic local diffeomorphisms, 
where the corresponding Silva--Lie algebras depend on certain 
parameters which are used to obtain a good topology. 

\Example VI.2.10. The formal analog of the Lie algebra $\gs_1(\R)_1$ is 
the Lie algebra $
\gf_1(\R)_1$ which is pro-nilpotent, 
hence in particular BCH. In contrast to this fact, {\smc Robart} observed 
that $\gs_1(\R)_1$ is not BCH. In fact, for the elements 
$\xi(x) = ax^2$, $\lambda \in \R$ and $\eta(x) = x^3$, we have 
$$\sum_{n = 0}^\infty \big((\ad \xi)^n \eta\big)(x) 
= x^3 \sum_{n = 0}^\infty a^n n! x^n, $$
which converges for no $x \not=0$ if $a \not=0$. 
With {\smc Floret}'s results from [Fl71, p.155], it follows that 
this series does not converge in the Silva space $\gs_1(\R)_1$, so that Theorem~IV.1.7 
shows that $\gs_1(\R)_1$ is not BCH. 
\qed

The following proposition is a variant of {\smc E.~Borel}'s theorem on 
the Taylor series of smooth functions. It provides an interesting connection 
between the smooth global and the formal perspective on diffeomorphism groups. 

\Proposition VI.2.11. {\rm(Gl\"ockner)} Let 
$M$ be a smooth finite-dimensional manifold, 
$m_0 \in M$ and $\Diff_c(M)_{m_0}$ the stabilizer of $m_0$. 
For each $\phi \in \Diff_c(M)_{m_0}$, let 
$T^\infty_{m_0}(\phi) \in \Gf_n(\R)$ denote the Taylor series of $\phi$ in $m_0$ 
with respect to some local chart. 
Then the map 
$$ T^\infty_{m_0} \: \Diff_c(M)_{m_0,0} \to \Gf_n(\R)_0 $$
is a surjective homomorphism of Lie groups, where 
$\Gf_n(\R)_0$ is the subgroup of index $2$, consisting of those formal 
diffeomorphisms $\psi$ with $\det(T_0(\psi)) > 0$. 
\qed

\Example VI.2.12.  Let $\gh_n(\C)$ denote the 
space of germs of holomorphic maps $f \: \C^n \to \C^n$ in $0$ satisfying $f(0) = 0$.  
We endow this space with the locally convex 
direct limit topology of the Banach spaces $E_k$ of holomorphic 
functions on the closed unit disc of radius ${1\over k}$ 
in $\C^n$ (with respect to any norm). Thinking of the elements of $\gh_n(\C)$ 
as germs of vector fields in $0$ leads to the Lie bracket 
$$ [f,g](z) := dg(z)f(z) - df(z) g(z), $$
which turns $\gh_n(\C)$ into a topological Lie algebra. 

The set $\Gh_n(\C)$ of all germs $[f]$ with $\det(f'(0)) \not=0$ is an open subset 
of $\gh_n(\C)$ which is a group with respect to composition 
$[f][g] := [f \circ g]$. In [Pis77], {\smc Pisanelli} 
shows that composition and inversion in 
$\Gh_n(\C)$ are holomorphic, so that $\Gh_n(\C)$ is a complex Lie group 
with respect to the manifold structure it inherits as an open subset of 
$\gh_n(\C)$. This Lie group has a holomorphic exponential function 
which is not locally surjective, where the latter fact can be obtained by 
adapting Sternberg's example $f(z) = e^{2\pi i \over m} z + p z^{m+1}$  
(Example~IV.1.14) ([Pis76]). 

Note that  $\Gh_n(\C) \cong \Gh_n(\C)_1 \rtimes \GL_n(\C)$, 
where $\Gh_n(\C)_1$ is the subgroup of all diffeomorphisms with linear term 
$\id_{\C^n}$. 
\qed

\Remark VI.2.13. Let $\g(A)$ be a symmetrizable Kac--Moody Lie algebra. 
In [Rod89], {\smc Rod\-ri\-guez-Carrington} describes certain 
Fr\'echet completions of $\g(A)$, 
including smooth $\g^\infty(A)$ and analytic versions 
$\g^\omega(A)$, which are BCH--Lie algebras ([Rod89,  Prop.~1]). 
Corresponding groups are constructed for the unitary real forms by unitary 
highest weight modules of $\g(A)$, as subgroups of the unitary groups of a 
Hilbert space (Corollary~IV.4.10). 
In [Su88], {\smc Suto} obtains closely related results, 
but no Lie group structures. 

In a different direction, {\smc Leslie} describes in [Les90] 
a certain completion $\oline\g(A)$ of $\g(A)$ which leads to a Lie group 
structure on the space $C^\infty([0,1],\oline\g(A))$, corresponding to the 
natural Lie algebra structure on this space. One thus obtains an integrable 
Lie algebra extension of $\oline\g(A)$ in the spirit of pre-integrable Lie algebras 
(Remark~IV.1.22). For an approach to Kac--Moody groups in the context of diffeological 
groups, we refer to [Les03] (cf.\ [So84]).
\qed

\subheadline{Open Problems for Section VI} 

\Problem VI.1. (Generators of smooth one-parameter groups) 
Let $E$ be a locally convex space and $D \: E \to E$ a continuous 
linear endomorphism. Characterize those linear operators $D$ 
for which there exists a homomorphism $\alpha \: \R \to \GL(E)$ 
defining a smooth action of $\R$ on $E$. In view of Theorem~VI.2.7, 
this is equivalent to the integrability of the $2$-step solvable 
Lie algebra 
$\g := E \rtimes_D \R$.

If $E$ is a Banach space, then each $D$ integrates to  
a homomorphism 
$\alpha$ which is continuous with respect to the norm topology on 
$\GL(E)$ and given by the convergent exponential series 
$\alpha(t) := \sum_{k = 0}^\infty {t^k \over k!} D^k.$

Since for each smooth linear $\R$-action on $E$, given by some $\alpha$ as above, 
the infinitesimal generator $\alpha'(0)$ is everywhere defined, 
this problem is not a problem about operators which are unbounded 
in the sense that they are only defined on dense subspaces. 
In some sense, the passage from Banach spaces to locally convex spaces 
takes care of this problem. If, e.g., $\alpha \: \R \to \GL(E)$ 
is a strongly continuous one-parameter group on a Banach space 
$E$, then the subspace $E^\infty \subeq E$ of smooth vectors 
carries a natural Fr\'echet topology inherited from the embedding 
$$E^\infty \into C^\infty(\R,E),\quad  v \mapsto \alpha^v, \quad 
\alpha^v(t) = \alpha(t)v, $$
and the induced one-parameter group 
$\alpha^\infty \: \R \to \GL(E^\infty)$ defines a smooth action. 
In this sense, each generator of a strongly 
continuous one-parameter group also generates a smooth one-parameter 
group on a suitable Fr\'echet space. 
\qed

\Problem VI.2. (Integrability of $2$-step solvable Lie algebras) Theorem VI.2.7 gives an integrability criterion 
for solvable Lie algebras of the type $\g = E \rtimes_D \R$. 

Since abelian Lie algebras are integrable for trivial reasons, it is natural 
to address the integrability problem for solvable Lie algebras by first 
restricting to algebras of {\it solvable class~$2$}, i.e., 
$D^1(\g) := \oline{[\g,\g]}$ is an abelian ideal of $\g$. Clearly, the adjoint 
action defines a natural topological module structure for the abelian Lie algebra 
$W := \g/D^1(\g)$ on $E := D^1(\g)$. Here are some problems concerning this situation: 
\litem{(1)} Does the integrability of $\g$ imply that the Lie algebra module  
structure of $W$ on $E$ integrates to a smooth action of the Lie group 
$(W,+)$ on $E$?  If $E$ is finite-dimensional, this can be proved by 
an argument similar to the proof of Theorem~VI.2.7. 
\litem{(2)} Assume that the Lie algebra module structure of $W$ on $E$ integrates 
to a smooth action of $(W,+)$. Does this imply that $\g$ is integrable? 

If $\g \cong V \rtimes W$ is a semidirect product, the latter is obvious, but 
if $\g$ is a non-trivial extension of $W$ by $V$, the situation is more complicated. 
Note that all solvable Banach--Lie algebras are integrable by Theorem~VI.1.21. 
\qed

\Problem VI.3. Is the group $\Gs_n(\R)_1$ of germs of diffeomorphisms $\phi$ of $\R^n$ fixing $0$, 
for which the linear term of $\phi - \id_{\R^n}$ vanishes, exponential? 
\qed

\Problem VI.4. Let $G$ be a regular Lie group. 
Is every finite codimensional closed subalgebra $\h \subeq \L(G)$ 
integrable to an integral subgroup? For $\mu$-regular groups this follows from 
Theorem~III.2.8. 
\qed

\Problem VI.5. Is the group $\Gh_n(\C)$ defined in Example VI.2.12 a regular Lie group? 
Is the subgroup $\Gh_n(\C)_1$ an exponential Lie group? (cf.\ Problem~VI.3) 
\qed

\Problem VI.6. Does Pestov's Theorem VI.1.24 
generalize to locally exponential Lie algebras? 
\qed

\Problem VI.7. For quotient maps $q \: E \to Q$ of Fr\'echet spaces, we may use 
[MicE59] to find a continuous linear cross section $\sigma \: Q \to E$, which implies in 
particular that $q$ defines a topologically trivial fiber bundle. For more general 
locally convex spaces, cross sections might not exist, but it would still be 
interesting if quotient maps of locally convex spaces are Serre fibrations, 
i.e., have the homotopy lifting property for cubes (cf.\ [Bre93]). If this is the 
case, the long exact homotopy sequence would also be available for quotient maps 
of locally exponential Lie groups, which would be an important tool to calculate homotopy 
groups of such Lie groups. 
\qed

\Problem VI.8. Prove an appropriate version of Theorem~VI.1.23 on the existence of a 
universal complexification for locally exponential Lie algebras. 

Note that this already 
becomes an interesting issue on the level of Lie algebras because the complexification 
of a locally exponential Lie algebra need not be locally exponential. 
In fact, in Example~IV.4.6 we have seen an exponential Lie algebra $\g$ 
containing an unstable closed subalgebra $\h$. If $\g_\C$ is locally exponential, 
as a complex Lie algebra, then the local multiplication in $\g_\C$ is holomorphic, 
so that $\g$ is BCH, contradicting the existence of unstable closed subalgebras. 
\qed

\sectionheadline{VII. Direct limits of Lie groups} 

\nin The systematic study of Lie group structures on direct limit 
Lie groups $G = \indlim G_n$ was started in the 1990s by {\smc J.~Wolf} 
and his coauthors ([NRW91/93]). 
They used certain conditions on the groups $G_n$ and the maps $G_n \to G_{n+1}$ 
to ensure that the direct limit 
groups are locally exponential.  
Since not all direct limit groups are locally exponential (Example~VII.1.4(c)), 
their approach does not cover all cases. The picture for countable 
direct limits of finite-dimensional Lie groups was nicely completed 
by {\smc Gl\"ockner} who showed that arbitrary countable limits of finite-dimensional 
Lie groups exist ([Gl03b/05a]). The key to these results are general 
construction principles for direct limits of finite-dimensional manifolds. 
These results are discussed in Section VII.1. 
In Section VII.2, we briefly turn to other types of direct limit constructions where 
the groups $G_n$ are infinite-dimensional Lie groups.

\subheadline{VII.1. Direct limits of finite-dimensional Lie groups} 

\Theorem VII.1.1. {\rm([Gl05])} {\rm(a)} For 
every sequence $(G_n)_{n \in \N}$ of finite-dimensional 
Lie groups $G_n$ with morphisms $\phi_n \: G_n \to G_{n+1}$, the direct limit 
group $G := \indlim G_n$ carries a regular Lie group structure. 
The model space $\L(G)  \cong \indlim \L(G_n)$ is countably dimensional 
and carries the finest locally convex topology, and $G$ has the universal 
property of a direct limit in  the category of Lie groups.  

{\rm(b)} Every countably dimensional locally finite 
Lie algebra $\g$, endowed with the finest locally convex topology, 
is integrable to a regular Lie group $G$. 

{\rm(c)} Every connected regular Lie group $G$ whose Lie algebra is countably 
dimensional, locally finite and carries the finest locally convex topology is a direct limit 
of finite-dimensional Lie groups.   
\qed

In the following, we shall call the class of Lie groups described by the preceding 
theorem {\it locally finite-dimensional} (regular) Lie groups. 

\Remark VII.1.2. (a) Beyond countable directed systems, several serious obstacles arise. 
First of all, for countably dimensional vector spaces, the finest locally 
convex topology coincides with the finest topology for which all 
inclusions of finite-dimensional subspaces are continuous. This is crucial 
for many arguments in this context. If $E$ is not of countable dimension, 
the addition on $E$ is not continuous for the latter topology. Similar problems 
occur for uncountable direct limits of topological groups: in many cases the direct 
limit topology does not lead to a continuous multiplication 
(cf.\ [Gl03b] for more details). 

(b) Any countably dimensional space $E$, endowed with the finest locally convex topology 
can be considered as a direct limit space of finite-dimensional subspaces $E_n$ of 
$\dim E_n = n$. Since each $E_n$ is a closed subspace which is Banach, and all 
inclusions $E_n \to E_{n+1}$ are compact operators, $E$ is an LF space and a Silva space 
at the same time. 
\qed

\Theorem VII.1.3. {\rm([Gl05/06d])} Let $G$ be a locally finite-dimensional Lie group. Then the 
following assertions hold: 
\litem{(1)} Every subalgebra $\h \subeq \L(G)$ integrates to an integral subgroup. 
\litem{(2)} Every closed subgroup $H$ is a split submanifold, so that 
$H$ is a locally finite-dimensional Lie  group, and 
the quotient space $G/H$ carries a natural manifold structure. 
\litem{(3)} Every locally compact subgroup $H \subeq G$ 
is a finite-dimensional Lie group. 
\litem{(4)} $G$ does not contain small subgroups. 
\qed

\Example VII.1.4. (a) One of the most famous examples of a direct limit Lie group 
is the group 
$$ \GL_\infty(\R) := \indlim \GL_n(\R) $$
with the connecting maps 
$$ \phi_n \: \GL_n(\R) \to \GL_{n+1}(\R), \quad a \mapsto \pmatrix{ a & 0 \cr 0 & 1\cr}. $$
Its Lie algebra is the Lie algebra 
$\gl_\infty(\R)$ of all $(\N \times \N)$-matrices with 
only finitely many non-zero entries (cf.\ [NRW91], [Gl03b]). 

In [KM97, Thm.~47.9], it is shown that every subalgebra $\h$ of $\gl_\infty(\R)$ is integrable 
to an integral subgroup, 
which is a special case of Theorem VII.1.3. Here $\h$ is even BCH. 

(b) In the context of $C^*$-algebras, direct limits of finite-dimensional ones are particularly 
interesting objects. On the level of unit groups one encounters in particular 
groups of the form 
$$ G := \indlim \GL_{2^n}(\C), \quad \phi_n(a) = \pmatrix{ a & 0 \cr 0 & a\cr}. $$

(c) Let $E := \C^{(\N)}$ be the free vector space with basis 
$(e_n)_{n \in \N}$ and $D \in {\cal L}(E)$ be defined by 
$D(e_n) = 2\pi i n e_n$ (cf.\ Example~II.5.9(a)). 
Then the Lie algebra $\g := E \rtimes_D \R$ is locally finite and 
we obtain a corresponding locally finite-dimensional Lie group 
$G = E \rtimes_\alpha \R$, where $\alpha(t) = e^{tD}$. 
Since the sequence $(0,{1\over n})_{n \in \N}$ consists of singular points 
for the exponential function, the Lie algebra $\g$ is not locally exponential 
(cf.\ Remark~II.5.8). 
\qed 

\Theorem VII.1.5. Every continuous homomorphism between locally 
finite-dimensional Lie groups is smooth. 
\qed

As the corresponding result for locally exponential Lie groups (Theorem~IV.1.18) did, 
the preceding theorem implies that 
locally finite-dimensional Lie groups form a full sub-category 
of topological groups. 
We even have the following stronger version of the preceding theorem: 

\Theorem VII.1.6. Let $G = \indlim G_n$ be a locally finite-dimensional 
Lie group and $H$ a Lie group. 
\litem{(a)} A group homomorphism $\phi \: G \to H$ is smooth if and only if the corresponding 
homomorphisms $\phi_n \: G_n \to H$ are smooth. 
\litem{(b)} If $H$ has a smooth exponential map, then each continuous homomorphism 
$\phi \: G \to H$ is smooth. 

\Proof. (a) is contained in [Gl05]. In view of (a), part (b) follows 
from the finite-dimensional case, which in turn follows from the existence of 
local coordinates of the second kind: 
$(t_1,\ldots, t_n)$ $\mapsto 
\prod_{i = 1}^n \exp_G(t_i x_i).$ 
\qed

\subheadline{VII.2. Direct limits of infinite-dimensional Lie groups} 

Direct limit constructions also play an important role when applied to sequences of 
infinite-dimensional Lie groups. On the level of Banach-, resp., Fr\'echet spaces, 
different types of directed systems lead to the important classes of 
LF spaces and Silva spaces  (cf.\ Definition I.1.2). 

If $M$ is a $\sigma$-compact finite-dimensional manifold 
and $K$ a Lie group, then the groups $C^\infty_c(M,K)$ of compactly supported 
smooth maps $M \to K$ are direct limits of the subgroups 
$C^\infty_X(M,K) := \{ f \in C^\infty(M,K) \: \supp(f) \subeq X\}$, 
which, for Banach--Lie groups  $K$, are Fr\'echet--Lie groups. 
On $C^\infty_c(M,K)$ this leads to the structure of an LF--Lie group 
if $K$ is Fr\'echet, but the construction of a Lie group structure works for general 
$K$ (Theorem~II.2.8). For $\dim K < \infty$, these groups are also discussed in 
[NRW94] as direct limit Lie groups which are BCH. 

Many interesting direct limits of mapping groups and other interesting 
classes embed naturally into certain direct sums, also called restricted direct products, 
often given by a nice atlas of a manifold. 
Therefore the following theorem turns out to be quite useful because it provides  
realizations as subgroups of a Lie group, and it usually is easier to verify that 
subgroups of Lie groups are Lie groups, than to construct the Lie group structures 
directly. 

\Theorem VII.2.1. {\rm([Gl03c])} If $(G_i)_{i \in I}$ is a family of locally exponential 
Lie groups, then their direct sum 
$$ G : =\bigoplus_{i \in I} G_i := \Big\{ (g_i)_{i \in I} \in \prod_{i \in I} G_i \: 
|\{i \: g_i \not=\1\}| < \infty\Big\} $$
carries a natural Lie group structure, where 
$\L(G) \cong \bigoplus_{i \in I} \L(G_i)$ carries the locally convex direct sum topology. 
\qed

\Theorem VII.2.2. {\rm([Gl06c])} For a $\sigma$-compact, non-compact manifold $M$ of 
positive dimension, 
the Lie group $\Diff_c(M)$ of compactly supported diffeomorphisms, 
endowed with the Lie group structure modeled on the LF space ${\cal V}_c(M)$  
is not a direct limit of the subgroups $\Diff_{M_n}(M)$, $(M_n)_{n \in \N}$ an exhaustion of $M$, 
in the category of smooth manifolds, but a homomorphism 
$\Diff_c(M) \to H$ to a Lie group $H$ is smooth if and only if it is smooth on 
each subgroup $\Diff_{M_n}(M)$. 
\qed

A crucial tool for the proof of the preceding theorem is the following lemma:

\Lemma VII.2.3. {\rm(Fragmentation Lemma)} Let $M$ be a $\sigma$-compact finite-dimensional 
manifold. Then there exists a locally finite cover $(K_n)_{n \in \N}$ of $M$ by compact 
sets, an open identity neighborhood $U \subeq \Diff_c(M)$ and a smooth mapping 
$\Phi \: U \to \bigoplus_{n \in \N} \Diff_{K_n}(M)$
which satisfies 
$\gamma = \Phi(\gamma)_1 \circ \ldots \circ \Phi(\gamma)_n$ 
for each $\gamma \in U$. 
\qed

\Theorem VII.2.4. {\rm([Gl06c])} For a $\sigma$-compact, non-compact manifold $M$ 
of positive dimension 
and a finite-dimensional Lie group $K$ of positive dimension, 
the Lie group $C^\infty_c(M,K)$ of compactly supported $K$-valued 
smooth functions, endowed with the Lie group structure modeled on the direct limit 
space $C^\infty_c(M,K) = \indlim C^\infty_{M_n}(M,K)$, $(M_n)_{n \in \N}$ an exhaustion of $M$, 
is not a direct limit of the subgroups $C^\infty_{M_n}(M,K)$ 
in the category of smooth manifolds, but a homomorphism 
$C^\infty_c(M,K) \to H$ to a Lie group $H$ is smooth if and only if it is smooth on 
each subgroup $C^\infty_{M_n}(M,K)$. 
\qed

\subheadline{Open Problems for Section VII} 

\Problem VII.1. Is every Lie group $G$ whose Lie algebra 
$\L(G)$ is countably dimensional, locally finite, and endowed with the finest locally convex 
topology regular? (cf.\ Problem II.2). 
\qed

\Problem VII.2. Are locally finite-dimensional Lie groups 
topological groups with Lie algebra? It is not clear that the compact open 
topology on ${\frak L}(G) \cong \Hom(\R,G)$ coincides with the given one on 
$\L(G)$ if the group is not locally exponential. 
\qed

\Problem VII.3. Does every subgroup $H$ of a locally finite-dimensional Lie group $G$ 
carry an initial Lie subgroup structure? (cf.\ (FP5))
\qed

\Problem VII.4. Let $M$ be a locally convex manifold and 
$\g \subeq {\cal V}(M)$ a countably dimensional locally finite-dimensional 
subalgebra consisting of complete vector fields. 
Does the inclusion $\g \to {\cal V}(M)$ integrate to a smooth action 
of a corresponding Lie group $G$ with $\L(G) = \g$? (cf.\ (FP7)) 

The first step should be to prove this for finite-dimensional Lie algebras $\g$, 
using local coordinates of the second kind and then to use that locally finite-dimensional 
Lie groups 
are direct limits in the category of smooth manifolds ([Gl05]). 
\qed

\Problem VII.5. The methods developed in [Gl03b] for the analysis of direct limit 
Lie groups seem to have potential to apply to more general classes of 
Lie groups $G$ which are direct limits of finite-dimensional manifolds $M_n$, $n \in \N$, 
with the property that for $n,m \in \N$ there exist $c(n,m)$ and $d(n)$ with 
$$ M_n \cdot M_m \subeq M_{c(n,m)} \quad \hbox{ and } \quad M_n^{-1} \subeq M_{d(n)}, $$
a situation which occurs in free constructions. 
Similar situations, with infinite-dimensional $M$, occur in the ind-variety 
description of Kac--Moody groups (cf.\ [Kum02], [BiPi02]). 
\qed

\sectionheadline{VIII. Linear Lie groups} 

\nin In this section, we take a closer look at linear Lie groups, i.e., 
Lie subgroups of CIAs. The  point of departure is that the unit 
group of a Mackey complete CIA $A$ is a BCH--Lie group 
(Theorem IV.1.11). This permits us to use the full machinery 
described Section IV for linear Lie groups. 

\Definition VIII.1. A {\it linear Lie group} is a Lie group  which can be realized 
as a locally exponential Lie subgroup of the unit group of some unital CIA. 
\qed

We collect some of the basic tools in the following theorem. 

\Theorem VIII.2. The following assertions hold: 
\litem{(1)} Linear Lie groups are BCH. 
\litem{(2)} Continuous homomorphisms of linear Lie groups are analytic.
\litem{(3)} If $G$ is a linear Lie group, then each closed Lie subalgebra 
$\h \subeq \L(G)$ integrates to a linear Lie group. 
\litem{(4)} For each morphism $\phi \: G \to H$ of linear Lie groups 
the kernel is a linear Lie group. 
\litem{(5)} For each $n \in \N$, the algebra $M_n(A)$ also is a CIA and 
$\GL_n(A) = M_n(A)^\times$ is a linear Lie group. 

\Proof. (1)-(4) follow from the fact that $A^\times$ is BCH 
(Theorem IV.1.11), Theorem IV.1.8, and the corresponding assertions on 
BCH--Lie groups in Section~IV. 

For (5), we refer to [Gl02b] (see also [Sw77]). 
\qed 

Linear Lie groups traditionally play an important role 
as groups of operators on Hilbert spaces, where they mostly occur 
as Banach--Lie subgroups (cf.\ [PS86], [Ne02b]). 
The connection between Lie theory and CIAs is more recent. 
The first systematic investigation of CIAs from a Lie theoretic 
perspective has been undertaken by {\smc Gl\"ockner} in [Gl02b]. 
Originally, complex CIAs came up in the 1950s as a natural 
class of locally convex associative algebras still permitting a powerful 
holomorphic functional calculus (cf.\ [Wae54a/b], [Al65]; see also [Hel93], 
and [Gram84] for Fr\'echet algebras of pseudo-differential operators). 

In K-theory, the condition on a topological ring $R$ that its unit group 
$R^\times$ is open and that the inversion map is continuous is quite natural 
because it is a crucial assumption for the analysis of idempotents in matrix 
algebras, resp., finitely generated projective modules, and the 
natural equivalence classes ([Swa62]; Section~VI.1). 

To get an impression of the variety of linear Lie groups, 
we describe some examples of CIAs: 

\Examples VIII.3. (a) Unital Banach algebras are CIAs. 

(b) If $M$ is a compact smooth manifold (with boundary) and $A$ is a CIA over 
$\K \in \{\R,\C\}$, 
then for each $r \in \N_0 \cup \{\infty\}$, the algebra 
$C^r(M,A)$ of $A$-valued $C^r$-functions on $M$ is a CIA ([Gl02b]). 
If $M$ is non-compact, but $\sigma$-compact, then 
$C^\infty_c(M,A)$, endowed with the direct limit topology of the 
subalgebras $C^\infty_X(M,A)$, is a non-unital CIA (Definition II.1.3(b)).  

(c) For $A = M_n(\C)$, the preceding construction leads in particular to the CIAs 
$$C^r(M,M_n(\C)) \cong M_n(C^r(M,\C)),$$ whose unit groups are the mapping groups 
$C^r(M,\GL_n(\C))$. 

(d) Let $X$ be a compact subset of $\C^n$ and 
$(U_n)_{n \in \N}$ a sequence of compact neighborhoods of $X$ with $\bigcap_n U_n = X$. 
In [Wae54b], {\smc Waelbroeck} shows that the algebra ${\cal O}(X,\C)$ of germs of holomorphic 
functions on $X$ is a CIA if it is endowed with the locally convex direct 
limit topology of the Banach algebras $C_{\cal O}(U_n,\C)$ of those continuous 
functions on $U_n$ which are holomorphic on the interior of $U_n$ (the {\it van Hove topology}). 
This defines on ${\cal O}(X,\C)$ the structure of a Silva space. 
The continuity of the multiplication and the completeness of this algebra 
is due to {\smc van Hove} ([vHo52a]). 

(e) If $A$ is a Banach algebra, $M$ a smooth manifold, 
$\alpha \: M \to \Aut(A)$ a map and $\alpha^a(m) := \alpha(m)(a)$, 
then the subspace 
$A^\infty := \{ a \in A \: \alpha^a \in C^\infty(M,A)\}$ 
is a CIA ([Gram84]). 
\qed

Part (d) of the preceding example shows in particular that for each compact subset $X \subeq \C^n$ 
the unit group 
${\cal O}(X,\C^\times)$ of the CIA ${\cal O}(X,\C)$ is a Lie group. 
In [Gl04b], {\smc Gl\"ockner} generalizes this Lie group construction as follows: 

\Theorem VIII.4. Let $X$ be a compact subset of a 
metrizable topological vector space, $\K \in \{\R,\C\}$ and $K$ a Banach--Lie group over 
$\K$. Then the group 
${\cal O}(X,K)$ of germs of $K$-valued analytic maps on open neighborhoods of $X$ 
is a $\K$-analytic BCH--Lie group. 
\qed

\Examples VIII.5. The following examples are Fr\'echet algebras with continuous inversion 
which are not CIAs because their unit groups are not open: 

(1) $A = C^\infty(M,\C)$, where $M$ is a non-compact 
$\sigma$-compact finite-dimensional manifold (cf.\ Remark~II.2.10). 

(2) $A = {\cal O}(M,\C)$, where $M$ is a complex submanifold of some $\C^n$, 
i.e., a Stein manifold. 

(3) $A = \R^\N$ with componentwise multiplication. 
\qed

In finite dimensions, a connected Lie group is called {\it linear} if 
it is isomorphic to a Lie subgroup of some $\GL_n(\R)$. Not all connected 
finite-dimensional Lie groups are linear. Typical examples of non-linear 
Lie groups are the universal covering $\tilde\SL_2(\R)$ of $\SL_2(\R)$ and 
the quotient $H/Z$, where 
$$ H = \pmatrix{ 1 & \R & \R \cr 0 & 1 & \R \cr 0 & 0 & 1\cr} $$
is the $3$-dimensional Heisenberg group (Example~V.3.5(c)) 
and $Z \subeq Z(H)$ is a non-trivial cyclic subgroup of its center ([Wie49]). 
It is a natural question whether the linearity condition on a connected finite-dimensional 
Lie group becomes weaker if we only require that it is a Lie subgroup of 
the unit group of some Banach algebra or even a CIA. 
According to the following theorem, this is not the case ([BelNe06]). 
Its Banach version is due to {\smc Luminet} and {\smc Valette} ([LV94]). 

\Theorem VIII.6. For a connected finite-dimensional Lie group $G$,  
the following are equivalent: 
\litem{(1)} The continuous homomorphisms $\eta \: G \to A^\times$ 
into the unit groups of  Mackey complete CIAs separate the points of $G$. 
\litem{(2)} $G$ is linear in the classical sense. 
\qed

\Remark VIII.7. Let us call a Banach--Lie algebra $\g$ 
{\it linear} if it has a faithful homomorphism into some Banach algebra~$A$. 

According to Ado's Theorem ([Ado36]), each finite-dimensional Lie algebra is linear, 
but the situation becomes more interesting, and also more complicated, 
for Banach--Lie algebras. 

In view of Corollary~IV.4.10, enlargeability is necessary for linearity, but it 
is not sufficient. In fact, if the Lie algebra $\g$ of a $1$-connected Banach--Lie group $G$ 
contains elements $p,q$ for which $[p,q]$ is a non-zero central element  
with $\exp_G([p,q]) = \1$, then $\g$ is not linear, because 
any morphism $\g \to A$ would lead to a linear representation of the 
quotient $H/Z$ of the $3$-dimensional Heisenberg group modulo a cyclic central subgroup $Z$. 
Such elements exist in the Lie algebra $\hat\g$ of the central 
extension of the Banach--Lie algebra $C^1(\SS^1, \su_2(\C))$ by $\R$ 
(Example~VI.1.16; [ES73]). 

Since for each Banach--Lie algebra $\g$ the quotient $\g_{\rm ad} = \g/\z(\g)$ 
is linear, the intersection $\n$ of all kernels of linear representations of $\g$ 
is a central ideal of $\g$. This links the linearity problem intimately with 
central extensions: When is a central extension of a linear Banach--Lie algebra 
linear? As the enlargeability is necessary, the discreteness of the corresponding 
period group is necessary (Theorem~VI.1.6), but what else?  

In [ES73], {\smc van Est} and {\smc \'Swierczkowski} 
describe a condition on the cohomology class of a central extension which 
is sufficient for linearity. They apply this in particular to show that, 
under some cohomological condition involving the center, 
for a Banach--Lie algebra $\g$, the Banach--Lie algebra $C^1_*([0,1],\g)$ 
of $C^1$-curves $\gamma$ in $\g$ with $\gamma(0) = 0$ is linear. It is remarkable 
that their argument does not work  for $C^0$-curves. 
Closely related to this  circle of ideas is {\smc van Est}'s proof of 
Ado's theorem, based on the vanishing of $\pi_2$ for each finite-dimensional 
Lie group ([Est66]). 

It is also interesting to note that for a real Banach--Lie algebra $\g$,  
linearity implies the linearity of the complexification $\g_\C$, which in turn implies that 
$\g_\C$ is enlargeable, which is crucial for the existence of universal 
complexifications of the corresponding groups (cf.\ Theorem~VI.1.23). 
In view of Corollary~IV.4.10, we thus have the implications 
$$ \g \quad \hbox{\rm linear} \quad 
\Rarrow \quad \g_\C \quad \hbox{\rm enlargeable} \quad \Rarrow \quad \g 
\quad \hbox{\rm enlargeable}. 
\qeddis

\Remark VIII.8. [GN07] (a) Let $A$ be a unital CIA and $n \in \N$. 
Further let $\sL_n(A) \trile \gl_n(A)$ denote the closed commutator algebra 
(cf.\ the end of Section VI.1). As this is a closed subalgebra, it generates some 
integral subgroup $S \to \GL_n(A)$ with $\L(S) = \sL_n(A)$. But in 
general $S$ will not be a Lie subgroup. This problem is caused by 
the fact that $\GL_n(A)$ need not be simply connected. 

Let $q \: \tilde\GL_n(A) \to \GL_n(A)_0$ denote the universal covering group 
of the identity component of $\GL_n(A)$. Then the Lie algebra morphism 
$$ \Tr \: \gl_n(A) \to A/\oline{[A,A]}, \quad (a_{ij}) \mapsto \Big[\sum_{i = 1}^n a_{ii}\Big] $$
satisfies $\ker \Tr = \sL_n(A)$. 
Let $HC_0(A)$ denote the completion of $A/\oline{[A,A]}$. Then $\Tr \: \gl_n(A) \to HC_0(A)$ 
integrates to a morphism of BCH--Lie groups 
$$ \tilde D \: \tilde\GL_n(A) \to HC_0(A), $$
and $\hat S := \ker D \trile \tilde\GL_n(A)$ is a BCH--Lie subgroup whose 
identity component $\hat S_0$ is a covering group of $S$. 
If the image of the induced period homomorphism 
$$ \per_{\Tr} \: \pi_1(\GL_n(A)) \to HC_0(A) \leqno(8.1) $$
is discrete, then $Z := HC_0(A)/\im(\per_{\Tr})$ is a Lie group 
and $D$ factors through a homomorphism 
$D \: \GL_n(A)_0 \to Z$, which can be considered as a generalization of the determinant. 
Now $\ker D$ is a BCH--Lie subgroup of $\GL_n(A)$ with  Lie algebra $\sL_n(A)$, 
which implies that $(\ker D)_0 = S$. It is interesting to compare this 
situation with the one in Remark V.2.14(c), where the group of Hamiltonian 
diffeomorphisms of a symplectic manifolds plays a similar role. 

If $A$ is commutative, then the determinant $\det \: \GL_n(A) \to A^\times$ 
is a morphism of Lie groups and $\SL_n(A) \trile \GL_n(A)$ is a 
normal BCH--Lie subgroup with Lie algebra $\sL_n(A)$. 

Since the period maps (8.1) are compatible for different $n$, they lead to a 
homomorphism 
$$ \per_A^0 \: K_2(A) = \indlim \pi_1(\GL_n(A)) \to HC_0(A) $$
(cf.\ the end of Section VI.1). If $A$ is complex, we may compose with the 
Bott isomorphism $\beta_A^0 \: K_0(A) \to K_2(A)$ to get a natural transformation 
$$ T_A := \per_A^0 \circ \beta_A^0 \: K_0(A) \to HC_0(A), $$
which is unique and therefore given by $T_A([e]) = \Tr(e)$. It follows that 
the image of $\per_A^0$ is discrete if and only if the image of the trace map 
$$\Tr \: \bigcup_{n = 1}^\infty \Idem(M_n(A)) \to HC_0(A)$$
generates a discrete subgroup. 

If $A$ is commutative, then $HC_0(A) = A$, and the image of the trace map 
lies in the discrete subgroup ${1\over 2\pi i} \ker(\exp_{A_\C})$ of $A_\C$. 
Hence the image of the trace map is discrete for each commutative CIA. 
\qed

\Remark VIII.9. The set $\Idem(A)$ of idempotents of a CIA plays a central 
role in (topological) $K$-theory. In [Gram84], {\smc Gramsch} shows that this set always 
carries a natural manifold structure, which implies in particular that its connected 
components are open subsets. The key point is to use rational methods to obtain charts 
on this set. 

In a similar spirit, it is explained in [BerN04/05] how Jordan methods can be used 
in an infinite-dimensional context to obtain manifold structures 
on geometrically defined manifolds generalizing symmetric spaces and Gra\3mann manifolds. 
\qed

\subheadline{Open Problems for Section VIII} 

\Problem VIII.1. Show that the completion of a CIA $A$ is again a CIA or give a 
counterexample. 
\qed

\Problem VIII.2. Characterize those Banach--Lie algebras which are linear in 
the sense that they have an injective homomorphism into some Banach algebra  
(Remark~VIII.7). 

Not much seems to be known about this problem, which is partly related to the 
non-existence of Lie's theorem on the representation of solvable Banach--Lie algebras. 
In view of this connection, the class of Banach--Lie algebras of the form 
$\g = E \rtimes_D \R$, where $D$ is a continuous linear operator on the Banach space 
$E$ should be a good testing ground. 

As each real Banach space $E$ is isomorphic to ${\cal L}(\R,E)$, which can be embedded as a 
Banach--Lie algebra into the Banach algebra ${\cal L}(E \oplus \R)$, each abelian 
Banach--Lie algebra is linear. What about nilpotent ones? 
\qed

\Problem VIII.3. Show that each linear Lie group is regular. We know that this is 
the case for unit groups of CIAs. If, in addition, $A$ is $\mu$-regular in 
the sense of Definition~III.2.4, then Theorem~III.2.10 implies the 
$\mu$-regularity of each Lie subgroup. 
\qed

\Problem VIII.4. Is the tensor product $A \otimes B$ of two CIAs, endowed 
with the projective tensor topology a CIA? 
This is true for $B = M_n(\K)$, $n \in \N$, where $A \otimes B 
\cong M_n(A)$. Is it also true if $B$ 
is the algebra of rapidly decreasing matrices or the direct limit 
algebra $M_\infty(\K) := \indlim M_n(\K)$? 
\qed

\Problem VIII.5. Let $\g$ be a locally convex Lie algebra. 
Does the enveloping algebra $U(\g)$ carry a natural  
topology for which the multiplication is continuous and the natural 
map $\g \to U(\g)$ is continuous? 

More generally, let $E$ be a locally convex space and endow its 
tensor algebra ${\cal T}(E) = \bigoplus_{n = 0}^\infty E^{\otimes n}$ 
with the locally convex direct limit topology, where the subspaces 
$E^{\otimes n}$ carry the 
projective tensor topology. Is the multiplication on ${\cal T}(E)$ continuous? 
\qed

\Problem VIII.6. (a) Does every locally convex Lie algebra $\g$ 
have a faithful topological module~$E$? If Problem VIII.5 has a positive solution, 
then we may simply take $E := U(\g)$. 

(b) Does every nilpotent locally convex Lie algebra have a 
faithful nilpotent topological module? Is this true in the Banach category? 
\qed

\sectionheadline{IX. Lie transformation groups} 

\nin One of the fundamental references on topological transformation groups 
is the monograph [MZ55] by {\smc Montgomery} and {\smc Zippin}. 
Since it also deals with differentiability properties of transformation groups on manifolds, 
some of the techniques described there have interesting applications 
in the context of infinite-dimensional Lie theory. 

\subheadline{IX.1. Smooth Lie group actions} 

\Theorem IX.1.1. {\rm([BoMo45], [MZ55, p.212])} Any continuous action $\sigma\: G \times M \to M$ 
of a finite-dimensional Lie group on a finite-dimensional 
smooth manifold $M$ by diffeomorphisms 
is smooth. 
\qed

For compact manifolds we obtain the following ``automatic smoothness'' result 
on homomorphisms of Lie groups (see also [CM70] for one-parameter groups; and 
[Gl02d] for the non-compact case). 

\Corollary IX.1.2. If $M$ is a $\sigma$-compact finite-dimensional 
manifold and $G$ a finite-di\-men\-sio\-nal Lie group, 
then any continuous homomorphism $\phi \: G \to \Diff_c(M)$ is smooth. 
\qed 

The following result provides a positive answer to (FP9) for diffeomorphism groups. 

\Theorem IX.1.3. {\rm([MZ55, Th.~5.2.2, p.~208])} If 
a locally compact group $G$ acts faithfully on a smooth finite-dimensional manifold $M$ 
by diffeomorphisms, then $G$ is a finite-dimensional Lie group. 
If $M$ is compact, then each locally compact subgroup of 
$\Diff(M)$ is a Lie group. 
\qed

The preceding results take care of the actions of locally compact groups 
on manifolds. As the work of {\smc de la Harpe} and {\smc Omori} ([OdH71/72]) 
shows, the situation for Banach--Lie groups is more subtle: 

\Theorem IX.1.4. {\rm([OdH72])} Let $G$ be a Banach--Lie group. 
If $\L(G)$ has no proper finite-codimensional closed ideals, then 
$\L(G)$ has no proper finite-codimensional closed subalgebra and each smooth action  
of $G$ on a finite-dimensional manifold is trivial. 
\qed

If $\alpha \: \g \to {\cal V}(M)$ is an injective map, then for each 
$p \in M$ the subspace 
$$\g_p := \{ x \in \g \: \alpha(x)(p) =0\}$$ 
is a finite-codimensional subalgebra 
with $\bigcap_p \g_p = \{0\}$. Therefore the existence of many finite-codimensional 
subalgebras is necessary for Lie algebras to be realizable by vector fields on 
a finite-dimensional manifold. 

\Theorem IX.1.5. {\rm([OdH72])} If a Banach--Lie group $G$ acts smoothly, effectively, 
amply (for each $m \in M$ the evaluation map $\g \to T_m(M)$ is surjective),  
and primitively (it leaves no foliation invariant) on a 
finite-dimensional manifold $M$, then it is finite-dimensional. 
\qed

\Theorem IX.1.6. {\rm([Omo78, Th.~B/C])} 
Let $G$ be a connected Banach--Lie group acting smoothly, effectively 
and transitively on a finite-dimensional manifold $M$. 
\litem{(1)} If $M$ is compact, then $G$  is finite-dimensional. 
\litem{(2)} If $M$ is non-compact, then $\L(G)$ contains a finite-codimensional 
closed solvable ideal. 
\qed

Since $\Diff(M)$ acts smoothly, effectively and transitively on $M$, this implies: 

\Corollary IX.1.7. $\Diff(M)$ cannot be 
given a Banach-Lie group structure for which the natural action on~$M$ is smooth. 
\qed

In Section 4 of [OdH72], {\smc Omori} and {\smc de la Harpe} construct 
an example of a Banach--Lie group $G$ acting smoothly and amply, but not primitively 
on $\R^2$. 

\msk 

The preceding discussion implies in particular that Banach--Lie groups 
rarely act on finite-dimensional manifolds. As the gauge groups of 
principal bundles $q \: P \to M$ over compact manifolds $M$ show, the 
situation is different for locally exponential Lie groups (cf.\ Theorem~IV.1.12). 
Therefore it is of some interest to have good criteria for 
the integrability of infinitesimal actions of locally exponential 
Lie algebras on finite-dimensional manifolds (cf.\ (FP7)).  

We start with a more general setup for infinite-dimensional manifolds which 
need extra smoothness assumptions: 

\Theorem IX.1.8. {\rm(Integration of locally exponential Lie algebras of vector fields; [AbNe06])} Let 
$M$ be a smooth manifold modeled on a locally convex space, 
$\g$ a locally exponential Lie algebra and 
$\alpha \: \g \to {\cal V}(M)$ a homomorphism of Lie algebras whose 
range consists of complete vector fields. 
Suppose further that the map 
$$ \Exp \: \g \to \Diff(M), \quad x \mapsto \Fl^{\alpha(x)}_1 $$
is smooth in the sense of {\rm Definition II.3.1} and that 
$0$ is isolated in $\z(\g) \cap \Exp^{-1}(\id_M)$. Then there exists a 
locally exponential Lie group $G$ and a smooth action 
$\sigma \: G \times M \to M$ whose derived action 
$\dot\sigma \: \g \to {\cal V}(M)$ coincides with $\alpha$. 
\qed

In the finite-dimensional case, the smoothness assumptions in Theorem~IX.1.8 follows from the 
smooth dependence of solutions of ODEs on parameters and initial values, and 
the condition on the exponential function can be verified with methods 
to be found in [MZ55]. This leads to the following less technical 
generalization of the Lie--Palais Theorem which subsumes in particular 
{\smc Omori}'s corresponding results for Banach--Lie algebras ([Omo80, Th.~A], 
[Pe95b, Th.4.4]).

\Theorem IX.1.9. Let $M$ be a smooth finite-dimensional 
manifold, $\g$ a locally exponential Lie algebra and 
$\alpha \: \g \to {\cal V}(M)$ a continuous homomorphism of Lie algebras 
whose range consists of complete vector fields. Then there exists a 
locally exponential Lie group $G$ and a smooth action 
$\sigma \: G \times M \to M$ with $\dot\sigma = \alpha$. 
\qed

The following result is a generalization of Palais' Theorem in another direction. 
Since $\Diff(M)$ is $\mu$-regular (Theorem~III.3.1), it also follows from Theorem~III.2.8. 

\Theorem IX.1.10. {\rm([Les68])} If $M$ is compact, then 
a subalgebra $\g \subeq {\cal V}(M)$ is integrable to an integral subgroup 
if $\g$ is finite-dimensional or closed and finite-codimensional. 
\qed

\subheadline{IX.2. Groups of diffeomorphisms as automorphism groups} 

In this subsection, we simply collect some results stating that 
automorphism groups of certain algebra, Lie algebra or groups associated to 
geometric structure on manifolds are what one expects. Most of the results 
formulated below for automorphisms of structures attached to a manifold $M$ 
generalize to results saying that if $M_1$ and $M_2$ are two manifolds and 
two objects of the same kind attached to $M_1$ and $M_2$ are isomorphic, 
then this isomorphism can be implemented by a diffeomorphism $M_1 \to M_2$, compatible 
with the geometric structures under consideration.

\Theorem IX.2.1. Let $M$ be a $\sigma$-compact finite-dimensional 
smooth manifold. Then the following assertions hold: 
\litem{(1)} For the Fr\'echet algebra $C^\infty(M,\R)$, each homomorphism 
to $\R$ is a point evaluation. 
\litem{(2)} $\Aut(C^\infty(M,\R)) \cong \Diff(M)$. 
\litem{(3)} $\Aut({\cal V}_c(M)) \cong \Aut({\cal V}(M)) \cong \Diff(M)$. 
\litem{(4)} If, in addition, $M$ is complex and ${\cal V}^{(1,0)}(M)\subeq {\cal V}(M)_\C$ 
is the Lie algebra of complex vector fields of type $(1,0)$, then 
$\Aut({\cal V}^{(1,0)}(M)) \cong \Aut_{\cal O}(M)$ is the group of biholomorphic 
automorphisms of $M$. 
\litem{(5)} For each finite-dimensional $\sigma$-compact manifold $M$ and each 
simple (real or complex) finite-dimensional Lie algebra $\k$, the natural homomorphism 
$$ C^\infty(M,\Aut(\k)) \rtimes \Diff(M) \to \Aut(C^\infty(M,\k)) $$
is surjective. 
\litem{(6)} If $M$ is a Stein manifold and $\k$ is a finite-dimensional 
complex simple Lie algebra, then 
$\Aut({\cal O}(M,\k)) \cong {\cal O}(M,\Aut(\k)) \rtimes \Aut_{\cal O}(M)$, where 
$\Aut_{{\cal O}}(M)$ denotes the group of biholomorphic diffeomorphisms of $M$. 
\litem{(7)} If $K \subeq \C^n$ is a polyhedral domain and 
${\cal O}(K,\C)$ the algebra of germs of holomorphic 
$\C$-valued functions in $K$, then the  
group $\Aut({\cal O}(K,\C))$ consists of the germs of biholomorphic maps 
of $K$ and $\der({\cal O}(K,\C))$ consists of the germs holomorphic vector fields on $K$. 

\Proof. (1) (cf.\ [My54] for the compact case; [Pu52]; [Co94]). 

(2) follows easily from (1) because each automorphism of the algebra $C^\infty(M,\R)$ 
acts on $\Hom(C^\infty(M,\R),\R) \cong M$. 

(3) The representability of each isomorphism of ${\cal V}_c(M)$ by a diffeomorphism 
is due to Pursell and Shanks ([PuSh54]), and the other assertion 
follows from Theorem 2 in [Ame75]. It is based on the fact that 
the maximal proper subalgebras of finite codimension are 
all of the form ${\cal V}(M)_m := \{ X \in {\cal V}(M) \: X(m) = 0\}$ 
for some $m \in M$, hence permuted by each automorphism; resp.\ the fact that 
all maximal ideals consist of all vector fields whose jet vanishes in some 
$m \in M$. 

(4) follows from Theorem 1 in [Ame75]. 

(5) [PS86, Prop.~3.4.2]. A central point is that every non-zero endomorphism of 
$\k$ is an automorphism. Further, it is used that 
$[\k, C^\infty(M,\k)] = C^\infty(M,\k)]$ and that distributions supported 
by one point are of finite order. 

(6) [NeWa06b]. 

(7) This is [vHo52b, Th.~III], where it is first shown that the maximal ideals 
in the Silva CIA ${\cal O}(K,\C)$ (Example VIII.3(d)) are the kernels of the point evaluations ([vHo52b, Th.~I]). 
\qed

\Remark IX.2.2. Let $K \subeq \C^n$ be a compact subset and 
$\Aut_{{\cal O}}(K)$ the group of germs of bihomolorphic maps, defined on some neighborhood 
of $K$, mapping $K$ onto itself. 
In [vHo52a], {\smc van Hove} introduces a group topology on 
this group as the topology for which the map 
$$ \Aut_{{\cal O}}(K) \to {\cal O}(K,\C^n) \times {\cal O}(K,\C^n), \quad 
g \mapsto (g,g^{-1}) $$
is an embedding. He shows that, under certain geometric conditions
 on the set $K$, this group is complete and contains no small subgroups. 
Moreover, its natural action on ${\cal O}(K,\C)$ is continuous. 
\qed

\Theorem IX.2.3. {\rm([Omo74], \S 10])} Let $M$ be a $\sigma$-compact finite-dimensional 
smooth manifold. For a differential form $\alpha$ on $M$ we put 
${\cal V}(M,\alpha) := \{ X \in {\cal V}(M) \: {\cal L}_X\alpha = 0\}.$
Then the following assertions hold: 
\litem{(1)} If $\mu$ is a volume form or a symplectic form on $M$, 
then every (algebraic) automorphism of 
${\cal V}(M,\mu)$ is induced by an element of the group 
$$\{ \phi \in \Diff(M) \: \phi^*\mu \in \R \mu\}.$$
\litem{(2)} If $\alpha$ is a contact $1$-form on $M$, then every (algebraic) automorphism of 
${\cal V}(M,\alpha)$ is induced by an element of the group 
$\{ \phi \in \Diff(M) \: \phi^*\alpha \in C^\infty(M,\R^\times)\cdot \alpha\}.$
\qed

In [Omo80], one finds another interesting result of this type. 
Let $V$ be a germ of an affine variety in $0 \in \C^n$. Two such germs $V$ and $V'$ 
are said to 
be {\it biholomorphically equivalent} if there exists an element 
$\phi \in \Gh_n(\C)$ of the group of germs of biholomorphic maps fixing $0$ 
(as in Example VI.2.12), such that 
$\phi(V) = \phi(V')$. On the infinitesimal level the automorphisms of a germ 
$V$ are given by the Lie algebra 
$$ \g(V) := \{ X \in \gh_n(\C) \: X.J(V) \subeq J(V)\}, $$
where $J(V) \subeq {\cal O}(0,\C)$ (the germs of holomorphic functions in $0$) 
is the annihilator ideal of $V$. Let $\g(V)_k \trile \g(V)$ denote the ideal 
consisting of all vector fields vanishing of order $k$ in $0$ and form the 
projective limit Lie algebra 
$$ \oline\g(V) := \prolim \g(V)/\g(V)_k, $$
which can be viewed as a Fr\'echet completion of $\g(V)$. 

An element $X \in \gh_n(\C)$ is called {\it semi-expansive} 
if it is $\Gh_n(\C)$-conjugate to a linear diagonalizable vector field for which all 
eigenvalues lie in some open halfplane. 
The germ $V$ is called an {\it expansive singularity} if 
$\g(V)$ contains an expansive vector field. 

\Theorem IX.2.4. Two expansive singularities $V$ and $V'$ are biholomorphically equivalent 
if and only if the pro-finite 
Lie algebras $\oline\g(V)$ and $\oline\g(V')$ are isomorphic. 
Moreover, $\Aut(\oline\g(V))$ can be identified with the stabilizer $\Gh_n(\C)_V$ of 
$V$ in the group $\Gh_n(\C)$. 
\qed

On the group level, we have the following analog of Theorem~IX.2.3 (cf.\ [Fil82] for (1) 
and [Ban97, Thms.~7.1.4/5/6] for (2)-(4)): 

\Theorem IX.2.5. Let $M$ be a $\sigma$-compact connected finite-dimensional 
smooth manifold. Then the following assertions hold: 
\litem{(1)} Every (algebraic) automorphism of $\Diff(M)$ is inner. 
\litem{(2)} If $\alpha$ is a contact $1$-form on $M$, then every (algebraic) automorphism of 
$\Diff(M,\alpha)$ is conjugation with an element of the group 
$\{ \phi \in \Diff(M) \: \phi^*\alpha \in C^\infty(M,\R^\times)\cdot \alpha\}.$
\litem{(3)} If $\omega$ is a symplectic form and $M$ is compact of dimension $\geq 2$, 
then every (algebraic) automorphism of $\Diff(M,\omega)$ is conjugation by an element of the group 
$$\{ \phi \in \Diff(M) \: \phi^*\omega \in \R \omega\}.$$
\litem{(4)} If $\mu$ is a volume form and $M$ is of dimension $\geq 2$, 
then every (algebraic) automorphism of $\Diff(M,\mu)$ is conjugation by an element of the group 
$\{ \phi \in \Diff(M) \: \phi^*\mu \in \R \mu\}.$
\qed

\subheadline{Open Problems for Section IX} 

\Problem IX.1. Let ${\cal V}(M)_{\rm cp}$ denote the set of complete 
vector fields on the finite-dimensional manifold $M$ (Remark~II.3.8). 
Then we have an exponential 
function 
$$ \Exp \: {\cal V}(M)_{\rm cp} \to \Diff(M), \quad X \mapsto \Fl^X_1. $$
Is it true that $0$ is isolated in $\Exp^{-1}(\id_M)$ with respect to the natural 
Fr\'echet topology on ${\cal V}(M)$ (cf.\ Definition~I.5.2)? 

That this is true for compact manifolds follows from  
Newman's Theorem ([Dr69, Th.~2]). 
For the proof of Theorem IX.1.9, we show 
for each continuous homomorphism $\alpha \: \g \to {\cal V}(M)$ of a 
locally exponential Lie algebra $\g$ to ${\cal V}(M)$ with 
range in ${\cal V}(M)_{\rm cp}$ that $0$ is isolated in 
$(\Exp \circ \alpha)^{-1}(\id_M)$, which is a weaker statement. 

Since the set $\Exp^{-1}(\id_M)$ is in one-to-one correspondence with the 
smooth $\T$-actions on $M$, the problem is to show that the trivial 
action is isolated in this ``space'' of all smooth $\T$-actions on $M$. 

If $M$ is the real Hilbert space $\ell^2(\N,\R)$ with the Hilbert basis $e_n$, $n \in \N$, 
then we have linear vector fields $X_n(v) := 2\pi i \la v, e_n \ra e_n$ with 
$\exp(X_n) = \id_M$ and $X_n \to 0$ uniformly on compact subsets of $E$. Hence 
the finite-dimensionality of $M$ is crucial. 
\qed

\Problem IX.2. (Banach symmetric spaces) 
Let $M$ be a smooth manifold. We say that $(M,\mu)$
is a {\it symmetric space} (in the sense of {\smc Loos}) (cf.\ [Lo69]) if 
$\mu \: M \times M \to M, (x,y) \mapsto x \cdot y$ 
is a smooth map with the following properties: 
\litem{(S1)} $x \cdot x$ for all $x \in M$. 
\litem{(S2)} $x \cdot (x \cdot y) =y$ for all $x,y \in M$. 
\litem{(S3)} $x \cdot (y \cdot z) = (x \cdot y) \cdot (x \cdot z)$ 
for all $x,y \in M$. 
\litem{(S4)} $T_x(\mu_x) = -\id_{T_x(M)}$ for $\mu_x(y) := \mu(x,y)$ and each $x \in M$. 

{\rm(a)} Is it true that the automorphism group $\Aut(M,\mu)$ of a Banach symmetric space 
$(M,\mu)$ is a 
Banach--Lie group? (cf.\ [Ne02c], [La99]) 

{\rm(b)} The tangent spaces $T_x(M)$ of a symmetric space carry natural structures of 
Lie triple systems. Develop a Lie theory for locally exponential, resp., 
Banach--Lie triple systems, including 
criteria for the integrability of morphisms and enlargeability (cf.\ Sections IV and VI). 
\qed

\Problem IX.3. A fundamental problem in the theory of Banach transformation 
groups is that we do not know if orbits carry natural manifold structures. As in finite dimensions, 
the main point is to find good criteria for a closed subgroup $H$ of a Banach--Lie group $G$ 
to ensure that the coset space $G/H$ has a natural manifold structure for which 
the action of $G$ on $G/H$ is smooth and the quotient map 
$q \: G \to G/H$ is a ``weak'' submersion in the sense that all its differentials are linear 
quotient maps. 
In view of Remark IV.4.13, this is true if 
\litem{(1)} $H$ is a split submanifold (the same proof as in 
finite dimensions works), 
\litem{(2)} $H$ is a normal Banach--Lie subgroup 
(without any splitting requirements) (Corollary IV.3.6), and 
\litem{(3)} $G$ is a Hilbert--Lie group, which implies the splitting condition (1). 

Here are some concrete problems: 

\litem{(a)} Suppose that $G/H$ is a smooth manifold with submersive $q$ and a smooth action 
of $G$. Does this imply that $H$ is a Lie subgroup of $G$? 

\litem{(b)} Are the stabilizer groups $G_m$ for a smooth action of a Banach--Lie group $G$ on a 
Banach manifold $M$ Lie subgroups? For linear actions this follows from 
Proposition IV.3.4 and Corollary IV.3.13. 

\litem{(c)} Characterize those Lie subgroups $H$ for which $G/H$ is a smooth manifold. 

\litem{(d)} Let $H \subeq G$ be a closed subgroup and $\h := \L^e(H)$ its Lie algebra. 
Then the normalizer $N_G(\h)$ of $\h$ is a Lie subgroup (Proposition IV.3.4,  Corollary IV.3.12). 
Is it true that $\Ad(G).\h \cong G/N_G(\h)$ carries a natural manifold structure? 

Note that, if $H$ is connected, it is a normal subgroup of $N_G(\h)$. If $H$ is a Lie subgroup, 
this implies that $N_G(\h)/H$ carries a Lie group structure and therefore a manifold structure.
\qed

\Problem IX.4. Show that for each compact subset $K \subeq \C^n$ the group 
$\Aut_{{\cal O}}(K)$ from Remark~IX.2.2 is a Lie group with respect to the manifold structure 
inherited from the embedding into ${\cal O}(K,\C^n)$ (cf.\ Remark IX.2.2).
\qed

\Problem IX.5. (Automorphisms of gauge algebras) 
Let $q \: P \to M$ be a smooth $K$-principal bundle over the (compact) manifold $M$. 
Determine the group $\Aut(\gau(P))$ of automorphisms of the gauge Lie algebra. 
Does it coincide with the automorphism group $\Aut(\ad(P))$ of the adjoint bundle,  
whose space of sections $\gau(P)$ is? If $K$ is a simple complex Lie group, 
then the results in [Lec80, Th.~16] provide a local description of the 
automorphisms of this Lie algebra in terms of diffeomorphisms of $M$ and 
sections of the automorphism bundle $\Aut(\ad(P))$. 
\qed

\Problem IX.6. Determine the automorphism groups of the Lie algebras 
$\gf_n(\K)$, $\gs_n(\K)$ and $\gh_n(\C)$. 

\Problem IX.7. Describe all connected Banach--Lie groups acting 
smoothly, effectively and transitively on a finite-dimensional manifold. 
In view of Theorem~IX.1.6, for each Banach--Lie group $G$, 
the Lie algebra $\L(G)$ contains a finite-codimensional closed solvable ideal. 
If, conversely, $\g$ is a Banach--Lie algebra with a finite-dimensional 
closed solvable ideal, then Theorem VI.1.19 implies that $\g$ is enlargeable. 
Under which conditions do the corresponding groups $G$ act effectively on 
some finite-dimensional homogeneous space? (see also the corresponding 
discussion in [Omo97]). 
\qed

\Problem IX.8. Let $G$ be a Banach--Lie group and $H \subeq G$ a closed subgroup for which 
$\L^e(H)$ has finite codimension. Does this imply that $G/H$ is a manifold? 
\qed

\sectionheadline{X. Projective limits of Lie groups}  

\nin Projective limits play an important role in several branches of 
Lie theory. Since complete locally convex spaces are nothing but closed subspaces 
of products of Banach spaces, on the level of the model spaces, the projective 
limit construction leads us from Banach spaces to the locally convex setting. 
On the group level, the situation is more involved because, although projective limits 
of Lie groups are often well-behaved topological groups, in general they are not Lie 
groups. In this section, we briefly report on some aspects of projective limit Lie theory 
and the recent theory of {\smc Hofmann} and {\smc Morris} of projective limits 
of finite-dimensional Lie groups. 

\subheadline{X.1. Projective limits of finite-dimensional Lie groups} 

In their recent monograph [HoMo06], {\smc Hofmann} and {\smc Morris} 
approach projective limits of finite-dimensional Lie groups, 
so-called {\it pro-Lie groups}, from a topological 
point of view. We refer to [HoMo06] for details on the results mentioned 
in this subsection. 

Clearly, arbitrary products of finite-dimensional Lie groups, such as 
$$ \R^J, \quad \Z^J, \quad \SL_2(\R)^J $$
for an arbitrary set $J$, are pro-Lie groups. 
The following theorem gives an abstract characterization of pro-Lie groups:

\Theorem X.1.1. A topological group $G$ is a pro--Lie group if and only if 
it is isomorphic to a closed subgroup of a product of finite-dimensional Lie 
groups. In particular, closed subgroups of pro-Lie groups are pro-Lie groups. 
\qed

A crucial observation is that the class of topological groups with 
Lie algebra (cf.\ Definition IV.1.23) is closed under projective limits and that 
$$ \fL(\prolim G_j) \cong \prolim \fL(G_j) $$
as locally convex Lie algebras. 
Let us call topological vector spaces isomorphic to $\R^J$ for some set $J$ 
{\it weakly complete}. These are the dual spaces of the vector spaces 
$\R^{(J)}$, endowed with the weak-$*$-topology. This provides a duality between 
real vector spaces and weakly complete locally convex spaces, which implies 
in particular that each closed subspace of a weakly complete space 
is weakly complete and complemented. In particular, weakly complete spaces are 
nothing but the projective limits of finite-dimensional vector spaces. 
These considerations lead to: 

\Theorem X.1.2. Every pro-Lie group has a Lie algebra which is a 
a projective limit of finite-dimensional Lie algebras, hence a weakly complete 
topological Lie algebra. The image of the exponential function 
$$ \exp_G \: \fL(G) \to G, \quad \gamma \mapsto \gamma(1) $$
generates a dense subgroup of the identity component $G_0$. 
\qed

In the following, we call projective limits of finite-dimensional 
Lie algebras {\it pro-finite Lie algebras} (called {\it pro-Lie algebras} in [HoMo06]). 

\Remark X.1.3. According to Yamabe's Theorem ([MZ55]), 
each locally compact group $G$ for which $G/G_0$ is compact 
is a pro-Lie group. Since the totally disconnected locally compact group 
$G/G_0$ contains an open compact subgroup, 
each locally compact group $G$ contains an open subgroup 
with a Lie algebra, hence is a topological 
group with a Lie algebra (cf.\ [Las57], [HoMo05, Prop.~3.5]). 
\qed

In view of Theorem X.1.1, the category of pro-Lie groups is closed under products 
and projective limits, which are remarkable closedness properties which in 
turn lead to the existence of an adjoint functor $\Gamma$ for the Lie functor 
$\fL$: 

\Theorem X.1.4. {\rm(Lie's Third Theorem for Pro-Lie Groups; [HoMo05])} 
The Lie functor $\fL$ from the category of pro-Lie groups 
to the category of pro-Lie algebras has a left adjoint $\Gamma$. It associates 
with each pro-finite Lie algebra $\g$ a connected pro-Lie group $\Gamma(\g)$ 
and a natural isomorphism $\eta_\g \: \g \to \fL(\Gamma(\g))$, such that for 
every morphism $\phi \: \g \to \fL(G)$ there exists a unique morphism 
$\phi' \: \Gamma(\g) \to G$ with 
$\fL(\phi') \circ \eta_\g = \phi$. 
\qed

The first part of the following structure theorem can be found in [HoMo06]. 
The second part follows from the fact that finite-dimensional tori are the only 
abelian connected compact Lie groups. 
\Theorem X.1.5. A connected abelian pro-Lie group is of the form 
$\R^J \times C$ for a compact connected abelian group $C$. 
It is a Lie group if and only if $C$ is a finite-dimensional torus. 
\qed

It is quite remarkable that the category of pro-finite Lie algebras 
permits to develop a structure theory which is almost as strong as in finite dimensions. 
In particular, there is a Levi decomposition. To describe it, we call a 
pro-finite Lie algebra $\g$ {\it pro-solvable} if it is a projective limit 
of finite-dimensional solvable Lie algebras: 

\Theorem X.1.6. {\rm(Levi decomposition of pro-finite Lie algebras and groups; [HoMo06])} 
Each pro-finite Lie algebra $\g$ contains a unique maximal pro-solvable ideal $\r$. 
There is a family $(\s_j)_{j \in J}$ of finite-dimensional simple Lie algebras 
such that $\s := \prod_{j \in J} \s_j$ satisfies 
$$ \g \cong \r \rtimes \s. $$

For the corresponding pro-finite Lie group $\Gamma(\g)$ we then have 
$$ \Gamma(\g) \cong R \rtimes S, \quad \hbox{ where } \quad 
S \cong \prod_{j \in J} S_j, $$
where $S_j$ is a $1$-connected Lie group with Lie algebra $\s_j$ and 
$R$ is diffeomorphic to $N \times \R^K$ for some set $K$ and some simply connected 
pro-nilpotent Lie group $N \cong (\L(N),*)$. 
\qed 

More concretely, the closed commutator algebra $\n := \oline{[\r,\r]}$ is pro-nilpotent, 
because all images of this subalgebra in finite-dimensional solvable quotients of $\r$ 
are nilpotent. If $\e \subeq \r$ is a closed vector space complement of $\n$ in $\r$ 
([HoMo06, 4.20/21]), 
then the map 
$$ \Phi \: \n \times \e \mapsto R, \quad (x,y) \mapsto \exp_R(x)\exp_R(y)\leqno(10.1.1)  $$
is a homeomorphism ([HoMo06, Th.~8.13]). 
The point of view of [HoMo06] is purely topological, so that infinite-dimensional 
Lie group structures are not discussed. We note that (10.1.1) can be viewed as a chart 
of the topological group $R$, and it is not hard to see that it defines on $R$ 
the structure of a smooth Lie group. 

Based on the preceding theorem, one can characterize those pro-finite Lie algebras 
which are integrable to locally convex Lie groups ([HoNe06]): 

\Theorem X.1.7. For a pro-finite Lie algebra $\g$ the following are equivalent: 
\litem{(1)} $\g$ is the Lie algebra of a locally convex Lie group $G$ 
with smooth exponential function. 
\litem{(2)} $\g$ is the Lie algebra of a regular locally convex Lie group $G$. 
\litem{(3)} $\g$ has a Levi decomposition 
$\g \cong \r \rtimes \s$, where only finitely many factors in 
$\s = \prod_{j \in J} \s_j$ are not isomorphic to $\sL_2(\R)$. 
\litem{(4)} The group $\Gamma(\g)$ is locally contractible. 
\litem{(5)} The group $\Gamma(\g)$ carries the structure of a regular 
Lie group, compatible with its topology. 

\Proof. (Sketch) Let $R$ denote the $1$-connected group with Lie algebra $\r$ 
constructed from the chart (10.1.1) and $S_j$ be the $1$-connected 
Lie group with Lie algebra $\s_j$. If $\s_j$ is isomorphic to $\sL_2(\R)$, 
then $S_j \cong \tilde \SL_2(\R)$ is diffeomorphic to $\R^3$, so that 
$S := \prod_{j\in J} S_j$ carries a natural manifold structure, turning it into a Lie group. 
One verifies that $S$ acts smoothly on $R$, so that $G := R \rtimes S$ is a Lie 
group with Lie algebra $\g$. 

For the converse, let $G$ be a Lie group with Lie algebra $\g$ and a smooth exponential 
function $\exp_G \: \g \to G$; put $J_0 := \{ j \in J \: \s_j \not\cong \sL_2(\R)\}.$ 
For each $j \in J_0$, we then have morphisms  
$\alpha_j \: S_j \to G$, $\beta_j \: G \to S_j$ with $\beta_j \circ \alpha_j = \id_{S_j}$. 
If $J_0$ is infinite, this contradicts the local contractibility of $G$. 
\qed

\Remark X.1.8. To describe all connected regular Lie groups $G$ with a pro-finite Lie 
algebra~$\g$, we have to describe the discrete central subgroups of the $1$-connected ones, 
which are isomorphic, as topological groups, to some $G := \Gamma(\g)$ 
(cf.\ [HoMo06]). In all these groups, the exponential function restricts to an 
isomorphism $\z(\g) \to Z(G)_0$ of topological groups, so that 
$Z(G)_0 \cong \R^X$ for some set $X$. Based on the information 
provided in the preceding theorem, it is shown in [HoMo06] that a subgroup 
$\Gamma \subeq Z(\Gamma(\g))_0$ is discrete if and only if it is finitely 
generated and its intersection with $Z(G)_0$ is discrete. 
This characterization provides a quite good description of all discrete 
subgroups of $Z(G)$, hence all non-simply connected regular Lie groups 
with Lie algebra $\g$. 

If $S \cong \prod_{j \in J} S_j$ and 
$J$ is infinite, then infinitely many factors $S_j$ are isomorphic to 
$\tilde\SL_2(\R)$, whose center is isomorphic to $\Z$. The subgroup of index $2$ 
acts trivially in each finite-dimensional representation, which leads to 
$Z(G) \cap S_j \cong \Z$. Hence $Z(G) \cap S$ contains a subgroup isomorphic 
to $\Z^{J \setminus J_0}$, which implies in particular that the adjoint group of $\g$ 
is {\sl not} a Lie group. 
\qed

We have already seen that all pro-nilpotent Lie algebras are exponential, 
which applies in particular to all pro-nilpotent pro-finite Lie algebras. 
The following theorem provides a characterization of the locally exponential 
pro-Lie algebras ([HoNe06]): 

\Theorem X.1.9. For a pro-Lie algebra $\g$, the following are equivalent: 
\litem{(1)} $\g$ is locally exponential. 
\litem{(2)} There exists a $0$-neighborhood 
$U \subeq \g$, consisting of $\exp$-regular points, i.e., 
$\kappa_\g(x)$ is invertible for each $x \in U$. 
\litem{(3)} $\Gamma(\g)$ is a locally exponential topological group. 

If these conditions are satisfied, then $\g$ contains a closed ideal of finite codimension 
which is exponential. In particular, $\g$ is virtually pro-solvable, 
i.e., $\g = \r \rtimes \s$ with a finite-dimensional semisimple Lie algebra $\s$,  
and $\g$ is enlargeable. 
\qed

Recall that we have seen in Example II.5.9(a) an example of a pro-finite Lie algebra 
$\g \cong \R^\N \rtimes_D \R$ which is not locally exponential. Since this Lie algebra 
has an abelian closed hyperplane, the existence of an exponential hyperplane 
ideal is not sufficient for local exponentiality. 

\subheadline{X.2. Projective limits of infinite-dimensional Lie groups} 

As we have seen in the preceding Subsection X.1, 
the extent to which the structure theory of finite-dimensional Lie groups can be 
carried forward to projective limits is quite surprising. There are also natural 
classes of topological groups which are natural projective limits 
of infinite-dimensional Lie groups. Therefore it would be of some interest to develop 
a systematic ``pro--Lie theory'' for such groups. 

One of the most natural classes of such groups are the mapping groups. 
Let $M$ be a $\sigma$-compact finite-dimensional smooth manifold, 
$r \in \N_0 \cup \{\infty\}$, $K$ a Lie group, and 
$G := C^r(M,K)$, endowed with the compact open $C^r$-topology 
(Definition~II.2.7). 

Then there exists a sequence $(M_n)_{n \in \N}$ of compact subsets of $M$ 
which is an exhaustion, in the sense that $M_n \subeq \Int(M_{n+1})$. 
Using the usual Morse theoretic arguments, we 
may assume that the subsets $M_n$ are compact manifolds with boundary. Then 
each compact subset of $M$ is contained in some $M_n$, which implies that 
$$ G = C^r(M,K) \cong \prolim C^r(M_n,K) $$
is a projective limit, where the projection maps are given by restriction. In view of 
Theorem II.2.8, the groups $C^r(M_n,K)$ are Lie groups, so that the topological group 
$C^r(M,K)$ is a projective limit of Lie groups. 

If $K$ is locally exponential, then each $C^r(M_n,K)$ inherits this property, 
so that $C^r(M_n,K)$ is a topological group with Lie algebra 
(Definition IV.1.23), and this implies that $C^r(M,K)$ also is a topological 
group with Lie algebra, where 
$$ {\frak L}(C^r(M,K)) 
\cong  \prolim {\frak L}(C^r(M_n,K))
=  \prolim \L(C^r(M_n,K))
=  \prolim C^r(M_n,\L(K)) = C^r(M,\L(K)). $$
As in the case of compact manifolds $M$, the exponential function of $C^r(M,K)$ is given by 
$$ \exp \: C^r(M,\L(K)) \to C^r(M,K), \quad \xi \mapsto \exp_K \circ \xi. $$

\msk 

If $M$ is a compact complex manifold and $K$ is a linear complex Lie group, 
then all holomorphic functions $M \to K$ are constant. Therefore the groups 
${\cal O}(M,K)$ are trivial in this case. If $M$ is non-compact and $K$ is a complex 
Lie group, we use the compact open topology to turn ${\cal O}(M,K)$ into a 
topological group and observe that we have a topological embedding 
${\cal O}(M,K) \into C^r(M,K)$ for each $r \in \N_0 \cup \{\infty\}$. 
Those cases for which we know to have honest Lie group structures on 
${\cal O}(M,K)$ are quite limited (cf.\ Theorem III.1.9), but it 
seems that projective limit theory is also a useful tool to study these groups 
of holomorphic maps. Let $(M_n)_{n \in \N}$, as above, be an exhaustion of $M$ by 
compact submanifolds with boundary. 
Then the groups ${\cal O}(M_n,K)$, defined appropriately, carry Lie group 
structures, for which ${\cal O}(M_n,\L(K))$ is the corresponding Lie algebra (cf.\ [Wo05b]), 
so that 
$${\cal O}(M,K) \cong\prolim {\cal O}(M_n,K) $$
is a projective limit of Lie groups.

It would be of considerable interest to find a good categorical framework 
for such classes of projective limits of Lie groups. Of particular relevance 
would be to understand the ``right class'' of central extensions of the groups 
${\cal O}(M,K)$ and $C^r(M,K)$ in the same spirit as for the groups 
$C^\infty_c(M,K)$ of compactly supported maps (cf.\ [Ne04c]). 

\subheadline{Open Problems for Section X} 

\Problem X.1. Are strong ILB--Lie groups, resp., $\mu$-regular Lie groups, 
topological groups with Lie algebra? What about diffeomorphism groups? 
Does it suffice that the Lie group $G$ has a smooth exponential function (cf.\ Problem~VII.2)?  
\qed

\Problem X.2. Let $\g$ be a pro-finite Lie algebra and $\n \trile \g$ a closed 
exponential ideal of finite codimension. Characterize the local exponentiality 
of $\g$ in terms of the spectra of the operators $\ad_\n x := \ad x\res_\n$ 
(cf.\ Proposition~X.1.9). Are all locally exponential pro-finite Lie algebras BCH? 
\qed

\Problem X.3. For the description of the non-simply connected Lie groups among the 
projective limits of finite-dimensional Lie groups, it is important to understand the 
structure of the center of the simply connected ones. 
Let $G$ be a such a $1$-connected group and $Z(G)$ its center. 
The Lie algebra of $Z(G)$ is $\z(\g)$, which lies in the pro-solvable radical, 
so that $Z(G)_0 \cong \z(\g)$. On the other hand, we have seen in 
Remark~X.1.8 that $Z(G)$ may contain non-discrete subgroups isomorphic to 
$\Z^\N$. Is it possible to determine the structure of $Z(G)$ as a topological group   
(see [HoMo06] for more details)? 
\qed

\Problem X.4. Determine the automorphism groups of pro-finite Lie algebras. Under 
which conditions are they Lie groups? An interesting situation where the automorphism 
group of a pro-finite Lie algebra is a closed subgroup of a Lie groups is described 
in Theorem~IX.2.4. 
\qed

\def\entries{

\[AbNe06 Abouqateb, A., and K.-H. Neeb, {\it Integration of locally exponential 
Lie algebras of vector fields}, in preparation 

\[ARS84 Adams, M., Ratiu, T., and R. Schmid, 
{\it The Lie group structure of diffeomorphism groups and 
invertible Fourier integral operators, with applications}, 
in ``Infinite-dimensional groups with applications'' 
(Berkeley, Calif., 1984), Math. Sci. Res. Inst. Publ. {\bf 4}, 
Springer, New York, 1985, 1--69

\[ARS86a ---, {\it A Lie group structure for pseudodifferential operators}, 
Math. Ann. {\bf 273:4} (1986), 529--551

\[ARS86b ---, {\it A Lie group structure for Fourier integral operators}, 
Math. Ann. {\bf 276:1} (1986), 19--41

\[Ado36 Ado, I., {\it \"Uber die Darstellung von Lieschen Gruppen durch lineare Substitutionen}, 
Bull. Soc. Phys.-Math. Kazan, III. Ser. {\bf 7} (1936), 3--43 

\[AHMTT93 Albeverio, S. A., H\o{}egh-Krohn, R. J., Marion, J. A., Testard, D. H., 
and \break B.~S. Torr\'esani, ``Noncommutative distributions. Unitary representation 
of \break Gauge Groups and Algebras,'' Monographs and Textbooks in Pure and Applied 
Mathematics {\bf 175}, Marcel Dekker, Inc., New York, 1993

\[Al65 Allan, G.R., {\it A spectral theory for locally convex algebras}, 
Proc. London Math. Soc. (3) {\bf 15} (1965), 399--421

\[AABGP97 Allison, B. N., Azam, S., Berman, S., Gao, Y., and
A. Pianzola, ``Extended Affine Lie Algebras and Their Root Systems,'' 
Memoirs of the Amer. Math. Soc. {\bf 603}, Providence R.I., 1997 

\[ABG00 Allison, B., Benkart, G., and  Y. Gao, {\it Central 
extensions of Lie algebras gra\-ded by finite-root
systems}, Math.\ Ann.\ {\bf 316} (2000), 499--527 

\[Ame75 Amemiya, I., {\it Lie algebra of vector fields and complex structure}, 
J. Math. Soc. Japan {\bf 27:4} (1975), 545--549 

\[Arn66 Arnold, V. I., {\it Sur la g\'eom\'etrie diff\'erentielle des groupes 
de Lie de dimension infinie et ses applications \`a l'hydrodynamique des 
fluides parfaits}, Ann. Inst. Fourier {\bf 16:1} (1966), 319--361 

\[AK98 Arnold, V. I., and B. A. Khesin, ``Topological Methods in Hydrodynamics,'' 
Springer-Verlag, 1998 

\[AI95 de Azcarraga, J. A., and J. M. Izquierdo, ``Lie Groups, Lie
Algebras, Cohomology and some Applications in Physics,'' Cambridge
Monographs on Math. Physics, 1995 

\[Bak01 Baker, H. F., {\it On the exponential theorem for a simply 
transitive continuous group, and the calculation of the finite equations from 
the constants of structure}, 
J. London Math. Soc. {\bf 34} (1901), 91--127 

\[Bak05 ---, {\it On the calculation of the finite equations of a 
continuous group}, London M. S. Proc. {\bf 35} (1903), 332--333

\[Ban97 Banyaga, A., ``The Structure of Classical Diffeomorphism
Groups,'' Kluwer Academic Publishers, 1997

\[Bas64 Bastiani, A., {\it Applications diff\'{e}rentiables
et vari\'{e}t\'{e}s diff\'{e}rentiables de dimension infinie},
J. Anal.\ Math.\ {\bf 13} (1964), 1--114

\[Beg87 Beggs, E. J., {\it The de Rham complex on infinite dimensional manifolds}, 
Quart. J. Math. Oxford (2) {\bf 38} (1987), 131--154 

\[Bel04 Belti\c{t}\u{a}, D., {\it Asymptotic 
products and enlargibility of Banach-Lie
algebras},
J. Lie Theory {\bf 14} (2004), 215--226

\[Bel06 ---, ``Smooth Homogeneous Structures 
in Operator Theory,''  Chapman and Hall, CRC Monographs and 
Surveys in Pure and Applied mathematics, 2006 

\[BelNe06 Belti\c t\u a, D., and Neeb, K.-H., 
{\it Finite-dimensional Lie subalgebras of algebras with continuous inversion}, 
Preprint, 2006 

\[BR05a Belti\c t\u a, D., and T. S. Ratiu, {\it 
Geometric representation theory for unitary groups of operator algebras}, 
Advances in Math., to appear 

\[BR05b ---, {\it Symplectic leaves in real Banach Lie-Poisson spaces}, 
Geom. Funct. Anal. {\bf 15:4} (2005), 753--779

\[BerN04 Bertram, W., and Neeb, K.-H., {\it Projective completions of Jordan pairs, Part I. 
The generalized projective geometry of a Lie algebra}, 
J. of Algebra {\bf 277:2} (2004), 474--519

\[BerN05 ---, {\it Projective completions of Jordan pairs, Part II}, 
Geom. Dedicata {\bf 112:1} (2005), 75--115 

\[BGN04 Bertram, W., H. Gl\"ockner, and K.-H. Neeb, {\it 
Differential Calculus over General Base Fields and Rings}, 
Expositiones Math. {\bf 22} (2004), 213--282

\[BiY03 Billig, Y., {\it Abelian extensions of the group of
  diffeomorphisms of a torus}, Lett. Math. Phys.  {\bf 64:2}  (2003),  155--169 

\[BiPi02 Billig, Y., and A. Pianzola, {\it Free Kac-Moody groups and their Lie algebras}, 
Algebr. Represent. Theory {\bf 5:2} (2002), 115--136

\[Bir36 Birkhoff, G., {\it Continuous groups and linear spaces}, 
Mat. Sbornik {\bf 1} (1936), 635--642 

\[Bir38 ---, {\it Analytic groups}, 
Trans. Amer. Math. Soc. {\bf 43} (1938), 61--101 

\[Bl98 Blackadar, B., ``K-theory for Operator Algebras,'' 2nd edition, 
Cambridge Univ. Press, 1998

\[BoSi71 Bochnak, J., and J. Siciak, {\it Analytic functions in topological vector 
spaces}, Studia Math. {\bf 39} (1971), 77--112 

\[BoMo45 Bochner, S., and D. Montgomery, {\it Groups 
of differentiable and real or complex analytic transformations}, 
Ann. of Math. (2) {\bf 46} (1945), 685--694 

\[BD73  Bonsall, F. F., and J. Duncan, ``Complete Normed Algebras,'' 
Ergeb. Mathem. und ihrer Grenzgebiete {\bf 80}, Springer-Verlag, New York, Heidelberg, 1973

\[BCR81 Boseck, H., Czichowski, G., and K.-P. Rudolph,
``Analysis on Topological \break Groups -- General Lie Theory,''
Teubner, Leipzig, 1981

\[Bos90 Bost, J.-B., {\it Principe d'Oka, $K$-theorie et syst\`emes dynamiques 
non-commu\-ta\-tifs}, 
Invent. Math. {\bf 101} (1990), 261--333 

\[Bo77 Bott, R., {\it On the characteristic classes of groups of diffeomorphisms}, 
Enseignement Math. (2) {\bf 23:3-4} (1977), 209--220 

\[Bou87 Bourbaki, N., ``Topological Vector Spaces, Chaps. 1--5,'' Springer-Verlag, Berlin, 1987 

\[Bou89 ---, ``Lie Groups and Lie Algebras, Chapter 1--3,'' Springer-Verlag, Berlin, 1989

\[Bre93 Bredon, G.\ E., ``Topology and Geometry,'' Graduate Texts in
Mathematics {\bf 139}, Springer-Verlag, Berlin, 1993 

\[Bry93 Brylinski, J.-L., ``Loop Spaces, Characteristic Classes and
Geometric Quantization,'' Progr. in Math. {\bf 107}, Birkh\"auser
Verlag, 1993 

\[BuGi02 Burgos Gil, J. I., ``The Regulators of Beilinson and Borel,'' 
CRM Monograph {\bf 15}, Amer. Math. Soc., 2002 

\[Cal70  Calabi, E., {\it On the group of automorphisms of a symplectic manifold}, 
in ``Problems in Analysis,´´ (Lectures at the Sympos. in honor of Salomon Bochner; 1--26, 
Princeton Univ. Press, Princeton, N.J., 1970

\[Cam97 Campbell, J. E., {\it On a law of combination of operators bearing on the theory 
of continuous transformation groups}, Proc. of the London Math. 
Soc. {\bf 28} (1897), 381--390 

\[Cam98 ---, {\it On a law of combination of operators. (Second paper)}, 
Proc. of the London Math. 
Soc. {\bf 28} (1897), 381--390 

\[CaE98 Cartan, E., {\it Les groupes bilin\'eaires et les syst\`emes 
de nombres complexes}, Ann. Fac. Sci. Toulouse Sci. Math. Sci. Phys. 
{\bf 12:1} (1898), B1--B64

\[CaE01 ---, {\it L'int\'egration des syst\`emes d'\'equations 
aux différentielles totales}, Ann. Sci. École Norm. Sup. (3) {\bf 18} 
(1901), 241--311

\[CaE04 ---, {\it Sur la structure des groups infinies des transformations}, 
Ann. Sci. Ecol. Norm. Sup. {\bf 21} (1904), 153--206; 
{\bf 22} (1905), 219--308 

%

\[CaE30 ---, {\it Le troisi\`eme th\'eor\`eme fondamental de Lie}, 
C. R. Acad. Sci. Paris {\bf 190} (1930), 914-916, 1005-1007 

\[CaE36 ---, {\it La topologie des groupes de Lie. (Expos\'es de g\'eom\'etrie Nr. 8.)}, 
Actual. sci. industr. {\bf 358} (1936), 28 p 


\[CaE52 ---, {\it La topologie des espaces repr\'esentifs de groupes de
Lie}, Oeuvres I, Gauthier--Villars, Paris, {\bf 2} (1952), 1307--1330 

\[CVLL98 Cassinelli, G., E. de Vito, P. Lahti and A. Levrero, {\it Symmetries of the quantum 
state space and group representations}, Reviews in Math. Physics {\bf 10:7} 
(1998), 893--924 

\[CM70 Chernoff, P., and J. Marsden, 
{\it On continuity and smoothness of group actions}, Bull. Amer. Math. Soc. {\bf 76} (1970), 
1044--1049

\[Ch46 Chevalley, C., ``Theory of Lie Groups I,'' Princeton Univ.\ Press, 1946

\[ChE48 Chevalley, C. and S. Eilenberg, {\it Cohomology theory of Lie groups and Lie 
algebras}, Transactions of the Amer. Math. Soc. {\bf 63} (1948), 85--124 
 
\[Co94 Connes, A., ``Non-commutative  Geometry,'' Academic Press, 1994 

\[CGM90 Cuenca Mira, J.\ A., A.\ Garcia Martin, and C.\ Martin
Gonzalez, {\it Structure theory of $L^*$-algebras}, Math.\ Proc.\
Camb.\ Phil.\ Soc.\ {\bf 107} (1990), 361--365

\[DP03 Dai, J., and D. Pickrell, {\it The orbit method and the 
Virasoro extension of \break $\Diff_+(\SS^1)$. I. Orbital integrals}, 
J. Geom. Phys. {\bf 44:4} (2003), 623--653 

\[Da94 Dazord, P., {\it Lie groups and algebras in infinite dimension: a new approach}, 
in ``Symplectic Geometry and Quantization''; 17--44, 
Contemp. Math. {\bf 179}, Amer. Math. Soc., Providence, RI, 1994

\[De32 Delsartes, J., ``Les groups de transformations lin\'eaires dans l'espace 
de Hilbert,'' M\'emoirs des Sciences Math\'ematiques, Fasc. {\bf 57}, Paris



\[DiPe99 Dimitrov, I., and I.\ Penkov, {\it Weight modules of direct
limit Lie algebras}, International Math. Res. Notices {\bf 5} (1999), 
223--249  

\[DoIg85 Donato, P., and P. Iglesias, {\it Examples de groupes diff\'eologiques: 
flots irrationnels sur le tore}, C. R. Acad. Sci. Paris, ser. 1  {\bf 301} (1985), 
127--130 


\[DL66 Douady, A., and M. Lazard, {\it Espaces fibr\'es en alg\`ebres
de Lie et en groupes}, Invent. math. {\bf 1} (1966),  133--151 

\[Dr69 Dress, A., {\it Newman's Theorem 
on transformation groups}, Topology {\bf 8} (1969), 203--207 

\[Dy47 Dynkin, E. B., {\it Calculation
of the coefficients in the Campbell--Hausdorff
formula} (Russian), Doklady Akad.\ Nauk.\
SSSR (N.S.) {\bf 57} (1947), 323--326 

\[Dy53 ---, ``Normed Lie Algebras and Analytic Groups,'' 
Amer. Math. Soc. Translation 1953, no. 97, 66 pp 

\[Eb70 Ebin, D. G., {\it The manifold of riemannian metrics}, 
in ``Global Analysis'' (Berkeley, Calif., 1968), Amer. Math. Soc. 
Proc. Symp. Pure Math. {\bf 15} (1970), 11--40 

\[EM69 Ebin, D. G., and J. E. Marsden, {\it Groups of diffeomorphisms and 
the solution of the classical Euler equations for a perfect fluid}, 
Bull. Amer. Math. Soc. {\bf 75} (1969), 962--967

\[EM70 ---, {\it Groups of diffeomorphisms and the motion of 
an incompressible fluid}, Ann. of Math. {\bf 92} (1970), 102--163 

\[EMi99 Ebin, D. G., and G. Misiolek, {\it The exponential map on 
${\cal D}^s_\mu$}, in ``The Arnoldfest (Toronto, ON, 1997),'' 
Fields Inst. Commun. {\bf 24},
Amer. Math. Soc., Providence, RI, 1999, 153--163

\[Ee58 Eells, J. Jr., {\it On the geometry of function spaces}, in 
``International Symposium on Algebraic Topology,''  pp. 303--308; 
Universidad Nacional Autonoma de M\'exico and UNESCO, Mexico City

\[Ee66 ---, {\it A setting for global analysis}, Bull.\ Amer.\
Math.\ Soc. {\bf 72} (1966), 751--807 

\[ES96 Eichhorn, J., and R. Schmid, {\it Form 
preserving diffeomorphisms on open manifolds}, 
Ann. Global Anal. Geom. {\bf 14:2} (1996), 147--176

\[ES01 ---, {\it Lie groups of Fourier integral operators on open manifolds}, 
Comm. Anal. Geom.  {\bf 9:5} (2001),  983--1040

\[Ei68 Eichler, M., {\it A new proof of the Baker--Campbell--Hausdorff formula}, 
J. Math. Soc. Japan {\bf 20} (1968), 23--25 





\[Est62 van Est, W. T., {\it Local and global groups}, 
Proc. Kon. Ned. Akad. v. Wet. Series A 65, Indag. math. {\bf 24} (1962), 391--425

\[Est66 ---, {\it On Ado's theorem}, 
Proc. Kon. Ned. Akad. v. Wet. Series A 69, Indag. Math. {\bf 28} (1966), 176--191 

\[Est84 ---, {\it Rapport sur les S-atlas}, Ast\'erisque {\bf 116} (1984), 235--292 

\[Est88 ---, {\it Une d\'emonstration de \'E.\ Cartan du
troisi\`eme th\'eor\`eme de Lie}, 
in P.\ Dazord et al eds., ``Seminaire Sud-Rhodanien de Geometrie
VIII: Actions Hamiltoniennes de Groupes; Troisi\`eme Th\'eor\`eme de
Lie,'' Hermann, Paris, 1988 

\[EK64 van Est, W.\ T., and Th.\ J.\ Korthagen, {\it Non enlargible Lie
algebras}, 
Proc.\ Kon.\ Ned. Acad. v. Wet. Series A, Indag. Math. {\bf 26} (1964), 15--31

\[ES73 van Est, W. T., and S. \'Swierczkowski, {\it The path functor and 
faithful representability of Banach Lie algebras}, 
in ``Collection of articles dedicated to the memory of Hannare Neumann, I.''
J. Austral. Math. Soc. {\bf 16} (1973), 54--69

\[EF94 Etinghof, P. I., and I. B. Frenkel, {\it Central extensions of current 
groups in two dimensions}, Commun. Math. Phys. {\bf 165} (1994), 429--444 

\[Fil82 Filipkiewicz, R. P., {\it Isomorphisms between diffeomorphism groups}, 
Ergodic Theory Dynamical Systems {\bf 2:2} (1983), 159--171 

\[Fl71 Floret, K., {\it Lokalkonvexe Sequenzen mit kompakten Abbildungen}, 
J. Reine \break Angew. Math. {\bf 247} (1971) 155--195 

\[Fre68 Freifeld, Ch., {\it One-parameter subgroups do not fill a 
neighborhood of the identity in an infinite-dimensional Lie (pseudo-) group}, 
1968, Battelle Rencontres. 1967 Lectures in Mathematics and Physics; 538--543, 
Benjamin, New York

\[FB66 Fr\"olicher, A., and W. Bucher, ``Calculus in Vector Spaces without Norm,'' 
Lecture Notes in Math. {\bf 30}, Springer-Verlag, Berlin, 1966 

\[FK88 Fr\"olicher, A., and A.\ Kriegl, ``Linear Spaces and
Differentiation Theory,'' J.~Wiley, Interscience, 1988 

\[Fu86 Fuks, D. B. ``Cohomology of Infinite-Dimensional Lie Algebras,'' 
Consultants Bureau, New York, London, 1986 

\[Ga96 Galanis, G., {\it Projective limits of Banach-Lie groups}, 
Period. Math. Hungar. {\bf 32:3} (1996), 179--191

\[Ga97 ---, {\it On a type of linear differential equations in Fr\'echet spaces}, 
Ann. Scuola Norm. Sup. Pisa Cl. Sci. (4) {\bf 24:3} (1997), 501--510

\[GG61 Glashow, S. L., and Gell-Mann, M., {\it Gauge theories of vector particles}, 
Ann. Phys. {\bf 15}  (1961), 437--460 

\[Gl02a Gl\"ockner, H., {\it Infinite-dimensional Lie groups without completeness 
restrictions}, in  ``Geometry and Analysis on Finite and Infinite-dimensional 
Lie Groups,'' 
A. Strasburger, W. Wojtynski, J.\ Hilgert and K.-H. Neeb (Eds.), 
Banach Center Publications {\bf 55} (2002), 43--59 

\[Gl02b ---, {\it Algebras whose groups of units are Lie groups}, 
Studia Math. {\bf 153} (2002), 147--177 

\[Gl02c ---,  {\it Lie group structures on quotient groups and universal 
complexifications for infinite-dimensional Lie groups}, 
J. Funct. Anal. {\bf 194} (2002), 347--409 

\[Gl02d ---, {\it Patched locally convex spaces, almost local mappings,
and diffeomorphism groups of non-compact manifolds}, Manuscript, TU
Darmstadt, 26.6.02 

\[Gl03a ---, {\it Implicit functions from topological vector spaces to 
Banach spaces}, Israel J. Math., to appear; math.GM/0303320 

\[Gl03b ---, {\it Direct limit Lie groups and manifolds}, J. Math. Kyoto Univ. {\bf 43} (2003), 1--26

\[Gl03c ---, {\it Lie groups of measurable mappings}, Canadian J. Math. {\bf 55:5} (2003), 969--999

\[Gl04a ---, {\it Tensor products in the category of topological vector spaces are 
not associative}, Comment. Math. Univ. Carol. {\bf 45:4} (2004), 607--614

\[Gl04b ---, {\it Lie groups of germs of analytic mappings}, 
pp.\ 1-16 in: Turaev, V. and T. Wurzbacher (Eds.), ``Infinite Dimensional Groups and 
Manifolds,'' IRMA Lecture Notes in Mathematics and Theoretical Physics, de Gruyter, 2004

\[Gl05 ---, {\it Fundamentals of direct limit Lie theory}, 
Compositio Math. {\bf 141} (2005), 1551--1577


\[Gl06a ---, {\it Discontinuous non-linear mappings on locally convex direct 
limits}, Publ. Math. Debrecen {\bf 68} (2006), 1--13 

\[Gl06b ---, {\it Fundamental problems in the theory of infinite-dimensional Lie groups}, J. Geom. Symm. Phys., to appear 

\[Gl06c ---, {\it Direct limits of infinite-dimensional Lie groups compared to direct 
limits in related categories}, in preparation 

\[Gl06d ---, {\it Direct limit groups do not have small subgroups}, 
Preprint;  math.GR/ \break 0602407 

\[GN03 Gl\"ockner, H., and K.-H. Neeb, {\it Banach--Lie quotients, enlargibility, 
and universal complexifications}, J. Reine Angew. Math. {\bf 560} (2003), 1--28 

\[GN06 ---, ``Infinite-dimensional Lie groups, Vol. I, Basic Theory and Main Examples,'' 
book in preparation 

\[GN07 ---, ``Infinite-dimensional Lie groups, Vol. II, Geometry and Topology,'' 
book in preparation 

\[Go04 Goldin, G. A., {\it Lectures on 
diffeomorphism groups in quantum physics}, in 
``Contemporary Problems in Mathematical Physics,'' Proc.\ of the third internat.\ 
workshop (Cotonue, 2003), 2004; 3--93

\[GW84 Goodman, R., and N. R. Wallach, {\it Structure and unitary
cocycle representations of loop groups and the group of
diffeomorphisms of the circle}, J. Reine Ang. Math. {\bf 347} (1984), 69--133

\[GW85 ---, {\it Projective unitary positive energy representations of
$\Diff(\SS^1)$}, J. Funct. Anal. {\bf 63} (1985), 299--312 

\[Go69 Goto, M., {\it On an arcwise 
connected subgroup of a Lie group}, Procedings of the Amer. 
Math. Soc. {\bf 20} (1969), 157--162  

\[Grab88 Grabowski, J., {\it Free subgroups of diffeomorphism groups}, 
Fund. Math. {\bf 131:2} (1988), 103--121

\[Grab93 ---, {\it Derivative of the exponential mapping for
infinite-dimensional Lie groups}, Ann. Global Anal. Geom. {\bf 11:3} (1993), 
213--220 

\[GVF01 Gracia-Bondia, J. M., Vasilly, J. C., and H. Figueroa,
``Elements of Non-com\-mu\-ta\-tive Geometry,'' Birkh\"auser Advanced Texts,
Birkh\"auser Verlag, Basel, 2001 

\[Gram84 Gramsch, B., {\it Relative Inversion in der
St\"{o}rungstheorie von Operatoren und $\Psi$-Algebren},
Math.\ Ann.\ {\bf 269} (1984), 22--71

\[GrNe06 Grundling, H., and K. - H. Neeb, {\it 
Lie group extensions associated to modules of continuous inverse algebras}, in preparation

\[Gu77 Gutknecht, J., {\it Die $C^\infty_\Gamma$-Struktur auf der 
Diffeomorphismengruppe einer kompakten Mannigfaltigkeit}, 
Ph.D. Thesis, Eidgen\"ossische Technische Hochschule Z\"urich, Diss. No. {\bf 5879}, 
Juris Druck + Verlag, Zurich, 1977 

\[Ham82 Hamilton, R., {\it The inverse function theorem of Nash and
  Moser}, Bull. Amer. Math. Soc. {\bf 7} (1982), 65--222

\[dlH72 de la Harpe, P., ``Classical Banach--Lie Algebras and
Banach--Lie Groups of Operators in Hilbert Space,'' Lecture Notes in
Math.\ {\bf 285}, Springer-Verlag, Berlin, 1972 

\[HK77 Harris, L. A., and W. Kaup, {\it Linear algebraic groups in infinite dimensions}, 
Illinois J. Math. {\bf 21} (1977), 666--674 

\[Hau06 Hausdorff, F., {\it Die symbolische Exponentialformel in der 
Gruppentheorie}, Leip\-zi\-ger Berichte {\bf 58} (1906), 19--48 

\[HS68 Hausner, M., and J. T. Schwartz, ``Lie Groups; Lie Algebras,''
Gordon and Breach, New York, London, Paris, 1968 

\[HeMa02 Hector, G., and E. Mac\'{\i}as-Virg\'{o}s, {\it Diffeological groups},
Research Exp.\ Math.\ {\bf 25} (2002), 247--260

\[Hel93 Helemskii, A.Ya., ``Banach and Locally Convex Algebras,'' 
Oxford Science Publications,  Oxford University Press, New York, 1993 

\[Hi99 Hiltunen, S., {\it Implicit functions from locally convex
spaces to Banach spaces}, Studia Math. {\bf 134:3} (1999), 235--250 

\[Ho51 Hochschild, G., {\it Group extensions of Lie groups I, II}, 
Annals of Math. {\bf 54:1} (1951), 96--109 and 
 {\bf 54:3} (1951), 537--551 

\[Ho65 ---, ``The Structure of Lie Groups,'' Holden Day, San 
Francisco, 1965

\[Hof68 Hofmann, K. H., ``Introduction to the Theory of Compact Groups. Part I,'' 
Dept. Math. Tulane Univ., New Orleans, LA, 1968 

\[Hof72 ---, {\it Die Formel von Campbell, Hausdorff und Dynkin 
und die Definition Liescher Gruppen}, in ``Theory 
of Sets and Topology'' (in honour of Felix Hausdorff, 1868--1942); 
 251--264. VEB Deutsch, Verlag Wissensch., Berlin, 1972 

\[Hof75 ---, {\it Analytic groups without analysis}, Symposia Mathematica {\bf 16} 
(Convegno sui Gruppi Topologici e Gruppi di Lie, INDAM, Rome, 1974); 357--374, 
Academic Press, London, 1975 

\[HoMo98 Hofmann, K.\ H., and S.\ A.\ Morris, ``The Structure of
Compact Groups,'' Studies in Math., de Gruyter, Berlin, 1998 

\[HoMo05 ---, {\it Sophus Lie's third fundamental theorem and the adjoint functor theorem}, 
J. Group Theory {\bf 8} (2005), 115--123 

\[HoMo06 ---, ``{\it The Lie Theory of Connected Pro-Lie Groups--A 
Structure Theory for Pro-Lie Algebras, Pro-Lie Groups and Connected 
Locally Compact Groups}, EMS Publishing House, Z\"urich 2006, to appear

\[HMP04 Hofmann, K. H., Morris, S. A., and D. Poguntke, {\it 
The exponential function of locally connected compact abelian groups}, 
Forum Math. {\bf 16:1} (2004), 1--16

\[HoNe06 Hofmann, K.~H., and K.-H.\ Neeb, 
{\it Pro-Lie groups as infinite-dimensional Lie groups}, in preparation

\[Hop42 Hopf, H., {\it Fundamentalgruppe und zweite Bettische Gruppe}, 
Comment. Math. Helv. {\bf 14} (1942), 257--309

\[vHo52a van Hove, L., {\it Topologie des espaces fonctionnels analytiques, et des groups 
infinis des transformations}, Bull. de la Classe de Sc., Acad. Roy. de Belgique, 
S\'er. 5, {\bf 38} (1952), 333--351 

\[vHo52b ---, {\it L'ensemble des fonctions analytiques sur un compact en tant 
qu'alg\`ebre topologique}, Bull. Soc. Math. Belgique 1952, (1953), 8--17

\[Is96 Ismagilov, R. S., ``Representations of Infinite-Dimensional
Groups,'' Translations of Math. Monographs {\bf 152},
Amer. Math. Soc., 1996 

\[Ka85 Kac, V. G., {\it Constructing groups associated to infinite-dimensional 
Lie algebras}, in  ``Infinite-Dimensional Groups with Applications,'' 
V. Kac., Ed., MSRI Publications {\bf 4}, 
Springer-Verlag, Berlin, Heidelberg, New York, 1985

\[Ka90 ---, ``Infinite-dimensional Lie Algebras", Cambridge University 
Press, 1990

\[KP83 Kac, V. G., and D. H. Peterson, {\it Regular functions on certain
in\-fi\-nite-di\-men\-sio\-nal groups}, in ``Arithmetic and Geometry,'' Vol. 2,
M.\ Artin, J.\ Tate, eds., Birk\-h\"auser, Boston, 1983 

\[KaRo01 Kamran, N., and T. Robart, {\it A manifold structure
for analytic Lie pseudogroups of infinite type}, 
J. Lie Theory {\bf 11} (2001), 57--80 

\[KaRo04 ---, {\it An infinite-dimensional manifold structure 
for analytic Lie pseudogroups of infinite type}, 
Internat. Math. Res. Notices {\bf 34} (2004), 1761--1783 

\[Ka81  Kaup, W., {\it \"Uber die Klassifikation der symmetrischen
her\-mi\-te\-schen Mannigfaltigkeiten unendlicher Dimension I},
Math. Annalen {\bf 257} (1981), 463--486 

\[Ka83a ---, {\it A Riemann mapping theorem for bounded symmetric
domains in complex Banach spaces}, Math. Z. {\bf 183} (1983), 503--529

\[Ka83b ---, {\it \"Uber die Klassifikation der symmetrischen
her\-mi\-te\-schen Mannigfaltigkeiten unendlicher Dimension II},
Math. Annalen {\bf 262} (1983), 57--75

\[KKM05 Kedra, J., Kotchick, D., and S. Morita, {\it Crossed flux homomorphisms 
and vanishing theorems for flux groups}, Preprint, 
math.AT/0503230, Aug. 2005

\[Ke74 Keller, H. H., ``Differential Calculus
in Locally Convex Spaces,'' Springer-Verlag, 1974


\[Ki05 Kirillov, A., {\it The orbit method beyond Lie groups. Infinite-dimensional 
groups}, Surveys in modern mathematics, 292--304; 
London Math. Soc. Lecture Note Ser. {\bf 321}, 
Cambridge Univ. Press, Cambridge, 2005 

\[KY87 Kirillov, A.~A., and D.~V.~Yuriev, {\it K\"ahler geometry of
the infinite-dimensional homogeneous space $M =
\Diff_+(\SS^1)/\Rot(\SS^1)$}, Funct.\ Anal.\ and Appl. {\bf 21}
(1987), 284--294 

\[KYMO85 Kobayashi, O., Yoshioka, A., Maeda, Y., and H. Omori, {\it The 
theory of infinite-dimensional Lie groups and its applications}, 
Acta Appl. Math.  {\bf 3:1} (1985), 71--106

\[K\"o69 K\"othe, G., ``Topological Vector Spaces I,'' Grundlehren der
Math. Wis\-sen\-schaf\-ten {\bf 159}, Springer-Verlag, Berlin etc., 1969 

\[Kop70 Kopell, N., {\it Commuting diffeomorphisms}, Proc. Symp. Pure Math. 
{\bf 14}, Amer. Math. Soc., 1970, 165--184 

\[Kos70 Kostant, B., {\it Quantization and unitary representations}, 
in ``Lectures in Modern Analysis and Applications III,''  
Springer Lecture Notes Math. {\bf 170} (1970), 87--208

\[KM97 Kriegl, A., and P.\ Michor, ``The Convenient Setting of
Global Analysis,'' Math.\ Surveys and Monographs {\bf 53}, Amer.\
Math.\ Soc., 1997 

\[KM97b ---, {\it Regular infinite-dimensional Lie groups}, 
J. Lie Theory {\bf 7} (1997), 61--99

\[Ku65 Kuiper, N. H., {\it The homotopy type of the unitary group of
Hilbert space}, Topology {\bf 3} (1965), 19--30 

\[Kum02 Kumar, S., ``Kac-Moody Groups, their Flag Varieties and Representation Theory,'' 
Progress in Math. {\bf 204}, Birkh\"auser Boston, Inc., Boston, MA, 2002

\[Kur59 Kuranishi, M., {\it On the local theory of continuous infinite pseudo groups I}, 
Nagoya Math. J. {\bf 15} (1959), 225--260 

\[LMP98 Lalonde, F., D. McDuff, and L. Polterovich, 
{\it On the flux conjectures}, in ``Geometry, topology, and dynamics'' 
(Montreal, PQ, 1995); 69--85, CRM Proc. Lecture Notes {\bf  15}, 
Amer. Math. Soc., Providence, RI, 1998 

\[La99 Lang, S., ``Fundamentals of Differential Geometry,'' Graduate
Texts in Math. {\bf 191}, Springer-Verlag, 1999 

\[Lar99 Laredo, V.\ T., {\it Integration of unitary representations of
infinite dimensional Lie groups}, J.\ Funct.\ Anal. {\bf 161:2}
(1999), 478--508 

\[Las57 Lashof, R. K., {\it Lie algebras of locally compact groups}, Pacific J. Math. {\bf 7} 
(1957), 1145--1162

\[Lau55 Laugwitz, D., {\it \"Uber unendliche kontinuierliche Gruppen. I. Grundlagen der 
Theorie; Untergruppen}, Math. Ann. {\bf 130} (1955), 337--350 

\[Lau56 ---, {\it \"Uber unendliche kontinuierliche Gruppen. II. 
Strukturtheorie lokal Banachscher Gruppen}, Bayer. Akad. Wiss. Math.-Nat. 
Kl. S.-B.  1956,  (1957), 261--286

\[LaTi66 Lazard, M. and J. Tits, {\it Domaines d'injectivit\'e de l'application exponentielle}, Topology {\bf 4} (1966), 315--322 


\[Lec80 Lecomte, P., {\it Sur 
l'alg\`ebre de Lie des sections d'un fibr\'e en alg\`ebre de Lie}, 
Ann. Inst. Fourier {\bf 30} (1980), 35--50 


\[Lec85 ---, {\it  Sur la suite exacte canonique associ\'ee \`a un fibr\'e 
principal}, Bull Soc. Math. Fr. {\bf 13} (1985), 259--271 


\[Lem95 Lempert, L., {\it The Virasoro group as a complex manifold}, 
Math. Res. Letters {\bf 2} (1995), 479--495 

\[Lem97 ---, {\it The problem of complexifying
a Lie group}, in: Cordaro, P.\ D.\ et al.\ (Eds.),
``Multidimensional Complex Analysis and Partial Differential
Equations,'' AMS, Contemp.\ Math.\ {\bf 205} (1997), 169--176

\[Les66 Leslie, J. A., {\it On a theorem of E. Cartan}, 
Ann. Mat. Pura Appl. (4) {\bf 74} (1966), 173--177 

\[Les67 ---, {\it On a differential structure for the group 
of diffeomorphisms}, Topology {\bf 6} (1967), 263--271 

\[Les68 ---, {\it Some Frobenius theorems in global analysis}, 
J. Diff. Geom. {\bf 2} (1968), 279--297 


\[Les82 ---, {\it On the group of real analytic
diffeomorphisms of a compact real analytic manifold}, 
Trans.\ Amer.\ Math.\ Soc.\ {\bf 274} (1982), 651--669 

\[Les83 ---, {\it A Lie group structure for the group of analytic diffeomorphisms}, 
Boll. Un. Mat. Ital. A (6) {\bf 2:1} (1983), 29--37

\[Les90 ---, {\it A path functor for Kac-Moody Lie algebras}, in ``Lie 
Theory, Differential Equations and Representation Theory (Montreal, PQ, 1989),'' 
265--270; Univ. Montréal, Montreal, QC, 1990 

\[Les92 ---, {\it Some integrable subalgebras of infinite-dimensional Lie groups}, 
Trans.\ Amer.\ Math.\ Soc.\ {\bf 333} (1992), 423--443 

\[Les93 ---, {\it On the integrability of some infinite dimensional Lie algebras}, 
Preprint, Howard University, 1993  

\[Les03 ---, {\it On a diffeological group realization of
certain generalized symmetrizable Kac-Moody Lie
algebras}, J. Lie Theory {\bf 13} (2003), 427--442 

\[Lew39 Lewis, D., {\it Formal power series transformations}, Duke Math. J. 
{\bf 5} (1939), 794--805 


\[Lie80 Lie, S., {\it Theorie der Transformationsgruppen I},  Math. Ann. 
{\bf 16:4} (1880), 441--528

\[Lie95 ---, {\it Unendliche kontinuierliche Gruppen}, Abhandlungen, S\"achsische 
Gesell\-schaft der Wissenschaften {\bf 21} (1895), 43--150 

\[Lo98 Loday, J.-L., ``Cyclic Homology,'' Grundlehren der
math. Wissenschaften {\bf 301}, Springer-Verlag, Berlin, 1998 

\[Lo69 Loos, O., ``Symmetric Spaces I: General Theory,'' Benjamin, New York, 
Amsterdam, 1969

\[Los92 Losik, M. V., {\it Fr\'{e}chet manifolds
as diffeologic spaces}, Russ.\ Math. {\bf 36} (1992), 31--37 

\[LV94 Luminet, D., and A. Valette, {\it Faithful uniformly continuous 
representations of Lie groups}, J. London Math. Soc. (2) {\bf 49} (1994), 100--108 

\[ML63 MacLane, S., ``Homology,'' Grundlehren der Math. Wiss. {\bf 114}, Springer-Verlag, 1963 

\[ML78 ---, {\it Origins of the cohomology of groups}, L'Enseignement math\'em. 
{\bf 24:2} (1978), 1--29 

\[MOKY85 Maeda, Y., Omori, H., Kobayashi, O., and A. Yoshioka, {\it On 
regular Fr\'echet-Lie groups. VIII. Primordial operators and Fourier integral operators}, 
Tokyo J. Math. {\bf 8:1} (1985), 1--47

\[Mai02 Maier, P., {\it Central extensions of topological
current algebras}, in 
``Geometry and Analysis on Finite-
and Infinite-Dimensional Lie Groups,'' Eds. A.~Strasburger et al Eds., 
Banach Center Publications {\bf 55}, Warszawa, 2002

\[MN03 Maier, P., and K.-H. Neeb,  {\it Central extensions of current groups}, 
Math. Annalen {\bf 326:2} (2003), 367--415 

\[Mau88 Maurer, L., {\it \"Uber allgemeinere Invarianten-Systeme}, 
M\"unchner Berichte {\bf 43} (1888), 103--150 

\[Mais62 Maissen, B., {\it Lie-Gruppen mit Banachr\"aumen als
Parameterr\"aume}, Acta Math. {\bf 108} (1962), 229--269 

\[Mais63 ---, {\it \"Uber Topologien im Endomorphismenraum eines topologischen 
Vektorraums}, Math. Ann. {\bf 151} (1963), 283--285 

\[MR95 Marion, J., and T. Robart, {\it Regular Fr\'{e}chet Lie groups of 
invertibe elements in some inverse limits of unital involutive Banach algebras},
Georgian Math.\ J. {\bf 2} (1995), 425--444 

\[Mar67 Marsden, J. E., {\it Hamiltonian one parameter groups: A mathematical 
exposition of infinite dimensional Hamiltonian systems with applications 
in classical and quantum mechanics}, Arch. Rational Mech. Anal. {\bf 28} (1968), 
362--396

\[MA70 Marsden, J. E., and R. Abraham, {\it Hamiltonian mechanics on Lie groups and 
Hydrodynamics}, in ``Global Analysis'', Proc. Symp. Pure Math. {\bf 16}, 1970,  
Eds. S.~S.~Chern and S.~Smale, Amer. Math. Soc., Providence, RI, 237--244  



\[MaTh35 Mayer, W., and T. Y. Thomas, {\it Foundations of the theory of Lie groups}, 
Ann. of Math. {\bf 36:3} (1935), 770--822 

\[MDS98 McDuff, D., and D. Salamon, ``Introduction to Symplectic
Topology,'' Oxford Math. Monographs, 1998 

\[MD05 McDuff, D., {\it Enlarging the Hamiltonian group}, 
Preprint, math.SG/0503268, May 2005 

\[MicE59 Michael, E., {\it Convex structures and continuous selections},
Can. J. Math. {\bf 11} (1959), 556--575 

\[MicA38 Michal, A. D. {\it Differential calculus in linear topological spaces}, 
Proc. Nat. Acad. Sci. USA {\bf 24} (1938), 340--342 

\[MicA40 ---, {\it Differential of functions with arguments and values in 
topological abelian groups}, 
Proc. Nat. Acad. Sci. USA {\bf 26} (1940), 356--359 

\[MicA45 ---, {\it The total differential equation for the exponential 
function in non-com\-mu\-ta\-tive normed linear rings}, 
Proc. Nat. Acad. Sci. U. S. A.  {\bf 31} (1945),  315--317

\[MicA48 ---, {\it Differentiable infinite continuous groups in abstract 
spaces}, Revista Ci., Lima  {\bf 50} (1948), 131--140

\[MiE37  Michal, A. D., and V. Elconin, {\it Differential 
properties of abstract transformation groups
with abstract parameters}, Amer.\ J. Math.\ {\bf 59} (1937), 129--143

\[Mi80 Michor, P. W., ``Manifolds of Differentiable
Mappings,'' Shiva Publishing, Or\-pington, Kent (U.K.), 1980 

\[Mi84 ---, {\it A convenient setting for differential geometry
and global analysis I, II}, Cah.\ Topologie G\'{e}om.\ Differ.\ {\bf 25} (1984), 
63--109, 113--178

\[Mi87 ---, {\it The cohomology of the diffeomorphism group of a manifold is a 
Gelfand-Fuks cohomology}, in ``Proc. of the 14th Winter School on Abstr. Analysis, 
Srni, 1986,'' Suppl. Rend. del Circ. Mat. Palermo II {\bf 14} (1987), 235--246 

\[Mi91 ---, ``Gauge Theory for Fiber Bundles,'' Bibliopolis, ed. di fil. sci., 
Napoli, 1991 

\[MT99 Michor, P., and J.~Teichmann, {\it Description of infinite
dimensional abelian regular Lie groups}, J. Lie Theory {\bf 9:2}
(1999), 487--489 

\[Mick87 Mickelsson, J., {\it Kac-Moody groups, topology of the Dirac determinant bundle, and fermionization}, Commun. Math. Phys. {\bf 110} (1987), 173--183 

\[Mick89 ---, ``Current algebras and groups,'' Plenum Press,
New York, 1989 

\[Mil82 Milnor, J., {\it On infinite-dimensional Lie groups},
Preprint, Institute of Adv. Stud. Princeton, 1982 

\[Mil84 ---, {\it Remarks on infinite-dimensional Lie groups}, 
in DeWitt, B., Stora, R. (eds), 
``Relativit\'{e}, groupes et topologie II'' (Les Houches, 1983), 
North Holland, Amsterdam, 1984; 1007--1057 

\[MZ55 Montgomery, D., and L. Zippin, ``Topological Transformation Groups,'' 
Interscience, New York, 1955 

\[Mo61 Moser, J., {\it A new technique for the construction of solutions of nonlinear 
differential equations}, Proc. Nat. Acad. Sci. USA {\bf 47} (1961), 1824--1831 

\[My54 Myers, S. B., {\it Algebras of differentiable functions}, 
Proc. Amer. Math. Soc. {\bf 5} (1954), 917--922

\[Naga66 Nagano, T., {\it Linear differential systems with singularities 
and an application to transitive Lie algebras}, J. Math. Soc. Japan {\bf 18} (1966), 
398--404

\[Nag36 Nagumo, M., {\it Einige analytische Untersuchungen in linearen metrischen Ringen}, 
Jap. J. Math. {\bf 13} (1936), 61--80 

\[NRW91 Natarajan, L., Rodriguez-Carrington, E., and J. A. Wolf, 
{\it Differentiable structure for direct limit groups}, 
Letters Math.\ Phys. {\bf 23} (1991), 99--109

\[NRW93 ---, {\it Locally convex Lie groups}, 
Nova Journal of Algebra and Geometry {\bf 2:1}(1993), 59--87 

\[NRW94 ---, {\it New classes of infinite dimensional Lie groups}, 
Proc.\ of Symp.\ in Pure Math.\ {\bf 56:2} (1994), 377--392 

\[NRW99 ---, {\it The Bott--Borel--Weil theorem for direct limit groups}, 
Trans.\ Amer.\ Math.\ Soc. {\bf 353} (2001), 4583--4622

\[Nat35 Nathan, D. S., {\it One-parameter groups of transformations in abstract vector spaces}, 
 Duke Math. J. {\bf 1:4}  (1935), 518--526

\[Ne98 Neeb, K.--H., {\it Holomorphic highest weight representations
of infinite dimensional complex classical groups}, 
J.\ Reine Angew. Math.\ {\bf 497} (1998), 171--222  

\[Ne99 ---, ``Holomorphy and Convexity in Lie Theory,'' 
Expositions in Mathematics {\bf 28}, de Gruyter Verlag, Berlin, 1999 

\[Ne01a ---, {\it Representations of infinite dimensional
groups}, in ``Infinite Dimensional K\"ahler Manifolds,'' 
A. Huckleberry, T. Wurzbacher (Eds.), DMV-Seminar {\bf 31}, 
Birkh\"auser Verlag, 2001; 131--178  

\[Ne01b ---, {\it Locally finite Lie algebras with unitary highest weight
representations}, manu\-scrip\-ta mathe\-matica {\bf 104:3} (2001), 343--358

\[Ne02a ---, {\it Central extensions of infinite-dimensional
Lie groups}, Annales de l'Inst. Fourier {\bf 52:5} (2002), 1365--1442 

\[Ne02b ---, {\it Classical Hilbert--Lie groups, their extensions and their
homotopy groups}, 
in  ``Geometry and Analysis on Finite and Infinite-dimensional Lie Groups,''
Eds. A. Strasburger, W. Wojtynski, J.\ Hilgert and K.-H. Neeb, 
Banach Center Publications {\bf 55}, Warszawa 2002; 87--151 

\[Ne02c ---, {\it A Cartan--Hadamard Theorem for
Banach--Finsler manifolds}, Geometriae Dedicata {\bf 95} (2002), 115--156 

\[Ne02d ---, {\it Universal central extensions of Lie groups},
Acta Appl. Math. {\bf 73:1,2} (2002), 175--219 

\[Ne03 ---, {\it Locally convex root graded Lie algebras}, 
Travaux Math. {\bf 14} (2003), 25--120  

\[Ne04a ---, {\it Abelian extensions of infinite-dimensional Lie
groups}, Travaux Math. {\bf 15} (2004), 69--194 

\[Ne04b ---, {\it Infinite-dimensional Lie groups and their
representations}, in ``Lie Theory: Lie Algebras and Representations,'' 
Progress in Math. {\bf 228}, Ed. J.~P.~Anker, B.~\O{}rsted, 
Birkh\"auser Verlag, 2004; 213--328 

\[Ne04c ---,  {\it Current groups for non-compact manifolds and their 
central extensions}, in "Infinite Dimensional Groups and Manifolds". 
T.~Wurzbacher (Ed.). IRMA Lectures in Mathematics and
Theoretical Physics {\bf 5}, de Gruyter Verlag, Berlin, 2004; 109--183

\[Ne05 ---, {\it Non-abelian extensions of infinite-dimensional Lie groups}, 
Ann. Inst. \break Fourier, to appear 
 
\[Ne06a ---, {\it Lie algebra extensions and higher order cocycles}, 
J. Geom. Symm. Phys. {\bf 5} (2006), 48--74 

\[Ne06b ---, {\it Non-abelian extensions of topological Lie algebras}, 
Communications in Algebra {\bf 34} (2006), 991--1041 

\[Ne06c ---, {\it On the period group of a continuous inverse algebra}, in preparation  

\[NS01 Neeb, K.-H., and N. Stumme, {\it On the classification of 
locally finite split simple Lie algebras}, J. Reine Angew. Math. {\bf 533} (2001),
25--53

\[NV03 Neeb, K.-H. and C. Vizman, {\it Flux homomorphisms and principal bundles over
infinite-dimensional manifolds}, Monatshefte Math. {\bf 139} (2003), 309--333

\[NeWa06a Neeb, K.-H., and F. Wagemann, {\it The second cohomology of current algebras 
of general Lie algebras}, Canadian J. Math., to appear 

\[NeWa06b ---, {\it Lie group structures on groups of maps on non-compact manifolds}, 
in preparation 

\[Neh93 Neher, E., {\it Generators and relations for $3$-graded Lie
algebras}, Journal of Algebra {\bf 155} (1993), 1--35 

\[Neh96 ---, {\it Lie algebras graded by  $3$-graded root systems and 
Jordan pairs covered by grids}, Amer. J. Math. {\bf 118} (1996), 439--491 

\[vN29 von Neumann, J., {\it \"Uber die analytischen Eigenschaften von Gruppen 
linearer Transformationen}, Math. Zeit. {\bf 30:1} (1929), 3--42 

\[Olv93 Olver, P. J., ``Applications of Lie Groups to Differential Equations.''  
Second edition,  Graduate Texts in Math. {\bf 107}, 
Springer-Verlag, New York, 1993

\[Omo70 Omori, H., {\it On the group of diffeomorphisms on a compact manifold}, 
in ``Global Analysis'' (Proc. Sympos. Pure Math., Vol. XV, Berkeley, Calif., 1968), 
Amer. Math. Soc., Providence, R.I., 1970, 167--183

\[Omo73 ---, {\it Groups of iffeomorphisms and their subgroups}, 
Trans. Amer. Math. Soc. {\bf 179} (1973), 85--122

\[Omo74 ---, ``Infinite Dimensional Lie Transformation Groups,'' 
Lect. Notes Math. {\bf 427}, Springer-Verlag, Berlin-New York, 1974

\[Omo78 ---, {\it On Banach-Lie groups acting on finite-dimensional 
manifolds}, T\^{o}hoku \break Math.\ J. {\bf 30} (1978), 223--250 

\[Omo80 ---, {\it A method of classifying expansive singularities}, 
J. Diff. Geom. {\bf 15} (1980), 493--512 

\[Omo81 ---, {\it A remark on non-enlargible Lie algebras}, 
J. Math. Soc. Japan {\bf 33:4} (1981), 707--710 

\[Omo97 ---, {\it Infinite-Dimensional Lie Groups}, Translations
of Math. Monographs {\bf 158}, Amer. Math. Soc., 1997 

\[OdH71 Omori, H., and P. de la Harpe, {\it Op\'eration de groupes de Lie 
banachiques sur les vari\'et\'es diff\'erentielles de dimension finie}, 
C. R. Acad. Sci. Paris S\'er. A-B  {\bf 273} (1971), A395--A397

\[OdH72 ---, {\it About interactions between 
Banach-Lie groups and finite dimensional mani\-folds}, 
J. Math. Kyoto Univ. {\bf 12} (1972), 543--570

\[OMY80 Omori, H., Maeda, Y. and A. Yoshioka, {\it On regular Fr\'echet-Lie groups. I. 
Some differential geometrical expressions of Fourier 
integral operators on a Riemannian manifold}, Tokyo J. Math.  {\bf 3:2} (1980), 
353--390

\[OMY81 ---, {\it On regular Fr\'echet-Lie groups. II. Composition 
rules of Fourier-integral operators on a Riemannian manifold}, 
Tokyo J. Math.  {\bf 4:2} (1981), 221--253

\[OMYK81 Omori, H., Maeda, Y., Yoshioka, A., and O. Kobayashi, 
{\it On regular Fr\'echet-Lie groups. III. A second cohomology class related 
to the Lie algebra of pseudodifferential operators of order one}, 
Tokyo J. Math.  {\bf 4:2}  (1981), 255--277

\[OMYK82 ---, {\it On 
regular Fr\'{e}chet-Lie groups IV. Definition and fundamental theorems},
Tokyo J. Math. {\bf 5} (1982), 365--398 

\[OMYK83a ---, {\it On regular Fr\'echet-Lie groups. V. Several basic properties}, 
Tokyo J. Math. {\bf 6:1} (1983), 39--64 

\[OMYK83b ---, {\it On regular Fr\'echet-Lie groups. VI. Infinite-dimensional 
Lie groups which appear in general relativity}, 
Tokyo J. Math.  {\bf 6:2}  (1983), 217--246

\[Ono04 Ono, K., {\it Floer-Novikov cohomology and the flux conjecture}, Preprint, 2004 

\[Ot95 Ottesen, J. T., ``Infinite Dimensional Groups and Algebras in Quantum 
\linebreak Physics'', Springer Verlag, Lecture Notes in 
Physics {\bf m 27}, 1995 

\[Pa57 Palais, R. S., ``A Global Formulation of the Lie Theory of Transformation Groups,'' 
Mem. Amer. Math. Soc. {\bf 22}, Amer. Math. Soc., 1957 


\[Pala71 Palamodov, V. P., {\it Homological methods in the theory of locally convex 
spaces}, Russian Math. Surveys {\bf 26:1} (1971), 1--64 

\[Pali68 Palis, J., {\it On Morse-Smale dynamical systems}, Topology {\bf 8} (1968), 385--404

\[Pali74 ---, {\it Vector fields generate few diffeomorphisms}, Bull Amer. Math. 
Soc. {\bf 80:3} (1974), 503--505 

\[Pe92 Pestov, V.\ G., {\it Nonstandard hulls
of Banach--Lie groups and algebras},
Nova J.\ Algebra Geom.\ {\bf 1} (1992),
371--381

\[Pe93a ---, {\it Free Banach--Lie algebras, couniversal Banach-Lie groups, 
and more},  Pacific J. Math. {\bf 157:1}  (1993),  137--144 

\[Pe93b ---, {\it Enlargible Banach--Lie algebras
and free topological groups},
Bull.\ Aust.\ Math.\ Soc.\ {\bf 48} (1993), 13--22

\[Pe95a  ---, {\it Correction to ``Free Banach--Lie algebras, couniversal Banach-Lie groups, 
and more''}, Pacific J. Math.  {\bf 171:2}  (1995), 585--588 

\[Pe95b ---, {\it Regular Lie groups and a theorem of Lie-Palais}, 
J. Lie Theory {\bf 5:2} (1995), 173--178

\[Pic00a Pickrell, D., ``Invariant Measures for Unitary Groups 
Associated to Kac-Moody Lie Algebras,'' 
Mem. Amer. Math. Soc.  {\bf 146/693}, 2000 

\[Pic00b ---, {\it On the action of the group of 
diffeomorphisms of a surface on sections of the determinant line bundle}, 
Pacific J. Math. {\bf 193:1}  (2000),  177--199 

\[Pis76 Pisanelli, D., {\it An extension of the exponential of a matrix 
and a counter example to the inversion theorem in a space $H(K)$}, 
Rendiconti di Matematica (6) (3) {\bf 9:3} (1976), 465--475

\[Pis77 ---, {\it An example of an infinite
Lie group}, Proc.\ Amer.\ Math.\ Soc.\ {\bf 62} (1977), 156--160

\[Pis79 ---, {\it The second Lie theorem in the group $\Gh(n,\C)$}, 
in ``Advances in Holomorphy,'' J. A. Barroso (ed.), North Holland Publ., 1979 


\[Pol01 Polterovich, L.,``The geometry of the group of 
symplectic diffeomorphisms,'' Lectures in Mathematics, ETH Z\"urich, 
Birkh\"auser Verlag, Basel, 2001

\[Po39 Pontrjagin, L., ``Topological Groups,'' 
Princeton Mathematical Series, v. {\bf 2}, 
Princeton University Press, Princeton, 1939

\[PS86 Pressley, A., and G. Segal, ``Loop Groups," Oxford University Press, 
Oxford, 1986

\[Pu52 Pursell, M. E., ``Algebraic structures associated with smooth manifolds,'' 
Thesis, Purdue Univ., 1952 

\[PuSh54 Pursell, M. E., and M. E. Shanks, {\it The Lie algebra of a smooth manifold}, 
Proc. Amer. Math. Soc. {\bf 5} (1954), 468--472 

\[PW52 Putnam, C.~R., and A.\ Winter, {\it The orthogonal group in
Hilbert space}, Amer.\ J.\ Math.\ {\bf 74} (1952), 52--78

\[RS81 Ratiu, T., and Schmid, R., {\it The differentiable structure
of three remarkable diffeomorphism groups}, Math. Z. {\bf 177} (1981),
81--100

\[RO03 Ratiu, T., and A. Odzijewicz, {\it Banach Lie--Poisson spaces and reduction}, 
Comm. Math. Phys. {\bf 243:1} (2003),  1--54

\[RO04 ---, {\it Extensions of Banach Lie--Poisson spaces}, 
 J. Funct. Anal. {\bf 217:1} (2004), 103--125


\[Ri50 Ritt, J. F., {\it Differential groups and formal Lie theory for an infinite number 
of parameters}, Ann. of Math. (2) {\bf 52} (1950), 708--726

\[Rob96 Robart, T., {\it Groupes de Lie de dimension infinie. Second et troisi\`eme 
th\'eor\`emes de Lie. I. Groupes de premi\`ere esp\`ece}, 
C. R. Acad. Sci. Paris S\'er. I Math. {\bf 322:11} (1996), 1071--1074

\[Rob97 ---, {\it Sur l'int\'{e}grabilit\'{e}
des sous-alg\`{e}bres de Lie en dimension infinie}, 
Canad.\ J. Math. {\bf 49:4} (1997), 820--839

\[Rob02 ---, {\it Around the exponential mapping}, in ``Infinite Dimensional 
Lie Groups in Geometry and Representation Theory,'' 11--30; World Sci. Publ., 
River Edge, NJ, 2002 

\[Rob04 ---, {\it On Milnor's regularity and the path-functor for the 
class of infinite dimensional Lie algebras of CBH type}, 
Algebras, Groups and Geometries {\bf 21} (2004), 367--386 

\[RK97 Robart, T., and N. Kamran, {\it Sur la th\'eorie locale des pseudogroupes 
de transformations continus infinis. I}, Math. Ann. {\bf 308:4} (1997), 593--613

\[Rod89 Rodriguez-Carrington, E., 
{\it Lie groups associated to Kac--Moody Lie algebras: an analytic approach}, 
in ``Infinite-dimensional Lie Algebras and Groups,'' (Luminy-Marseille, 1988), Adv. 
Ser. Math. Phys. {\bf 7}, World Scientific Publ., Teaneck, NJ, 1989; 
57--69 

\[Rog95 Roger, C., {\it Extensions centrales d'alg\`ebres et de groupes
de Lie de dimension infinie, alg\`ebres de Virasoro et
g\'en\'eralisations}, Reports on Math. Phys. {\bf 35} (1995), 225--266

\[Ros94 Rosenberg, J., ``Algebraic K-theory and its Applications,''
Graduate Texts in Math. {\bf 147}, Springer-Verlag, 1994 

\[Ru73 Rudin, W., ``Functional Analysis,'' McGraw Hill, 1973

\[Sch87 Schmid, R., ``Infinite-dimensional Hamiltonian systems,'' 
Monographs and Textbooks in Physical Science, Lecture Notes {\bf 3}, 
Bibliopolis, Naples, 1987 

\[Sch04 ---, {\it Infinite dimensional Lie groups with applications 
to mathematical physics}, J. Geom. Symm. Phys. {\bf  1}  (2004), 54--120

\[SAR84 Schmid, R., Adams, M. and T. Ratiu, {\it The group of Fourier 
integral operators as symmetry group}, in ``XIIIth international colloquium 
on group theoretical methods in physics'' (College Park, Md., 1984), 
World Sci. Publishing, Singapore, 1984, 246--249 

\[Sc60 Schue, J.\ R., {\it Hilbert space methods in the theory of Lie
algebras}, Transactions of the Amer.\ Math.\ Soc. {\bf 95} (1960),
69--80 

\[Sc61 ---, {\it Cartan decompositions for $L^{*}$-algebras}, Trans.\
Amer.\ Math.\ Soc.\ {\bf 98} (1961), 334--349

\[SchF90a Schur, F., {\it Neue Begr\"undung der 
Theorie der endlichen Transformationsgruppen}, 
Math. Ann. {\bf 35} (1980), 161--197 

\[SchF90b ---,  {\it Beweis f\"ur die Darstellbarkeit 
der infinitesimalen Transformationen aller transitiven endlichen Gruppen durch 
Quotienten best\"andig convergenter Potenz\-rei\-hen}, Leipz. Ber. {\bf 42} (1890), 1--7 

\[Se81 Segal, G., {\it Unitary representations of some
infinite-dimensional groups}, Comm.\ Math.\ Phys. {\bf 80} (1981), 301--342 

\[Ser65 Serre, J.-P., ``Lie algebras and Lie groups,'' 
Lect. Notes. Math. {\bf 1500}, Springer-Verlag, 1965 (1st ed.) 

\[SiSt65 Singer, I. M., and S. Sternberg, {\it The infinite groups of Lie and Cartan. I. The transitive 
groups}, J. d'Anal. Math. {\bf 15} (1965), 1--114 

\[So84 Souriau, J.-M., {\it  Groupes diff\'{e}rentiels de physique math\'{e}matique},
in ``Feuilletages 
et Quantification G\'{e}ome\-trique,'' Dazord, P. and N. Desolneux-Moulis (Eds.), 
Journ. lyonnaises Soc.\ math.\ 
France 1983, S\'{e}min.\ sud-rhodanien de g\'{e}om.\ II, Hermann, Paris, 1984; 
73--119  

\[So85 ---, {\it Un algorithme g\'en\'erateur de structures quantiques}, 
Soc. Math. Fr., Ast\'e\-risque, hors s\'erie, 1985, 341--399 

\[St61 Sternberg, S., {\it Infinite Lie groups and the formal aspects 
of dynamical systems}, J. Math. Mech. {\bf 10} (1961), 451--474 

\[St99 Stumme, N., ``The Structure of Locally Finite Split Lie 
algebras,'' Ph. D.\ thesis, Darmstadt University of Technology, 1999

\[Sus73 Sussmann, H. J., {\it Orbits of families of vector fields and 
integrability of distributions}, Trans. Amer. Math. Soc. {\bf 180} (1973), 171--188

\[Su88 Suto, K., {\it Groups associated with unitary forms of
Kac--Moody algebras}, J.\ Math.\ Soc.\ Japan {\bf 40:1} (1988),
85--104 

\[Su97 ---, {\it Borel--Weil type theorem for the flag manifold of a 
generalized Kac--Moody algebra}, J. Algebra {\bf 193:2} (1997), 529--551

\[Swa62 Swan, R. G., {\it Vector bundles and projective modules}, 
Trans. Amer. Math. Soc. {\bf 105} (1962), 264--277 


\[Swi65 Swierczkowski, S., {\it  Embedding theorems for local analytic groups}, 
Acta Math. {\bf 114} (1965), 207--235 

\[Swi67 ---,  {\it Cohomology of local group extensions}, 
Transactions of the Amer. Math. Soc. {\bf 128} (1967), 291--320 

\[Swi71 ---, {\it  The path-functor on Banach Lie
algebras}, Nederl.\ Akad.\ Wet., Proc.\ A 74;  Indag.\ Math. {\bf 33} (1971), 235--239



\[Ta68 Tagnoli, A., {\it La variet\`a analitiche reali come spazi omogenei}, 
Boll. Un. Mat. Ital. (s4) {\bf 1} (1968), 422--426 




\[Ti83 Tits, J., ``Liesche Gruppen und Algebren", Springer, New York, 
Heidelberg, 1983

\[Tr67 Treves, F., ``Topological Vector Spaces, Distributions, and
Kernels,'' Academic Press, New York, 1967 


\[Up85 Upmeier, H., ``Symmetric Banach Manifolds and Jordan 
$C^*$-algebras,'' North Holland Mathematics Studies, 1985 

\[Va84 Varadarajan, V. S., ``Lie Groups, Lie Algebras, and Their Representations,'' 
Graduate Texts in Math. {\bf 102}, Springer--Verlag, 1984 

\[Vo87 Vogt, D., {\it On the functors $\Ext^1(E,F)$ for Fr\'echet spaces}, Studia 
Math. {\bf 85} (1987), 163--197 

\[Wae54a Waelbroeck, L., {\it Les alg\`{e}bres \`{a} inverse continu}, 
C. R. Acad.\ Sci.\ Paris {\bf 238} (1954), 640--641 
  
\[Wae54b ---,  {\it Le calcul symbolique dans les alg\`ebres 
commutatives}, J. Math. Pures Appl. {\bf 33} (1954), 147--186

\[Wae54c ---, {\it Structure des alg\`{e}bres \`{a} inverse continu}, 
C. R. Acad.\ Sci.\ Paris {\bf 238} (1954), 762--764 
  

\[Wae71 ---, ``Topological Vector Spaces and Algebras,'' 
Springer-Verlag, Berlin, Heidelberg, New York, 1971 



\[Wa72 Warner, G., ``Harmonic Analysis on Semisimple Lie Groups I,'' Springer, 
Berlin, Heidelberg, New York, 1972

\[Wei69 Weinstein, A., {\it  Symplectic structures on Banach manifolds}, 
Bull. Amer. Math. Soc. {\bf 75} (1969), 1040--1041


\[Wer95 Werner, D., ``Funktionalanalysis,'' Springer-Verlag, Berlin,
Heidelberg, 1995 

\[Wie49 Wielandt, H., {\it \"Uber die Unbeschr\"anktheit der Operatoren 
der Quantenmechanik},  Math. Ann.  {\bf 121} (1949), 21

\[Wo05a Wockel, Chr., {\it The Topology of Gauge Groups}, 
submitted; math-ph/0504076 

\[Wo05b ---, {\it Smooth Extensions and Spaces of Smooth and Holomorphic 
Mappings}, J. Geom. Symm. Phys., to appear; math.DG/0511064 

\[Wol05 Wolf, J. A., {\it Principal series representations
of direct limit groups},  Compos. Math. {\bf 141:6}  (2005), 1504--1530

\[Woj06 Wojty\'nski, W., {\it Effective integration of Lie algebras}, 
J. Lie Theory {\bf 16} (2006), to appear 

\[W\"u03 W\"ustner, M., {\it Supplements on the theory of exponential Lie groups}, 
J. Algebra {\bf 265:1}  (2003), 148--170

\[W\"u05 ---, {\it The classification of all simple Lie groups with surjective 
exponential map},  J. Lie Theory  {\bf 15:1} (2005), 269--278


\[Yo36 Yosida, K., {\it On the groups embedded in the metrical complete ring}, 
Japanese J. Math. {\bf 13} (1936), 7-26 

}

\def\address
{Karl-Hermann Neeb

Technische Universit\"at Darmstadt 

Schlossgartenstrasse 7

D-64289 Darmstadt 

Deutschland

neeb@mathematik.tu-darmstadt.de}

\references
\lastpage 
\bye